# NEW CLASSES OF NEUTROSOPHIC LINEAR ALGEBRAS

**W. B. Vasantha Kandasamy**
e-mail: **vasanthakandasamy@gmail.com**
web: **http://mat.iitm.ac.in/home/wbv/public_html/**
**www.vasantha.in**

**Florentin Smarandache**
e-mail: **smarand@unm.edu**

**K. Ilanthenral**
e-mail: **ilanthenral@gmail.com**

**2010**

# NEW CLASSES OF NEUTROSOPHIC LINEAR ALGEBRAS

**W. B. Vasantha Kandasamy**
**Florentin Smarandache**
**K. Ilanthenral**

**2010**

# CONTENTS

# PREFACE

In this book we introduce mainly three new classes of linear algebras; neutrosophic group linear algebras, neutrosophic semigroup linear algebras and neutrosophic set linear algebras. The authors also define the fuzzy analogue of these three structures.

This book is organized into seven chapters. Chapter one is introductory in content. The notion of neutrosophic set linear algebras and neutrosophic neutrosophic set linear algebras are introduced and their properties analysed in chapter two. Chapter three introduces the notion of neutrosophic semigroup linear algebras and neutrosophic group linear algebras. A study of their substructures are systematically carried out in this chapter.

The fuzzy analogue of neutrosophic group linear algebras, neutrosophic semigroup linear algebras and neutrosophic set linear algebras are introduced in chapter four of this book. Chapter five introduces the concept of neutrosophic group bivector spaces, neutrosophic bigroup linear algebras, neutrosophic semigroup (bisemigroup) linear algebras and neutrosophic biset bivector spaces. The fuzzy analogue of these concepts are given in chapter six. An interesting feature of this book is it contains nearly 424 examples of these new notions. The final chapter suggests over 160 problems which is another interesting feature of this book.

Finally it is an immense pleasure to thank Dr. K. Kandasamy for proof-reading and Kama and Meena without whose help the book would have been impossibility.

W.B.VASANTHA KANDASAMY
FLORENTIN SMARANDACHE
K.ILANTHENRAL

**Chapter One**

# INTRODUCTION

In this chapter we assume fields to be of any desired characteristic and vector spaces are taken over any field. We denote the indeterminacy by 'I' as i will make a confusion, as it denotes the imaginary value, viz. $i^2 = -1$ that is $\sqrt{-1} = i$. The indeterminacy I is such that $I \cdot I = I^2 = I$.

In this chapter we just recall some of the basic neutrosophic structures used in this book.

In this chapter we recall the notion of neutrosophic groups. Neutrosophic groups in general do not have group structure. We also define yet another notion called pseudo neutrosophic groups which have group structure. As neutrosophic groups do not have group structure the classical theorems viz. Sylow, Lagrange or Cauchy are not true in general which forces us to define notions like Lagrange neutrosophic groups, Sylow neutrosophic groups and Cauchy elements.

We just give the basic definition alone as we use it only in the construction of neutrosophic group vector spaces and neutrosophic group linear algebras which are analogous structures of neutrosophic set vector spaces, neutrosophic semi group vector spaces, neutrosophic set linear algebras and neutrosophic semigroup linear algebras.

**DEFINITION 1.1:** *Let $(G, *)$ be any group, the neutrosophic group is generated by I and G under $*$ denoted by $N(G) = \{\langle G \cup I \rangle, *\}$.*

*Example 1.1:* Let $Z_7 = \{0, 1, 2, \ldots, 6\}$ be a group under addition modulo 7. $N(G) = \{\langle Z_7 \cup I \rangle, \text{`+'} \text{ modulo } 7\}$ is a neutrosophic group which is in fact a group. For $N(G) = \{a + bI \mid a, b \in Z_7\}$ is a group under '+' modulo 7. Thus this neutrosophic group is also a group.

*Example 1.2:* Consider the set $G = Z_5 \setminus \{0\}$, G is a group under multiplication modulo 5. $N(G) = \{\langle G \cup I \rangle$, under the binary operation, multiplication modulo 5$\}$. $N(G)$ is called the neutrosophic group generated by $G \cup I$. Clearly $N(G)$ is not a group, for $I^2 = I$ and I is not the identity but only an indeterminate, but $N(G)$ is defined as the neutrosophic group.

Thus based on this we have the following theorem:

**THEOREM 1.1:** *Let $(G, *)$ be a group, $N(G) = \{\langle G \cup I \rangle, *\}$ be the neutrosophic group.*

1. *$N(G)$ in general is not a group.*
2. *$N(G)$ always contains a group.*

*Proof:* To prove $N(G)$ in general is not a group it is sufficient we give an example; consider $\langle Z_5 \setminus \{0\} \cup I \rangle = G = \{1, 2, 4, 3, I, 2I, 4I, 3I\}$; G is not a group under multiplication modulo 5. In fact $\{1, 2, 3, 4\}$ is a group under multiplication modulo 5.

N(G) the neutrosophic group will always contain a group because we generate the neutrosophic group N(G) using the group G and I. So $G \subsetneq N(G)$; hence N(G) will always contain a group.

Now we proceed onto define the notion of neutrosophic subgroup of a neutrosophic group.

**DEFINITION 1.2:** *Let $N(G) = \langle G \cup I \rangle$ be a neutrosophic group generated by G and I. A proper subset P(G) is said to be a neutrosophic subgroup if P(G) is a neutrosophic group i.e. P(G) must contain a (sub) group.*

*Example 1.3*: Let $N(Z_2) = \langle Z_2 \cup I \rangle$ be a neutrosophic group under addition. $N(Z_2) = \{0, 1, I, 1 + I\}$. Now we see $\{0, I\}$ is a group under + in fact a neutrosophic group $\{0, 1 + I\}$ is a group under '+' but we call $\{0, I\}$ or $\{0, 1 + I\}$ only as pseudo neutrosophic groups for they do not have a proper subset which is a group. So $\{0, I\}$ and $\{0, 1 + I\}$ will be only called as pseudo neutrosophic groups (subgroups).

We can thus define a pseudo neutrosophic group as a neutrosophic group, which does not contain a proper subset which is a group. Pseudo neutrosophic subgroups can be found as a substructure of neutrosophic groups. Thus a pseudo neutrosophic group though has a group structure is not a neutrosophic group and a neutrosophic group cannot be a pseudo neutrosophic group. Both the concepts are different.

Now we see a neutrosophic group can have substructures which are pseudo neutrosophic groups which is evident from the following example:

*Example 1.4:* Let $N(Z_4) = \langle Z_4 \cup I \rangle$ be a neutrosophic group under addition modulo 4. $\langle Z_4 \cup I \rangle = \{0, 1, 2, 3, I, 1 + I, 2I, 3I, 1 + 2I, 1 + 3I, 2 + I, 2 + 2I, 2 + 3I, 3 + I, 3 + 2I, 3 + 3I\}$. $o(\langle Z_4 \cup I \rangle) = 4^2$.

Thus neutrosophic group has both neutrosophic subgroups and pseudo neutrosophic subgroups. For $T = \{0, 2, 2 + 2I, 2I\}$ is a neutrosophic subgroup as $\{0\ 2\}$ is a subgroup of $Z_4$ under addition modulo 4. $P = \{0, 2I\}$ is a pseudo neutrosophic group under '+' modulo 4.

**DEFINITION 1.3:** *Let N(G) be a neutrosophic group. The number of distinct elements in N(G) is called the order of N(G). If the number of elements in N(G) is finite we call N(G) a finite neutrosophic group; otherwise we call N(G) an infinite neutrosophic group, we denote the order of N(G) by o(N(G)) or |N(G)|.*

**DEFINITION 1.4:** *Let S be a semigroup, the semigroup generated by S and I i.e. S ∪ I denoted by ⟨S ∪ I⟩ is defined to be a neutrosophic semigroup.*

It is interesting to note that all neutrosophic semigroups contain a proper subset, which is a semigroup.

*Example 1.5:* Let $Z_{12}$ = {0, 1, 2, …, 11} be a semigroup under multiplication modulo 12. Let N(S) = ⟨$Z_{12}$ ∪ I⟩ be the neutrosophic semigroup. Clearly $Z_{12}$ ⊂ ⟨$Z_{12}$ ∪ I⟩ and $Z_{12}$ is a semigroup under multiplication modulo 12.

*Example 1.6:* Let Z = {the set of positive and negative integers with zero}, Z is only a semigroup under multiplication. Let N(S) = {⟨Z ∪ I⟩} be the neutrosophic semigroup under multiplication. Clearly Z ⊂ N(S) is a semigroup.

Now we proceed on to define the notion of the order of a neutrosophic semigroup.

**DEFINITION 1.5:** *Let N(S) be a neutrosophic semigroup. The number of distinct elements in N(S) is called the order of N(S), denoted by o(N(S)).*

If the number of elements in the neutrosophic semigroup N(S) is finite we call the neutrosophic semigroup to be finite otherwise infinite. The neutrosophic semigroup given in example 1.5 is finite where as the neutrosophic semigroup given in example 1.6 is of infinite order.

Now we proceed on to define the notion of neutrosophic subsemigroup of a neutrosophic semigroup N(S).

**DEFINITION 1.6:** *Let N(S) be a neutrosophic semigroup. A proper subset P of N(S) is said to be a neutrosophic subsemigroup, if P is a neutrosophic semigroup under the operations of N (S). A neutrosophic semigroup N(S) is said to have a subsemigroup if N(S) has a proper subset, which is a semigroup under the operations of N(S).*

It is interesting to note a neutrosophic semigroup may or may not have a neutrosophic subsemigroup but it will always have a subsemigroup.

Now we proceed on to illustrate these by the following examples:

*Example 1.7:* Let $Z^+ \cup \{0\}$ denote the set of positive integers together with zero. $\{Z^+ \cup \{0\}, +\}$ is a semigroup under the binary operation '+'. Now let N(S) = $\langle Z^+ \cup \{0\}^+ \cup \{I\}\rangle$. N(S) is a neutrosophic semigroup under '+'. Consider $\langle 2Z^+ \cup I\rangle$ = P, P is a neutrosophic subsemigroup of N(S). Take R = $\langle 3Z^+ \cup I\rangle$; R is also a neutrosophic subsemigroup of N(S).

**DEFINITION 1.7:** *Let K be the field of reals. We call the field generated by $K \cup I$ to be the neutrosophic field for it involves the indeterminacy factor in it. We define $I^2 = I$, $I + I = 2I$ i.e., $I +...+ I = nI$, and if $k \in K$ then $k.I = kI$, $0I = 0$. We denote the neutrosophic field by K(I) which is generated by $K \cup I$ that is $K(I) = \langle K \cup I\rangle$. $\langle K \cup I\rangle$ denotes the field generated by K and I.*

*Example 1.8:* Let R be the field of reals. The neutrosophic field of reals is generated by R and I denoted by $\langle R \cup I\rangle$ i.e. R(I) clearly $R \subset \langle R \cup I\rangle$.

*Example 1.9:* Let Q be the field of rationals. The neutrosophic field of rationals is generated by Q and I denoted by Q(I).

**DEFINITION 1.8:** *Let K(I) be a neutrosophic field we say K(I) is a prime neutrosophic field if K(I) has no proper subfield, which is a neutrosophic field.*

***Example 1.10:*** Q(I) is a prime neutrosophic field where as R(I) is not a prime neutrosophic field for Q(I) ⊂ R(I).

Likewise we can define neutrosophic subfield.

**DEFINITION 1.9:** *Let K(I) be a neutrosophic field, P ⊂ K(I) is a neutrosophic subfield of P if P itself is a neutrosophic field. K(I) will also be called as the extension neutrosophic field of the neutrosophic field P.*

We can also define neutrosophic fields of prime characteristic p (p is a prime).

**DEFINITION 1.10:** *Let $Z_p = \{0, 1, 2, \ldots, p-1\}$ be the prime field of characteristic p. $\langle Z_p \cup I \rangle$ is defined to be the neutrosophic field of characteristic p. Infact $\langle Z_p \cup I \rangle$ is generated by $Z_p$ and I and $\langle Z_p \cup I \rangle$ is a prime neutrosophic field of characteristic p.*

***Example 1.11:*** $Z_7 = \{0, 1, 2, 3, \ldots, 6\}$ be the prime field of characteristic 7. $\langle Z_7 \cup I \rangle = \{0, 1, 2, \ldots, 6, I, 2I, \ldots, 6I, 1+I, 1+2I, \ldots, 6+6I\}$ is the prime field of characteristic 7.

**DEFINITION 1.11:** *Let G(I) by an additive abelian neutrosophic group and K any field. If G(I) is a vector space over K then we call G(I) a neutrosophic vector space over K.*

Elements of K(I) or $\langle Z_p \cup I \rangle$ or Q(I) or R(I) will also be known as neutrosophic numbers.

For more about neutrosophy please refer [10-7, 19, 25-6]. These concepts have good relevance in research, notably Smarandache's neutrosophic method which is a generalization of Hegel's dialectic, and suggests that scientific research will progress via studying the opposite ideas and the neutral ideas related to them in order to have a bigger picture.

**Chapter Two**

# SET NEUTROSOPHIC LINEAR ALGEBRA

In this chapter we for the first time introduce the notion of set neutrosophic vector spaces and set neutrosophic linear algebras and study their properties.

This chapter has four sections. Section one introduces the concept of neutrosophic sets. The notion of set neutrosophic vector spaces are introduced in section two. Section three introduces the concept of set neutrosophic linear algebras. Mixed set neutrosophic rational vector spaces and their properties are discussed in section four.

## 2.1 Types of Neutrosophic Sets

In this section we introduce a few types of neutrosophic sets essential to define the notions of set neutrosophic vector spaces and set neutrosophic linear algebras. Throughout this book set implies subset of integers or subset of rationals or subset of complex numbers or subset of modulo integers modulo n (n $\in$ N) or subset of reals. By finite set we mean a set S with finite number of distinct elements in them. If a set S has infinite number of elements in them, then we say S is of infinite cardinality.

**DEFINITION 2.1.1:** *Let $S = \{x_1, \ldots, x_n\}$; $n \in N$; if each $x_i$ is a neutrosophic number say of the form $a_i + b_i I$, $b_i \neq 0$, $a_i, b_i \in Z$ then we call S a pure neutrosophic set of integers or pure integer neutrosophic subset of pure integer neutrosophic set.*

*Example 2.1.1:* Let $S = \{5 + 2I, 7 - 3I, 15 + 8I, -9 + 3I, 8 + 27I, 12 - 43I, 43I, -50I\}$, S is a pure integer neutrosophic subset of the pure integer neutrosophic set.

*Note:* Let $PN(Z) = \{a + bI \mid a, b \in Z \text{ and } b \neq 0\}$, we call $PN(Z)$ to be the pure integer neutrosophic set or pure neutrosophic integer set. $N(Z) = \{a + bI \mid a, b \in Z\}$ is the mixed set of neutrosophic integers or mixed neutrosophic set of integers.

Clearly $PN(Z) \subseteq N(Z)$, i.e., pure neutrosophic set of integers is always a proper subset of mixed neutrosophic set of integers.
    $P(N(Z)) \cup \{0\}$ is called the pure neutrosophic set of integers with zero.

*Example 2.1.2:* Let $P = \{2I, 0, 3I + 1, 9, 4I - 5, 8 - 9I, -14, 10I\}$, P is a mixed neutrosophic subset of $N(Z)$. Clearly P is not a pure neutrosophic subset of $N(Z)$.

**THEOREM 2.1.1:** *Every pure neutrosophic subset of N(Z) is a subset of mixed subset of neutrosophic set and not conversely.*

*Proof:* Since $PN(Z)$ is a proper subset of $N(Z)$, every subset of $PN(Z)$ is also a subset of $N(Z)$. Hence the claim.

Now consider $T = \{9 + 3I, 2I, 3 - 5I, 19I, 20 - 31I, 0, 5, 7, 4I + 2\} \subseteq N(Z)$. Since $\{5, 7, 0\} \not\subseteq PN(Z)$ we see T is not a pure neutrosophic subset of $PN(Z)$, it is only a mixed neutrosophic subset of $N(Z)$.

Likewise we define N(Q) to be the mixed neutrosophic set of rationals, i.e., N(Q) = {a + bI | a, b ∈ Q} and PN(Q) = {a + bI | a ∈ Q, b ≠ 0 ∈ Q} is defined as the pure neutrosophic set of rationals. Clearly PN(Q) is a proper subset of N(Q).

*Example 2.1.3:* Let
$$T = \{\frac{17}{5}I, 5 + \frac{3}{4}I, \frac{7}{3} + I, \frac{-3}{4} + 2I, \frac{-18}{7}I, 8 + 25I\}.$$
T is clearly a pure neutrosophic rational subset of PN(Q).

*Example 2.1.4:* Let B = {0, 3, 5I + 1, 7I + 8, 14I, –17I, –26} ⊆ N(Q). Clearly B is a mixed neutrosophic rational subset of N(Q). We see B ⊊ PN(Q).

It is however interesting to note that T ⊆ N(Q).

In view of this we have the following theorem.

**THEOREM 2.1.2:** *Every pure neutrosophic rational subset is a subset of N(Q). However a mixed neutrosophic rational subset is not a subset of PN(Q).*

The proof is left as an exercise for the reader.

**DEFINITION 2.1.2:** *Let N(C) = {a + bI | a, b ∈ C} (C, the field of complex numbers); N(C) is defined as the mixed neutrosophic complex number set or mixed neutrosophic set of complex numbers. Let PN(C) = {a + bI | a ∈ C, b ≠ 0, b ∈ C}, PN(C) is defined as the pure neutrosophic set of complex numbers or pure neutrosophic complex number set. Thus we see PN(C) ⊆ N(C).*
We have following relations;

$$N(Z) \subseteq N(Q) \subseteq N(C) \text{ and}$$
$$PN(Z) \subseteq PN(Q) \subseteq PN(C).$$

*So the relation $Z \subseteq Q \subseteq C$ is preserved under the mixed neutrosophy and pure neutrosophy. Now if $Z_n$ denotes the set of integers modulo n i.e., $Z_n = \{0, 1, 2, ..., n - 1\}$. $N(Z_n) = \{a + bI \mid a, b \in Z_n\}$; $N(Z_n)$ is defined as the mixed neutrosophic set of modulo integers $Z_n$ or mixed neutrosophic modulo set of integers.*

*$PN(Z_n) = \{a + bI \mid a \in Z_n, b \neq 0, b \in Z_n\}$ is defined as the pure neutrosophic set of modulo integers $Z_n$ or the pure neutrosophic modulo integers set.*

*$N(Z_n) = \{a + bI \mid a, b \in Z_n\}$ is defined as the mixed neutrosophic set of modulo integers. Clearly $PN(Z_n) \subseteq N(Z_n)$.*

***Example 2.1.5:*** *Let $Z_3 = \{0, 1, 2\}$ be the ring of integers modulo 3. $N(Z_3) = \{0, 1, 2, I, 2I, 1 + I, 1 + 2I, 2 + 2I, 2 + I\}$ is the mixed neutrosophic set of integers modulo 3.*

*$PN(Z_3) = \{I, 2I, 1 + I, 1 + 2I, 2 + I, 2 + 2I\}$ is the pure neutrosophic set of integers modulo 3 and $PN(Z_3) \cup \{0\} = \{0, I, 2I, 1 + I, 2 + I, 1 + 2I, 2 + 2I\}$ is the pure neutrosophic set of integers with zero modulo 3. $T = \{I, 2 + 2I, 1 + 2I, 1 + I, 2I\} \subseteq PN(Z_3)$ is a pure neutrosophic subset of integers modulo 3. $P = \{0, I, 2, 2I, 1 + I\} \subseteq N(Z_3)$ is a mixed neutrosophic subset of integers modulo 3. $S = \{0, 2I, I + 2, 1 + I\} \subseteq PN(Z_3) \cup \{0\}$ is a pure neutrosophic subset of modulo integers with zero.*

Thus we have given the basic concepts of the types of neutrosophic subsets and sets which will be used in this book.

2.2 Set Neutrosophic Vector Spaces

In this section we proceed onto define the new notion of set neutrosophic vector spaces and discuss a few of the properties associated with them.

**DEFINITION 2.2.1:** *Let $S = \{x_1, ..., x_n\}$ be a mixed neutrosophic subset of integers. Let $P \subseteq Z$ be the subset of integers. If for* *every $x_i \in S$ and for every $p \in P$ we have $px_i$, $x_i p \in S$ then we call S to be a mixed neutrosophic set vector space of integers over the set P ($|P| \geq 2$).*

We shall now illustrate this by some examples.

***Example 2.2.1:*** *Let $S = \{0, I, 2, 3, 4, I + 5, -3 + I, 14 + 4I, 3 + 2I\}$ be the mixed subset of neutrosophic integers. Let $P = \{0, 1\}$ be a subset of Z. S is a mixed neutrosophic set vector space of integers over the set P.*

***Example 2.2.2:*** *Let $T = \{0, 2 + I, 7 - I, 81I, 40 - 51I, -64, 640I + 1\}$ be a proper subset of mixed neutrosophic set N(Z). T is a mixed neutrosophic set vector space over the set $P = \{0, 1\}$.*

***Example 2.2.3:*** *Let $D = \{0, \pm(2 - 5I), \pm(17 + 3I), \pm(16 - I), \pm 7I, \pm(-8I), \pm 9, \pm(-7 + 3I)\}$ be a mixed neutrosophic subset of N(Z). Take $P = \{-1, 0, 1\} \subseteq Z$. D is a mixed neutrosophic set vector space over P.*

***Example 2.2.4:*** *Let $P = \{0, I, 2I, 3I, 4I, 5I, 3 + 3I, 5 + 5I, 8 + 8I, 6I, 10I, 16I\} \subseteq PN(Z)$ be the pure neutrosophic subset of PN(Z) with $\{0\}$. We see P is a not a mixed neutrosophic set vector space over the set $T = \{0, 1, I\} \subseteq N(Z)$, further $T \not\subseteq Z$, so this is not even neutrosophic set vector space as it is not defined over a proper subset of Z. This will be dealt later.*

We proceed onto define more new concepts.

**DEFINITION 2.2.2:** *Let $S = \{y_1, \ldots, y_m\} \subseteq PN(Z) \cup \{0\}$; $y_i \in PN(Z) \cup \{0\}$; $1 \leq i \leq m$; $m \in N$ be a proper subset of $PN(Z) \cup \{0\}$. Take $\phi \neq P \subseteq Z$ to be a subset of Z with $|P| \geq 2$. If for every $y_i \in S$ and $t \in P$; $y_i t$, $ty_i \in S$ then we call S to be a pure neutrosophic integer set vector space over $P \subseteq Z$ with zero or simply pure neutrosophic integers set vector space over S.*

We now illustrate this situation by an example.

*Example 2.2.5:* Let T = {0, 1 + I, 2 + I, 3 + I, I, 2I, 3I, 9I + 4, 20I – 5, 6I – 71, 8I – 351} ⊆ PN(Z) be a proper subset of PN(Z). Take S = {0, 1} ⊆ Z. T is clearly a pure neutrosophic integer set vector space over S.

*Example 2.2.6:* Let V = {0, 1 + I, 2 + 2I, 3 + 3I, 42, 80, 4I, 8I, 6 + 6I} ⊆ PN(Z). V is a pure neutrosophic integer set vector space over the set S = {0, 1}. We see V is not a pure neutrosophic integer set vector space over the set T = {0, 2} or any P = {0, n}; n ∈ N \ {1}.

*Example 2.2.7:* Let V = {0, mI, m | m ∈ N} ⊆ N(Z); V is a mixed neutrosophic integer set vector space over the set S = $Z^+$ ∪ {0}; $Z^+$ the set of positive integers. However it is easily verified that V is a mixed neutrosophic integer set vector space over any proper subset of S. But it can also be verified that V is not a mixed neutrosophic integer set vector space over P = $Z^-$ ∪ {0}; $Z^-$ is the set of negative integers or any subset of P.

*Example 2.2.8:* Let V = {0, 3nI | n ∈ Z} ⊆ PN(Q) ∪ {0}. V is a pure neutrosophic integer set vector space over every proper subset of Z.

*Example 2.2.9:* Let W = {0, mI, m | m ∈ Z} ⊆ N(Z). W is a mixed neutrosophic integer set vector space over any subset of Z.

*Example 2.2.10:* Let W = {m, 0, nI | m, n ∈ $2Z^+$} ⊆ N(Z). W is a mixed neutrosophic integer set vector space over every subset of $Z^+$. However W is not a mixed neutrosophic integer set vector space over any subset of $Z^-$ ∪ {0} or on Z.

Now we proceed onto define the notion of mixed neutrosophic integer set vector subspace and pure neutrosophic integer set vector subspace of a mixed neutrosophic integer set vector space.

**DEFINITION 2.2.3:** *Let V be a mixed neutrosophic integer set vector space over a subset S of Z. Let W ≠ φ be a proper subset*

*of V. If W itself is a mixed neutrosophic integer set vector space over S then we call W to be a mixed neutrosophic integer set vector subspace of V over the set S.*

We first illustrate this situation by some examples.

***Example 2.2.11:*** Let $V = \{n + nI \mid n \in Z\} \subseteq N(Z)$ be a mixed neutrosophic integer set vector space over $S = \{Z^+ \cup \{0\}\} \subseteq Z$. Let $W = \{m + mI \mid m \in 5Z\} \subseteq V \subseteq N(Z)$. We see W is a mixed neutrosophic integer set vector space over S. Hence W is a mixed neutrosophic integer set vector subspace of V over S.

***Example 2.2.12:*** Let $V = \{0, 3, 2I, 3I, 5 + 2I, 16, 3 – 3I, 14, 17 + 2I, –15, 15I + 1\} \subseteq N(Z)$ be a mixed neutrosophic integer set vector space over $S = \{0, 1\} \subseteq Z$.

Take $W = \{0, 2I, 16, 17 + 2I\} \subseteq V$; W is a mixed neutrosophic integer set vector subspace of V over S. It is interesting to see that every subset of V which contains 16 or –15 or 3 or 14 is a mixed neutrosophic integer set vector subspace of V over $S = \{0, 1\}$.

***Example 2.2.13:*** Let $V = \{0, 3I + 1, 3 + 5I, 2 – 5I, 17I + 3, 15I + 30, 21, 17I – 153, 15, 412, 317I\} \subseteq N(Z)$; V is a mixed neutrosophic integer set vector space over $S = \{0, 1\}$. Every subset W which contains 15 or 21 or 412 is a mixed neutrosophic integer set vector subspace of V over S. Take $W_1 = \{0, 17I + 3\} \subseteq V$; $W_1$ is a pure neutrosophic integer set vector space contained in V over the set $S = \{0, 1\} \subseteq Z$.

***Example 2.2.14:*** Let $V = \{np, mpI, np + mpI \mid n, m \in Z^+ \text{ and } p = 3\} \subseteq N(Z)$ be a mixed neutrosophic integer set vector space over the set $S = Z^+ \subseteq Z$. Take $W = \{3 + 3I, 6 + 6I, 21 + 21I, 15 + 15I, 27 + 27I, 279 + 279I\} \subseteq V$. W is not a mixed neutrosophic integer set vector subspace over the set $S = Z^+$. We see no finite proper subset of V is a mixed neutrosophic integer set vector subspace of V over the set $S = Z^+$.

***Example 2.2.15:*** Let V = {0, 1, 1 + I} ⊆ N(Z); V is a mixed neutrosophic integer set vector space over S = {0, 1} ⊆ Z. V has no proper subset, hence V has no mixed neutrosophic integer set vector subspace over S = {0, 1} ⊆ Z.

***Example 2.2.16:*** Let V = {m + mI, m, mI | m ∈ $Z^+$} ⊆ N(Z); V is a mixed neutrosophic integer set vector space over the set $3Z^+$ ⊆ Z. W = {n + nI | n ∈ $2Z^+$} ⊆ V ⊆ N (Z); W is not a mixed neutrosophic integer set vector subspace of V over the set $3Z^+$.

It is interesting to note that 0 ∉ V so V cannot have the zero element as the usual set vector space.

In view of this we have the following example, definition and result.

***Example 2.2.17:*** Let V = {9, 9I, 0, 2 + 3I, 4 + 5I, 7, –81, 27, 51I, 91I} ⊆ N(Z). V is a mixed neutrosophic integer set vector space over the set S = {0, 1} ⊆ Z.

Take P = {0, 9, 7, –81, 27} ⊆ V, P is a set vector space over the set S = {0, 1}. We call P to be a pseudo neutrosophic integer set vector subspace of V over S.

**DEFINITION 2.2.4:** *Let V = {$x_1$, ..., $x_n$} ⊆ N(Z) be a mixed neutrosophic integer set vector space over the set S ⊆ Z. Suppose P ⊆ V such that (0) ≠ P ⊆ Z and if P is a set vector space over S then we call P to be a pseudo mixed neutrosophic integer set vector subspace of V over S.*

We will illustrate this by some examples.

***Example 2.2.18:*** Let V = {n + nI, $Z^+$ | n ∈ Z} ⊆ N(Z) be a mixed neutrosophic integer set vector space over $Z^+$ = S. Take P = $3Z^+$ ⊆ V, P is a pseudo neutrosophic integer set vector subspace of V over S.

***Example 2.2.19:*** Let $V = \{nI, mZ^+ \mid n, m \in Z^+\} \subseteq N(Z)$ be a mixed neutrosophic integer set vector space over $S = Z^+$. Take $P = \{mZ^+ \mid m \in Z^+\} \subseteq V$, P is a pseudo neutrosophic integer set vector subspace of V over $Z^+$.

***Example 2.2.20:*** Let $V = \{0, 3I, 4 + 3I, 2 + I, –I, 3 – 8I, 8 + 5I, 7, 250, 49, –560, 2069, 42I + 3\} \subseteq N(Z)$ be a mixed neutrosophic integer set vector space over the set $S = \{0, 1\} \subseteq Z$. Take $P = \{0, 7, 2069, –560\} \subseteq V$, P is a pseudo neutrosophic integer set vector subspace over the set $S = \{0, 1\} \subseteq Z$.

Now we proceed onto define pseudo pure neutrosophic integer set vector subspace.

**DEFINITION 2.2.5:** *Let $W = \{x_1, ..., x_n\} \subseteq N(Z)$ be a mixed neutrosophic integer set vector space over the set $S \subseteq Z$. Let $V \subseteq W$, where $V \subseteq PN(Z)$ (V a subset of W containing only elements from PN(Z)). If V is a pure neutrosophic integer set vector space over the set $S \subseteq Z$, then we define V to be a pseudo pure neutrosophic integer set vector subspace of W.*

We illustrate this by some examples.

***Example 2.2.21:*** Let $V = \{nI, 3Z^+ \mid n \in N\}$ be a mixed neutrosophic integer set vector space over the set $Z^+ \subseteq Z$. Take $P = \{nI \mid n \in N\} \subseteq V$; P is a pure neutrosophic integer set vector space over the set $Z^+$. P is clearly a pseudo pure neutrosophic integer set vector subspace of V over the set $Z^+$.

***Example 2.2.22:*** Let $W = \{3I, 0, 46 + 7I, 91I + 27, 5I, 7, 982, 47I, 61, -257, 96 + 2I\} \subseteq N(Z)$ be a mixed neutrosophic integer set vector space over the set $S = \{0, 1\} \subseteq Z$. Take $V = \{0, 7, 982, 61, -257\} \subseteq W$. V is a pseudo neutrosophic integer set vector subspace of W over S.

***Example 2.2.23:*** Let $V = \{mI, mZ^+ \mid m \in N\} \subseteq N(Z)$ be a mixed neutrosophic integer set vector space defined over the set $P = Z^+ \subseteq Z$. Take $P = \{mI \mid m \in Z^+\} \subseteq V$. P is a pseudo pure neutrosophic integer set vector over the set $Z^+ = P$.

Now these concepts cannot be even imitated in case of pure neutrosophic integer set vector space V as it cannot contain a proper subset which is mixed(i.e., integers). Now we proceed on to define yet another substructure of both pure and mixed neutrosophic integer set vector spaces.

**DEFINITION 2.2.6:** *Let $V = \{x_1, \ldots, x_n\}$ be a mixed neutrosophic set vector space over the set $S \subseteq Z$. Let $W \subseteq V$ be a proper subset of V, if there exists a proper subset $T \subset S$ such that W is a mixed neutrosophic integer set vector space over T then we call W to be a mixed neutrosophic integer subset vector subspace of V over the subset T of S. If $W \subseteq V$ is such that W is a proper subset of PN(Z) and W is a pseudo mixed neutrosophic vector space over T, then we call W to be a pseudo mixed neutrosophic integer subset vector subspace of V over the subset T of S.*

*If $W \subseteq V$ is such that W is a subset of Z and S has a proper subset of integers say $B \subset Z$ (for $B \subseteq S$) and if W is a set vector space over B then we call W to be a pseudo set integer subset vector subspace of V over the subset B of S.*

We will illustrate this by the following examples.

***Example 2.2.24:*** $V = \{0, \pm I, \pm 4I, \pm(3 + 2I), \pm(-5, 7I), \pm(8 - 25I), \pm 20, \pm 246, \pm 284I, \pm 261, \pm 85, \pm 98\}$ be a mixed neutrosophic integer set vector space over the set $S = \{0, -1, 1\} \subseteq Z$. Consider $W = \{0, 4I, 3 + 2I, 246, 20\} \subseteq V$. Let $T = \{0, 1\} \subseteq S \subseteq Z$.

W is a mixed neutrosophic integer subset vector subspace of V over the subset T of S. Take $X = \{\pm I, \pm 4I, \pm(3 + 2I), \pm(8 - 25I)\} \subseteq V$ and the subset $A = \{-1, 1\} \subseteq \{0, 1, -1\}$, X is a pseudo pure neutrosophic integer subset vector subspace of V over the subset X of S. Consider $Y = \{0, \pm 261, \pm 85, \pm 20, \pm 98, \pm 246\} \subseteq V$. Take $D = \{-1, 1\} \subseteq \{0, -1, + 1\}$, Y is a pseudo set integer subset vector subspace of V over the set of integers $\{-1, 1\}$, a subset of S.

***Example 2.2.25:*** Let $V = \{3mI, n(2I + 2), Z^+ \mid m, n \in Z^+\}$ be a mixed neutrosophic integer set vector space over the set $S = Z^+ \subseteq Z$.

Take $W = \{n(2I+ 2), Z^+ \mid n \in Z^+\}$; W is a mixed neutrosophic integer subset set vector subspace over $3Z^+ \subseteq Z^+$. Take $X = \{5Z^+\} \subseteq V$, X is a pseudo set integer subset vector subspace of V over the subset $2Z^+ \subseteq Z^+$. Consider $Y = \{3mI \mid m \in Z^+\} \subseteq V$; V is a pseudo pure neutrosophic integer subset of the vector subspace V over the subset $6Z^+$ of $Z^+$.

Now having seen several new substructures of a mixed neutrosophic integer set vector space we now proceed on to define the notion of mixed neutrosophic integer set linear algebra over a subset of integers.

**DEFINITION 2.2.7:** *Let $V = \{x_1, …, x_n\}$ be a mixed neutrosophic integer set vector space over the set $S \subseteq Z^+$. If on V we can define a closed binary operation '+' such that for all $x_i, x_j \in V$, $x_i + x_j \in V$ then we define V to be a mixed neutrosophic integer set linear algebra over S.*

We illustrate this situation by some examples.

***Example 2.2.26:*** Let $V = \{3mI, 0, 3m, 3m + 3mI \mid m \in Z^+\}$, V is a mixed neutrosophic integer set linear algebra over the set $Z^+ \subseteq Z$.

***Example 2.2.27:*** Let
$$V = \left\{ \begin{pmatrix} mI & 0 \\ 0 & mI \end{pmatrix}, \begin{pmatrix} 0 & 0 \\ 0 & 0 \end{pmatrix} \middle| m \in Z^+ \right\}$$
be a mixed neutrosophic integer set linear algebra over the set $Z^+ \subseteq Z$.

***Example 2.2.28:*** Let

$$V = \left\{ \begin{pmatrix} mI & mI & mI \\ 0 & 0 & 0 \end{pmatrix}, \begin{pmatrix} 0 & 0 & 0 \\ m & m & m \end{pmatrix}, \begin{pmatrix} mI & mI & mI \\ m & m & m \end{pmatrix} \middle| m \in Z^+ \right\}.$$

V is a mixed neutrosophic integer set linear algebra over the set $3Z^+ \subseteq Z$.

***Example 2.2.29:*** Let $V = \{0, \pm nI, \pm n \mid n \in Z\}$, V is not a mixed neutrosophic integer set linear algebra over any subset of Z. It is only a mixed neutrosophic set vector space over every proper subset of Z and including Z.

In view of this we have the following theorem.

**THEOREM 2.2.1:** *Let V be a mixed neutrosophic integer set linear algebra over a subset S of Z. Then V is a mixed neutrosophic set vector space over the subset S of Z.*

*Proof:* Clear from the definition of mixed neutrosophic integer set linear algebra over S.

**COROLLARY 2.2.1:** *A mixed neutrosophic integer set vector space V over a set S ($S \subseteq Z$) in general is not a mixed neutrosophic integer set linear algebra over S.*

*Proof:* The proof is given using an example.
  Consider $V = \{2, 0, 2I, 4I, 3I + 1, 2I - 27, 28, 41I - 38, 1\} \subseteq N(Z)$. V is a mixed neutrosophic integer set vector space over the set $S = \{0, 1\} \subseteq Z$. Clearly V is not closed under addition for $2 + 2I \notin V$, $2I - 27 + 28 \notin V$ and so on.
  Thus a mixed neutrosophic integer set vector space in general is not a mixed neutrosophic integer set linear algebra over S. Hence the claim.

Now we proceed onto define substructures in mixed neutrosophic integer set linear algebra.

**DEFINITION 2.2.8:** *Let V be a mixed neutrosophic integer set linear algebra over the set S ($S \subseteq Z$). Let W be a subset of V such that W itself is a mixed neutrosophic integer set linear* *algebra over S; then we call W to be a mixed neutrosophic integer set linear subalgebra of V over the set $S \subseteq Z$.*

We now illustrate this situation by some examples.

***Example 2.2.30:*** Let $V = \{5nI, 0, 5n, 5nI + 5n \mid n \in N\}$ be a mixed neutrosophic integer set linear algebra over the set $S = \{0, 1\} \subseteq Z$. Take $W = \{0, 25n, 25nI, 25n + 25nI \mid n \in N\} \subseteq V$; W is clearly a mixed neutrosophic integer set linear subalgebra of V over the set $S = \{0, 1\}$.

***Example 2.2.31:*** Let $V = \{0, m + mI \mid m \in Z^+\}$ be a mixed neutrosophic integer set linear algebra of V over $S = mZ^+ \subseteq Z$. Take $W = \{0, 3m + 3mI \mid m \in Z^+\} \subseteq V$; clearly W is a mixed neutrosophic integer set linear subalgebra of V over $S = mZ^+$.

Now we proceed onto define the notion of yet another new substructure.

**DEFINITION 2.2.9:** *Let V be a mixed neutrosophic integer set linear algebra over the subset $S \subseteq Z$. Let $W \subseteq V$ be such that $W \subseteq Z$. If W itself is a set linear algebra over the set S then we call W to be a pseudo integer set linear subalgebra over S; $S \subseteq Z$.*

We shall illustrate this by some examples.

***Example 2.2.32:*** Let $V = \{3nI + 3n, 3nI, 3n \mid n \in Z^+\} \subseteq N(Z)$ be a mixed neutrosophic integer set linear algebra over the set $S = Z^+ \subseteq Z$. $P = \{3n \mid n \in Z^+\} \subseteq V$ is a pseudo integer set linear subalgebra of V over $S \subseteq Z$.

***Example 2.2.33:*** Let $V = \{2Z + 2ZI, 2Z, 2ZI\} \subseteq N(Z)$ be a mixed neutrosophic integer set linear algebra over the set $S = Z^+$. Take $W = \{2Z^+\} \subseteq V$; W is a pseudo integer set linear subalgebra of V over $S = Z^+$.

***Example 2.2.34:*** Let $V = \{2ZI, 2Z, 2Z + 2ZI\} \subseteq N(Z)$ be a mixed neutrosophic integer set linear algebra over the set $S = $

$7Z^+ \subseteq Z$. Take $P = \{16Z\} \subseteq V$; P is a pseudo integer set linear subalgebra of V over $S = 7Z^+ \subseteq Z$.

**DEFINITION 2.2.10:** *Let V be a pure neutrosophic integer set vector space over the set $S \subseteq Z$. If in addition for every x, y $\in$ V; x + y $\in$ V we call V to be a pure neutrosophic integer set linear algebra over the set $S \subseteq Z$.*

We illustrate this situation by some examples.

*Example 2.2.35:* Let $V = \{2ZI\} \subseteq PN(Z)$ be pure neutrosophic integer set linear algebra over the set $S = Z^+ \subseteq Z$.

*Example 2.2.36:* Let $V = \{3ZI\} \subseteq PN(Z)$. V is a pure neutrosophic integer set linear algebra over the set $3Z^+ \subseteq Z$.

*Example 2.2.37:* Let
$$V = \left\{ \begin{pmatrix} 2ZI & 4Z \\ 4Z & 8ZI \end{pmatrix} \right\}.$$
V is a pure neutrosophic integer set linear algebra over $2Z^+ \subseteq Z$.

*Example 2.2.38:* Let $V = \{mZ + mZI \mid m \in Z\} \subseteq PN(Z)$. V is a pure neutrosophic integer set linear algebra over $5Z^+ \subseteq Z$.

Now we proceed onto describe the substructure of pure neutrosophic integer set linear algebra.

**DEFINITION 2.2.11:** *Let V be a pure neutrosophic integer set linear algebra over a set $S \subset Z$. Let $W \subseteq V$ be a proper subset of V; if W is pure neutrosophic integer set linear algebra over S; then we call W to be a pure neutrosophic integer set linear subalgebra of V over the set S.*

We shall illustrate this situation by some examples.

*Example 2.2.39:* Let $V = \{3Z + 3ZI\} \subseteq PN(Z)$ be a pure neutrosophic integer set linear algebra over the set $S = Z^+ \subseteq Z$.

Take $W = \{6Z + 6ZI\} \subseteq V \subseteq PN(Z)$; W is a pure neutrosophic integer set linear subalgebra of V over the set S.

*Example 2.2.40:* Let

$$V = \left\{ \begin{pmatrix} 2Z & 3ZI \\ 3ZI & 2Z \end{pmatrix} \right\}$$

be a pure neutrosophic integer set linear set linear algebra over the set $S = Z^+$.
Take

$$W = \left\{ \begin{pmatrix} 8Z & 12ZI \\ 12ZI & 8Z \end{pmatrix} \right\}$$

W is a pure neutrosophic integer set linear subalgebra of V over the set $S = Z^+$.

*Example 2.2.41:* Let

$$V = \left\{ \begin{pmatrix} 3Z & 3Z & 3Z \\ 4ZI & 4ZI & 4ZI \end{pmatrix} \right\} \subseteq PN(Z),$$

be a pure neutrosophic integer set linear algebra over the set $S = 3Z$.
Take

$$W = \left\{ \begin{pmatrix} 6Z & 6Z & 6Z \\ 8ZI & 8ZI & 8ZI \end{pmatrix} \right\} \subseteq V$$

to be a proper subset of V. W is a pure neutrosophic integer set linear subalgebra of V over 3Z.

Now we proceed onto define yet another type of substructure in pure neutrosophic integer set linear algebra.

**DEFINITION 2.2.12:** *Let V be a pure neutrosophic integer set linear algebra over the set $S \subseteq Z$. Let W be a proper subset of V* and T ⊆ S be a proper subset of S. If W is a pure neutrosophic integer set linear algebra over the set T, (T ⊆ S) then we call W to be a pure neutrosophic integer subset linear subalgebra of V over the subset T of S.

We illustrate this definition by some examples.

**Example 2.2.42:** Let V = {3Z + 3ZI} ⊆ PN(Z) ∪ {0} be a pure neutrosophic integer set linear algebra over S = $Z^+$. Take W = {27Z + 27ZI} ⊆ V and T = $3Z^+$ ⊆ $Z^+$ = S. W is a pure neutrosophic integer subset linear subalgebra of V over the subset T of S.

**Example 2.2.43:** Let

$$V = \left\{ \begin{pmatrix} 2ZI & Z \\ Z & 4ZI \end{pmatrix} \right\}$$

V is a pure neutrosophic integer set linear algebra over S = Z. Take

$$W = \left\{ \begin{pmatrix} 10ZI & 5Z \\ 5Z & 20ZI \end{pmatrix} \right\} \subseteq V,$$

W is a pure neutrosophic integer subset linear subalgebra of V over the subset T = $Z^+$ ⊆ Z = S.

It is pertinent to mention here that pure neutrosophic integer set linear algebras do not have proper pseudo neutrosophic integer substructures. It is also important to mention here that every pure neutrosophic integer set linear algebra is a pure neutrosophic integer set linear algebra.
We prove the later part of the claim by the following example:

**Example 2.2.44:** Let V = {2I + 2, 0, 5 + 5I, 7I, -28I, 19-21I} ⊆ PN (Z) be a pure neutrosophic integer set vector space over the set S = {0, 1} ⊆ Z. We see V is not a pure neutrosophic integer set linear algebra as 5 + 5I + 7I = 5 + 12I ∉ V, –28I + 7I = –21I

∉ V and so V is not an integer set linear algebra of S = {0, 1} ⊆ Z.

We now proceed onto define two more new concepts on the integer set Z.

**DEFINITION 2.2.13:** *Let $V \subseteq NP(Z)$ be a pure neutrosophic integer set linear algebra over a set $S \subseteq Z$. If V has no pure neutrosophic integer subset linear subalgebra then we call V to be a pure neutrosophic integer set simple linear algebra.*

We first illustrate this by some examples.

*Example 2.2.45:* Let V = {nI | n ∈ Z \ {0}} be the pure neutrosophic integer set linear algebra over S = {0, 1}. V is not a pure neutrosophic integer set simple linear algebra.

*Example 2.2.46:* Let
$$V = \left\{ \begin{pmatrix} nI & mI \\ mI & nI \end{pmatrix} \middle| m, n \in Z^+ \right\}.$$

V is a pure neutrosophic integer set simple linear algebra over S = {0, 1}.

**DEFINITION 2.2.14:** *Let $V \subseteq NP(Z)$ be a pure neutrosophic integer set linear algebra over the set $S \subseteq Z$. If V has no proper subset $W \subseteq V \subseteq NP(Z)$ such that W is a pure neutrosophic integer set linear subalgebra or W is not a pure neutrosophic integer subset linear subalgebra for any subset $T \subseteq S \subseteq Z$ over any subset $T \subseteq S \subseteq Z$; then we call V to be pure neutrosophic integer set weakly simple linear algebra.*

We illustrate this by some simple examples.

*Example 2.2.47:* Let
$$V = \left\{ \begin{pmatrix} pI & 0 \\ 0 & pI \end{pmatrix} \middle| p \in Z^+ \right\}$$

be a pure neutrosophic integer set linear algebra over $\{0, 1\} = S \subseteq Z$. Since S has no proper subsets, V is a pure neutrosophic integer set weakly simple linear algebra over $S = \{0, 1\} \subseteq Z$.

*Example 2.2.48:* Let

$$V = \left\{ \begin{pmatrix} nI & 0 & tI \\ mI & qI & pI \end{pmatrix} \middle| n, t, m, p, q \in Z^+ \right\}$$

be a pure neutrosophic integer set linear algebra over the set $S = \{0, 1\} \subseteq Z$. It is easily verified V is a pure neutrosophic integer set weakly simple linear algebra.

Now we proceed onto define the notion of set neutrosophic integer set linear transformation and set neutrosophic integer set linear operator.

**DEFINITION 2.2.15:** *Let V and W be any two mixed neutrosophic integer set vector spaces over the same set $S \subseteq Z$. Let T be a map from V into W satisfying the following conditions:*
  *(1) T(I) = I*
  *(2) T(sυ) = sT(υ) for all $s \in S$ and for all $υ \in V$ and $T(υ) \in W$.*

*We define T to be a set neutrosophic integer linear transformation of V into W. If V = W then we call the set neutrosophic integer set linear transformation to be the set neutrosophic integer set linear operator on V.*

We illustrate this by simple examples.

*Example 2.2.49:* Let V = {8I, 0, 5I, 22I, 46, 3 + 25I} and W = {0, 46 + I, 8I 22I + 3, 7I, 21 5I 25I, 63I} be mixed neutrosophic integer set vector space over the set $S = \{0, 1\} \subseteq Z$. Let T: V → W be a map such that
$$T(I) = I$$

$$T(8I) = 8I$$
$$T(0) = 0$$
$$T(5I) = 5I$$
$$T(22I) = 22I + 3$$
$$T(46) = 21$$
$$T(3 + 25I) = 25I$$
$$T(I) = I.$$

Clearly T(0.x) = 0.T(x) i.e.,
$$T(0) = 0 \text{ and } T(1x) = 1\ T(x) = T(x);$$
$$\text{as } 0.x = 0 \text{ and } 1.x = x \text{ for all } x \in V.$$

*Example 2.2.50:* Let V = {0, 2nI + 2n | n ∈ Z$^+$} be a pure neutrosophic integer set vector space over Z$^+$. Define the map T : V → V as
$$T(0) = 0$$
$$T(I) = I$$
$$T(2nI + 2n) = 2(n + 2)I + 2(n + 2),$$
i.e., T(x) = 2x for all x ∈ V. Clearly T is a linear operator on V.

*Note:* It is interesting and important to note that V and W can be both mixed neutrosophic integer set vector spaces or both can be pure mixed neutrosophic integer set vector spaces or one can be mixed neutrosophic integer set vector space and other can be pure neutrosophic integer set vector space; still the definition of the set neutrosophic integer set linear transformation remains the same.

The only main criteria is that T(I) = I for any set neutrosophic integer set linear transformation T from V to W except in case of the zero linear transformation 0(I) = 0; but however this special transformation is of no use to the real world problems or applications.

Now we proceed onto give the definition of set neutrosophic integer set linear transformation of mixed neutrosophic integer set linear algebras and pure neutrosophic integer set linear algebras.

**DEFINITION 2.2.16:** *Let V and W be any two mixed neutrosophic integer set linear algebras defined over the same*

*set S. Let T be a set neutrosophic integer set linear transformation from V to W such that T: V → W is a semigroup homomorphism with respect to addition then we define T to be a set neutrosophic integer set linear algebra transformation from V to W. If V = W then we define T to be a set neutrosophic integer set linear algebra operator from V to W.*

We illustrate this by some examples.

***Example 2.2.51:*** Let V = {2nI, 2n, 2m + 2tI | m, t, n ∈ $Z^+$ {m, t and n need not be taking always distinct values}} and

$$W = \left\{ \begin{pmatrix} 2nI & 0 \\ 0 & 0 \end{pmatrix}, \begin{pmatrix} 0 & 0 \\ 0 & 2n \end{pmatrix}, \begin{pmatrix} 2mI & 0 \\ 0 & 2t \end{pmatrix} \middle| n, m, t \in Z^+ \right\}$$

be two mixed neutrosophic integer set linear algebras over S = {0,1}. Let T : V → W be defined as

$$T(2nI) = \begin{pmatrix} 2nI & 0 \\ 0 & 0 \end{pmatrix}$$

$$T(2n) = \begin{pmatrix} 0 & 0 \\ 0 & 2n \end{pmatrix}$$

$$T(2m + 2tI) = \begin{pmatrix} 2tI & 0 \\ 0 & 2m \end{pmatrix}$$

for 2nI, 2n, 2m + 2tI ∈ V we see T is a set neutrosophic integer set linear algebra transformation of V to W.

***Example 2.2.52:*** Let

$$V = \left\{ \begin{pmatrix} 2n & 0 \\ 0 & 2n \end{pmatrix}, \begin{pmatrix} 0 & 2mI \\ 2mI & 0 \end{pmatrix}, \begin{pmatrix} 2n & 2mI \\ 2mI & 2n \end{pmatrix} \middle| m, n \in Z^+ \right\}$$

be a mixed neutrosophic integer set linear algebra over the set S = {0, 1}. Let T: V → V be defined by

$$T\begin{pmatrix} 2n & 0 \\ 0 & 2n \end{pmatrix} = \begin{pmatrix} 2n+2 & 0 \\ 0 & 2n+2 \end{pmatrix}$$

$$T\begin{pmatrix} 0 & 2mI \\ 2mI & 0 \end{pmatrix} = \begin{pmatrix} 0 & 2mI+2I \\ 2mI+2I & 0 \end{pmatrix}.$$

Clearly T is a set neutrosophic integer set linear algebra operator over the set S = {0,1}.

*Example 2.2.53:* Let V = {(2n + 2nI) | n ∈ $Z^+$} and

$$W = \left\{ \begin{pmatrix} 2n+2nI & 0 \\ 0 & 2n+2nI \end{pmatrix} \middle| n \in Z^+ \right\}$$

be two pure neutrosophic integer set linear algebras over the set S = {0,1}. Define T: V → W by

$$T(2n + 2nI) = \begin{pmatrix} 2n+2nI & 0 \\ 0 & 2n+2nI \end{pmatrix};$$

T is a set neutrosophic integer set linear algebra transformation of V into W.

*Example 2.2.54:* Let

$$V = \left\{ \begin{pmatrix} 2n & 2nI \\ 2nI & 2n \end{pmatrix} \middle| n \in Z^+ \right\}$$

be a pure neutrosophic integer set linear algebra over the set S = {0,1} ⊆ Z. Let T : V → V be defined by

$$T = \begin{pmatrix} 2n & 2nI \\ 2nI & 2n \end{pmatrix} = \begin{pmatrix} 2n+2 & 2nI+2 \\ 2nI+2 & 2n+2 \end{pmatrix}.$$

T is easily verified to be a set neutrosophic integer set linear operator on V.

It is left as a problem for the reader to prove $\text{Hom}_S(V, W)$ is a set of all neutrosophic integer set vector spaces over $S \subseteq Z$ where V and W are mixed (pure) neutrosophic integer set vector spaces over the set $S \subseteq Z$. Let $\text{Hom}_S(V, V)$ denote the collection of all set neutrosophic integer set linear operators of V over the set $S \subseteq Z$; where V is the pure (mixed) neutrosophic set vector space over S. What is the structure of $\text{Hom}_S(V, V)$?

If V and W are mixed (pure) neutrosophic integer set linear algebras defined over S, will $\text{Hom}_S(V, W)$ be a mixed (pure) neutrosophic integer set linear algebras?

Let us define some more properties of set neutrosophic integer set linear algebra.

**DEFINITION 2.2.17:** *Let T be a set neutrosophic integer set linear transformation from V to W. If atleast one vector subspace P of V is mapped into a vector subspace of W then we say T weakly preserves subspaces i.e., T(P) is a vector subspace of W for atleast one vector subspace P of V; we define T to be a set neutrosophic integer set weak subspace preserving linear transformation of V to W.*

*If every subspace P of V is preserved by a set neutrosophic integer set linear transformation then we call T to a set neutrosophic integer set strong subspace preserving linear transformation of V to W.*

*In an analogous way one can define these two concepts for set neutrosophic integer set linear operator on V.*

We illustrate these definitions by some examples.

***Example 2.2.55:*** Let V = {0, 3I + 2, 7I + 4, –2 + I, 80I, 92 – 8I, –47, –6I, –4, 50 – 2I} and W = {0, 6I + 4, –35 14I + 8, –4 + 2I,

–12I, –8, 25 – 9I, – 48 + I, 97– 4I, –52I, –2I, 40I} be two mixed neutrosophic integer set vector spaces over the set S = {0,1} ⊆ Z. Let T: V → W given by

$$T(0) = 0$$
$$T(3I + 2) = 6I + 4$$
$$T(7I + 4) = 14I + 8$$
$$T(-2 + I) = -4 + 2I$$
$$T(81I) = 40I$$
$$T(92 - 8I) = 97 - 4I$$
$$T(-47) = -35$$
$$T(-6I) = -12I$$
$$T(-4) = -8$$
$$T(50 - 2I) = -2I.$$

We see P = {0, 3I + 2, 7I + 4, –2 + I, –6I, –4} is a mixed neutrosophic integer set vector subspace of V. Also we see T(P) = {0, 6I + 4, 14I + 8, – 4 + 2I, –12I, – 8} ⊆ W is a mixed neutrosophic integer set vector subspace of W. So T is a set neutrosophic integer set weak linear transformation of V into W.

*Example 2.2.56:* Let V = {20, I, 0, 10I, 26 + I} ⊆ N(Z) be a mixed neutrosophic integer set vector space over S = {0, 1} ⊆ Z. The mixed neutrosophic subspaces of V are

$$P_1 = \{0, I, 20\},$$
$$P_2 = \{0, I, 20, 10I\}$$
$$P_3 = \{0, I, 20\ 26 + I\}$$
$$P_4 = \{0, 10I, 20\}$$
$$P_5 = \{0, 10I, 20\ 26 + I\}$$
$$P_6 = \{0, 26 + I, 20\}$$

Define a set neutrosophic integer set linear operator T on such that T(0) = 0 and T(20) = 0 and others in any compatible way then T is not a set neutrosophic integer set weak linear operator on V or T is not a set neutrosophic integer set strong linear operator on V. Now define $T_1 : V \to V$ as follows.

$$T_1(0) = (0)$$
$$T_1(20) = 20$$

$$T_1(I) = I$$
$$T(10I) = 26 + I$$

$T_1$ is a set neutrosophic integer set weak linear operator on V.

We shall discuss about the generating subset of a mixed (pure) neutrosophic integer set vector spaces.

**DEFINITION 2.2.18:** *Let $V = \{x_1, …, x_n\}$ be a mixed (pure) neutrosophic integer set vector space over the set $S \subseteq Z$. Suppose $T = \{x_1, …, x_m \mid m \leq n\} \subseteq V$ is such every $x_i \in V$ can be represented as $sx_j$ for some $s \in S$ and $x_j \in T$ i.e., $x_i = sx_j$, then we say T generates V over S and T is called the mixed (pure) neutrosophic integer set generator subset of V over S.*

*Note:* It may at times so happen that $T = V$.

We shall illustrate this situation by some examples.

*Example 2.2.57:* Let $V = \{0, 2I, 24, 2 + 3I, 41 - I, 37 + 44I\} \subseteq N(Z)$ be a mixed neutrosophic integer set vector space over $S = \{0, 1\} \subseteq Z$. $T = \{2I, 24, 2 + 3I, 41 - I, 37 + 44I\} \subseteq V$ is the mixed neutrosophic integer set generator subset of V over the set $S = \{0, 1\}$. We see $V \neq T$.

*Example 2.2.58:* Let $V = \{\pm 3I, \pm 1, \pm (22 + I), \pm (5I - 20), \pm 70I, \pm (8I + 4)\} \subseteq N(Z)$ be a mixed neutrosophic integer set vector space over the set $S = \{-1, 1\}$. We see $T = \{3I, 1, 22 + I, 5I - 20, 70I, 8I + 4\} \subseteq V$ is a mixed neutrosophic integer set generator subset of V over the set $S = \{-1, 1\}$. We see $V \neq T$. In fact $|V| = 12$ and $|T| = 6$.

*Example 2.2.59:* Let $V = \{3ZI, 8Z\} \subseteq N(Z)$ be a mixed neutrosophic integer set vector space over Z. The mixed neutrosophic integer set generator subset of V over the set Z is given by $T = \{3I, 8\} \subseteq V$. We see $|T| = 2$ where as $|V| = \infty$.

*Example 2.2.60:* Let $V = \{2nI \mid n \in Z\}$ be a pure neutrosophic integer set vector space over $S = \{0, 1\} \subseteq Z$. $T = \{2nI \mid n \in Z \setminus$

{0}} ⊆ V is the pure neutrosophic integer set generator of V over S = {0, 1}.

*Example 2.2.61:* Let V = {2I, 9I, –8I, 14I, 27I + 4, 44 – 2I, 0, 14I – 9} ⊆ PN(Z) be a pure neutrosophic integer set vector space with zero over S = {0,1} ⊆ Z. T = {2I, 9I, –8I, 14I, 27I + 4, 44 – 2I, 14I – 9} ⊆ V i.e., V \ {0} = T is the pure neutrosophic integer set generator of V over S.

*Example 2.2.62:* Let V = {3ZI} be the pure neutrosophic integer set vector space over the set S = Z. T = {3} ⊆ V is the pure neutrosophic integer set generator of V over S = Z. Thus | T | = 1.

*Remark:* We see the pure (mixed) neutrosophic integer set generator of a vector space may be finite or infinite. The cardinality in some case depends on the set over which they are defined.

This is proved by the following examples:

*Example 2.2.63:* Let V = {5ZI} be a pure neutrosophic integer set vector space over S = {0, 1} ⊆ Z. T = 5ZI \ {0} is the pure neutrosophic integer set generator of V over S. Clearly | T | = ∞ infact | T | = | V \ {0}|.

Now we consider the same pure neutrosophic integer set vector space V over a different set S ⊆ Z and find the cardinality of the pure neutrosophic integer set generator of V.

*Example 2.2.64:* Let V = {5ZI} be a pure neutrosophic integer set vector space over the set S = Z. Now T = {5I} is the pure neutrosophic integer generator of V over S = Z. Clearly | T | = 1. So we see depending on the set S over which V is defined the cardinality may be one or ∞.

Now we define the pure (mixed) neutrosophic integer set generator of a pure (mixed) neutrosophic integer set linear algebra over the set S ⊆ Z.

Let $T \subseteq V$; if every $v \in V$ can be represented as $v = st$ or

$$v = \sum_{i} s_i t_i,$$

for some $s_i$, $s \in S$ and $t_i$, $t \in T$ then we define T to be a pure (mixed) neutrosophic integer set generator of V over S.

We shall now illustrate this situation by some examples.

***Example 2.2.65:*** Let $V = \{n + nI \mid n \in Z^+ \cup \{0\}\} \subseteq PN(Z)$ be the pure neutrosophic integer set linear algebra over $S = \{0, 1\}$. $T = \{1 + 1I\}$ is the pure neutrosophic integer set generator of V over S. Thus $|T| = 1$.

***Example 2.2.66:*** Let $V = \{n, nI, m + tI \mid m, n, t \in Z\} \subseteq N(Z)$ be a mixed neutrosophic integer set linear algebra over $Z^+ \subseteq Z$. Take $T = \{\pm 1, \pm I, 0\} \subseteq V$; T is a mixed neutrosophic integer set generator of V over $Z^+$. $|T| = 5$.

***Example 2.2.67:*** Let

$$V = \left\{ \begin{pmatrix} nI & 0 \\ p & mI \end{pmatrix} \middle| n, m, p \in Z \right\}$$

be a pure neutrosophic integer set linear algebra over Z. Take

$$T = \left\{ \begin{pmatrix} I & 0 \\ 0 & 0 \end{pmatrix}, \begin{pmatrix} 0 & 0 \\ 0 & I \end{pmatrix}, \begin{pmatrix} 0 & 0 \\ 1 & 0 \end{pmatrix} \right\} \subseteq V;$$

T is the pure neutrosophic integer set generator of the linear algebra over Z.

If we change the set over which these spaces are defined then their generating set are also different. This is described by the following examples:

***Example 2.2.68:*** Let $V = \{n + nI \mid n \in Z^+ \cup \{0\}\}$ be the pure neutrosophic integer set linear algebra over $S = \{0, 1\}$. $T = \{1 + I\} \subseteq V$ is the pure neutrosophic integer set generator of $V$ over $S = \{0, 1\}$.

If we replace $S = \{0, 1\}$ by $M = Z^+ \cup \{0\}$. We see $T = \{1 + I\} \subseteq V$ is the pure neutrosophic integer set generator of $V$ over $S = \{0, 1\}$.

Suppose $V = \{n + nI \mid n \in Z^+ \cup \{0\}\}$ is a pure neutrosophic integer set linear algebra over any subset $S \subseteq Z^+ \cup \{0\} \subseteq Z$ have same $T = \{1 + I\}$ to be the pure neutrosophic integer set generator of $V$ over $S$.

***Example 2.2.69:*** Let $V = \{2Z^+, 0, mZ + nZI \mid m, n \in Z^+\} \subseteq N(Z)$ be a mixed neutrosophic integer set linear algebra over $Z^+ \cup \{0\}$; $T = \{2, m + nI \mid m, n \in Z^+\} \subseteq V$; $T$ is the mixed neutrosophic integer set generator of $V$ over $Z^+ \cup \{0\}$.

We just wish to show that a $V$ treated as a mixed (pure) neutrosophic integer set vector space in general has a distinct generator set from the same $V$ treated as a mixed (pure) neutrosophic integer set linear algebra.

The following example shows the above claim.

***Example 2.2.70:*** Take $V = \{m + mI \mid m \in Z^+ \cup \{0\}\} \subseteq PN(Z)$ a pure neutrosophic integer set vector space over $S = \{0, 1\} \subseteq Z$. $T = V \setminus \{0\} \subseteq V$ is the pure neutrosophic integer set generator of $V$.

Clearly $|T| = \infty$. Now $V = \{m + mI \mid m \in Z^+ \cup \{0\}\}$ is a pure neutrosophic integer set linear algebra over $S = \{0, 1\}$. $T = \{1 + I\}$ is the pure neutrosophic set generator of $V$. We see $|T| = 1$.

From this example the reader can understand the vast difference between the pure (mixed) neutrosophic integer set linear algebra $V$ and pure (mixed) neutrosophic integer set vector space $V$ (same $V$) over the same set $S$.

## 2.3 Neutrosophic-Neutrosophic Integer Set Vector Spaces

In this section we introduce yet another new type of vector spaces called neutrosophic-neutrosophic integer set vector spaces. Here we cannot have two types of vector spaces. Recall $N(Z) = \{a + bI \mid a, b \in Z\}$ is the set of neutrosophic integers.

**DEFINITION 2.3.1:** *Let $V = \{v_1, \ldots, v_n\}$ where $v_i \in N(Z)$; $1 \leq i \leq n$. We say V is a neutrosophic-neutrosophic integer set vector space over $S \subseteq N(Z)$ ($S \nsubseteq Z$) if $s v_i = v_i s \in V$ for every $v_i \in V$ and $s \in S$. We shall for easy representation write neutrosophic-neutrosophic integer vector space as n-n integer set vector space.*

We now illustrate this new structure by some examples.

***Example 2.3.1:*** Let $V = \{0, 1 + (2^n - 1)I \mid n = 1, 2, \ldots, \infty\} \subseteq N(Z)$. V is a n-n integer set vector space over $S = \{0, 1 + I, 1\} \subseteq N(Z)$.

***Example 2.3.2:*** Let $V = \{I, 2I, 5I, 7I, 0, 8I, 27I\} \subseteq N(Z)$. V is a n-n integer set vector space over $S = \{0, I\} \subseteq N(Z)$.

***Example 2.3.3:*** Let $V = ZI \subseteq N(Z)$, V is a n-n integer set vector space over $Z = \{0, I\}$.

***Example 2.3.4:*** Let $V = \{ZI\} \subseteq N(Z)$; V is a n-n integer set vector space over $ZI \subseteq N(Z)$.

***Example 2.3.5:*** Let $V = \{0, 1 - I\} \subseteq PN(Z) \cup \{0\}$ be the n-n integer set vector space over the set $S = \{0, 1, 1+2I\} \subseteq N(Z)$.

*Note:* It is important and interesting to note that if V is a n-n integer set vector space then V cannot contain any integer from Z; i.e., if any $a \in V$ then $a \notin Z$.

**DEFINITION 2.3.2:** *Let V be a n-n set vector space over a set $S \subseteq N(Z)$. Suppose W is a proper subset of V and W is itself a n-n* *set vector space over the same set $S \subseteq N(Z)$; then we call W to be a n-n set vector subspace of V over S.*

We will illustrate this situation by some examples.

*Example 2.3.6:* Let V = {ZI} be a n-n set vector space over the set $S = Z^+I$. Take W = {2ZI} $\subseteq$ V; V is a n-n set vector subspace of V over $Z^+I$.

*Example 2.3.7:* Let V = {0, m ± mI | m ∈ $Z^+$} be the n-n set vector space over the set S = {0, 1, 1 – I} $\subseteq$ N(Z). Take W = {2m ± 2mI | m ∈ $Z^+$} be the n-n set vector subspace of V over $S \subseteq N(Z)$.

*Example 2.3.8:* Let V = {ZI, m ± mI | m ∈ $Z^+$} be a n-n set vector space over the set S = {0, 1, 1 – I} $\subseteq$ N(Z). Take W = {ZI} $\subseteq$ V; W is a n-n set subvector space of V over S.

Now we will proceed onto define the notion to neutrosophic-neutrosophic set linear algebra (n-n set linear algebra).

**DEFINITION 2.3.3:** *Let V be a n-n set vector space over $S \subseteq N(Z)$, a subset of N(Z). If V is such that for every a, b ∈ V, a + b, b + a ∈ V then we call V to be a neutrosophic-neutrosophic set linear algebra over S (n-n set linear algebra over S).*

We first illustrate this definition by some examples before we prove some properties about them.

*Example 2.3.9:* Let V = {m – mI | m ∈ $Z^+$} be a n-n set linear algebra over the set S = {1, 1 – I} $\subseteq$ N(Z).

*Example 2.3.10:* Let V = {ZI} $\subseteq$ PN(Z). V is a n-n set linear algebra over the set S = {0, 1, 1 – I} $\subseteq$ N(Z).

*Example 2.3.11:* Let

$$V = \left\{ \begin{pmatrix} ZI & ZI \\ ZI & ZI \end{pmatrix} \right\},$$

V is a n-n set linear algebra over the set $S = \{0, 1, 1 - I\}$.

It is important and interesting to note that as in case of linear algebra, n-n linear algebras is a set vector space but in general a n-n set vector space is not a n-n linear algebra. The following examples show that a n-n set vector space is not a n-n set linear algebra.

***Example 2.3.12:*** Let $V = \{3I, 24I, 41I, 26I, 0, -13I, 48I\} \subseteq PN(Z)$. V is a n-n set vector space over the set $S = \{0, 1, I\} \subseteq N(Z)$. We see V is not a n-n set linear algebra over S; as $3I + 24I = 27I \notin V$ and so on.

Thus in general a n-n set vector space is not a n-n set linear algebra over S. Now we proceed onto define the new concept of n-n set linear subalgebra.

**DEFINITION 2.3.4:** *Let V be a n-n set linear algebra over the set $S \subseteq N(Z)$. ($S \not\subseteq Z$). Suppose W is a subset of V such that W is a n-n set linear algebra over the set $S \subseteq N(Z)$ the we call W to be a n-n set linear subalgebra of V over the set S.*

We will illustrate this situation by some examples.

***Example 2.3.13:*** Let $V = \{3ZI\}$ be a n-n set linear algebra over the set $S = Z^+I$. Take $W = \{9ZI\} \subseteq V$; W is a n-n set linear algebra of V over S.

***Example 2.3.14:*** Let $V = \{m - mI \mid m \in Z^+\}$ be a n-n set linear algebra over the set $S = \{0, 1, 1 - I\} \subseteq N(Z)$. $W = \{3m - 3mI \mid m \in Z^+\} \subseteq V$ is a n-n set linear subalgebra of V over S.

***Example 2.3.15:*** Let
$$V = \left\{ \begin{pmatrix} m - mI & 0 \\ 3Z^+I & Z^+I \end{pmatrix}, \begin{pmatrix} m - mI & 0 \\ 0 & 0 \end{pmatrix} \middle| m \in Z^+ \right\}$$

be a n-n linear algebra over the set $S = \{0, 1, 1 - I\} \subseteq N(Z)$.
Take

$$W = \left\{ \begin{pmatrix} m - mI & 0 \\ 0 & 0 \end{pmatrix} \middle| m \in Z^+ \right\} \subseteq V;$$

W is a n-n linear subalgebra over the set $S = \{0, 1, 1 - I\}$.

*Example 2.3.16:* Let

$$V = \left\{ \begin{pmatrix} 0 & n - nI & p - pI \\ t - tI & 0 & 0 \end{pmatrix} \middle| n, p, t \in Z^+ \right\}$$

be a n-n set linear algebra over the set $S = \{0, 1, 1 - I\} \subseteq N(Z)$.

$$W = \left\{ \begin{pmatrix} 0 & 2n - 2nI & 2p - 2pI \\ 2t - 2tI & 0 & 0 \end{pmatrix} \middle| n, p, t \in Z^+ \right\} \subseteq V$$

is a n-n set linear subalgebra of V over the set S.

Now we proceed onto define yet another new substructure.

**DEFINITION 2.3.5:** *Let V be a n-n set linear algebra over S. Suppose W is a proper subset of V and W is only a n-n set vector space over S then we call W to be a pseudo n-n set vector subspace of V over S.*

We will illustrate this by some simple examples.

*Example 2.3.17:* Let

$$V = \left\{ \begin{pmatrix} m - mI & 0 \\ 0 & 0 \end{pmatrix}, \begin{pmatrix} 0 & 0 \\ 0 & m - mI \end{pmatrix}, \begin{pmatrix} m - mI & 0 \\ 0 & m - mI \end{pmatrix} \middle| m \in Z^+ \right\}$$

be a n-n set linear algebra over the set $S = \{0, 1, 1 - I\}$.

Take

$$W = \left\{ \begin{pmatrix} m-mI & 0 \\ 0 & 0 \end{pmatrix}, \begin{pmatrix} 0 & 0 \\ 0 & m-mI \end{pmatrix} \middle| m \in 2Z^+ \right\} \subseteq V.$$

W is a pseudo n-n set vector subspace of V over S.

***Example 2.3.18:*** Let $V = \{m \pm mI \mid m \in 2Z^+\}$ be a n-n set linear algebra over the set $S = \{0, 1, 1 - I\} \subseteq N(Z)$. Take $W = \{3 + 3I, 3 - 3I, 5 - 5I, 5 + 5I, 27 + 27I, 0\} \subseteq V$. W is only a pseudo n-n set vector space over the set $S = \{0, 1, 1 - I\} \subseteq N(Z)$.

***Example 2.3.19:*** Let $V = \{2ZI, 2Z, 2nZ + 2mZI \mid m, n \in Z^+\} \subseteq N(Z)$ be a n-n set linear algebra over the set $S = \{0, 1, 1 - I\}$. Take $W = \{2ZI, 2Z\} \subseteq V$; W is a pseudo n-n set vector subspace of V over the set S.

Now we define yet another new substructure.

**DEFINITION 2.3.6:** *Let V be a n-n set linear algebra over the set $S \subseteq N(Z)$. Let $W \subseteq V$ be a proper subset of V and $T \subseteq S$ be a proper subset of S. If W is a n-n subset linear algebra over T then we call W to be a n-n subset linear subalgebra of V over the subset T of S.*

We will illustrate this by some simple examples.

***Example 2.3.20:*** Let

$$V = \left\{ \begin{pmatrix} m-mI & 0 & 0 \\ 0 & m-mI & 0 \end{pmatrix}, \begin{pmatrix} 0 & m-mI & 0 \\ m-mI & 0 & m-mI \end{pmatrix}, \right.$$

$$\left. \begin{pmatrix} m_1-m_1I & m_2-m_2I & 0 \\ m_3-m_3I & m_4-m_4I & m_5-m_5I \end{pmatrix} \middle| m_i, m \in Z^+ \cup \{0\}; 1 \le i \le 5 \right\}$$

be a n-n set linear algebra over the set $S = \{0, 1, 1 - I\} \subseteq N(Z)$.

Take
$$W = \left\{ \begin{pmatrix} 0 & m-mI & 0 \\ m-mI & 0 & m-mI \end{pmatrix} \middle| m \in Z^+ \right\} \subseteq V.$$

W is a n-n set linear algebra over the subset $T = \{0, 1 - I\}$ of S. Hence W is a n-n subset linear subalgebra of V over T.

*Example 2.3.21:* Let $V = \{2Z^+I, 7Z^+I, 3Z^+I, 5Z^+I\} \subseteq N(Z)$ be a n-n set linear algebra over $S = Z^+$ I. $W = \{2Z^+I\} \subseteq V$ is a n-n subset linear subalgebra of V over the subset $P = 3Z^+$ I of S.

We now proceed onto define the notion of pseudo n-n subset vector subspace of a n-n set linear algebra.

**DEFINITION 2.3.7:** *Let V be a n-n set linear algebra over the set $S \subseteq N(Z)$. Let W be a subset of V and if W is a n-n set vector space over a subset T of S then we define W to be a pseudo n-n subset vector subspace of V over the subset T of S.*

We illustrate this definition by some examples.

*Example 2.3.22:* Let $V = \{m - mI \mid m \in Z^+ \cup \{0\}\}$ be a n-n set linear algebra over the set $S = \{0, 1, 1 - I\} \subseteq N(Z)$. Take $W = \{5 - 5I, 28 - 28I, 40 - 40I, 0, 18 - 18I\} \subseteq V$; W is a pseudo n-n subset vector subspace of V over the subset $T = \{0, 1 - I\}$ of S.

*Example 2.3.23:* Let
$$V = \left\{ \begin{pmatrix} m-mI & 0 \\ m-mI & 0 \end{pmatrix} \right\}$$

such that $m \in Z^+ \cup \{0\}$ be a n-n set linear algebra over the set $S = \{0, 1, 1 - I\}$.
Take
$$W = \left\{ \begin{pmatrix} 5-5I & 0 \\ 5-5I & 0 \end{pmatrix}, \begin{pmatrix} 8-8I & 0 \\ 8-8I & 0 \end{pmatrix}, \begin{pmatrix} 7-7I & 0 \\ 7-7I & 0 \end{pmatrix}, \right.$$

$$\left. \begin{pmatrix} 25-25I & 0 \\ 25-25I & 0 \end{pmatrix}, \begin{pmatrix} 0 & 0 \\ 0 & 0 \end{pmatrix} \right\} \subseteq V;$$

W is a pseudo n-n subset vector subspace of V over the subset T = {0, 1 – I} of S.

Now we proceed onto define the notion of linear transformation and linear operator of n-n set vector spaces defined over a set S ⊆ N(Z).

**DEFINITION 2.3.8:** *Let V and W be any two n-n set vector spaces over the same set S ⊆ N(Z). A map T: V → W is said to be a n-n set linear transformation of V into W if T (I) = I and T (αυ) = αT (υ) for every α ∈ S and for every υ ∈ V. If W = V in this definition we call T from V to V to be a n-n set linear operator of V.*

We illustrate this by some simple examples.

*Example 2.3.24:* Let V = {2 ± 2I, 8 ± 8I, 27 ± 27I, 45 ± 45I, 0, 35 ± 35I} and W = {0, 1 ± I, 6 ± 6I, 20 ± 20I, 49 ± 49I, 26 ± 26I, 11 ± 11I, 8 ± 8I, 17 ± 17I} be two n-n set vector spaces defined over the set S = {0, 1, 1 –I} ⊆ N(Z).

Define T: V → W by T (0) = 0;
$$T(2 \pm 2I) = 6 \pm 6I$$
$$T(8 \pm 8I) = 20 \pm 20I,$$
$$T(27 \pm 27I) = 26 \pm 26I$$
$$T(45 \pm 45I) = 49 \pm 49I \text{ and}$$
$$T(35 \pm 35I) = 11 \pm 11I.$$

T is a n-n set linear transformation of V to W.

*Example 2.3.25:* Let V = {0, 2m ± 2mI, 5m ± 5mI | m ∈ $Z^+$} be a n-n set neutrosophic vector space over S = {0, 1, 1 –I} and W = {27m + 27mI, 0, 8m + 8mI | m ∈ $Z^+$} be a set neutrosophic vector space over te set S = {0, 1, 1 – I}.

Define T: V → W as T (0) = 0.
$$T(2m \pm 2mI) = 8m \pm 8mI.$$
$$T(5m \pm 5mI) = 27m + 27mI$$
then T is a n-n set linear transformation of V into W.

Now we give some examples of n-n set linear operators on V, V a n-n set vector space defined over the set S.

*Example 2.3.26:* Let V = {7 ± 7I, 0, 21 ± 21I, 63 ± 63I, 3I, 9I, 27I, 45I, 63I, 15I} be a n-n set vector space over the set S = {0, 1, 1 – I}.
Define T: V → V by
$$T(0) = 0$$
$$T(I) = I$$
$$T(7 \pm 7I) = 21 \pm 21I$$
$$T(21 \pm 21I) = 63 \pm 63I$$
$$T(63 \pm 63I) = 63I$$
$$T(3I) = 9I$$
$$T(9I) = 27I$$
$$T(45I) = 15I$$
$$T(27I) = 27I$$
$$T(15I) = 45I \text{ and}$$
$$T(63I) = 63 \pm 63I.$$

T is a n-n set linear operator on V.

*Example 2.3.27:* Let V = {2ZI, 81ZI, 47ZI, 0} be a n-n set vector space over the set S = $2Z^+ I$.
T : V → V as follows:
Define T (0) = 0
$$T(2ZI) = 47ZI$$
$$T(81ZI) = 2ZI$$
$$T(47ZI) = 81ZI,$$
T is a n-n set linear operator on V.

**DEFINITION 2.3.9:** *Let V and W be n-n set linear algebra defined on the same set S ⊆ N(Z). A map T: V → W is said to a n-n set linear transformation from V to W if the following conditions are satisfied:*

$$T(I) = I$$
$$T(v_1 + v_2) = T(v_1) + T(v_2)$$
$$T(\alpha v) = \alpha T(v)$$

*for all $v_1$, $v_2$, $v \in V$ and $\alpha \in S$. If $W = V$ then the n-n set linear transformation is defined to be a n-n set linear operator on V.*

We shall illustrate this situation by the following examples:

*Example 2.3.28:* Let

$$V = \left\{ \begin{pmatrix} m-mI & 0 \\ 0 & 0 \end{pmatrix} \middle| m \in Z^+ \cup \{0\} \right\}$$

and $W = \{m - mI \mid m \in Z^+ \cup \{0\}\}$ be two n-n set linear algebras over the set $S = \{0, 1, 1 - I\} \subseteq N(Z)$.
Define $T: V \to W$ by

$$T\begin{pmatrix} m-mI & 0 \\ 0 & 0 \end{pmatrix} = 2m - 2mI$$

for every

$$\begin{pmatrix} m-mI & 0 \\ 0 & 0 \end{pmatrix} \in V$$

and $T(0) = 0$. Thus T is a n-n set linear transformation of V into W.

*Example 2.3.29:* Let $V = \{2ZI\}$ and

$$W = \left\{ \begin{pmatrix} 2ZI & 0 \\ 0 & m-mI \end{pmatrix} \middle| m \in Z^+ \cup \{0\} \right\}$$

be two n-n set linear algebras over the set $S = \{0, 1, 1 - I\} \subseteq N(Z)$.

Define $T: V \to W$ by

$$T(2ZI) = \begin{pmatrix} 2ZI & 0 \\ 0 & 0 \end{pmatrix} \text{ and } T(0) = \begin{pmatrix} 0 & 0 \\ 0 & 0 \end{pmatrix},$$

T is a n-n set linear transformation of V into W.

*Example 2.3.30:* Let

$$V = \left\{ \begin{pmatrix} ZI & 2ZI & 3ZI \\ 2ZI & 0 & 0 \\ 3ZI & 0 & 0 \end{pmatrix} \right\}$$

be a n-n set linear algebra over the set $S = Z^+ \cup \{0\}$.
Define a map $T : V \to V$

$$T \begin{pmatrix} 0 & 0 & 0 \\ 0 & 0 & 0 \\ 0 & 0 & 0 \end{pmatrix} = \begin{pmatrix} 0 & 0 & 0 \\ 0 & 0 & 0 \\ 0 & 0 & 0 \end{pmatrix}$$

$$T \begin{pmatrix} ZI & 2ZI & 3ZI \\ 2ZI & 0 & 0 \\ 3ZI & 0 & 0 \end{pmatrix} = \begin{pmatrix} 2ZI & 4ZI & 6ZI \\ 4ZI & 0 & 0 \\ 6ZI & 0 & 0 \end{pmatrix}.$$

It is easily verified that T is a n-n set linear operator on V.

*Example 2.3.31:* Let $V = \{12ZI\}$ be a n-n set linear algebra over the set $S = 2Z^+I \subseteq N(Z)$.
Define T: $V \to V$ by $T(12ZI) = 24ZI$; it is easily verified to be a n-n set linear operator on V.

**DEFINITION 2.3.10:** *Let V and W be n-n set linear algebras over the set $S \subseteq N(Z)$. If T: $V \to W$ be a n-n set linear transformation of V into W such that T preserves atleast one n-n set linear subalgebra of V then we define T to be a weak n-n set subalgebra preserving linear transformation. If T preserves every n-n set linear subalgebra of V then T we define T to be a strong n-n set subalgebra preserving linear transformation. If*

*W = V then we call T to be a strong (weak) n-n set subalgebra preserving operator on V.*

We shall illustrate these concepts by some examples.

***Example 2.3.32:*** Let $V = \{m - mI \mid m \in Z^+ \cup \{0\}\}$ and

$$W = \left\{ \begin{pmatrix} m - mI & 0 \\ m - mI & m - mI \end{pmatrix} \middle| m \in Z^+ \cup \{0\} \right\}$$

be two n-n set linear algebras over the set $S = \{0, 1, 1 - I\} \subseteq N(Z)$. If $T: V \to W$ be such that

$$T(m - mI) = \begin{pmatrix} m - mI & 0 \\ m - mI & m - mI \end{pmatrix}$$

and

$$T(0) = \begin{pmatrix} 0 & 0 \\ 0 & 0 \end{pmatrix}$$

then T is a strong n-n set subalgebra preserving linear transformation of V into W.

***Example 2.3.33:*** Let $V = \{ZI\}$ and $W = \{5ZI\}$ be two n-n set linear algebras over the set $Z^+ I \cup \{0\}$. The map $T: V \to W$ given by $T(0) = 0$, $T(I) = 5I$ is a strong n-n set subalgebra preserving linear transformation of V into W.

***Example 2.3.34:*** Let $V = \{m - mI, 0 \mid m \in Z^+\}$ be a n-n set linear algebra over the set $S = \{0, 1, 1 - 2I\}$. The map $T: V \to V$ such that $T(m - mI) = 2m - 2mI$, for every $m - mI \in V$ is a strong n-n set subalgebra preserving linear operator on V.

Now we proceed onto define the notion of generator for n-n set vector space and n-n set linear algebra.

**DEFINITION 2.3.11:** *Let V be a n-n set vector space over the set S $\subseteq$ N(Z). Suppose B $\subseteq$ V is a non empty subset of V such that every v $\in$ V can be represented by v = sb for some b $\in$ B and s $\in$ S then we call B to be the n-n generator set (or n-n generating set) of V over the set S.*

We illustrate this by some simple examples.

*Example 2.3.35:* Let V = {3 ± 3I, 0, 7 ± 7I, 15 ± 15I, 20 ± 20I} $\subseteq$ N(Z) be a n-n set vector space over the set S = {0, 1, 1 – I} $\subseteq$ N(Z). Take B = {3 + 3I, 7 + 7I, 15 + 15I, 20 + 20I} $\subseteq$ V. It is easily verified B is a n-n generator of V over the set S.

*Example 2.3.36:* Let V = {2ZI, 0, 15ZI} $\subseteq$ N(Z) be the n-n set vector space over the set S = $Z^+$I $\subseteq$ N(Z). Take B = {±2I, ±15I, 0} $\subseteq$ V, B is n-n generator set V over the set S.

**DEFINITION 2.3.12:** *Let V be a n-n set linear algebra over the set S $\subseteq$ N(Z). Let C $\subseteq$ V be a subset of V such that every element v $\in$ V can be represented as v = sc or v = $\sum s_i c_i$ for some s, $s_i$ $\in$ S and c, $c_i$ $\in$ C then we call C to be a n-n generator set (generating set) of the n-n set linear algebra V over the set S.*

We illustrate this by some simple examples.

*Example 2.3.37:* Let V = {m – mI | m $\in$ $Z^+$ $\cup$ {0}} be a n-n set linear algebra over the set S = {0, 1, 1 – I}. Take C = {m – mI | m $\in$ $Z^+$} $\subseteq$ V is the n-n generator set of V over S.

*Example 2.3.38:* Let V = {2ZI} $\subseteq$ N(Z) be a n-n set linear algebra over the set S = $Z^+$ $\subseteq$ N(Z). Let C = {± 2I} $\subseteq$ V, C is a n-n generator set of V over the set S = $Z^+$ $\subseteq$ N(Z). When the n-n generator set C of V (V; n-n set linear algebra or V a n-n set vector space) has finite number of elements in C then we say V has n-n finite generator C. If C has infinite cardinality then we say V has n-n infinite generator set.

## 2.4 Mixed Set Neutrosophic Rational Vector Spaces and their Properties

In the section we introduce the notion of mixed (pure) set neutrosophic rational vector space and describe some of their properties. Throughout this book $Q(I)$ denotes the neutrosophic rational field; i.e., $Q = Q(I) = \{a + bI \mid a, b \in Q\}$. $PN(Q)$ contains only neutrosophic rational numbers of the form $\{a + bI \mid b \neq 0 \text{ and } a, b \in Q\}$. So $PN(Q)$ hereafter will be known as pure set neutrosophic rational numbers.

**DEFINITION 2.4.1:** *Let $V \subseteq N(Q)$ ($PN(Q)$) be a proper subset of $N(Q)$ or $V$ contains elements from $N(Q)$ ($PN(Q)$) ($V \nsubseteq Q$). Let $S \subseteq N(Q)$ be a proper subset of $N(Q)$. We say $V$ is a mixed (pure) set neutrosophic rational vector space over $S$ if $sv \in V$ for every $s \in S$ and $v \in V$.*

*Example 2.4.1:* Let
$$V = \left\{\frac{2}{7} - \frac{2I}{7}, 0, \frac{19}{2} - \frac{19I}{2}, \frac{27}{5} - \frac{27I}{5},\right.$$

$$\left. 17I, 48 - 48I, \frac{28}{13} - \frac{28I}{13}, \frac{47I}{5}\right\} \subseteq N(Q).$$

Take
$$S = \{0, 1, \frac{11}{7} - \frac{11I}{7}, 1 - I\} \subseteq N(Q).$$

It is easily verified $V$ is a pure set neutrosophic rational vector space over the set $S$.

*Example 2.4.2:* Let
$$V = \left\{0, \frac{27}{5}, \frac{27}{5} - \frac{27I}{5}, \frac{48}{7}, \frac{48}{7} - \frac{48I}{7}, \frac{19I}{7}, \frac{8I}{23}\right\} \subseteq N(Q)$$

be a mixed set neutrosophic rational vector space over the set $S = \{0, 1, 1 - I\}$.

Now we proceed to define mixed (pure) set neutrosophic rational linear algebra over the set $S \subseteq N(Q)$.

**DEFINITION 2.4.2:** *Let $V \subseteq N(Q)$ (or V contains entries from N(Q)) be a mixed (pure) set neutrosophic rational vector space over a set $S \subseteq N(Q)$. If in V we have for every $v, u, \in V$ $u + v$ and $v + u \in V$ then we call V to be a mixed (pure) set neutrosophic rational linear algebra over the set S.*

We shall illustrate this by the following examples:

***Example 2.4.3:*** Let $V = \{m - mI \mid m \in Q\} \subseteq PN(Q) \cup \{0\}$. V is a pure set neutrosophic rational linear algebra over the set $S = \{0, 1, 1 - I, I\}$.

***Example 2.4.4:*** Let

$$V = \left\{ \begin{pmatrix} QI & m - mI \\ m - mI & QI \end{pmatrix} \middle| m \in Q \right\};$$

V is a pure set neutrosophic rational linear algebra over $S = \{0, 1, 1 - I\}$.

***Example 2.4.5:*** Let $V = \{QI, Q, Q + QI\} = \{mI, n, t + pI \mid m, n, t, p \in Q\} \subseteq N(Q)$. V is a mixed set neutrosophic rational linear algebra over $S = \{1, 0, 1 - I\}$.

*Remark:* The following facts are interesting about these new structures.

(1) We can define them over Q or N(Q) or PN(Q); we do not want to distinguish it by different names. We have to show that they are different whenever the set over which they are defined are different. However the reader can know the difference by the context.

(2) We have just said S ⊆ N(Q) so S ⊆ Z or S ⊆ Q or S ⊆ N(Q) and S ⊆ Q. But while studying the reader can understand over which they are defined.

(3) Further we understand all mixed (pure) set neutrosophic integer vector spaces (linear algebras) are mixed (pure) set neutrosophic rational vector spaces (linear algebras). However the mixed (pure) set neutrosophic rational vector space (or linear algebra) in general is not a mixed (pure) set neutrosophic integer vector spaces (or linear algebras).

(4) The notion of n-n set integer vector space is merged in case of neutrosophic set rational vector space.

However the difference in the structure is evident to any reader.

Now we proceed onto define the substructures of these new structures.

**DEFINITION 2.4.3:** *Let V be any mixed set neutrosophic rational vector space over the set S ⊆ N(Q). Let W be any proper subset of V. If W is a mixed set neutrosophic rational vector space over S, then we define W to be a mixed set neutrosophic rational vector subspace of V over S.*

We will illustrate this by some simple examples.

*Example 2.4.6:* Let

$$V = \left\{\frac{5}{3} - \frac{5I}{3}, \frac{20}{7} - \frac{20I}{7}, \frac{20}{7}, \frac{21}{8} - \frac{21I}{8}, \frac{5}{3}, 0, \frac{48}{19} - \frac{48I}{19}, \frac{122}{31} - \frac{122I}{31}, \frac{122}{31}\right\} \subseteq N(Q)$$

be a mixed set neutrosophic rational vector space over S = {0, 1, 1 – I}. Take

$$W = \left\{0, \frac{20}{7}, \frac{20}{7} - \frac{20I}{7}, \frac{5}{3}, \frac{5}{3} - \frac{5I}{3}\right\} \subseteq N(Q).$$

W is a mixed set neutrosophic rational vector subspace of V over S. It is interesting to note that every subset of V need not in general be a mixed set neutrosophic rational vector subspace of V over S.

For take

$$W_1 = \left\{ \frac{20}{7}, \frac{5}{3} - \frac{5I}{3} \right\} \subseteq V;$$

$W_1$ is not a mixed set neutrosophic rational vector subspace of V over S as

$$\frac{20}{7}(1 - I) = \frac{20}{7} - \frac{20I}{7} \notin W_1 \text{ also } 0 \cdot \frac{5}{3} - \frac{5I}{3} = 0 \notin W_1.$$

***Example 2.4.7:*** Let $V = \{m - mI, mI \mid m \in Q^+ \cup \{0\}\}$ be the pure set neutrosophic rational vector space over the set $S = \{0, 1, I, 1 - I\} \subseteq N(Q)$. Take $W = \{mI \mid m \in Q^+ \cup \{0\}\} \subseteq V$; W is a pure set neutrosophic rational vector subspace of V over S.

**DEFINITION 2.4.4:** *Let V be a set neutrosophic rational linear algebra over the set $S \subseteq N(Q)$. Let W be a subset of V; if W is a set neutrosophic rational linear algebra over the set $S \subseteq N(Q)$ then we call W to be a set neutrosophic linear subalgebra of V over the set S.*

We will illustrate this by some examples.

***Example 2.4.8:*** Let $V = \{0, m - mI \mid m \in Q^+\}$ be a set neutrosophic set linear algebra over the set $S = \{0, 1, 1 - I\} \subseteq N(Q)$.

Take $W = \{0, 3m - 3mI \mid m \in Z^+\} \subseteq V$; W is a set neutrosophic rational linear subalgebra of V over S. Infact W is a set neutrosophic integer set linear subalgebra of V over S. Take $W_1 = \{0, m - mI \mid m \in \{1/2n\} \mid n \in Z^+\} \subseteq V$. $W_1$ is a set neutrosophic rational linear subalgebra of V over S.

Infact $W_1$ is not a set neutrosophic integer set linear subalgebra of V over S.

*Example 2.4.9:* Let

$$V = \left\{ \begin{pmatrix} m - mI & 0 \\ 0 & pI \end{pmatrix} \middle| p, m \in Q^+ \cup \{0\} \right\}$$

be a set neutrosophic rational linear algebra over the set $S = \{0, 1, I, 1 - I\}$.
Take

$$W = \left\{ \begin{pmatrix} m - mI & 0 \\ 0 & 0 \end{pmatrix} \middle| m \in Q^+ \cup \{0\} \right\} \subseteq V;$$

W is a set neutrosophic rational linear subalgebra of V over S.

$$W_1 = \left\{ \begin{pmatrix} 0 & 0 \\ 0 & pI \end{pmatrix} \middle| p \in Q^+ \cup \{0\} \right\}$$

is also a set neutrosophic rational linear subalgebra of V over S.

**DEFINITION 2.4.5:** *Let V be a set neutrosophic rational vector space over the set S. Suppose $W \subseteq V$ is such that W is a set neutrosophic rational vector space over a proper subset $T \subseteq S$; then we call W to be a subset neutrosophic rational vector subspace of V over the subset T of S.*

We will illustrate this by some examples.

*Example 2.4.10:* Let $V = \{5 - 5I, 0, 25I, 41 - 41I, 60 - 60I, 60, 60I\}$ be a set neutrosophic rational vector space over the set $S = \{0, 1, I, 1 - I\}$. $W = \{0, 60, 60I, 60 - 60I, 25I\} \subseteq V$ is a subset neutrosophic rational vector subspace of V over $\{0, I, 1 - I\} = T \subseteq S$.

*Example 2.4.11:* Let

$$V = \left\{ \begin{pmatrix} 25I & 0 \\ 0 & 4I \end{pmatrix}, \begin{pmatrix} 0 & 9-9I \\ 8-8I & 0 \end{pmatrix}, \begin{pmatrix} 5I & 21-21I \\ 48I & 51-51I \end{pmatrix}, \right.$$

$$\left. \begin{pmatrix} 0 & 0 \\ 0 & 0 \end{pmatrix}, \begin{pmatrix} 5I & 0 \\ 48I & 0 \end{pmatrix}, \begin{pmatrix} 0 & 21-21I \\ 0 & 51-51I \end{pmatrix} \right\}$$

be a set neutrosophic rational vector space over the set S = {0, 1, I, 1-I}.

Take

$$W = \left\{ \begin{pmatrix} 0 & 0 \\ 0 & 0 \end{pmatrix}, \begin{pmatrix} 5I & 0 \\ 48I & 0 \end{pmatrix}, \begin{pmatrix} 25I & 0 \\ 0 & 4I \end{pmatrix} \right\} \subseteq V$$

is a subset neutrosophic rational vector subspace over T = {0, I} ⊆ S.

Now we proceed onto define yet another new substructure.

**DEFINITION 2.4.6:** *Let V be a set neutrosophic rational vector space over the set S. Let W ⊆ V be a proper subset of V if W is a set neutrosophic rational linear algebra over the set S; then we call W to be a pseudo set neutrosophic rational linear subalgebra of V over the set S.*

We illustrate this situation by some examples.

*Example 2.4.12:* Let

$$V = \{25p - 25pI, 25, 25I, \frac{17}{8} - \frac{17I}{8},$$

$$0, \frac{8}{7} - \frac{8I}{7}, \frac{3I}{7}, \frac{4I}{3}, 42I \mid p \in Z^+\}$$

be a set neutrosophic rational vector space over the set S={0, 1, I, 1 − I}. Let W = {25p − 25pI, 0 | p ∈ Z$^+$} ⊆ V. W is a pseudo set neutrosophic rational linear subalgebra of V over the set S.

*Example 2.4.13:* Let

$$V = \left\{ \begin{pmatrix} m - mI & 0 \\ 0 & m - mI \end{pmatrix}, \frac{22}{83} - \frac{22I}{83}, \frac{42I}{13}, \frac{88I}{9} \,\middle|\, m \in Q^+ \cup \{0\} \right\}$$

be a set neutrosophic rational vector space over the set $S = \{0, 1, I, 1-I\}$.
Take

$$W = \left\{ \begin{pmatrix} m - mI & 0 \\ 0 & m - mI \end{pmatrix} \,\middle|\, m \in Q^+ \cup \{0\} \right\} \subseteq V,$$

W is a pseudo set neutrosophic rational linear subalgebra of V over the set S.

**DEFINITION 2.4.7:** *Let V be a set neutrosophic rational linear algebra over the set S. Let $W \subseteq V$ be a set neutrosophic rational linear algebra over a proper subset T of S. We define W to be a subset neutrosophic rational linear subalgebra of V over the subset T of S.*

We will illustrate this situation by some examples.

*Example 2.4.14:* Let $V = \{m - mI \mid m \in Q^+ \cup \{0\}\}$ be a set neutrosophic rational linear algebra over the set $S = \{Z^+I, 0, 1 - I, 1\} \subseteq N(Q)$. Choose $W = \{m - mI \mid m \in Z^+ \cup \{0\}\} \subseteq V$; W is a subset neutrosophic rational sublinear algebra over the subset $T = \{2Z^+ I, 0\} \subseteq S$.

*Example 2.4.15:* Let

$$V = \left\{ \begin{pmatrix} 2ZI & m - mI \\ m - mI & 3ZI \end{pmatrix} \,\middle|\, m \in Q^+ \right\}$$

be a set neutrosophic rational linear algebra over the set $\{1 - I, 1, 2 - 2I, 5 - 5I, -19I + 19\} = S$.

Take

$$W = \left\{ \begin{pmatrix} 16ZI & 2m - 2mI \\ 5m - 5mI & 15ZI \end{pmatrix} \middle| m \in Q^+ \right\} \subseteq V;$$

W is a subset neutrosophic rational linear subalgebra over the set $T = \{1 - I, 1, 19 - 19I\} \subseteq S$.

As in case of neutrosophic integer set vector spaces (linear algebras) we can define the notion of set neutrosophic linear transformation, set neutrosophic linear operator, set neutrosophic subspace (sublinear algebra) preserving linear transformation and linear operator. Interested reader can construct examples of them as the analogous definition are identical. Now we proceed on to give the notion of set neutrosophic real set vector spaces and set neutrosophic real set linear algebra.

Throughout this book N(R) will denote the set of all neutrosophic reals, i.e., $N(R) = \{a + bI \mid a, b \in R\}$; PN(R) denote the pure neutrosophic reals i.e., elements of the form $\{a + bI \mid b \neq 0; a, b \in R\}$. Clearly $PN(R) \subset N(R)$, $R \subset N(R)$ but $R \not\subseteq PN(R)$ and $\langle R \cup I \rangle = N(R)$ in our usual notation.

**DEFINITION 2.4.8:** *Let $V \subseteq N(R)$ (or PN(R)) we say V is a set neutrosophic real vector space over a set $S \subseteq N(R)$ if for every $v \in V$ and for every $s \in S$, sv, vs $\in V$.*

We will illustrate this by some examples.

***Example 2.4.16:*** Let $V = \{r - rI \mid r \in R^+ \cup \{0\}\} \subseteq N(R)$; V is a set neutrosophic real vector space over the set $S = \{R^+I\} \subseteq PN(R)$.

***Example 2.4.17:*** Let

$$V = \left\{ \begin{pmatrix} RI & 0 \\ 0 & m - mI \end{pmatrix} \middle| m \in R \right\}$$

take S = {R$^+$I}. V is a set neutrosophic real vector space over the set S.

It is interesting to note that all neutrosophic integer set vector spaces and set neutrosophic rational vector spaces are also set neutrosophic real vector spaces. However when we work with real world problems according to need we can choose any of these set neutrosophic vector spaces.

Further, as we do not need all the axioms of a neutrosophic vector space to be satisfied by these structures these structures can be realized as the most generalized structures in neutrosophic vector spaces.

**DEFINITION 2.4.9:** *Let V be a subset of N(R) or PN(R) (or has entries from N(R) or PN(R)). We say V is a set neutrosophic real linear algebra if*
  1. *for every u, v ∈ V, u + v and v + u ∈ V and*
  2. *for every s ∈ S and v ∈ V; v s and sv ∈ V.*

All set neutrosophic real linear algebras are set neutrosophic real vector spaces but a set neutrosophic real vector space in general need not be a set neutrosophic real linear algebra.

We shall illustrate this by some simple examples.

*Example 2.4.17:* Let

$$V = \left\{ \begin{pmatrix} \sqrt{2} - \sqrt{2}I & 0 \\ -\sqrt{19}I + \sqrt{19} & 0.14 - 0.14I \end{pmatrix}, \right.$$

$$\begin{pmatrix} 20 - 20I & 0 \\ \sqrt{7} - \sqrt{7}I & \sqrt{11} - \sqrt{11}I \\ 0 & \dfrac{23}{\sqrt{3}} - \dfrac{23I}{\sqrt{3}} \\ 3I & 41 - 41I \end{pmatrix},$$

$$\left\{ \frac{46}{\sqrt{5}} - \frac{46I}{\sqrt{5}}, -11\sqrt{2}\,I + 11\sqrt{2}, 4I, \frac{\sqrt{17}}{5} - \frac{\sqrt{17}I}{5}, 0 \right\}$$

be a set neutrosophic real vector space over the set $S = \{1, 1 - I, 5 - 5I\}$.

We see clearly V is not a set neutrosophic real linear algebra over S.

*Example 2.4.18:* Let

$$V = \{(mI, 0, 4m - 4mI, \frac{8m}{\sqrt{3}} - \frac{8mI}{\sqrt{3}}) \mid m \in Z^+\} \subseteq N(R)$$

be a set neutrosophic real linear algebra over the set

$$S = \{1, 1 - I, \frac{2}{\sqrt{7}} - \frac{2I}{\sqrt{7}}, \frac{5}{\sqrt{3}} - \frac{5I}{\sqrt{3}} \}.$$

*Example 2.4.19:* Let

$$V = \left\{ \begin{pmatrix} mI & 0 \\ \frac{3m}{\sqrt{2}} - \frac{3m}{\sqrt{2}} & m - mI \end{pmatrix} \middle| m \in R^+ \right\}$$

be a set neutrosophic real linear algebra over the set

$$S = \{1, 1 - I, \frac{41}{\sqrt{7}} - \frac{41I}{\sqrt{7}}, \frac{2}{\sqrt{7}} - \frac{2I}{\sqrt{7}}, \frac{21}{\sqrt{41}} - \frac{21I}{\sqrt{41}} \}.$$

Now all substructures pseudo substructures and set neutrosophic real transformation, set neutrosophic real operator can be defined as in case of neutrosophic integer set vector space and set neutrosophic rational vector space.

The following results can be easily proved by the reader.

(1) Set neutrosophic real vector space is not a set neutrosophic rational vector space or a neutrosophic integer set vector space.
(2) {The class of neutrosophic integer set vector space} $\subseteq$ {The class of set neutrosophic rational vector space} $\subseteq$ {The class of set neutrosophic real vector space} – prove.
(3) Study set neutrosophic real linear transformation of set neutrosophic real vector spaces V and W over a set S.
(4) Study and obtain some interesting properties about set neutrosophic real linear operator of the set neutrosophic real vector space V over the set S.
(5) Obtain interesting results about set neutrosophic real linear transformation(operator) which preserves subspaces.

Next we proceed onto define the new notion of neutrosophic modulo integers. $\langle Z_n \cup I \rangle = \{a + bI \mid a, b \in Z_n\} = N(Z_n)$ will denote the neutrosophic modulo integers. $PN(Z_n) = \{a + bI \mid b \in Z_n \setminus \{0\}\}$ denotes the pure neutrosophic modulo integers.

**DEFINITION 2.4.10:** *Let $V \subseteq N(Z_n)$ or a subset of $PN(Z_n)$ (V can be a set with has entries from $N(Z_n)$ or $PN(Z_n)$). Let $S \subseteq N(Z_n)$. We say V is a set neutrosophic modulo integer vector space if for every $v \in V$ and $s \in S$ vs for sv is in V.*

We will illustrate this by some examples.

*Example 2.4.20:* Let $V = \{0, 2, 6, 4, 6I, 8I, 10I\} \subseteq N(Z_{12})$ and $S = \{0, 3\} \subseteq N_{12}$. V is a set neutrosophic modulo integer vector space over the set S.

*Example 2.4.21:* Let $V = \{0, 2I, 8I, 6I, 4I, 2, 4, 6, 8\} \subseteq N(Z_{10})$ and $S \subseteq \{0, 1, 2, I, 2I\} \subseteq N(Z_{10})$. V is a set neutrosophic modulo integer vector space over the set S.

*Example 2.4.22:* Let V = {0, 1I, 2I, 3I, 4I, 5I, 6I} $\subseteq$ N($Z_7$). S = {0, 1, I} $\subseteq$ N($Z_7$). V is a set neutrosophic modulo integer vector space over the set S.

**DEFINITION 2.4.11:** *Let V $\subseteq$ N($Z_n$) (or V has entries from N($Z_n$)) and S $\subseteq$ N($Z_n$). If V in addition being a set neutrosophic modulo integer vector space over S satisfies the condition, that for every pair v, u $\in$ V, u + v and v + u $\in$ V; then we call V to be a set neutrosophic modulo integer linear algebra over S.*

We illustrate this by some simple examples.

*Example 2.4.23:* Let V = {0, I, 2I, 3I, 4I, 5I, 6I, 7I, 8I, 9I, 10I} $\subseteq$ N($Z_{11}$) and S = {0, 1, I, 5, 3, 2I, 6I, 8I} $\subseteq$ N($Z_{11}$). V is a set neutrosophic modulo integer linear algebra over S.

*Example 2.4.24:* Let V = {0, 2I, 4I, 6I, 8I, 10I, 12I, 14I, 16I} $\subseteq$ N($Z_{18}$), S = {0, 1, 2, 4, 8, 2I, 6I, 10I}. V is a set neutrosophic modulo integer linear algebra over S.

*Example 2.4.25:* Let V = {0, 1 + 9I, 2 + 8I, 3 + 7I, 4 + 6I, 5I, 8 + 2I, 7 + 3I, 6 + 4I, I + 9} $\subseteq$ N($Z_{10}$). V is a set neutrosophic modulo integer linear algebra over S = {0, 1, 1 + 9I, I+ 9} $\subseteq$ N($Z_{10}$).

It is left as an exercise for the reader to prove "A set neutrosophic modulo integer vector space in general is not a set neutrosophic modulo integer linear algebra and every set neutrosophic modulo integer linear algebra is a set neutrosophic modulo integer vector space".

*Example 2.4.26:* Let V = {I, 0, 2I, 8I, 24I, 17I, 17, 22I} $\subseteq$ N($Z_{25}$) and S = {0, 1, I} $\subseteq$ N($Z_{25}$). V is a set neutrosophic modulo integer vector space over S only and is not a set neutrosophic modulo integer linear algebra over S.

Now we proceed on to define set neutrosophic modulo integer vector space.

**DEFINITION 2.4.12:** *Let $V \subseteq N(Z_n)$ be a set neutrosophic modulo integer vector space over the set $S \subseteq N(Z_n)$. Suppose $W \subseteq V$ is a proper subset of V such that W is a set neutrosophic modulo integer vector space over S then we define W to be a set neutrosophic modulo integer vector subspace of V over S.*

We will illustrate this by some examples.

*Example 2.4.27:* Let $V = \{2, 0, 4I, 26I, 14I, 4, 12I, 20I, 20, 24, 24I\} \subseteq N(Z_{28})$ be a set neutrosophic modulo integer vector space over the set $S = \{0, 1, I\} \subseteq N(Z_{28})$. Take $W = \{0, 4I, 12I, 20I, 26I, 14I\} \subseteq V$; W is a set neutrosophic modulo integer set vector subspace of V over S.

*Example 2.4.28:* Let $V = \{3I, 0, 6I, 21I, 9I, 3, 9, 6, 18I\} \subseteq N(Z_{27})$ and $S = \{0, 1, 3I, 3\} \subseteq N(Z_{27})$. V is a set neutrosophic modulo integer vector space over the set S. Take $W = \{0, 3, 3I, 6I, 18I, 9I\} \subseteq V$; W is a set neutrosophic modulo integer vector subspace of V over the set S.

Now we proceed onto define the notion of subset neutrosophic modulo integer vector subspace of V over a subset T of S.

**DEFINITION 2.4.13:** *Let $V \subseteq N(Z_n)$ be a set neutrosophic modulo integer vector space over the set $S \subseteq N(Z_n)$. Let $W \subseteq V$; W is said to be a subset neutrosophic modulo integer vector space of V over the subset T of S if W is a set neutrosophic modulo integer vector space over the set T.*

We will illustrate this by some examples.

*Example 2.4.29:* Let $V = \{0, 3I, 12, 4I, 5I, 4, 12I, 10I, 5, 3\} \subseteq N(Z_{15})$ be a set neutrosophic modulo integer vector space over the set $S = \{0, 3, 3I, 5, 5I\} \subseteq N(Z_{15})$. Take $W = \{0, 3I, 3, 5I, 5\} \subseteq V$ and $T = \{0, 3, 3I\} \subseteq S$. W is a subset neutrosophic modulo integer vector subspace of V over the subset T of S.

***Example 2.4.30:*** Let V = {0, 1, I, 18, 18I, 16, 16I, 13, 13I, 20, 20I} $\subseteq$ N($Z_{23}$) be a set neutrosophic modulo integer vector space over the set S = {0, 1, I} $\subseteq$ N($Z_{23}$). Take W = {0, 20, 20I, 16, 16I, 13, 13I} $\subseteq$ V and T = {0, I} $\subseteq$ S, W is a subset neutrosophic modulo integer vector subspace of V over the subset T of S.

It is important and interesting to note that all set neutrosophic modulo integer vector spaces have only finite number of elements in them. It also implies the set neutrosophic generator subset of a set neutrosophic modulo integer vector space is always finite.

**DEFINITION 2.4.14:** *Let V $\subseteq$ N($Z_n$) be a set neutrosophic modulo integer vector space over the set S $\subseteq$ N($Z_n$). Let W $\subseteq$ V be such that W is a set neutrosophic modulo integer linear algebra over S then we call (or define) W to be a pseudo set neutrosophic modulo integer linear subalgebra of V over S.*

We will illustrate this situation by some examples.

***Example 2.4.31:*** Let V = {0, 2I, 4I, 6I, 8I, 10I, 12I, 14I, 16I, 18I, 17, 17I, 11, 11I, 19I} $\subseteq$ N($Z_{20}$) and S = {0, 1, 2I, I} $\subseteq$ N($Z_{20}$). V is a set neutrosophic modulo integer vector space over S.
    Take W = {0, 2I, 4I, 6I, 8I, 10I, 12I, 14I, 16I, 18I} $\subseteq$ V; W is a pseudo set neutrosophic modulo integer linear subalgebra of V over the set S.
    Take $W_1$ = {0, 4I, 8I, 12I, 16I} $\subseteq$ V; $W_1$ is a pseudo set neutrosophic modulo integer linear subalgebra of V over S.

***Example 2.4.32:*** Let V = {0, 3I, 3, 6I, 9I, 12I, 15I, 18I, 21I, 24I, 27I, 12, 18, 21} $\subseteq$ N($Z_{30}$) be a set neutrosophic modulo integer vector space over the set S = {0, 1, I} $\subseteq$ N($Z_{30}$).
    Take W = {0, 3I, 6I 9I, 12I, 15I, 18I, 21I, 24I, 27I} $\subseteq$ V. W is a pseudo set neutrosophic modulo integer linear subalgebra of V over S.

Now we proceed on to define the notion of substructures in set neutrosophic modulo integer linear algebras.

**DEFINITION 2.4.15:** *Let $V \subseteq N(Z_n)$ be a set neutrosophic modulo integer linear algebra over the set $S \subseteq N(Z_n)$. Suppose $W \subseteq V$ and if W is a set neutrosophic modulo integer linear algebra over the same set $S \subseteq N(Z_n)$; the we define W to be a set neutrosophic modulo integer linear subalgebra of V over the set S.*

We shall illustrate this situation by some examples.

*Example 2.4.33:* Let $V = \{0, I, 2I, 3I, 4I, …, 23I\} \subseteq N(Z_{24})$ be a set neutrosophic modulo integer linear algebra over the set $S = \{0, 1, I, 2I\} \subseteq N(Z_{24})$. Take $W = \{0, 2I, 4I, 6I, 8I, 10I, 12I, 14I, 16I, 18I, 20I, 22I\} \subseteq V$; W is a set neutrosophic modulo integer linear subalgebra of V over the set $S \subseteq N(Z_{24})$.

*Example 2.4.34:* Let $V = \{7nI, 0 \mid n = 1, 2, …, 6\} \subseteq N(Z_{49})$; i.e., $V = \{0, 7I, 14I, 21I, 28I, 35I, 42I\} \subseteq N(Z_{49})$ over the set $S = \{0, 3, 1, I, 2I\} \subseteq N(Z_{49})$. V is a set neutrosophic modulo integer linear algebra. But V has no set neutrosophic modulo integer sublinear algebra over S.

*Example 2.4.35:* Let $V = \{0, I, 2I, …, (p – I)I\} \subseteq N(Z_p)$ be a set neutrosophic modulo integer linear algebra over $S = \{0, 1, 2I, 3I, I, 5I\} \subseteq N(Z_p)$. V has no set neutrosophic modulo integer linear subalgebras.

Now we proceed onto define pseudo set neutrosophic modulo integer vector subspace of a set neutrosophic modulo integer vector linear algebra.

**DEFINITION 2.4.16:** *Let $V \subseteq N(Z_n)$ be a set neutrosophic modulo integer linear algebra over the set $S \subseteq N(Z_n)$. Suppose $W \subseteq V$; W is a subset of V such that W is a set neutrosophic modulo integer vector space over S then we define W to be a* *pseudo set neutrosophic integer vector subspace of V over the set S.*

*Example 2.4.36:* Let $V = \{0, I, 2I, 4I, 3I, \ldots, 24I\} \subseteq N(Z_{25})$ be a set neutrosophic modulo integer linear algebra over the set $S = \{0, 1, I, 1-I\} \subseteq N(Z_{25})$. Take $W = \{0, I, 4I, 2I, 24I\} \subseteq V$; W is a pseudo set neutrosophic modulo integer vector subspace of V over the set S.

*Example 2.4.37:* Let $V = \{0, 1 + 6I, I + 6, 2 + 5I, 2I + 5, 3 + 4I, 4 + 3I\} \subseteq N(Z_7)$ be a set neutrosophic modulo integer set linear algebra over the set $S = \{0, 1, I, 1 + 6I\} \subseteq N(Z_7)$. We see $W = \{0, 2 + 5I, 3 + 4I\} \subseteq V$ is a pseudo set neutrosophic modulo integer vector subspace of V over the set S.

It is left as a research problem for the reader to characterize those set neutrosophic modulo integer linear algebra V over a set S such that every subset of V is a pseudo set neutrosophic modulo integer vector subspace of V over the set S. Now we proceed onto define the notions of subset neutrosophic modulo integer linear subalgebra and pseudo subset neutrosophic modulo integer vector subspace of a set neutrosophic modulo integer linear algebra over the set S.

**DEFINITION 2.4.17:** *Let $V \subseteq N(Z_n)$ be a set neutrosophic modulo integer linear algebra over the set $S \subseteq N(Z_n)$. Let $W \subseteq V$ be a set neutrosophic modulo integer vector space over a subset T of S. We call W to be a pseudo subset neutrosophic modulo integer vector subspace of V over the subset $T \subseteq S$.*

We will illustrate this by some examples.

*Example 2.4.38:* Let $V = \{0, 1 + 28I, 28 + I, 2 + 27I, 27 + 2I, 3 + 26I, 26 + 3I, 4 + 25I, 25 + 4I, 5 + 24I, 24 + 5I, 6 + 23I, 23 + 6I, 7 + 22I, 22 + 7I, 21 + 8I, 8 + 21I, 9 + 20I, 20 + 9I, 19 + 10I, 10 + 19I, 11 + 18I, 18 + 11I, 17 + 12I, 12 + 17I, 13 + 16I, 16 + 13I, 14 + 15I, 15 + 14I\} \subseteq N(Z_{29})$ be a set neutrosophic modulo integer linear algebra over the set $S = \{0, 1, I, 1 + 28I\} \subseteq N(Z_{29})$.

Take W = {0, 1 + 28I, 9 + 20I, 15 + 14I, 12 + 17I, 7 + 22I, 3 + 26I} ⊆ V and T = {0, 1, I} ⊆ S; W is a pseudo subset neutrosophic modulo integer vector subspace of V over the subset T of S.

*Example 2.4.39:* Let V = {0, I, 2I, 3I, …, 26I} ⊆ $N(Z_{27})$ be a set neutrosophic modulo integer linear algebra over the set S = {0, 1, I, 1 + 26I} ⊆ $N(Z_{27})$. Take W = {0, I, 3I, 6I, 14I, 11I, 10I} ⊆ V and T = {0, 1, 1+26I} ⊆ S. W is a pseudo subset neutrosophic modulo integer vector subspace over the subset T of S.

**DEFINITION 2.4.18:** *Let V ⊆ $N(Z_n)$ be a set neutrosophic modulo integer linear algebra over a set S ⊆ $N(Z_n)$. Suppose W ⊆ V be such that W is a set neutrosophic modulo integer linear algebra over a subset T of S with cardinality of T greater than one i.e., |T| > 1, then we call W to be a subset neutrosophic modulo integer linear subalgebra of V over the subset T of S.*

We will illustrate this by some simple examples.

*Example 2.4.40:* Let V = {0, I, 2I, 3I, 4I, … , 48I} ⊆ $N(Z_{49})$ be a set neutrosophic modulo integer linear algebra over the set S = {0, 1, I, 7 + 24I} ⊆ $N(Z_{49})$. Take W = {0, 7I, 14I, 21I, 28I, 35I, 42I} ⊆ V; and T = {0, 7 + 42I} ⊆ S. W is a subset neutrosophic modulo integer linear subalgebra of V over the set T ⊆ S.

*Example 2.4.41:* Let V = {0, 1 + 14I, 14 + I, 2 + 13I, 13 + 2I, 12I + 3, 3I + 12, 4I + 11, 11I + 4, 10I + 5, 5I + 10, 6 + 9I, 6I + 9, 7I + 8, 8 + 7I} ⊆ $N(Z_{15})$ be a set neutrosophic modulo integer linear algebra over the set S = {0, 1, I, 1 + 14I, 14 + I} ⊆ $N(Z_{15})$. Take W = {1 + 14I, 14 + I, 0} ⊆ V and T = {0, 1, I}; W is a subset neutrosophic modulo integer linear subalgebra of V over the subset T ⊆ S.

Now we study some more properties of set neutrosophic modulo integer linear algebra.

**DEFINITION 2.4.19:** *Let $V \subseteq N(Z_n)$ be a set neutrosophic modulo integer linear algebra over the set $S \subseteq N(Z_n)$. If V has no proper subset W such that W is a set neutrosophic modulo integer linear subalgebra of V over the set S then we call W to be a simple set neutrosophic modulo integer linear algebra or set neutrosophic modulo integer simple linear algebra.*

We illustrate this situation by some simple examples.

*Example 2.4.42:* Let $V = \{0, I, 2I, 3I, 4I, 6I, 7I, 8I, 9I, 10I\} \subseteq N(Z_{11})$ be a set neutrosophic modulo integer linear algebra over the set $S = \{0, 1, I, 1 + 10I, 10 + 1I, 5 + 6I\} \subseteq N(Z_{11})$. V has no proper subset W such that W is a set neutrosophic modulo integer linear subalgebra of V. Thus V is a simple set neutrosophic modulo integer linear algebra over the set S.

*Example 2.4.43:* Let $V = \{0, 1 + 40I, 40 + I\} \subseteq N(Z_{41})$ be a set neutrosophic modulo integer linear algebra over the set $S = \{0, 1, 4 + 37I, 37 + 4I, I\} \subseteq N(Z_{41})$. V is a simple set neutrosophic modulo integer linear algebra over S.

We now proceed onto define the notion of weakly simple set neutrosophic modulo integer linear algebra.

**DEFINITION 2.4.20:** *Let $V \subseteq N(Z_n)$ be a set neutrosophic modulo integer linear algebra over a set $S \subseteq N(Z_n)$. If V has no proper subset W such that W is a subset neutrosophic modulo integer linear subalgebra over any proper subset T of S then we call W to be a weakly simple set neutrosophic modulo integer linear algebra over the set S.*

We will illustrate this situation by some simple examples.

*Example 2.4.44:* Let $V = \{0, I, 2I, …, 30I\} \subseteq N(Z_{31})$ be a set neutrosophic modulo integer linear algebra over the set $S = \{0, I\}$. V is a weakly simple set neutrosophic modulo integer linear algebra over the set S.

***Example 2.4.45:*** Let $V = \{0, 3 + 12I\} \subseteq N(Z_{15})$ be a set neutrosophic modulo integer set linear algebra over the set $S = \{0, I\} \subseteq N(Z_{15})$. V is a weakly simple set neutrosophic modulo integer linear algebra over S.

**THEOREM 2.4.1:** *Let $V \subseteq N(Z_n)$ be a set neutrosophic modulo integer linear algebra over a set $S \subseteq N(Z_n)$, $|S| = 2$ then V is a weakly simple set neutrosophic modulo integer linear algebra over S.*

*Proof:* Since $|S| = 2$ even if $W \subseteq V$ is such that W is a set neutrosophic modulo integer linear subalgebra over S it cannot be a subset neutrosophic modulo integer linear subalgebra over a subset T of S as $|S| = 2$ and $|T| > 1$ is not possible unless $|T| = |S|$. Hence the claim.

**THEOREM 2.4.2:** *Let $V = \{0, I, ..., (p-1)I\} \subseteq N(Z_p)$ where p is a prime be a set neutrosophic modulo integer linear algebra over a set $S \subseteq N(Z_p)$. V is both a simple set neutrosophic modulo integer linear algebra and a weakly simple set neutrosophic linear algebra.*

*Proof:* We see V has no subset W such that W is a set neutrosophic modulo integer linear subalgebra over S. So V is a simple set neutrosophic modulo integer linear algebra.

Now since V has no set neutrosophic modulo integer linear subalgebras even if $|S| > 2$ still we cannot find W in V such that W is a set neutrosophic modulo integer inear algebra over any proper subset of S. Hence the claim.

Now we proceed onto define the notion of set neutrosophic linear transformation of V, $W \subseteq N(Z_n)$ V and W neutrosophic modulo integer vector space and set neutrosophic linear operator on when W = V.

**DEFINITION 2.4.21:** *Let V and W be two set neutrosophic modulo integer vector spaces over a set $S \subseteq N(Z_n)$. A map $T: V \to W$ such that $T(sv) = sT(v)$ for all $s \in S$ and $v \in V$ is defined to be a set neutrosophic modulo integer linear transformation V to W.*

Let $N(Hom_S(V, W))$ be the collection of all set neutrosophic modulo integer linear transformations from V to W then $|N_S(Hom (V, W)| < \infty$. If $W = V$ then we call $T: V \to V$ to be a set neutrosophic modulo integer linear operator on V.

Clearly $|N_S(Hom(V, V))| < \infty$ where $N_S(Hom(V, V))$ denotes the collection of all set neutrosophic modulo integer linear operators from V to V.

The interested reader is expected to construct examples.

**DEFINITION 2.4.22:** *Let V and W be two set neutrosophic modulo integer linear algebras over the set $S \subseteq N(Z_n)$. A map $T: V \to W$ is said to be a set neutrosophic modulo integer linear transformation of these linear algebras V and W if*

*(1) $T(x + y) = T(x) + T(y)$*
*(2) $T(\alpha x) = \alpha T(x)$;*

*for all $x, y \in V$ and for all $\alpha \in S$.*

*If $V = W$, we call the map T to be a set neutrosophic modulo integer linear operator of V.*

**DEFINITION 2.4.23:** *Let V and W be two set neutrosophic modulo integer linear algebras over the set $S \subseteq N(Z_n)$.*

*A set neutrosophic modulo integer linear transformation T is said to be a strong set neutrosophic modulo integer subalgebra preserving linear transformation if T preserves set neutrosophic modulo integer linear subalgebras of V, i.e., if P is a set neutrosophic modulo integer linear subalgebra of V and if $T(P) = Q$ then Q is a set neutrosophic modulo integer linear subalgebra of W and this is true for every set neutrosophic modulo integer linear subalgebra P of V.*

*If atleast one set neutrosophic modulo integer linear subalgebra P of V is preserved under a set neutrosophic modulo integer linear transformation $T_1$ then we call $T_1$ to be a set neutrosophic modulo integer subalgebra preserving linear transformation of V to W.*

It is left as an exercise for the reader to prove the following theorem:

**THEOREM 2.4.3:** *Every strong set neutrosophic modulo integer linear subalgebra preserving linear transformation of V to W where V and W are set neutrosophic modulo integer vector spaces over a set $S \subseteq N(Z_n)$ is a set neutrosophic modulo integer linear subalgebra preserving linear transformation of V and W.*

Prove the converse is not true.

Several interesting results on set neutrosophic modulo integer linear algebras over a set $S \subseteq N(Z_n)$. The task of obtaining interesting and innovative results is left as an exercise for the reader.

We give an example of a different type of set neutrosophic modulo integer vector spaces.

***Example 2.4.46:*** Let

$$V = \left\{ \begin{pmatrix} 3I & 4I & 7I \\ 0 & 2I & 9I \\ 0 & 0 & 10I \end{pmatrix}, \begin{pmatrix} 2I & 5I & 9I \\ 0 & 10I & 7I \\ 0 & 0 & I \end{pmatrix}, \begin{pmatrix} 10I & 0 & 0 \\ 9I & 7I & 0 \\ 2I & I & 4I \end{pmatrix}, \right.$$

$$\left. \begin{pmatrix} 3I & 0 & 0 \\ 0 & 8I & 0 \\ 0 & 0 & 10I \end{pmatrix}, \begin{pmatrix} 3I & 7I & 9I \\ 0 & 0 & 0 \\ 0 & 0 & 0 \end{pmatrix}, \begin{pmatrix} 0 & 0 & 0 \\ 8I & I & 7I \\ 0 & 0 & 0 \end{pmatrix}, \begin{pmatrix} 0 & 0 & 0 \\ 0 & 0 & 0 \\ 6I & 7I & I \end{pmatrix} \right\}$$

where the entries of the neutrosophic matrices are from $N(Z_{11})$. Take $S = \{0, 1, I, 1 + 10I, 10 + I\} \subseteq N(Z_{11})$. V is a set neutrosophic modulo integer vector space over the set $S \subseteq N(Z_{11})$.

***Example 2.4.47:*** Let

$$V = \left\{ 0, 3I, \begin{pmatrix} 2I & 0 \\ 7I & 9I \end{pmatrix}, \begin{pmatrix} 12I & 3I & 0 \\ 0 & 4I & 5I \end{pmatrix}, \begin{pmatrix} 7I & 8I \\ 9I & 0 \\ 10I & 11I \end{pmatrix}, 4I, \right.$$

$$\left. \begin{pmatrix} 3I+10 & 2I+10 & 0 \\ 4I & 5I & 0 \\ 6I & 7I & 0 \\ 8I & 0 & 0 \end{pmatrix} \right\}$$

where the entries in V are from $N(Z_{13})$. Take S = {0, 1, I, 1 + 12I, I + 12, 10I + 3, 3I + 10} $\subseteq N(Z_{13})$. V is a set neutrosophic modulo integer vector spaces over the set S $\subseteq N(Z_{13})$.

It can also so happen that S is also a collection of neutrosophic matrices.

In view of this we define the new concept of set neutrosophic modulo integer matrix vector spaces and set neutrosophic modulo integer linear algebras.

**DEFINITION 2.4.24:** *Let V = {Any collection of m × n matrices with entries from N(Z)} and S $\subseteq$ N(Z); if V is such that for every M $\in$ V and s $\in$ S, sM and Ms $\in$ V then we define V to be a set neutrosophic integer matrix vector space over the set S. ( m and n can vary depending on S).*

We will exhibit this by some examples.

*Example 2.4.48:* Let

$$V = \left\{ \begin{pmatrix} 3I & 0 \\ 4I & 7I \end{pmatrix}, \begin{pmatrix} 14I & 7I \\ 4I & 2I \end{pmatrix}, \begin{pmatrix} I & -42I \\ 23I & 55I \end{pmatrix}, \begin{pmatrix} 16I & 12I \\ -98I & 9I \end{pmatrix}, \right.$$

$$(2I, 4I, 7I), (0\ 0\ 7I\ 14I), (0\ 0\ 0\ 11I\ 14I\ -I),$$

$$\left. \begin{pmatrix} 2I & 0 & -3I \\ -4I & I & -14I \\ 0 & 5I & -7I \end{pmatrix}, \begin{pmatrix} 16I & 2I & 0 & 0 \\ -4I & -I & I & -2I \\ -33I & 11I & 0 & -7I \end{pmatrix} \right\}$$

the entries of the matrices are from N(Z). Take S = {1, I} $\subseteq$ N(Z), V is a set neutrosophic integer matrix vector space; which we can also call as set neutrosophic integer mixed vector space over S.

*Example 2.4.49:* Let

$$V = \left\{ \begin{pmatrix} I & 4I & 0 \\ 2 & 3I & 0 \\ 1 & 0 & 5I \end{pmatrix}, \begin{pmatrix} 0 & 0 & 0 \\ 0 & 0 & 0 \\ 0 & 0 & 0 \end{pmatrix}, \begin{pmatrix} I & 0 & 9I \\ 7 & 4I & 0 \\ 1 & 0 & 2 \end{pmatrix}, \right.$$

$$\left. \begin{pmatrix} 0 & I & 2 \\ 3 & 4 & 5 \\ 6 & 7 & 8I \end{pmatrix}, \begin{pmatrix} 110I & 43I & 470 \\ 81 & 73I & 1 \\ 0 & 0 & 97 \end{pmatrix} \right\}$$

be the set of $3 \times 3$ matrices with entries from N(Z). Let

$$S = \left\{ \begin{pmatrix} 1 & 0 & 0 \\ 0 & 1 & 0 \\ 0 & 0 & 1 \end{pmatrix}, \begin{pmatrix} I & 0 & 0 \\ 0 & I & 0 \\ 0 & 0 & I \end{pmatrix}, \begin{pmatrix} 0 & 0 & 0 \\ 0 & 0 & 0 \\ 0 & 0 & 0 \end{pmatrix}, 0, 1 \right\}.$$

V is also a set neutrosophic integer matrix vector space over S which is clearly not a mixed one.

On similar lines we can define set neutrosophic rational matrix vector space over a set S or set neutrosophic real matrix vector space over S or a set neutrosophic complex matrix vector space over S or set neutrosophic modulo integer matrix vector space over S. Here we do not demand S to be a collection of neutrosophic numbers it can be matrices or row vectors.

We shall illustrate each one of these by examples.

*Example 2.4.50:* Let

$$V = \left\{ \left( \frac{5}{7}, 0, \frac{9I}{17}, 0 \right), \left( 0 - 9I, 25, \frac{41I}{7}, \frac{15I}{8} \right), \right.$$

$$(19I, 225, 0, I, -47, 5521),$$

$$\begin{pmatrix} 0 & 2I \\ 1 & 0 \\ I & 1 \end{pmatrix}, (0,0,0,0), (0,0,0,0,0), (0,0,0,0,0,0),$$

$$\begin{pmatrix} 7I & 2I/15 & 1/7 \\ 3I & 7/19I & 2I/11 \end{pmatrix},$$

$$\begin{pmatrix} 3I & 0 \\ I/11 & 2I/9 \end{pmatrix}, \begin{pmatrix} 0 & 0 \\ 0 & 0 \\ 0 & 0 \end{pmatrix}, \begin{pmatrix} 0 & 0 & 0 \\ 0 & 0 & 0 \end{pmatrix}, \begin{pmatrix} 0 & 0 \\ 0 & 0 \end{pmatrix} \right\}$$

be a collection of some matrices with entries from $N(Q) = \langle Q \cup I \rangle$. V is a set neutrosophic rational matrix vector space over the set $S = \{0, 1\} \subseteq N(Q)$. Infact V is a set neutrosophic rational mixed matrix vector space over S.

*Example 2.4.51:* Let

$$V = \left\{ \begin{pmatrix} \frac{3I}{7} & 0 & I & 0 \\ 0 & 4I & 0 & \frac{7I}{9} \\ -5I & 7I & 8I & -\frac{9I}{7} \\ \frac{4I}{5} & \frac{11I}{3} & \frac{12I}{11} & I \end{pmatrix}, \begin{pmatrix} 0 & 0 & 0 & 0 \\ 0 & 0 & 0 & 0 \\ 0 & 0 & 0 & 0 \\ 0 & 0 & 0 & 0 \end{pmatrix}, \right.$$

$$\begin{pmatrix} I & 0 & -I & \frac{4I}{7} \\ 7I & 8I & 0 & \frac{-12I}{5} \\ 14I & 0 & 9I & 0 \\ 0 & \frac{12I}{7} & 0 & \frac{140I}{23} \end{pmatrix}, \begin{pmatrix} 17I & 0 & 0 & 0 \\ 0 & 15I/4 & 0 & 0 \\ -5I & 7I & 8I & -9I/7 \\ \frac{-15I}{4} & 0 & 0 & \frac{-21I}{10} \end{pmatrix},$$

$$\left. \begin{pmatrix} -15I & 0 & \dfrac{71}{4} & 0 \\ 0 & 19I & 0 & \dfrac{5I}{14} \\ 0 & 0 & -11I & \dfrac{I}{7} \\ 0 & 0 & 0 & 21I \end{pmatrix} \right\}$$

be the 4 × 4 neutrosophic matrices with entries from N(Q). Take

$$S = \left\{ 0, 1, I, \begin{pmatrix} 1 & 0 & 0 & 0 \\ 0 & 1 & 0 & 0 \\ 0 & 0 & 1 & 0 \\ 0 & 0 & 0 & 1 \end{pmatrix}, \begin{pmatrix} I & 0 & 0 & 0 \\ 0 & I & 0 & 0 \\ 0 & 0 & I & 0 \\ 0 & 0 & 0 & I \end{pmatrix} \right\}.$$

V is a set neutrosophic rational matrix vector space which is clearly not a mixed one.

*Example 2.4.52:* Let

$$V = \{(\sqrt{2}\,I, 3I, 0, I), (\sqrt{17}\,I, \sqrt{13}\,I, \sqrt{14}\,I, 0, 0, 0, I), (0, 0, 0, 0), (0, 0, 0, 0, 0, 0, 0), (\sqrt{8}\,I, \sqrt{7}\,I, 0, I, \sqrt{5}\,I), (0, 0, 0, 0, 0),$$

$$(I, 6I, \sqrt{5}\,I, \sqrt{6}\,x, -\sqrt{3}\,I), \begin{pmatrix} I & \sqrt{2}I & 0 \\ 0 & -\sqrt{7}I & \sqrt{3}I \\ \sqrt{5}I & 6I & \sqrt{5}I \end{pmatrix},$$

$$\begin{pmatrix} \sqrt{10}I & 3I \\ 0 & -\sqrt{14}I \end{pmatrix}, \begin{pmatrix} 3I & \sqrt{2}I \\ 0 & \sqrt{7}I \\ \sqrt{91}I & \sqrt{5}I \end{pmatrix} \}$$

be a set of neutrosophic matrices with entries from the set N(R). Take $S = \{I, 1, 1 - I, \sqrt{3} - \sqrt{3}\,I\} \subseteq N(R)$. Clearly V is a set neutrosophic real matrix vector space over the set S.

Clearly V is a set neutrosophic real mixed matrix vector space over the set S.

*Example 2.4.53:* Let

$$V = \left\{ \begin{pmatrix} \sqrt{3}I & 0 & 0 \\ \sqrt{7}I & \sqrt{5}I & 0 \\ \sqrt{11}I & I & -10I \end{pmatrix}, \begin{pmatrix} \sqrt{2}I & 0 & \sqrt{7}I \\ -8I & \sqrt{3}I & 0 \\ -14I & 0 & \sqrt{17}I \end{pmatrix}, \right.$$

$$\left. \begin{pmatrix} 0 & 0 & 0 \\ 0 & 0 & 0 \\ 0 & 0 & 0 \end{pmatrix}, \begin{pmatrix} I & 4I & I \\ \sqrt{3}I & I & -I \\ I & \sqrt{5}I & \sqrt{7}I \end{pmatrix}, \begin{pmatrix} I & \sqrt{15}I & -4I \\ 0 & \sqrt{5}I & \sqrt{13}I \\ 0 & 0 & \sqrt{26}I \end{pmatrix} \right\}$$

be the set of $3 \times 3$ set neutrosophic matrix from the set $N(R)$. Take

$$S = \left\{ \begin{pmatrix} 0 & 0 & 0 \\ 0 & 0 & 0 \\ 0 & 0 & 0 \end{pmatrix}, \begin{pmatrix} I & 0 & 0 \\ 0 & I & 0 \\ 0 & 0 & I \end{pmatrix}, \begin{pmatrix} 1 & 0 & 0 \\ 0 & 1 & 0 \\ 0 & 0 & 1 \end{pmatrix}, \right.$$

$$\left. \begin{pmatrix} 1-I & 0 & 0 \\ 0 & 1-I & 0 \\ 0 & 0 & 1-I \end{pmatrix}, 0, I, 1 \right\}.$$

Clearly V is set neutrosophic real matrix vector space over S.

*Example 2.4.54:* Let

V = {(3 + 4i, a + bI, 1 – I, 5i + 2I, (6 + 5i) + (3 – 4i)I, 0),

$$\begin{pmatrix} 0 & 0 \\ 0 & 0 \end{pmatrix}, (c+dI, I, e+fI, 0, 7+4i), \begin{pmatrix} a+bI & 0 \\ 0 & c+dI \end{pmatrix}, (0, 0, 0,$$

$$0, 0), \begin{pmatrix} 0 & 0 & 0 & 0 \\ 0 & 0 & 0 & 0 \end{pmatrix}, \begin{pmatrix} 3+4i & x+yI & 0 & 10+I \\ 14I & 10-51U & -17I & 0 \end{pmatrix}\}$$

be a set of mixed matrices with entries from N(C). Take S = {0, 1} ⊆ N(C). V is a set neutrosophic complex matrix vector space over the set S. Infact V is a set neutrosophic complex mixed matrix vector space over the set S.

*Example 2.4.55:* Let

$$V = \left\{ \begin{pmatrix} (5+3i)I & 0 \\ 0 & 7I \end{pmatrix}, \begin{pmatrix} 0 & 0 \\ 0 & 0 \end{pmatrix}, \begin{pmatrix} (12+7i)I & 17I \\ 0 & (3-5i)I \end{pmatrix}, \right.$$

$$\left. \begin{pmatrix} 3I & 0 \\ 0 & (7+i)I \end{pmatrix}, \begin{pmatrix} 5I & 4I \\ (5-i)I & I \end{pmatrix}, \begin{pmatrix} 7iI & 18iI \\ 14I & (-3+4i)I \end{pmatrix} \right\}$$

be a set of 2 × 2 neutrosophic matrices from the set N(C). Take

$$S = \left\{ 0, 1, I, \begin{pmatrix} 0 & 0 \\ 0 & 0 \end{pmatrix}, \begin{pmatrix} 1 & 0 \\ 0 & 1 \end{pmatrix}, \begin{pmatrix} I & 0 \\ 0 & I \end{pmatrix} \right\};$$

clearly V is a set neutrosophic complex matrix vector space over the set S.

Now we proceed on to describe the set neutrosophic modulo integer vector spaces.

*Example 2.4.56:* Let

$$V = \{(I, 3I, 10I), (0, 4I, 0, 0, 5I), (I, I, I),$$

$$(9I, 5I, 3I, 4I, 6I, I, 0), \begin{bmatrix} 3I & 0 \\ I & 2I \end{bmatrix},$$

$$\begin{bmatrix} 3I & I & 3I & I \\ I & 3I & I & 3I \end{bmatrix}, \begin{bmatrix} 2I & I \\ 0 & 3I \\ 4I & 5I \\ 0 & 7I \end{bmatrix}$$

be a collection of mixed neutrosophic matrices with entries from the set $N(Z_{11})$. Take $S = \{1, I, 1 + 10I, 10 + I\} \subseteq N(Z_{11})$. V is a set neutrosophic modulo integer mixed matrix vector space over the set S.

*Example 2.4.57:* Let

$$V = \left\{ \begin{bmatrix} 3I & 4I \\ 6I & 2I \\ I & 5I \\ 2I & I \end{bmatrix}, \begin{bmatrix} 4I & 3I \\ I & 5I \\ 6I & 2I \\ 5I & 6I \end{bmatrix}, \begin{bmatrix} 0 & 0 \\ 0 & 0 \\ 0 & 0 \\ 0 & 0 \end{bmatrix}, \begin{bmatrix} 6I & I \\ 5I & 4I \\ 2I & 3I \\ 4I & 2I \end{bmatrix}, \right.$$

$$\left. \begin{bmatrix} I & 6I \\ 2I & 3I \\ 5I & 4I \\ 3I & 5I \end{bmatrix}, \begin{bmatrix} 5I & 2I \\ 3I & I \\ 4I & 6I \\ I & 4I \end{bmatrix}, \begin{bmatrix} 2I & 5I \\ 4I & 6I \\ 3I & I \\ 6I & 3I \end{bmatrix} \right\}$$

be a $4 \times 2$ neutrosophic matrices from the set $N(Z_7)$. Take $S = \{0, 1, I, 1 + 6I, I + 6, 2I + 5, 5 + 2I\} \subseteq N(Z_7)$. V is a set neutrosophic modulo integer matrix vector space over S. Clearly V is not a mixed one.

Now we proceed onto define the simple notion of set neutrosophic real (integer, modulo integer, rational, complex) matrix linear algebra and set neutrosophic real (integer, modulo integer, rational, complex) mixed matrix linear algebra does not exist as in case of vector spaces.

We will define for one type of set say reals and the same definition holds good for all types of sets.

**DEFINITION 2.4.25:** *Let V be a collection of n × m real neutrosophic matrices $1 \leq m, n < \infty$ with entries from N(R) once m and n are chosen they are fixed such that V is closed under matrix addition. Let $S \subseteq N(R)$. If for every $v \in V$ and $s \in S$, vs, $sv \in V$ then we call V to be a set neutrosophic real matrix linear algebra over the set S.*

On similar lines, integer, rational, complex and modulo integer linear algebras can be defined.

We will illustrate these types by some examples.

*Example 2.4.58:* Let

$$V = \left\{ \begin{pmatrix} mI & nI \\ tI & sI \end{pmatrix} \middle| m, n, s, t \in Z^+ \right\}$$

be a set of 2 × 2 neutrosophic matrices with real entries. Take

$$S = \left\{ (1, I, \begin{pmatrix} 1 & 0 \\ 0 & 1 \end{pmatrix}, \begin{pmatrix} I & 0 \\ 0 & I \end{pmatrix} \right\}.$$

V is a set neutrosophic integer matrix linear algebra over S.

*Example 2.4.59:* Let

$$V = \left\{ \begin{pmatrix} 0 & 0 \\ 0 & 0 \\ 0 & 0 \\ 0 & 0 \end{pmatrix}, \begin{pmatrix} nI & nI \\ nI & nI \\ nI & nI \\ nI & nI \end{pmatrix} \middle| n \in Z \right\}$$

be a collection of all 4 × 2 neutrosophic matrices with entries from the set of neutrosophic integers N(Z). Take S = {0, 1, I, 5I, 4, 8, 43, 7I, – 14I, 47I, – 2I} ⊆ N(Z).

Clearly V is a set neutrosophic integer matrix linear algebra over the set S.

*Example 2.4.60:* Let V = {(m ± mI, m ± mI, m ± mI, …, m ± mI) | m ∈ $Q^+$ ∪ {0}} be a 1 × n neutrosophic rational row vector. Take S = {0, $Q^+$ I, $Q^+$} ⊆ N(Q). V is a set neutrosophic rational matrix linear algebra over the set S.

*Example 2.4.61:* Let

$$V = \left\{ \begin{pmatrix} mI & 0 & 0 \\ mI & mI & 0 \\ mI & mI & mI \end{pmatrix} \middle| m \in Q^+ \cup \{0\} \right\}$$

be the collection of 3 × 3 lower triangular neutrosophic matrices with entries from the set N(Q). Set

$$S = \left\{ 0, 1, \begin{pmatrix} I & 0 & 0 \\ 0 & I & 0 \\ 0 & 0 & I \end{pmatrix}, \begin{pmatrix} 1 & 0 & 0 \\ 0 & 1 & 0 \\ 0 & 0 & 1 \end{pmatrix}, I \right\}$$

V is a set neutrosophic rational matrix linear algebra over set S.

*Example 2.4.62:* Let

$$V = \left\{ \begin{pmatrix} mI \\ mI \\ mI \\ mI \\ mI \\ mI \end{pmatrix} \middle| m \in R^+ \cup \{0\} \right\}$$

be the collection of $5 \times 1$ real neutrosophic column matrix. Let $S = \{0, 1, I, R^+, N(R^+)\} \subseteq N(R)$. V is a set neutrosophic real matrix linear algebra over the set S.

*Example 2.4.63:* Let

$$V = \left\{ \begin{pmatrix} mI & mI \\ 0 & 0 \end{pmatrix}, \begin{pmatrix} 0 & 0 \\ mI & mI \end{pmatrix}, \begin{pmatrix} mI & mI \\ nI & nI \end{pmatrix} \middle| m, n \in R^+ \cup \{0\} \right\}$$

be the special collection of $2 \times 2$ real neutrosophic matrices. Take

$$S = \left\{ 0, 1, I, (1-I), (r-rI) \mid r \in R^+, \begin{pmatrix} 1 & 0 \\ 0 & 1 \end{pmatrix}, \begin{pmatrix} I & 0 \\ 0 & I \end{pmatrix} \right\}.$$

We see V is a set neutrosophic real matrix linear algebra over the set S.

Next we proceed onto give some examples of set neutrosophic modulo integer matrix linear algebra.

*Example 2.4.64:* Let

$$V = \left\{ \begin{bmatrix} I & 2I \\ 3I & 4I \\ 5I & 6I \\ 7I & 8I \end{bmatrix}, \begin{bmatrix} 2I & 4I \\ 6I & 8I \\ I & 3I \\ 5I & 7I \end{bmatrix}, \begin{bmatrix} 3I & 6I \\ 0 & 3I \\ 6I & 0 \\ 3I & 6I \end{bmatrix}, \begin{bmatrix} 4I & 8I \\ 3I & 7I \\ 2I & 6I \\ I & 5I \end{bmatrix}, \begin{bmatrix} 5I & I \\ 6I & 2I \\ 7I & 3I \\ 8I & 4I \end{bmatrix}, \right.$$

$$\left. \begin{bmatrix} 6I & 3I \\ 0 & 6I \\ 3I & 0 \\ 6I & 3I \end{bmatrix}, \begin{bmatrix} 7I & 5I \\ 3I & I \\ 8I & 6I \\ 4I & 2I \end{bmatrix}, \begin{bmatrix} 8I & 7I \\ 6I & 5I \\ 4I & 3I \\ 2I & I \end{bmatrix}, \begin{bmatrix} 0 & 0 \\ 0 & 0 \\ 0 & 0 \\ 0 & 0 \end{bmatrix} \right\}$$

be a special collection of neutrosophic matrices with entries $N(Z_9)$. Take $S = \{0, 1, I, 1 + 8I, 8 + I, 3 + 6I, 6 + 3I\} \subseteq N(Z_9)$. V is a set neutrosophic modulo integer matrix linear algebra over the set S.

*Example 2.4.65:* Let V = {(2I, 3I, 4I, 6I, 7I, 10I, I), (4I, 6I, 8I, I, 3I, 9I, 2I), (6I, 9I, I, 7I, 10I, 8I, 3I), (8I, I, 5I, 2I, 6I, 7I, 4I), (10I, 4I, 9I, 8I, 2I, 6I, 5I), (I, 7I, 2I, 3I, 9I, 5I, 6I), (3I, 10I, 6I, 9I, 5I, 4I, 7I), (5I, 2I, 10I, 4I, I, 3I, 8I), (7I, 5I, 3I, 10I, 8I, 2I, 9I) (9I, 8I, 7I, 5I, 4I, I, 10I), (0, 0, 0, 0, 0, 0, 0) be a $1 \times 7$ neutrosophic row vector with entries from $N(Z_{11})$}. Take $S = \{0, 1, I, 1 + 8I, 8 + I, 3 + 8I, 8 + 3I\} \subseteq N(Z_{11})$. V is a set neutrosophic modulo integer matrix linear algebra over the set S.

Now we give some examples of set neutrosophic complex matrix linear algebra over a set $S \subseteq N(C)$.

*Example 2.4.66:* Let V = {(mI, mI, ..., mI) | m ∈ C} be a $1 \times n$ complex neutrosophic row vector with entries from the set N(C). Take $S = \{0, 1, I\} \subseteq N(C)$. V is a set neutrosophic complex matrix linear algebra over the set S.

*Example 2.4.67:* Let

$$V = \left\{ \begin{pmatrix} aI & 0 & aI \\ 0 & aI & 0 \end{pmatrix}, \begin{pmatrix} 0 & aI & 0 \\ aI & 0 & aI \end{pmatrix}, \begin{pmatrix} aI & bI & aI \\ bI & aI & bI \end{pmatrix} \middle| a, b \in C \right\}$$

be a collection of $2 \times 3$ neutrosophic complex matrices with entries from N(C). Let $S = \{0, I, 1, (1 - I)\} \subseteq N(C)$. Clearly V is a set neutrosophic complex matrix linear algebra over the set S.
 Having defined set neutrosophic matrix vector spaces and linear algebras we can define their substructures, linear transformations and linear operators analogously. This can be taken up as a simple exercise by the reader.

These matrix structures have applications in fields which include economic models, and neutrosophic bidirectional associative memories.

Next we proceed onto define the notion of set neutrosophic integer polynomial vector spaces and linear algebras.

**DEFINITION 2.4.26:** *Let N(Z) denote the set of neutrosophic integers. Let*

$$N(Z)[x] = \left\{ \sum_{i=0}^{m} n_i x^i \,\middle|\, n_i \in N(Z) \right\}$$

*and x an indeterminate and $m \in Z^+ \cup \{0\}$, N(Z)[x] denotes the neutrosophic integer coefficient polynomials in the variable x. Similarly N(Q)[x] denotes the neutrosophic rational coefficient polynomials in the variable x and N(R)[x], the neutrosophic real coefficient polynomials in the variable x, N(C)[x] the neutrosophic complex coefficient polynomial in the variable x and $N(Z_n)[x]$ the neutrosophic modulo integer coefficient polynomial in the variable x.*

Further we see
$$N(Z)[x] \subseteq N(Q)[x] \subseteq N(R)[x] \subseteq N(C)[x]$$

the containment is strict; we will define only for one; viz. neutrosophic integers coefficient polynomials vector space and linear algebra. The reader can take the simple exercise of defining using other coefficients. However we will give examples for all cases which will make the situation simple and easy to understand.

**DEFINITION 2.4.27:** *Let $V \subseteq N(Z)[x]$, and $S \subseteq N(Z)$. We say V is a set neutrosophic integer coefficient polynomial vector space in the variable x if sv and $vs \in V$ for every $v \in V$ and $s \in S$.*

We will illustrate this by some simple examples.

*Example 2.4.68:* Let

$$V = \{(3 + I)x^3 + (7 - 3I)x^2 - 5Ix + 12 - 17I,$$
$$Ix^6 + (1 - I)x^5 - 7Ix^3 + (15 - 21I)x^2 - 10Ix + 3, 0, Ix^7, 21x^3,$$

$(3 + 5I)x^{21}, 27, –3 + 29I, 49I\} \subseteq N(Z)[x]$.

Take $S = \{0, 1\} \subseteq N(Z)$. V is a set neutrosophic integer coefficient polynomial vector space over the set S.

The main advantage of using these structures is when we get a solution for certain equations which is not in the vector space V we can include it in V arbitrarily provided its scalar multiplication with S alone is compatible and nothing more. This sort of flexibility cannot be enjoyed by any of the algebraic structures. Hence these structures have advantage over them.

*Example 2.4.69:* Let $V = \{Z^+I[x]$; that is $Z^+I[x]$ consists of all polynomials in the variable x with coefficients from $Z^+I\} \subseteq N(Z)$.
   Take $S = 3Z^+I \subseteq N(Z)$. V is a set neutrosophic integer coefficient polynomials vector space over the set S.

The reader is given the task of defining set neutrosophic rational coefficient (real coefficient, complex coefficient and modulo integer coefficient) vector spaces over a suitable set S analogously.

We will give examples of these four types of vector spaces.

*Example 2.4.70:* Let $V = \{Q^+I[x] \cup \{0\}\} \subseteq N(Q)[x]$. Take $S = \{Q^+, Q^+I\} \subseteq N(Q)$. V is a set neutrosophic rational coefficients polynomial vector space over the set S.

*Example 2.4.71:* Let

$$V = \left\{\sum_{i=0}^{n}(m_i + m_iI)x^i \,\middle|\, m_i \in Q\right\} \subseteq N(Q)[x].$$

Take $S = Q^+ \subseteq N(Q)$. V is a set neutrosophic rational coefficient polynomial vector space over the set S.

*Example 2.4.72:* Let

$$V = \{(\frac{7}{5} - 8I)x^3, (\frac{3}{7} + \frac{14}{9}II)x^{14},$$
$$(\frac{27}{13} + \frac{43}{7}I)x, \frac{801}{7}Ix^7, 43x^6, 2011Ix^{12},$$
$$27I + 3(3+7I)x^8 - 4Ix^7 +$$
$$48x^5 - (50-I)x^3 + Ix^2 - (27-I)x + 48\} \subseteq N(Q)[x].$$

Take $S = \{0, 1\} \subseteq N(Q)$, V is a set neutrosophic rational coefficient polynomials vector space over the set S.

*Example 2.4.73:* Let

$$V = \left\{\sum_{i=1}^{n}(a_i + a_iI)x^i \,\middle|\, a_i \in C\right\} \subseteq N(C)[x].$$

Set $S = \{0, 1\} \subseteq N(C)$. V is a set neutrosophic complex coefficient polynomial vector space over S.

*Example 2.4.74:* Let $V = \{0, 1, (5 - 3I)x, [(5 + 2i) + (14 + 5i)I]x^3 - [(-3 + 2i) + (4 - i)I]x^2 + [(11 + 8i) + (3 - 11i)I]x + (3 - i) + 17 + 4i, [(11 + 4i) + (21 - 5i)I]x^{27}, 28, 4I, (11 + 48i)I\} \subseteq [N(C)[x]$. Set $S = \{0, 1\} \subseteq N(C)$. V is a set neutrosophic complex coefficient polynomial vector space over S.

*Example 2.4.75:* Let

$$V = \left\{\sum_{i=1}^{n}(a_i + a_iI)x^i \,\middle|\, n \in I \text{ and } a_i \in R^+\right\} \subseteq N(R)[x].$$

Set $S = \{0, 1, I, b - bI \mid b \in R^+\} \subseteq N(R)$. V is a set neutrosophic real coefficient polynomial vector space over the set S.

*Example 2.4.76:* Let

$$V = \{a_0x + a_1x^2 + \ldots + a_nx^n, 0, 1,$$

$$\sqrt{5}\,x + 17Ix^2 + (3I + 27)x^3 - (4I - \sqrt{3}\,)\,x^4$$
$$+ 48I\,x^5 - 27,\,x^3,\,\sqrt{21}\,x^7,\,(48I - 4)\,x^3$$
$$+ (\sqrt{41} - 7I)\,x^2 - \sqrt{20} \text{ where } a_i \in R^+ \} \subseteq N(R)\,[x].$$

Take $S = \{0, 1\} \subseteq N(R)$. V is a set neutrosophic real coefficient polynomial vector space over S.

*Example 2.4.77:* Let

$$V = \{0, I + x, Ix + 1,$$
$$2x + Ix^2 + 3x + 3Ix\,2x - 3Ix\,I + 1 +$$
$$(3 + 2I)x + 3I - 2)x^2 + 3x^3 + Ix^4 + x^5\}$$
$$\subseteq N(Z_4)\,[x].$$

Set $S = \{0, 1\} \subseteq N(Z_4)$. V is a neutrosophic modulo integer coefficient polynomial vector space over the set S.

*Example 2.4.78:* Let

$$V = \left\{ \sum_{i=0}^{9} (m_i - m_i I)x^i \,\middle|\, m_i \in Z_5 \right\}.$$

Let $S = \{0, 1, I, 1 + 4I, I + 4\} \subseteq N(Z_5)$. V is a set neutrosophic modulo integer coefficient polynomial vector space over the set S.

Now we proceed onto define the notion of set neutrosophic polynomial linear algebra over a set S.

**DEFINITION 2.4.28:** *Let $V \subseteq N(Z)[x]$ and $S \subseteq N(Z)$. We say V is a set neutrosophic integer coefficient polynomial linear algebra over the set S if V is a set neutrosophic integer coefficient polynomial vector space and V is a semigroup with respect to addition; that is for $a, b \in V$, $a + b$ and $b + a \in V$.*

We can define analogously set neutrosophic real coefficient or rational coefficient or complex coefficient or modulo integer coefficient linear algebra over a set S.

This task is left for the reader. However to make the concept clear we give examples of all these types of set neutrosophic polynomial linear algebras.

*Example 2.4.79:* Let

$$V = \left\{ \sum_{i=0}^{n}(m_i - m_i I)x^i \,\middle|\, m_i \in Z^+, n \in N \right\} \subseteq N(Z)[x].$$

Set $S = \{0, 1, I, m - mI\} \subseteq N(Z)$. V is a set neutrosophic integer coefficient polynomial linear algebra over the set S.

*Example 2.4.80:* Let

$$V = \left\{ \sum_{i=0}^{n} m_i I x^i \,\middle|\, m_i \in Z^+ \cup \{0\} \right\} \subseteq N(Z)[x]$$

and let $S = \{3Z^+I, 13Z^+\} \subseteq N(Z)$. V is a set neutrosophic integer coefficient polynomial linear algebra over the set S.

*Example 2.4.81:* Let

$$V = \left\{ \sum_{i=1}^{m}(a_i - a_i I)x^i \,\middle|\, a_i \in Q^+ \cup \{0\} \right\} \subseteq N(Q)[x]$$

and $S = \{0, 1, 1 - I, 5/7 - 5/7I\} \subseteq N(Q)$. V is a set neutrosophic rational coefficient polynomial linear algebra over the set S.

*Example 2.4.82:* Let

$$V = \left\{ \sum_{i=0}^{n} a_i I x^i \,\middle|\, a_i \in Z^+ \cup \{0\} \right\} \subseteq N(Q)[x]$$

and $S = \{0, 1, 3Z^+ I\} \subseteq N(Q)$. V is a set neutrosophic rational coefficient polynomial linear algebra over the set S.

*Example 2.4.83:* Let

$$V = \left\{ \sum_{i=0}^{n} \left( \frac{1}{\left(\sqrt{2}\right)^{n_i}} - \frac{I}{\left(\sqrt{2}\right)^{n_i}} \right) x^i \,\bigg|\, n_i \in Z^+ \cup \{0\} \right\} \subseteq N(R)\,[x]$$

and let $S = \{m - mI \mid m \in R^+\} \subseteq N(R)$. V is a set neutrosophic real coefficient polynomial linear algebra over the set S.

*Example 2.4.84:* Let

$$V = \left\{ \sum_{i=0}^{n} \left( \frac{\sqrt{p}}{\sqrt{q}} \right)^i Ix^i \,\bigg|\, \sqrt{p} \in R^+ \cup \{0\}; q \in R^+ \right\} \subseteq N(R)[x]$$

and $S = \{pI \mid p \in R^+ \cup \{0\}\} \subseteq N(R)$. V is a set neutrosophic real coefficient linear algebra over the set S.

*Example 2.4.85:* Let

$$V = \left\{ \sum_{i=0}^{n} (m - mI) x^i \,\bigg|\, m \in Z_7 \right\} \subseteq N(Z_7)\,[x]$$

and $S = \{Z_7 I\} \subseteq N(Z_7)$. V is a set neutrosophic modulo integer coefficient polynomial linear algebra over the set S.

*Example 2.4.86:* Let

$$V = \left\{ \sum_{i=0}^{n} mIx^i \,\bigg|\, m \in Z_{23} \right\} \subseteq N(Z_{23})\,[x]$$

and S = {– mI + m | m ∈ $Z_{23}$} ⊆ N($Z_{23}$). V is a set neutrosophic modulo integer coefficient polynomial linear algebra over S.

The interested reader is expected to define the properties like substructures, pseudo substructures, linear transformations and linear operators, analogously for set neutrosophic polynomial linear algebras of all types.

**Chapter Three**

# NEUTROSOPHIC SEMIGROUP LINEAR ALGEBRA

In this chapter we for the first time introduce the notion of neutrosophic semigroup vector space, neutrosophic group vector space and the analogous neutrosophic linear algebras and describe a few of its properties. This chapter has two sections. Section one introduces the notion of neutrosophic semigroup vector spaces and section two introduces the notion of neutrosophic group vector spaces and their properties. Notions about neutrosophic semigroups and neutrosophic groups are given in chapter one of this book.

## 3.1 Neutrosophic Semigroup Linear Algebras

In this section we introduce the notion of neutrosophic semigroup linear algebras and neutrosophic semigroup vector spaces. Several interesting properties about them are derived.

**DEFINITION 3.1.1:** *Let V be a neutrosophic set of integers (V ⊆ N(Z)). S any additive semigroup with 0. We call V to be a semigroup neutrosophic vector space or neutrosophic semigroup vector space over S if the following conditions hold:*

*(1) s v ∈ V for all s ∈ S and v ∈ V.*
*(2) 0 . v = 0 ∈ V for all v ∈ V and 0 ∈ S; 0 is the zero vector.*
*(3) $(s_1 + s_2) v = s_1 v + s_2 v$ for all $s_1, s_2 ∈ S$ and $v ∈ V$.*

*Note:* Even if S is just a semigroup without zero then also V is a neutrosophic semigroup vector space. The condition (2) of the definition will become superfluous.

We will illustrate this situation by some examples.

*Example 3.1.1:* Let V = {1-I, 0, 25I, 37, 8 + 8I, 47I, 52 – 3I, 46 + 23I, 3} ⊆ N(Z). Take S = {0, 1, 1 – I} ⊆ N(Z). Clearly S is not a semigroup neutrosophic vector space over S.

Thus we see here in example 3.1.1, S is not a semigroup.

*Example 3.1.2:* Let V = {$3Z^+ I$} ⊆ N(Z) and S {$Z^+ ∪ \{0\}$}, a semigroup under addition with 0. V is a neutrosophic semigroup vector space over the semigroup S.

*Example 3.1.3:* Let V = {3ZI, 5ZI, 2ZI} ⊆ N(Z) and S = Z the semigroup under addition. V is a neutrosophic semigroup vector space over the semigroup S.

*Example 3.1.4:* Let V = {m – mI | m ∈ $Q^+ ∪ \{0\}$} and S = {QI} be the semigroup under addition. V is a neutrosophic semigroup vector space over the semigroup S.

*Example 3.1.5:* Let V = {QI} ⊆ N(Q) and S = {m – mI | m ∈ $Q^+ ∪ \{0\}$} ⊆ N(Q). V is a neutrosophic semigroup vector space over the semigroup S.

***Example 3.1.6:*** Let $V = \{RI\} \subseteq N(R)$ and $S = R^+ \cup \{0\} \subseteq R$ the semigroup under addition. V is a neutrosophic semigroup vector space over the semigroup S.

***Example 3.1.7:*** Let $V = \{m - mI \mid m \in R^+ \cup \{0\}\}$ and $S = Q^+ \cup \{0\}$ the semigroup under addition. V is a neutrosophic semigroup vector space over the semigroup S or semigroup neutrosophic vector space over the semigroup S.

***Example 3.1.8:*** Let $V = \{0, 1 + 4I, I + 4, 2 + 3I, 3 + 2I\} \subseteq N(Z_5)$ and $S = \{Z_5I\} \subseteq N(Z_5)$ be a semigroup, V is a neutrosophic semigroup vector space over the semigroup S.

***Example 3.1.9:*** Let $V = \{Q^+I, ZI\} \subseteq N(R)$ and $S = Z$ the semigroup under addition. V is a neutrosophic semigroup vector space over the semigroup S.

Now we proceed onto define substructures.

**DEFINITION 3.1.2:** *Let $V \subseteq N(R)$ and $S \subseteq N(R)$ where S is a semigroup under addition such that V is a neutrosophic semigroup vector space over the semigroup S.*
  *Suppose $W \subseteq V$ is a proper subset of V such that W is itself a neutrosophic semigroup vector space over the semigroup S then we call W to be a neutrosophic semigroup vector subspace of V or neutrosophic semigroup subvector space of V over the semigroup S.*

We will illustrate this situation by some examples.

***Example 3.1.10:*** Let $V = \{Q^+I, ZI\} \subseteq N(Q)$ and $S = \{Z^+I\} \subseteq N(Q)$, a semigroup under addition. V is a neutrosophic semigroup vector space over S. Take $W = \{Z^+ I\} \subseteq V$.

It is easily verified W is a neutrosophic semigroup subvector space of V over S.

***Example 3.1.11:*** Let

$$V = \left\{ \begin{pmatrix} 3ZI & 0 \\ 0 & 2ZI \end{pmatrix}, m - mI, 0 \,\middle|\, m \in Q^+ \right\}$$

and 2 × 2 matrices with diagonal elements from 3ZI and 2ZI} and $S = Z^+I \cup \{0\}$ be a semigroup under addition. V is a neutrosophic semigroup vector space over S. Take W = {m – mI, 0 | m ∈ $Q^+$} ⊆ V. W is a neutrosophic semigroup vector subspace of V over S.

**DEFINITION 3.1.3:** *Let V be a neutrosophic semigroup vector space over the semigroup S. Suppose W ⊆ V be a proper subset of V and T ⊆ S be a proper subsemigroup of S with | T | > 1 such that W is a neutrosophic semigroup vector space over the semigroup S then we call W to be a neutrosophic subsemigroup vector subspace of V over the subsemigroup T of S.*

We will illustrate this situation by some examples.

***Example 3.1.12:*** Let V = {ZI, $R^+I$} ⊆ N(R) and S = {$Z^+I$} be a neutrosophic semigroup vector space over S. Take W = {$R^+$ I} ⊆ V and T = {$3Z^+$ I} ⊆ {$Z^+$ I} ⊆ S be a subsemigroup of S. W is a neutrosophic subsemigroup vector subspace of V over the subsemigroup T of S.

***Example 3.1.13:*** Let V = {$R^+I$, QI, m – mI | m ∈ $Z^+$} ⊆ N(R) and S = {$Q^+I$} be a semigroup. Take W = {QI, m – mI | m ∈ $Z^+$} ⊆ V and T = {$Z^+$ I} ⊆ S. W is a neutrosophic subsemigroup vector subspace of V over the subsemigroup T of S.

We see if a neutrosophic semigroup vector space has no neutrosophic subsemigroup vector subspace over the semigroup S then we define them to be simple.

**DEFINITION 3.1.4:** *Let V be a neutrosophic semigroup vector space over the semigroup S. If V has no neutrosophic subsemigroup vector subspace over any subsemigroup T of the semigroup S then we call V to be a neutrosophic semigroup simple vector space over S.*

We will illustrate this situation by some examples.

*Example 3.1.14:* Let V = {(0 0 0), (1 I 0), (I I I), (0 I I ), (1 0 1), (I I 0), ( I I I 0), (0 0 0 I), (1 1 0 I), (1 I 0 1), (1 1 1 1), (0 0 0 0)} be a neutrosophic semigroup vector space over the semigroup S = {0, 1} where (1 + 1) = 1 . It is easily verified S has no subsemigroup T such that | T | > 1. So V is a neutrosophic semigroup simple vector space over S.

*Example 3.1.15:* Let

$$V = \left\{ \begin{pmatrix} 1 & 1 \\ 1 & I \end{pmatrix}, \begin{pmatrix} 0 & 1 \\ I & 0 \end{pmatrix}, \begin{pmatrix} I & I \\ 0 & 1 \end{pmatrix}, \begin{pmatrix} 1 & 1 & 1 \\ I & I & 0 \end{pmatrix}, \begin{pmatrix} 1 & 1 & 0 \\ 0 & I & I \end{pmatrix}, \right.$$

$$\left. \begin{pmatrix} 0 & 0 \\ 0 & 0 \end{pmatrix}, \begin{pmatrix} 0 & 0 & 0 \\ 0 & 0 & 0 \end{pmatrix}, \begin{pmatrix} 1 & 1 & 1 \\ 0 & 0 & 0 \end{pmatrix} \right\}$$

be a neutrosophic semigroup over the semigroup S = {0, 1} with 1 + 1 = 1. V is also a neutrosophic semigroup simple vector space over S.

*Example 3.1.16:* Let V = {(I I I), (0 0 0)} be a neutrosophic semigroup vector space over the semigroup S = {0, 1} with 1 + 1 = 1. V is a neutrosophic semigroup simple vector space over S.

**DEFINITION 3.1.5:** *Let V be a neutrosophic semigroup vector space over the semigroup S. If V itself is a neutrosophic semigroup under addition, then we call V to be a neutrosophic semigroup linear algebra over the semigroup S, if $s(v_1 + v_2) = sv_1 + sv_2$ for all $v_1, v_2$ in V and for all $s \in S$.*

We will illustrate this by some examples.

*Example 3.1.17:* Let V = {$Z^+I \cup \{0\}$} × {$Z^+I \cup \{0\}$} × {$Z^+I \cup \{0\}$} × {$Z^+I \cup \{0\}$}. Take S = {$3Z^+I \cup \{0\}$} $\subseteq$ N(Z). Both V and S are semigroups under addition. So V is a neutrosophic semigroup linear algebra over the semigroup S.

**Example 3.1.18:** Let

$$V = \left\{ \begin{pmatrix} a & b & c & g \\ d & e & f & h \end{pmatrix} \middle| a,b,c,d,e,f,g,h \in N(Q) \right\}$$

and $S = \{Q^+I \cup \{0\}\}$. Both V and S are semigroups. Thus V is a neutrosophic semigroup linear algebra over the semigroup S.

The following theorem is left as an exercise for the reader to prove.

**THEOREM 3.1.1:** *Every neutrosophic semigroup linear algebra is a neutrosophic semigroup vector space; but in general a neutrosophic semigroup vector space is not a neutrosophic semigroup linear algebra.*

Now we proceed onto define the substructures in neutrosophic semigroup linear algebras.

**DEFINITION 3.1.6:** *Let V be a neutrosophic semigroup linear algebra over the semigroup S.*
*Suppose $W \subseteq V$ be a proper subset of V and W is a neutrosophic semigroup linear algebra over S then we call W to be a neutrosophic semigroup linear subalgebra of V over the semigroup S.*

We will illustrate this situation by some simple examples.

**Example 3.1.19:** Let

$$V = \left\{ \begin{pmatrix} a_1 & a_2 & a_3 & a_4 \\ a_5 & a_6 & a_7 & a_8 \end{pmatrix} \middle| a_i \in N(Q); 1 \le i \le 8 \right\}$$

be a neutrosophic semigroup under matrix addition and $S = Z^+ I \cup \{0\} \subseteq N(Q)$ be a semigroup, V is a neutrosophic semigroup linear algebra over the semigroup S.
Take

$$W = \left\{ \begin{pmatrix} a_1 & a_2 & 0 & a_3 \\ a_4 & 0 & 0 & 0 \end{pmatrix} \middle| a_i \in N(Q); 1 \leq i \leq 4 \right\} \subseteq V;$$

W is a neutrosophic semigroup linear subalgebra of V over the semigorup S.

***Example 3.1.20:*** Let $V = \{(Q^+I \cup \{0\}) \times (QI) \times (Q^+I \cup \{0\}) \times QI\}$ be a neutrosophic semigroup under component wise addition. Take $S = (Q^+I \cup \{0\}) \subseteq N(Q)$, a semigroup. V is a neutrosophic semigroup linear algebra over the semigroup S.

Choose $W = \{(0) \times QI \times (Q^+I \cup \{0\}) \times \{0\}\} \subseteq V$, W is a neutrosophic semigroup linear subalgebra of V over the semigroup S.

**DEFINITION 3.1.7:** *Let V be a neutrosophic semigroup linear algebra over the semigroup S and $W \subseteq V$ be such that W is a proper subset of V and W is a neutrosophic semigroup linear algebra over the proper subsemigroup T of S, then we call W to be a neutrosophic subsemigroup linear subalgebra of V over the subsemigroup T of the semigroup S.*

We shall illustrate this by some examples.

***Example 3.1.21:*** Let

$$V = \left\{ \begin{pmatrix} a & b \\ c & d \end{pmatrix} \middle| a, b, c, d \in N(Q) \right\}$$

and $S = (Q^+I \cup \{0\})$ be neutrosophic semigroups. V is a neutrosophic semigroup linear algebra over the semigroup S. Take

$$W = \left\{ \begin{pmatrix} a & b \\ c & d \end{pmatrix} \middle| a,b,c,d \in Q^+I \cup \{0\} \right\} \subseteq V$$

and $T = (Z^+I \cup \{0\}) \subseteq S$. W is a neutrosophic subsemigroup linear subalgebra of V over the subsemigroup T of the semigroup S.

*Example 3.1.22:* Let $V = N(Q)[x]$ that is the collection of polynomials in the variable x with coefficients from $N(Q)$. Clearly V is a neutrosophic semigroup under addition. Take $S = (Q^+I \cup \{0\}) \subseteq N(Q)$, S is also a semigroup. We see V is a neutrosophic semigroup linear algebra over the semigroup S. Suppose

$$W = \left\{ \text{All polynomials in V of the form} \sum_{i=0}^{n} a_i x^{2i} \middle| a_i \in N(Q); n \in N \right\} \subseteq V.$$

W is also a neutrosophic subsemigroup of V. Take $T = (Z^+I \cup \{0\}) \subseteq S$. T is also a semigroup. We see W is a neutrosophic subsemigroup linear subalgebra of V over the subsemigroup T or the semigroup S.

We shall define simple neutrosophic semigroup linear algebra of neutrosophic semigroup simple linear algebras.

**DEFINITION 3.1.8:** *Let V be a neutrosophic semigroup linear algebra over the semigroup S. If V has no proper neutrosophic subsemigroup linear subalgebras then we define V to be a neutrosophic semigroup simple linear algebra (or simple neutrosophic semigroup linear algebra) over the semigroup S.*

We will illustrate this by some simple examples.

*Example 3.1.23:* Let
$$V = \left\{ \begin{pmatrix} a & b & e \\ c & d & f \end{pmatrix} \middle| a,b,c,d,e,f \in N(Z_3) \right\}$$

be a neutrosophic semigroup under matrix addition modulo 3. Take $S = Z_3$, a semigroup under addition. V is a simple neutrosophic semigroup linear algebra over the semigroup S. For $S = Z_3$ has no proper subsemigroups.

***Example 3.1.24:*** Let $V = \{N(Z_5) \times N(Z_5) \times N(Z_5) \times N(Z_5) \times N(Z_5)\}$ be a neutrosophic simple semigroup linear algebra over the semigroup $Z_5$.

It is easily verified V is a simple neutrosophic semigroup linear algebra over the semigroup $Z_5$.

**DEFINITION 3.1.9:** *Let V be a simple neutrosophic semigroup linear algebra over the semigroup S. If V has no neutrosophic semigroup linear subalgebras then we call V to be doubly simple neutrosophic semigroup linear algebra or neutrosophic semigroup doubly simple linear algebra.*

We shall illustrate these situations by some examples.

***Example 3.1.25:*** Let $V = \{N(Z_7) \times N(Z_7) \times N(Z_7) \times N(Z_7)\}$ be a neutrosophic semigroup linear algebra over the semigroup $S = Z_7$. V is a neutrosophic semigroup simple linear algebra over S but V is not a doubly simple neutrosophic semigroup linear algebra over S as $W = \{0\} \times N(Z_7) \times N(Z_7) \times \{0\} \subseteq V$ is a neutrosophic semigroup linear subalgebra of V over S.

***Example 3.1.26:*** Let $V = \{0, 1 + I\} \subseteq N(Z_2)$ be a neutrosophic semigroup linear algebra over $Z_2 = S$. V is a doubly simple neutrosophic semigroup linear algebra over $S = Z_2$.

***Example 3.1.27:*** Let $V = \{Z_7I\} \subseteq N(Z_7)$ and $S = Z_7$. V is doubly simple neutrosophic semigroup linear algebra over $S = Z_7$.

Now we proceed on to define for neutrosophic semigroup vector spaces we can define the basis and generating set.

**DEFINITION 3.1.10:** *Let V be a neutrosophic semigroup vector space over the semigroup S under addition. Let $T = \{v_1, ..., v_n\} \subseteq V$ be a subset of V; we say T generates the neutrosophic semigroup vector space V over S if every element $v \in V$ can be got as $v = s v_i$; $v_i \in T$ and $s \in S$.*

***Example 3.1.28:*** Let $V = \{3Z^+I \cup \{0\}\}$ be a neutrosophic semigroup vector space over the semigroup $Z^+ \cup \{0\}$. Take $T = \{3I\} \subseteq V$, T is generates V over S.

***Example 3.1.29:*** Let $V = Z_{25}I$ modulo integers 25. $S = \{0, 5, 10, 15, 20\} \subseteq Z_{25}$ is a semigroup under addition modulo 25. V is a neutrosophic semigroup vector space over the semigroup S.
  $T = \{1I, 2I, 3I, 4I, 6I, 7I, 8I, 9I, 11I, 12I, 13I, 14I, 16I, 17I, 18I, 19I, 21I, 22I, 23I, 24I\}$ is a generating set of V over S.

We will illustrate by some examples that the generating set of V is dependent in general on the semigroup over which V is defined.

***Example 3.1.30:*** Let $V = Z_{20}I$ be a neutrosophic semigroup vector space over the semigroup $S = \{0,10\}$. The generating set of V over S is $T = \{1I, 2I, 3I, 4I, 5I, 6I, 7I, 8I, 9I, 11I, 12I, 13I, 14I, 15I, 16I, 17I, 18I, 19I\}$.
  If we take $S_1 = \{0, 5, 10, 15\}$ to be the semigroup over which the same V is defined we see the generating set of V over $S_1$ is $T_1 = \{1I, 2I, 3I, 4I, 5I, 6I, 7I, 8I, 9I, 11I, 12I, 13I, 14I, 15I, 16I, 17I, 18I, 19I\}$. We see $T_1 \neq T$.
  Thus in general the generating set of a neutrosophic semigroup vector space is dependent on the semigroup over which it is defined.

Now we proceed onto define the generating set of a neutrosophic semigroup linear algebra for which we also need the concept of independent set.

**DEFINITION 3.1.11:** *Let V be a neutrosophic semigroup linear algebra over the semigroup S. $T = \{v_1, ..., v_n\} \subseteq V$ is an independent set if $v_i \neq s v_j$, $i \neq j$ for some $s \in S$ and*

$$v_k \neq \sum_{i=1}^{m} s_i v_i \ ;$$

$i \leq k \leq n$, $m < n$, and $s_i \in S$. We say T is a generating subset of V if T is a linearly independent set and every element $v \in V$ can be represented as

$$v = \sum_{i=1}^{n} s_i v_i \ ;$$

$s_i \in S$; $1 \leq i \leq n$.

We will illustrate this by some examples.

***Example 3.1.31:*** Let $V = Q^+I \cup \{0\}$ be a neutrosophic semigroup linear algebra over the semigroup $S = Z^+I \cup \{0\}$. V is an infinite generating set over the semigroup S.

***Example 3.1.32:*** Let $V = Q^+I \cup \{0\}$ be the neutrosophic semigroup linear algebra over the semigroup $S = Q^+I \cup \{0\}$. V is generated by the set $B = \{1\}$ over $S = Q^+I \cup \{0\}$.

**DEFINITION 3.1.12:** *Let V be a neutrosophic semigroup vector space over the semigroup S. Let $W \subseteq V$ be such that W is a neutrosophic semigroup linear algebra over S, then we call W to be a pseudo neutrosophic semigroup linear subalgebra of V over S.*

*If V has no pseudo neutrosophic semigroup linear subalgebras then we call V to be a pseudo simple neutrosophic semigroup vector space or pseudo neutrosophic semigroup simple vector space.*

We shall illustrate these situations by some simple examples.

***Example 3.1.33:*** Let

$$V = \left\{ \begin{pmatrix} a & b \\ c & d \end{pmatrix}, (a_1 \ a_2 \ a_3 \ a_4) \middle| a,b,c,d,a_i \in Q^+I \cup \{0\}; 1 \leq i \leq 4 \right\}$$

be a neutrosophic semigroup vector space over the semigroup $S = Z^+I \cup \{0\}$.

Take

$$W = \left\{ \begin{pmatrix} a & b \\ c & d \end{pmatrix} \middle| a_i \in Q^+I \cup \{0\}; 1 \leq i \leq 4 \right\} \subseteq V.$$

$W_1$ is also a pseudo neutrosophic linear subalgebra of V over the semigroup S.

***Example 3.1.34:*** Let $V = \{(Q^+I \cup \{0\})[x], m - mI \mid m \in Z^+\}$ be a neutrosophic semigroup vector space over the semigroup $S = Z^+I \cup \{0\}$. $W = \{m - mI \mid m \in Z^+\} \subseteq V$ is a pseudo neutrosophic semigroup linear subalgebra of V over the semigroup S.

***Example 3.1.35:*** Let $V = \{I, 0, 1\} \subseteq N(Z_2)$ be a neutrosophic semigroup vector space over the semigroup $S = Z_2$. V has pseudo neutrosophic semigroup linear subalgebra over $S = Z_2$ so V is a pseudo neutrosophic semigroup simple vector space over the semigroup S.

Now we proceed onto define the notion of semigroup linear transformation.

**DEFINITION 3.1.13:** *Let V and W be any two neutrosophic semigroup vector spaces over the same semigroup S. We say a map T from V to W is a neutrosophic semigroup linear transformation if $T(c\alpha) = cT(\alpha)$ for all $c \in S$ and $\alpha \in V$.*

***Example 3.1.36:*** Let

$$V = \left\{ \begin{pmatrix} a & b \\ c & d \end{pmatrix} \middle| a, b, c, d \in Q^+I \cup \{0\} \right\}$$

and $W = \{Q^+I \cup \{0\} \times Q^+I \cup \{0\}\}$ be two neutrosophic semigroup vector spaces over the semigroup $Z^+I \cup \{0\}$.

Define $T : V \to W$ by

$$T\left\{\begin{pmatrix} a & b \\ c & d \end{pmatrix}\right\} = (a, d)$$

T is a neutrosophic semigroup linear transformation of V to W.

*Example 3.1.37:* Let

$$V = \{Q^+I \cup \{0\} \times Q^+I \cup \{0\} \times Q^+I \cup \{0\} \times Q^+I \cup \{0\}\}$$

and

$$W = \left\{\begin{pmatrix} a & b & e \\ c & d & f \end{pmatrix} \middle| a,b,c,d,e,f \in Q^+I \cup \{0\}\right\}$$

be two neutrosophic semigroup vector spaces define over the semigroup $S = Q^+I \cup \{0\}$.
Define $T: V \to W$ by

$$T(a, b, c, d) = \left\{\begin{pmatrix} a & 0 & b \\ c & d & 0 \end{pmatrix}\right\},$$

T is a neutrosophic semigroup linear transformation of V to W.

When the domain space and the range space are the same that is V = W then we call the neutrosophic semigroup linear transformation as neutrosophic semigroup linear operator on V.

We will illustrate this by some examples.

*Example 3.1.38:* Let

$$V = \left\{\begin{pmatrix} a_1 & a_2 & a_3 \\ a_4 & a_5 & a_6 \\ a_7 & a_8 & a_9 \end{pmatrix} \middle| a_i \in R^+I \cup \{0\}; 1 \le i \le 9\right\}$$

be a neutrosophic semigroup vector space over the semigroup $Q^+I \cup \{0\}$.
Define $T : V \to V$ by

$$T\begin{pmatrix} a_1 & a_2 & a_3 \\ a_4 & a_5 & a_6 \\ a_7 & a_8 & a_9 \end{pmatrix} = \begin{pmatrix} a_1 & a_2 & a_3 \\ 0 & a_5 & a_6 \\ 0 & 0 & a_9 \end{pmatrix}$$

It is easily verified T is a neutrosophic semigroup linear operator on V.

*Example 3.1.39:* Let $V = \{Q^+I \cup \{0\} \times Q^+I \cup \{0\} \times Q^+I \cup \{0\} \times Q^+I \cup \{0\} \times Q^+I \cup \{0\}\} = \{(a_1, a_2, a_3, a_4, a_5) \mid a_i \in Q^+I \cup \{0\}; 1 \leq i \leq 5\}$ be a neutrosophic semigroup linear algebra over the semigroup $S = Z^+I \cup \{0\}$.
Define $T : V \to V$ by
$$T(a_1, a_2, a_3, a_4, a_5) = (0, a_2, a_3, 0, a_5).$$
It is easily verified that T is a neutrosophic semigroup linear operator on V.

We see in case of neutrosophic semigroup linear algebras V and W over a semigroup S we need an additional condition to be satisfied by the neutrosophic semigroup linear transformation $T: V \to W$; $T(cu + v) = cT(u) + T(v)$; $c \in S$ and $u, v \in V$.

We will illustrate this by some examples.

*Example 3.1.40:* Let

$$V = \left\{ \begin{pmatrix} a_1 & a_2 & a_3 \\ a_4 & a_5 & a_6 \end{pmatrix} \middle| a_i \in R^+I \cup \{0\}; 1 \leq i \leq 6 \right\}$$

and

$$W = \left\{ \begin{pmatrix} a_1 & a_2 \\ a_3 & a_4 \\ a_5 & a_6 \end{pmatrix} \middle| a_i \in R^+I \cup \{0\}; 1 \leq i \leq 6 \right\}$$

be two neutrosophic semigroup linear algebras over the semigroup $S = R^+I \cup \{0\}$. Let $T: V \to W$ be defined as

$$T\left\{\begin{pmatrix} a_1 & a_2 & a_3 \\ a_4 & a_5 & a_6 \end{pmatrix}\right\} = \begin{pmatrix} a_1 & a_2 \\ a_3 & a_4 \\ a_5 & a_6 \end{pmatrix}.$$

It is easily verified T is a neutrosophic semigroup linear transformation of V to W.

*Example 3.1.41:* Let

$$V = \left\{ \begin{pmatrix} a_1 & a_2 & a_3 \\ 0 & a_4 & a_5 \\ 0 & 0 & a_6 \end{pmatrix} \middle| a_i \in Z^+I \cup \{0\}; 1 \le i \le 6 \right\}$$

and $W = \{(Z^+I \cup \{0\}) \times (Z^+I \cup \{0\}) \times (Z^+I \cup \{0\})\} = \{(a_1, a_2, a_3) \mid a_i \in Z^+I \cup \{0\}; 1 \le i \le 3\}$ be a two neutrosophic semigroup linear algebras over the semigroup $S = 2Z^+I \cup \{0\}$.
Define a map $T : V \to W$ by

$$T = \left\{ \begin{pmatrix} a_1 & a_2 & a_3 \\ 0 & a_4 & a_5 \\ 0 & 0 & a_6 \end{pmatrix} \right\} = (a_1, a_4, a_6).$$

It is easily verified T is a neutrosophic semigroup linear transformation from V to W.

*Example 3.1.42:* Let

$$V = \left\{ \begin{pmatrix} a_1 & a_2 \\ a_3 & a_4 \\ a_5 & a_6 \\ a_7 & a_8 \end{pmatrix} \middle| a_i \in Q^+I \cup \{0\}; 1 \le i \le 8 \right\}$$

be a neutrosophic semigroup linear algebra over the semigroup $S = Z^+I \cup \{0\}$.

Define T: V → V by

$$T\left\{\begin{pmatrix} a_1 & a_2 \\ a_3 & a_4 \\ a_5 & a_6 \\ a_7 & a_8 \end{pmatrix}\right\} = \begin{pmatrix} a_1 & 0 \\ a_3 & a_4 \\ 0 & a_6 \\ a_7 & 0 \end{pmatrix}.$$

It is easily verified T is a neutrosophic semigroup linear operator on V.

*Example 3.1.43:* Let

$$V = \left\{\sum_{i=0}^{n} a_i x^i \,\middle|\, n \in I \right.$$

and $a_i \in Q^+I \cup \{0\}$; i.e., all polynomials in the variable x with coefficient from $Q^+I \cup \{0\}\}$ be a neutrosophic semigroup linear algebra over the semigroup $S = Z^+I \cup \{0\}$.
Define T: V → V by

$$T \sum_{i=1}^{n} a_i x^i \to \sum_{i=1}^{n} a_{2i} x^{2i}$$

It is easily verified T is a neutrosophic semigroup linear operator on V.
Let T be a neutrosophic semigroup linear transformation from V into W. We say T is set invertible if there exist a neutrosophic semigroup linear transformation U from W into V such that U.T and T.U are neutrosophic semigroup identity maps on V and W respectively. If T is neutrosophic semigroup invertible, the map U is called the neutrosophic semigroup inverse of T and is unique and is denoted by $T^{-1}$.

The following theorem is left as an exercise for the reader to prove.

**THEOREM 3.1.2:** *Let V and W be two neutrosophic semigroup vector spaces over the semigroup S and T be a neutrosophic semigroup linear transformation from V into W. If T is invertible the inverse map $T^{-1}$ is a neutrosophic semigroup linear transformation from W onto V.*

**DEFINITION 3.1.14:** *Let V be a neutrosophic semigroup linear algebra over the semigroup S. Let $W \subseteq V$ be a neutrosophic subsemigroup linear subalgebra over the subsemigroup P of S. P a proper subsemigroup of the semigroup S. Let $T: V \to W$ be a map such that $T(\alpha v + u) = T(\alpha) T(v) + T(u)$ for all $u, v \in V$ and $T(\alpha) \in P$. We call T a pseudo neutrosophic semigroup linear operator on V.*

Interested reader is requested to construct examples.

**DEFINITION 3.1.15:** *Let V be a neutrosophic semigroup linear algebra over the semigroup S. Let W be a neutrosophic semigroup linear subalgebra of V over S. Let T be a neutrosophic linear operator on V. T is said to be a neutrosophic semigroup linear projection on W if $T(v) = w$; $w \in W$ and $T(\alpha u + v) = \alpha T(u) + T(v)$, $T(u)$ and $T(v) \in W$ for all $\alpha \in S$ and $u, v \in V$.*

We will illustrate this by a simple example.

***Example 3.1.44:*** Let $V = \{(Z^+I \cup \{0\}) \times (Z^+I \cup \{0\}) \times (Z^+I \cup \{0\}) \times (Z^+I \cup \{0\})\}$. V is a neutrosophic semigroup linear algebra over the semigroup $Z^+I \cup \{0\}$. Let $W = (2Z^+I \cup \{0\}) \times (2Z^+I \cup \{0\}) \times \{0\} \times \{0\} \subseteq V$ be a neutrosophic semigroup linear subalgebra of V over $Z^+I \cup \{0\}$.
Define $T: V \to V$ by
$$T(x, y, z, w) = (2x, 2y, 0, 0)$$

It is easily verified that T is a neutrosophic semigroup linear projection of V onto W.

**DEFINITION 3.1.16:** *Let V be a neutrosophic semigroup vector space over the semigroup S. Let $W \subseteq V$ be a neutrosophic semigroup vector subspace of V over the semigroup S. A neutrosophic linear operator T on V is said to be a neutrosophic semigroup projection operator of a subspace of V onto W if for $T : V \to W$; $T(V) \subseteq W$ that is $T(v) = w$ for every $v \in V$ and $w \in W$.*

We will illustrate this by a simple example.

*Example 3.1.45:* Let

$$V = \left\{ \begin{pmatrix} a_1 & a_2 & a_3 & a_4 \\ 0 & 0 & 0 & 0 \end{pmatrix}, \begin{pmatrix} 0 & 0 & 0 & 0 \\ b_1 & b_2 & b_3 & b_4 \end{pmatrix} \middle| \begin{array}{l} a_i, b_i \in Q^+I \cup \{0\} \\ 1 \le i \le 4 \end{array} \right\}$$

be a neutrosophic semigroup vector space over the semigroup S $= Z^+I \cup \{0\}$.
Let

$$W = \left\{ \begin{pmatrix} a_1 & a_2 & a_3 & a_4 \\ 0 & 0 & 0 & 0 \end{pmatrix} \middle| a_i \in Q^+I \cup \{0\}; 1 \le i \le 4 \right\} \subseteq V$$

be a neutrosophic semigroup vector subspace of V over the semigroup S. Let $T : V \to V$ be defined by

$$T\left\{ \begin{pmatrix} a_1 & a_2 & a_3 & a_4 \\ 0 & 0 & 0 & 0 \end{pmatrix} \right\} = \begin{pmatrix} a_1 & a_2 & a_3 & a_4 \\ 0 & 0 & 0 & 0 \end{pmatrix}$$

and

$$T\left\{ \begin{pmatrix} 0 & 0 & 0 & 0 \\ b_1 & b_2 & b_3 & b_4 \end{pmatrix} \right\} = \begin{pmatrix} 0 & 0 & 0 & 0 \\ 0 & 0 & 0 & 0 \end{pmatrix};$$

then T is a neutrosophic semigroup projection of V on W.

We will now define the concept of direct union of neutrosophic semigroup vector subspaces of a neutrosophic semigroup vector space.

**DEFINITION 3.1.17:** *Let V be a neutrosophic semigroup vector space over the semigroup S. Let $W_1, W_2, \ldots, W_n$ be a collection of neutrosophic semigroup vector subspaces of V; if $V = \cup W_i$ and $W_i \cap W_j = \phi$ or $\{0\}$ if $i \neq j$ then we say V is the direct union of the neutrosophic semigroup vector subspaces of the neutrosophic semigroup vector space V over S.*

We will illustrate this by some examples.

*Example 3.1.46:* Let

$$V = \left\{ \begin{pmatrix} a_1 & a_2 & a_3 \\ 0 & 0 & 0 \\ 0 & 0 & 0 \end{pmatrix}, \begin{pmatrix} 0 & 0 & 0 \\ b_1 & b_2 & b_3 \\ 0 & 0 & 0 \end{pmatrix}, \begin{pmatrix} 0 & 0 & 0 \\ 0 & 0 & 0 \\ c_1 & c_2 & c_3 \end{pmatrix} \middle| \begin{array}{l} a_i b_i c_i \in Z^+I \cup \{0\} \\ 1 \leq i \leq 3 \end{array} \right\}$$

be a neutrosophic semigroup vector space over the semigroup S $= Z^+I \cup \{0\}$.
Take

$$W_1 = \left\{ \begin{pmatrix} a_1 & a_2 & a_3 \\ 0 & 0 & 0 \\ 0 & 0 & 0 \end{pmatrix} \middle| a_i \in Z^+I \cup \{0\}; 1 \leq i \leq 3 \right\}$$

$$W_2 = \left\{ \begin{pmatrix} 0 & 0 & 0 \\ b_1 & b_2 & b_3 \\ 0 & 0 & 0 \end{pmatrix} \middle| b_i \in Z^+I \cup \{0\}; 1 \leq i \leq 3 \right\}$$

and

$$W_3 = \left\{ \begin{pmatrix} 0 & 0 & 0 \\ 0 & 0 & 0 \\ c_1 & c_2 & c_3 \end{pmatrix} \middle| c_i \in Z^+I \cup \{0\}; 1 \leq i \leq 3 \right\}$$

be a neutrosophic semigroup vector subspaces of V over the semigroup S. Clearly $V = W_1 \cup W_2 \cup W_3$ and

$$W_i \cap W_j = \begin{pmatrix} 0 & 0 & 0 \\ 0 & 0 & 0 \\ 0 & 0 & 0 \end{pmatrix}, \text{ if } i \neq j; \ 1 \leq i, j \leq 3.$$

Thus V is a direct union of neutrosophic semigroup vector subspaces of V over the semigroup S.

*Example 3.1.47:* Let

$$V = \left\{ \begin{pmatrix} a_1 \\ a_2 \\ a_3 \\ a_4 \end{pmatrix}, \begin{pmatrix} x_1 & x_2 \\ x_3 & x_4 \\ x_5 & x_6 \\ x_7 & x_8 \\ x_9 & x_{10} \end{pmatrix}, (b_1 \ b_2 \ b_3), \begin{pmatrix} c_1 & c_2 \\ c_3 & c_4 \end{pmatrix} \ \middle| \ \begin{array}{l} a_i, b_j, c_k, x_t \in Q^+I \cup \{0\} \\ 1 \leq i \leq 4, \\ 1 \leq j \leq 3, \\ 1 \leq k \leq 4, \\ 1 \leq t \leq 10. \end{array} \right\}$$

be a neutrosophic semigroup vector space over the semigroup S $= Z^+I \cup \{0\}$.
Take

$$W_1 = \left\{ \begin{pmatrix} a_1 \\ a_2 \\ a_3 \\ a_4 \end{pmatrix} \ \middle| \ a_i \in Q^+I \cup \{0\}; 1 \leq i \leq 4 \right\},$$

$$W_2 = \{(b_1 \ b_2 \ b_3) \mid b_j \in Q^+I \cup \{0\}; 1 \leq i \leq 3\},$$

$$W_3 = \left\{ \begin{pmatrix} x_1 & x_2 \\ x_3 & x_4 \\ x_5 & x_6 \\ x_7 & x_8 \\ x_9 & x_{10} \end{pmatrix} \ \middle| \ x_i \in Q^+I \cup \{0\}; 1 \leq i \leq 10 \right\}$$

and

$$W_4 = \left\{ \begin{pmatrix} c_1 & c_2 \\ c_3 & c_4 \end{pmatrix} \middle| c_i \in Q^+I \cup \{0\}; 1 \le i \le 4 \right\}$$

be neutrosophic semigroup vector subspaces of V over the semigroup S. Clearly $V = W_1 \cup W_2 \cup W_3 \cup W_4$ and $W_i \cap W_j = \phi$ if $i \ne j$; $1 \le i, j \le 4$.

Now we proceed onto define the analogous notion for neutrosophic semigroup linear algebras.

**DEFINITION 3.1.18:** *Let V be a neutrosophic semigroup linear algebra over the semigroup S. We say V is the direct sum of neutrosophic semigroup linear subalgebras $W_1, W_2, …, W_n$ of V if*
  *(1) $V = W_1 + … + W_n$*
  *(2) $W_i \cap W_j = \{0\}$ or $\phi$ if $i \ne j$ $1 \le i, j \le n$.*

We will illustrate this situation by some simple examples.

*Example 3.1.48:* Let

$$V = \left\{ \begin{pmatrix} a_1 & a_2 \\ a_3 & a_4 \\ a_5 & a_6 \end{pmatrix} \middle| a_i \in Q^+I \cup \{0\}; 1 \le i \le 6 \right\}$$

be a neutrosophic semigroup linear algebra over the semigroup $S = Z^+I \cup \{0\}$.
Take

$$W_1 = \left\{ \begin{pmatrix} a_1 & 0 \\ 0 & 0 \\ 0 & a_6 \end{pmatrix} \middle| a_1, a_6 \in Q^+I \cup \{0\} \right\},$$

$$W_2 = \left\{ \begin{pmatrix} 0 & a_2 \\ 0 & 0 \\ 0 & 0 \end{pmatrix} \middle| a_2 \in Q^+I \cup \{0\} \right\},$$

$$W_3 = \left\{ \begin{pmatrix} 0 & 0 \\ a_3 & 0 \\ a_5 & 0 \end{pmatrix} \middle| a_3, a_6 \in Q^+I \cup \{0\} \right\}$$

and

$$W_4 = \left\{ \begin{pmatrix} 0 & 0 \\ 0 & a_4 \\ 0 & 0 \end{pmatrix} \middle| a_4 \in Q^+I \cup \{0\} \right\}$$

be neutrosophic semigroup linear subalgebras of V over the semigroup S.
$$V = W_1 + W_2 + W_3 + W_4$$
and
$$W_i \cap W_j = \begin{pmatrix} 0 & 0 \\ 0 & 0 \\ 0 & 0 \end{pmatrix}$$

if $i \neq j$; $1 \leq i, j \leq 4$. This V is a direct sum of neutrosophic semigroup linear subalgebras.

A neutrosophic semigroup linear algebra is strongly simple if it cannot be written as a direct sum of neutrosophic semigroup linear subalgebras and has no proper neutrosophic semigroup linear subalgebras.

***Example 3.1.49:*** Let $V = \{0, I, 2I, 3I, \ldots, 10I\} \subseteq N(Z_{11})$ be a neutrosophic semigroup linear algebra $Z_{11}I$. Clearly V is a strongly simple neutrosophic linear algebra.

In view of this we have a nice theorem which guarantees a class of strongly simple neutrosophic linear algebras.

**THEOREM 3.1.3:** *Let $V = \{0, I, 2I, ..., (p – 1)I \mid p$ is any prime$\} \subseteq N(Z_p)$ and $S = Z_pI$ be the semigroup. Clearly V is a strongly simple neutrosophic semigroup linear algebra.*

*Proof:* Since V has no proper neutrosophic semigroup linear subalgebras we see V cannot be written as a direct sum of neutrosophic semigroup linear subalgebras.
Hence the claim.

In the next section we proceed on to define neutrosophic group vector spaces and neutrosophic group linear algebras.

### 3.2 Neutrosophic Group Linear Algebras

In this section we introduce the notion of neutrosophic group linear algebras. Already the notion of neutrosophic groups have been introduced in the chapter one of this book. We give several interesting properties about them. Infact we illustrate these new concepts by examples so that the reader can follow them easily.

**DEFINITION 3.2.1:** *Let V be a non empty subset (say N(R), N(C) or $N(Z_n)$ or N(Q) or N(Z)). Let G be a group under addition. We call V to be neutrosophic group vector space over G if the following conditions are true.*

*(1) For every $v \in V$ and $g \in G$, gv and vg is in V.*
*(2) $0.v = 0$ for every $v \in V$, 0 the additive identity of G.*

*Example 3.2.1:* Let $V = \{0, 2I, 4I, 6I, 8I, 10I\}$ be a subset of $N(Z_{12})$ and $G = \{0, I, 2I, 3I, 4I, 5I, 6I, 7I, 8I, 10I, 9I, 11I\} \subseteq N(Z_{12})$ be a group under addition modulo 12. V is a neutrosophic group vector space over the group G.

*Example 3.2.2:* Let

$$V = \left\{ \begin{pmatrix} 0 & 0 & 0 & 0 \\ a_1 & a_2 & a_3 & a_4 \end{pmatrix}, \begin{pmatrix} b_1 & b_2 & b_3 & b_4 \\ 0 & 0 & 0 & 0 \end{pmatrix}, \right.$$

$$\left. \left( \begin{array}{ccc} 0 & 0 & 0 \\ c_1 & c_2 & c_3 \end{array} \right), \left( \begin{array}{cc} a & b \\ c & d \end{array} \right) \right| a_i, b_j, c_k, a, b, c, d \in QI\,;$$

$$1 \leq i \leq 4,\ 1 \leq j \leq 4 \text{ and } 1 \leq k \leq 3\}$$

and $G = QI \subseteq N(Q)$ be a group under addition. V is a neutrosophic group vector space over the group G.

*Example 3.2.3:* Let $V = \{(0, a_1, 0, a_2, 0, a_3), (b_1, b_2, b_3, b_4), (c_1, c_2, c_3) \mid a_i, b_j, c_k \in ZI;\ 1 \leq i \leq 3,\ 1 \leq j \leq 4 \text{ and } 1 \leq k \leq 3\}$ and $G = ZI$ be a group under addition. V is a neutrosophic group vector space over G.

*Example 3.2.4:* Let

$$V = \left\{ \left( \begin{array}{cc} a_1 & b_1 \\ 0 & 0 \end{array} \right), \left( \begin{array}{cc} 0 & 0 \\ a_2 & b_2 \\ 0 & 0 \end{array} \right), \left( \begin{array}{cc} 0 & 0 \\ 0 & 0 \\ a_3 & b_3 \end{array} \right), \left( \begin{array}{ccc} x_1 & x_2 & x_3 \\ 0 & 0 & 0 \\ 0 & 0 & 0 \end{array} \right), \right.$$

$$\left. \left. \left( \begin{array}{ccc} 0 & 0 & 0 \\ 0 & x_4 & x_5 \\ x_6 & 0 & 0 \end{array} \right) \right| a_i, b_i, x_i \in RI; 1 \leq i, j \leq 3 \text{ and } 1 \leq k \leq 6 \right\}$$

and $G = QI$ a group under addition. V is a neutrosophic group vector space over the group G.

**DEFINITION 3.2.2:** *Let V be a neutrosophic group vector space over the group G. Suppose $W \subseteq V$ be a proper subset of V. We say W is a neutrosophic group vector subspace of V if W is itself a neutrosophic group vector space over G.*

We shall illustrate this situation by some examples.

*Example 3.2.5:* Let $V = \{(2ZI)[X] \text{ and } (5ZI)[x];$ be polynomials with coefficients from 2ZI and 5ZI respectively$\}$ and $G = ZI$ a group under addition. V is a neutrosophic group vector space over G. Take W = {(2ZI)[x]} ⊆ V; W is a neutrosophic group vector subspace of V over the group G.

*Example 3.2.6:* Let

$$V = \left\{ \begin{pmatrix} a_1 & 0 & 0 \\ a_2 & 0 & 0 \\ a_3 & 0 & 0 \end{pmatrix}, \begin{pmatrix} b_1 & 0 & 0 \\ 0 & b_2 & 0 \\ 0 & 0 & b_3 \end{pmatrix}, \begin{pmatrix} 0 & 0 & b_4 \\ b_5 & 0 & b_6 \\ 0 & b_7 & 0 \end{pmatrix} \middle| \begin{array}{l} a_i, b_j \in RI; \\ 1 \le i \le 3; \\ 1 \le j \le 7; \end{array} \right\};$$

be a neutrosophic group vector space over the group G = ZI. Take

$$W = \left\{ \begin{pmatrix} 0 & 0 & b_4 \\ b_5 & 0 & b_6 \\ 0 & b_7 & 0 \end{pmatrix} \middle| b_4, b_5, b_6, b_7 \in RI \right\} \subseteq V;$$

W is a neutrosophic group vector subspace of V over the group G.

Now we proceed onto define the concept of linearly independent subset of a neutrosophic group vector space.

**DEFINITION 3.2.3:** *Let V be a neutrosophic group vector space over the group G.*
  *We say a proper subset P of V to be a linearly dependent neutrosophic subset of V if for any $p_1$, $p_2$ in P ($p_1 \ne p_2$) $p_1 = ap_2$ or $p_2 = a_1 p_1$ for some a, $a_i \in G$. If for no distinct pair of elements $p_1, p_2 \in P$ we have $a_1, a_2 \in G$ such that $p_1 = a_1 p_2$ or $p_2 = a_2 p_1$ then we say P is a linearly independent neutrosophic subset of V.*

We will illustrate this by some examples.

*Example 3.3.7:* Let

$$V = \left\{ \begin{pmatrix} a_1 & a_2 \\ 0 & 0 \\ 0 & 0 \end{pmatrix}, \begin{pmatrix} 0 & 0 \\ a_1 & a_2 \\ 0 & 0 \end{pmatrix}, \begin{pmatrix} 0 & 0 \\ 0 & 0 \\ a_1 & a_2 \end{pmatrix} \middle| a_1, a_2 \in ZI \right\}$$

be a neutrosophic group vector space over the group ZI.
Take

$$P = \left\{ \begin{pmatrix} 2I & 4I \\ 0 & 0 \\ 0 & 0 \end{pmatrix}, \begin{pmatrix} I & 2I \\ 0 & 0 \\ 0 & 0 \end{pmatrix}, \begin{pmatrix} 6I & 12I \\ 0 & 0 \\ 0 & 0 \end{pmatrix}, \begin{pmatrix} 5I & 10I \\ 0 & 0 \\ 0 & 0 \end{pmatrix} \right\} \subseteq V.$$

P is a linearly dependent neutrosophic subset of V.
Take

$$Q = \left\{ \begin{pmatrix} 0 & 0 \\ I & I \\ 0 & 0 \end{pmatrix}, \begin{pmatrix} I & 4I \\ 0 & 0 \\ 0 & 0 \end{pmatrix}, \begin{pmatrix} 0 & 0 \\ 0 & 0 \\ 5I & 6I \end{pmatrix} \right\} \subseteq V.$$

Q is a linearly independent neutrosophic subset of V over G.
Take

$$T = \left\{ \begin{pmatrix} 8I & I \\ 0 & 0 \\ 0 & 0 \end{pmatrix}, \begin{pmatrix} 0 & 0 \\ 7I & 8I \\ 0 & 0 \end{pmatrix}, \begin{pmatrix} 0 & 0 \\ 0 & 0 \\ I & 7I \end{pmatrix} \right\} \subseteq V,$$

T is also a linearly independent neutrosophic subset of V over G.

Now we will define the notion of generating neutrosophic subset of a neutrosophic group vector space over a group G.

**DEFINITION 3.2.4:** *Let V be a neutrosophic group vector space over the group G. Suppose T is a subset of V which is a linearly independent neutrosophic subset of T and if T generates V that is using $t \in T$ and $g \in G$ we can get every $v \in V$ as $v = gt$ then we call T to be a generating neutrosophic subset of V over G. The number of elements in G gives the dimension of V. If T is of* *finite cardinality we say V is of finite dimension. If T is of infinite cardinality we say V is of infinite dimension.*

**Example 3.2.8:** Let

$$V = \left\{ \begin{pmatrix} a_1 & 0 \\ a_1 & 0 \end{pmatrix}, \begin{pmatrix} 0 & b_1 \\ 0 & b_1 \end{pmatrix} \middle| a_i, b_i \in ZI \right\}$$

be a neutrosophic group vector space over the group $ZI = G$. Take

$$T = \left\{ \begin{pmatrix} I & 0 \\ I & 0 \end{pmatrix}, \begin{pmatrix} 0 & I \\ 0 & I \end{pmatrix} \right\} \subseteq V;$$

T is generating neutrosophic subset of V. Clearly V is finite dimensional we can have independent neutrosophic subsets of V but they may not be generating subsets of V which will be illustrated.

**Example 3.2.9:** Let

$$V = \left\{ \begin{pmatrix} a_1 & a_2 \\ 0 & 0 \end{pmatrix}, \begin{pmatrix} 0 & 0 \\ b_1 & b_2 \end{pmatrix} \middle| a_i, b_i \in ZI; 1 \leq i \leq 2 \right\}$$

be a neutrosophic group vector space over the group $G = ZI$. We have several linearly independent neutrosophic subsets of V but V cannot be finitely generated over G. Thus dimension of V over G is infinite.
Take

$$T = \left\{ \begin{pmatrix} 0 & 0 \\ I & 0 \end{pmatrix}, \begin{pmatrix} 0 & 0 \\ 0 & I \end{pmatrix}, \begin{pmatrix} 0 & 0 \\ I & I \end{pmatrix}, \begin{pmatrix} I & 0 \\ 0 & 0 \end{pmatrix}, \begin{pmatrix} 0 & I \\ 0 & 0 \end{pmatrix}, \begin{pmatrix} I & I \\ 0 & 0 \end{pmatrix} \right\}$$

is a linearly independent neutrosophic subset of V but T cannot generate V over G.

*Example 3.2.10:* Let V = {(a a a a a a) | a ∈ QI} be a neutrosophic group vector space over the group QI = G. T = (I I I I I I)} ⊆ V is the generating subset of V over QI.

*Example 3.2.11:* Let V = {(a a a a a a) be such that a ∈ QI} be a neutrosophic group vector space over the group G = ZI. Clearly V is of infinite dimension over ZI.

However we have several finite linearly independent subsets of V.

**DEFINITION 3.2.5:** *Let V be a neutrosophic group vector space over the group G. Let W ⊆ V be a proper subset of V. H ⊆ G be a proper subgroup of G. If W is a neutrosophic group vector space over H then we call W to be a neutrosophic subgroup vector subspace of V over the subgroup H of G.*

We will illustrate this by some simple examples.

*Example 3.2.12:* Let

$$V = \left\{ (a\ a\ a), \begin{bmatrix} a & b \\ c & d \end{bmatrix}, \begin{bmatrix} x \\ y \\ z \\ w \end{bmatrix} \,\middle|\, a, b, c, d, x, y, z, w, \in QI \right\}$$

be a neutrosophic group vector space over the group G = QI. Take

$$W = \left\{ \begin{bmatrix} a & b \\ c & d \end{bmatrix} \,\middle|\, a, b, c, d \in QI \right\}$$

contained in V and H = ZI ⊆ G be a subgroup of G; clearly W is a neutrosophic subgroup vector subspace of V over the subgroup H = ZI of G.

*Example 3.2.13:* Let V = {QI × QI × QI} = {(x, y, z) | x, y, z ∈ QI} be a neutrosophic group vector space over the group G = ZI. Let W = {QI × QI × {0}} = {(x, y, 0) | x, y ∈ QI} ⊆ V; W is

a neutrosophic subgroup vector subspace of V over the subgroup H = 2ZI of G.

*Example 3.2.14:* Let

$$V = \left\{ \begin{pmatrix} a_1 & a_2 & a_3 \\ 0 & a_4 & a_5 \\ 0 & 0 & a_6 \end{pmatrix}, \begin{pmatrix} b_6 & 0 & 0 \\ b_1 & b_2 & 0 \\ b_3 & b_4 & b_5 \end{pmatrix} \middle| \begin{array}{l} a_i, b_j \in QI \\ 1 \le i \le 6 \\ 1 \le j \le 6 \end{array} \right\}$$

be a neutrosophic group vector space over the group G = QI. Take

$$W = \left\{ \begin{pmatrix} a_1 & a_2 & a_3 \\ 0 & a_4 & a_5 \\ 0 & 0 & a_6 \end{pmatrix} \middle| a_i \in ZI; 1 \le i \le 6 \right\}$$

contained in V. Let H = ZI ⊆ G be a subgroup of G. Clearly W is a neutrosophic subgroup vector subspace of V over the subgroup H of G.

*Example 3.2.15:* Let V = {$Z_6I \times Z_6I \times Z_6I \times Z_6I$} be a neutrosophic group vector space over the group G = $Z_6I$. Take W = {(0 0 0), (I I I), (3I, 3I, 3I)} ⊆ V; W is a neutrosophic subgroup vector subspace of V over the subgroup H = {0, 3I} ⊆ G. Clearly W is not a neutrosophic group vector subspace of V over the group G.

In the view of this we have the following theorem.

**THEOREM 3.2.1:** *Let V be a neutrosophic group vector space over the group G. If W is a neutrosophic group vector subspace of V then W need not be a neutrosophic subgroup vector subspace of V.*

The proof is left as an exercise for the reader.

**THEOREM 3.2.2:** *Let V be a neutrosophic group vector space over a group G. Suppose S $\subseteq$ V is a neutrosophic subgroup vector subspace of V then S need not in general be a neutrosophic group vector subspace of V over G.*

This proof is also left as an exercise for the reader.

Next we proceed onto define a neutrosophic duo subgroup vector subspace.

**DEFINITION 3.2.6:** *Let V be a neutrosophic group vector space over the group G. Let W $\subseteq$ V. If W is a neutrosophic subgroup vector subspace over a proper subgroup H of G as well as W is a neutrosophic group vector subspace of V over G then we call W to be a neutrosophic duo subgroup vector subspace of V.*

We will illustrate this by some examples.

*Example 3.2.16:* Let

$$V = \left\{ \begin{pmatrix} a_1 & a_2 \\ 0 & a_3 \end{pmatrix}, \begin{pmatrix} a_1 & 0 & 0 \\ a_2 & a_3 & 0 \\ a_4 & a_5 & a_6 \end{pmatrix} \middle| a_i \in QI; 1 \le i \le 6 \right\}$$

be a neutrosophic group vector space over the group G = QI. Let

$$W = \left\{ \begin{pmatrix} a_1 & a_2 \\ 0 & a_3 \end{pmatrix} \middle| a_i \in QI; 1 \le i \le 3 \right\} \subseteq V$$

be a neutrosophic group vector subspace of V over the group G. It is easy to verify W is also a neutrosophic subgroup vector subspace of V over the subgroup H = ZI $\subseteq$ QI $\subseteq$ G. Thus W is a neutrosophic due subgroup vector subspace of V.

*Example 3.2.17:* Let V = {($x_1$, $x_2$, $x_3$, $x_4$) | $x_i \in$ QI; 1 $\le$ i $\le$ 4} be a neutrosophic group vector space over the group G = ZI. Take W = {(0, $x_2$, $x_3$, 0)| $x_2$, $x_3 \in$ QI} $\subseteq$ V. W is a neutrosophic group

vector subspace of V over G and W is also a neutrosophic subgroup vector subspace of V over the subgroup $H = ZI \subseteq QI = G$. Thus W is a neutrosophic duo subgroup vector subspace of V.

The following theorem is evident from the very definition.

**THEOREM 3.2.3:** *Let V be a neutrosophic group vector space over the group G; if W is a neutrosophic duo subgroup vector subspace of V then W is both a neutrosophic group vector subspace of V as well as W is a neutrosophic subgroup vector subspace of V.*

**DEFINITION 3.2.7:** *Let V be a neutrosophic group vector space over the group G. Suppose V has no neutrosophic subgroup vector subspace then we call V to be a neutrosophic simple group vector space.*

We will illustrate this situation by some examples.

***Example 3.2.18:*** Let $V = \{Z_{11}I \times Z_{11}I \times Z_{11}I \times Z_{11}I \times Z_{11}I\} = \{(x, y, z, w, t) \,/\, x, y, z, w, t, \in Z_{11}I\}$. V is a neutrosophic simple group vector space over the group $G = Z_{11}I$.

***Example 3.2.19:*** Let

$$V = \left\{ \begin{pmatrix} a & b \\ c & d \end{pmatrix} \middle| a, b, c, d \in Z_{23}I \right\}$$

be a neutrosophic group vector space over the group $G = Z_{23}I$. V is a neutrosophic simple group vector space over G as G has no proper subgroups.

We have the following nice theorem which guarantees the existence of neutrosophic simple group vector spaces.

**THEOREM 3.2.4:** *Let V be a neutrosophic group vector space over a group G, which has no proper subgroups other than G and {0}, then V is a neutrosophic simple group vector space over G.*

*Proof:* Follows from the fact that G has no proper subgroup for a proper subset W to be a neutrosophic subgroup vector subspace; we need a proper subgroup in G over which W is a group vector space.

If G has no proper subgroup the existence of neutrosophic subgroup vector subspace is impossible.

We will now give a large class of neutrosophic simple group vector spaces.

**THEOREM 3.2.5:** *Let $V = Z_pI \times Z_pI \times \ldots \times Z_pI$ – n times be a neutrosophic group vector space over the group $Z_pI = G$, where p is a prime. V is a neutrosophic simple group vector space over $G = Z_pI$.*

*Proof:* Clear from the fact that $Z_pI$ has no proper subgroups.

**DEFINITION 3.2.8:** *Let V be a neutrosophic group vector space over the group G. Let $W \subseteq V$ and $S \subseteq G$ where S is a semigroup under +. If W is a neutrosophic semigroup vector subspace of V over S then we call W to be a neutrosophic pseudo semigroup vector subspace of V over S.*

We will illustrate this by some simple examples.

*Example 3.2.20:* Let $V = ZI \times ZI \times ZI$ be a neutrosophic group vector space over the group $G = ZI$. $W = ZI \times ZI \times \{0\} \subseteq V$. W is a neutrosophic pseudo semigroup vector subspace of V over the subsemigroup $S = Z^+I \cup \{0\} \subseteq ZI$.

*Example 3.2.21:* Let
$$V = \{(QI \times QI \times QI \times QI), \begin{pmatrix} a & b \\ c & d \end{pmatrix} \mid a, b, c, d \in QI\}$$
be a neutrosophic group vector space over the group $G = QI$. Let

$$W = \left\{ \begin{pmatrix} a & b \\ c & d \end{pmatrix} \middle| a,b,c,d \in QI \right\} \subseteq V$$

is a neutrosophic pseudo semigroup vector subspace of V over the semigroup $Z^+I \cup \{0\} \subseteq QI = G$.

*Example 3.2.22:* Let

$$V = \left\{ \begin{bmatrix} a \\ b \\ c \end{bmatrix}, [a\ b\ c\ d\ e]\ |\ a, b, c, d, e \in RI \right\}$$

be a neutrosophic group vector space over the group $G = QI$. Take W = {(a, b, c, d, e) | a, b, c, d, e $\in$ RI} $\subseteq$ V be a neutrosophic pseudo semigroup vector subspace of V over the semigroup $Z^+I \cup \{0\} \subseteq QI = G$.

We now proceed onto define the notion of neutrosophic group linear transformation.

**DEFINITION 3.2.9:** *Let V and W be two neutrosophic group vector spaces defined over the same group G. A map T from V to W will be called as the neutrosophic group linear transformation; if $T(\alpha v) = \alpha T(v)$ for all $\alpha \in G$ and for all $v \in V$.*

We will illustrate this by some simple examples.

*Example 3.2.23:* Let V = ZI × ZI × ZI and W = QI × QI × QI × QI × QI be two neutrosophic group vector spaces over the group G = ZI. Let T : V → W be defined by T (x, y, z) = (z, y, x, y, z) for all (x, y, z) $\in$ V. Clearly T is a neutrosophic group linear transformation of V into W.

*Example 3.2.24:* Let

$$V = \left\{ \begin{pmatrix} a & b \\ c & d \end{pmatrix} \middle| a, b, c, d \in QI \right\}$$

and

$$W = \left\{ \begin{bmatrix} a \\ b \\ c \\ d \end{bmatrix} \middle| a, b, c, d \in QI \right\}$$

be a neutrosophic group vector spaces over the group G = QI. Define T : V → W by

$$T \left\{ \begin{pmatrix} a & b \\ c & d \end{pmatrix} \right\} = \begin{bmatrix} a \\ b \\ c \\ d \end{bmatrix};$$

T is a neutrosophic group linear transformation of V into W.

**DEFINITION 3.2.10:** *Let V be a neutrosophic group vector space over the group G. Let T from V to V be a neutrosophic linear transformation then we call T to be a neutrosophic group linear operator on V.*

We will illustrate this by some examples.

*Example 3.2.25:* Let V = {(a, b, c, d) | a, b, c, d ∈ QI} be a neutrosophic group vector space over G = QI. Define T from V to V by T(a, b, c, d) = (d, c, b, a). Clearly T is a neutrosophic group linear operator on V.

*Example 3.2.26:* Let

$$V = \left\{ \begin{pmatrix} a & b & c \\ d & e & f \\ g & h & i \end{pmatrix} \middle| a,b,c,d,e,f,g,h,i \in RI \right\}$$

be a neutrosophic group vector space over the group G = ZI. Define T : V → V by

$$T\begin{pmatrix} a & b & c \\ d & e & f \\ g & h & i \end{pmatrix} = \begin{pmatrix} a & 0 & c \\ d & e & 0 \\ 0 & 0 & i \end{pmatrix}.$$

It is easily verified T is a neutrosophic group linear operator on V. Define N($M_G$ (V,W)) = {collection of all neutrosophic group linear transformations from V to W; V and W neutrosophic group vector spaces over the group G} and N($M_G$(V,V)) = {set of all neutrosophic group linear operators from V to V, V a neutrosophic group vector space over G}. The reader is expected to study the algebraic structure of N($M_G$(V,W)) and N($M_G$(V,V)).

We now proceed onto the notion of neutrosophic group linear algebra over a group.

**DEFINITION 3.2.11:** *Let V be a neutrosophic group vector space over the group G. If V is again a neutrosophic group under the operation of addition, then we call V to be a neutrosophic group linear algebra over G.*

*Example 3.2.27:* Let

$$V = \left\{ \begin{pmatrix} a & a & a \\ a & a & a \\ a & a & a \end{pmatrix} \middle| a \in QI \right\}.$$

V is a neutrosophic group linear algebra over the group ZI = G.

*Example 3.2.28:* Let V = {(x, y, z) | x, y, z ∈ QI} be a neutrosophic group linear algebra over the group G = QI.

It is important to mention at this juncture that every neutrosophic group linear algebra is a neutrosophic group vector space over a group G but however a neutrosophic group vector space over a group G in general is not a neutrosophic group linear algebra over a group G.

We will illustrate this by an example.

*Example 3.2.29:* Let

$$V = \left\{ \begin{pmatrix} a & b \\ 0 & 0 \end{pmatrix}, \begin{pmatrix} 0 & 0 \\ c & d \end{pmatrix}, \begin{pmatrix} a & 0 \\ 0 & c \end{pmatrix}, \begin{pmatrix} 0 & d \\ e & 0 \end{pmatrix} \middle| a,b,c,d,e \in QI \right\}$$

be a neutrosophic group vector space over the group $G = QI$. We see V is not a group under matrix addition. Thus V is only a neutrosophic group vector space over the group G and V is not a neutrosophic group linear algebra over the group G.

We proceed onto define the notion of dimension of a neutrosophic group linear algebra.

**DEFINITION 3.2.12:** *Let V be a neutrosophic group linear algebra over the group G. $X \subset V$ be a proper subset of V, we say X is a linearly independent subset of V if $X = \{x_1, ..., x_n\}$, for some $x_i \in G$; $1 \leq i \leq n$; $\sum_{i=1}^{n} \alpha_i x_i = 0$ if and only if each $\alpha_i = 0$. A linearly independent subset X of V is said to be a generator of V if every element v of V can be represented as*

$$v = \sum_{i=1}^{n} \alpha_i x_i \ ; \ \alpha_i \in G \ (1 \leq i \leq n).$$

We will illustrate this situation by some examples.

*Example 3.2.30:* Let $V = \{(x, y, z) \mid x, y, z \in Z_2I; Z_2 = \{0, I\}\}$ be a neutrosophic group linear algebra over the group $Z_2I = G$. V is generated by the set $X = \{(I\ 0\ 0)\ (0\ I\ 0), (0\ 0\ I)\}$. Clearly X is a linearly independent subset of V over the group $G = Z_2I$.

*Example 3.2.31:* Let

$$V = \left\{ \begin{pmatrix} a_1 & a_2 \\ a_3 & a_4 \\ a_5 & a_6 \end{pmatrix} \middle| a_i \in ZI;\ 1 \le i \le 6 \right\}$$

be the neutrosophic group linear algebra over the group $G = ZI$. Let

$$X = \left\{ \begin{pmatrix} I & 0 \\ 0 & 0 \\ 0 & 0 \end{pmatrix}, \begin{pmatrix} 0 & I \\ 0 & 0 \\ 0 & 0 \end{pmatrix}, \begin{pmatrix} 0 & 0 \\ I & 0 \\ 0 & 0 \end{pmatrix}, \begin{pmatrix} 0 & 0 \\ 0 & I \\ 0 & 0 \end{pmatrix}, \begin{pmatrix} 0 & 0 \\ 0 & 0 \\ I & 0 \end{pmatrix}, \begin{pmatrix} 0 & 0 \\ 0 & 0 \\ 0 & I \end{pmatrix} \right\} \subseteq V$$

is the generating subset of V over the group $G = ZI$.

We now proceed into define substructures of neutrosophic group linear algebras.

**DEFINITION 3.2.13:** *Let V be a neutrosophic group linear algebra over the group G. Let $W \subseteq V$ be a proper subset of V. We say W is a neutrosophic group linear subalgebra of V over G if W is itself a neutrosophic group linear algebra over G.*

We illustrate this situation by some examples.

*Example 3.2.32:* Let

$$V = \left\{ \begin{pmatrix} a & b & c \\ d & e & f \\ g & h & i \end{pmatrix} \middle| a, b, \ldots, i \in QI \right\}$$

be a neutrosophic group linear algebra over the group $G = ZI$. Take

$$W = \left\{ \begin{pmatrix} a & a & a \\ b & b & b \\ c & c & c \end{pmatrix} \middle| a, b, c \in QI \right\} \subseteq V;$$

W is a neutrosophic group linear subalgebra of V over the group G.

*Example 3.2.33:* Let V = {(a, b, c) | a, b, c ∈ QI} be a neutrosophic group linear algebra over the group G = ZI. Take W = {(a, b, c) | a, b, c ∈ ZI ⊆ QI} ⊆ V; W is a neutrosophic group linear subalgebra of V over the group G.

We will now proceed onto define the notion of direct sum of neutrosophic group linear subalgebra of a neutrosophic group linear algebra.

**DEFINITION 3.2.14:** *Let V be a neutrosophic group linear algebra over the group G. Let $W_1$, $W_2$, ..., $W_n$ be neutrosophic group linear subalgebras of V over the group G.*
   *We say V is a direct sum of the neutrosophic group linear subalgebras $W_1$, $W_2$, ..., $W_n$ if*

   *(1) $V = W_1 + ... + W_n$*
   *(2) $W_i \cap W_j = \{0\}$ if $i \neq j$; $1 \leq i, j \leq n$.*

We will illustrate this by some examples.

*Example 3.2.34:* Let $V = Z_{14}I \times Z_{14}I \times Z_{14}I$ be a neutrosophic group linear algebra over the group $G = Z_{14}I$, the group under addition modulo 14. Let $W_1 = Z_{14}I \times \{0\} \times \{0\}$, $W_2 = \{0\} \times Z_{14}I \times Z_{14}I$ be neutrosophic group linear subalgebras of V over the group G.
   We see $V = W_1 + W_2$ and $W_1 \cap W_2 = (0\ 0\ 0)$; if $i \neq j$; $1 \leq i, j \leq 3$.

*Example 3.2.35:* Let

$$V = \left\{ \begin{pmatrix} a & 0 & 0 \\ b & c & 0 \\ d & e & f \end{pmatrix} \middle| a,b,c,d,e,f \in QI \right\}$$

be a neutrosophic group linear algebra over $G = ZI$, the group under addition.
Take
$$W_1 = \left\{ \begin{pmatrix} a & 0 & 0 \\ 0 & 0 & 0 \\ 0 & 0 & f \end{pmatrix} \middle| a, f \in QI \right\},$$

$$W_2 = \left\{ \begin{pmatrix} 0 & 0 & 0 \\ b & c & 0 \\ 0 & 0 & 0 \end{pmatrix} \middle| b, c \in QI \right\},$$

$$W_3 = \left\{ \begin{pmatrix} 0 & 0 & 0 \\ 0 & 0 & 0 \\ d & 0 & 0 \end{pmatrix} \middle| d \in QI \right\}$$

and

$$W_4 = \left\{ \begin{pmatrix} 0 & 0 & 0 \\ 0 & 0 & 0 \\ 0 & e & 0 \end{pmatrix} \middle| e \in QI \right\}$$

be neutrosophic group linear subalgebras of V over the group G. We see $V = W_1 + W_2 + W_3 + W_4$ and

$$W_i \cap W_j = \begin{pmatrix} 0 & 0 & 0 \\ 0 & 0 & 0 \\ 0 & 0 & 0 \end{pmatrix} \; i \neq j; \; 1 \leq i, j \leq 4.$$

Now we proceed onto define the notion of pseudo direct sum of a neutrosophic linear algebra over the group G.

**DEFINITION 3.2.15:** *Let V be a neutrosophic group linear algebra over the group G. Suppose $W_1, W_2, \ldots, W_n$ are distinct neutrosophic group linear subalgebras of V over G. We say V is a pseudo direct sum if*

(1) $V = W_1 + \ldots + W_n$
(2) $W_i \cap W_j \neq \{0\}$ or $\phi$ in general even if $i \neq j$
(3) We need $W_i$'s to be distinct i.e., $W_i \cap W_j \neq W_i$ or $W_j$ if $i \neq j$; $1 \leq i, j \leq n$.

We will illustrate this situation by some examples.

*Example 3.2.36:* Let

$$V = \left\{ \begin{pmatrix} a_1 & a_2 & a_3 \\ a_4 & a_5 & a_6 \\ a_7 & a_8 & a_9 \end{pmatrix} \middle| a_i \in QI; 1 \leq i \leq 9 \right\}$$

be a neutrosophic group linear algebra over the group $G = ZI$. Take

$$W_1 = \left\{ \begin{pmatrix} a_1 & a_2 & 0 \\ 0 & a_5 & 0 \\ 0 & 0 & 0 \end{pmatrix} \middle| a_1, a_2, a_5 \in QI \right\}$$

$$W_2 = \left\{ \begin{pmatrix} 0 & a_2 & 0 \\ a_4 & 0 & 0 \\ 0 & 0 & 0 \end{pmatrix} \middle| a_2, a_4 \in QI \right\}$$

$$W_3 = \left\{ \begin{pmatrix} a_1 & 0 & 0 \\ a_4 & 0 & 0 \\ a_7 & a_8 & a_9 \end{pmatrix} \middle| a_1, a_4, a_7, a_8, a_9 \in QI \right\}$$

and

$$W_4 = \left\{ \begin{pmatrix} a_1 & 0 & 0 \\ a_4 & a_5 & a_6 \\ 0 & a_8 & 0 \end{pmatrix} \middle| a_1, a_4, a_5, a_6, a_8 \in QI \right\}$$

to be neutrosophic group linear subalgebras of V over the group $G = ZI$. We see $V = W_1 + W_2 + W_3 + W_4$

But

$$W_1 \cap W_2 = \begin{pmatrix} 0 & a_2 & 0 \\ 0 & 0 & 0 \\ 0 & 0 & 0 \end{pmatrix}, W_1 \cap W_3 = \begin{pmatrix} a_1 & 0 & 0 \\ 0 & 0 & 0 \\ 0 & 0 & 0 \end{pmatrix},$$

$$W_1 \cap W_4 = \begin{pmatrix} a_1 & 0 & 0 \\ 0 & 0 & 0 \\ 0 & 0 & 0 \end{pmatrix}, W_2 \cap W_3 = \begin{pmatrix} 0 & 0 & 0 \\ a_4 & 0 & 0 \\ 0 & 0 & 0 \end{pmatrix},$$

$$W_2 \cap W_4 = \begin{pmatrix} 0 & 0 & 0 \\ a_4 & 0 & 0 \\ 0 & 0 & 0 \end{pmatrix} \text{ and } W_3 \cap W_4 = \begin{pmatrix} a_1 & 0 & 0 \\ a_4 & 0 & 0 \\ 0 & a_8 & 0 \end{pmatrix}$$

and $W_i \not\subseteq W_j$ for $i \neq j$; $1 \leq i, j \leq 4$.

Thus V is a pseudo direct sum of neutrosophic group linear subalgebras over the same group G.

*Example 3.2.37:* Let $V = \{Z_{24}I \times Z_{24}I \times Z_{24}I \times Z_{24}I\}$ be a neutrosophic group linear algebra over the group $G = Z_{24}I$.
Take
$$W_1 = \{Z_{24}I \times Z_{24}I \times \{0\} \times \{0\}\},$$
$$W_2 = \{0 \times \{Z_{24}I\} \times Z_{24}I \times \{0\}\} \text{ and}$$
$$W_3 = \{\{0\} \times Z_{24}I \times Z_{24}I \times Z_{24}I\}$$

to be neutrosophic group linear subalgebras of V over G. We see $V = W_1 + W_2 + W_3$

$$W_1 \cap W_2 = \{0\} \times Z_{24}I \times \{0\} \times \{0\}$$
$$W_1 \cap W_3 = \{\{0\} \times Z_{24}I \times \{0\} \times \{0\}\}$$
$$W_2 \cap W_3 = \{0\} \times \{0\} \times Z_{24}I \times \{0\}$$

Thus V is a pseudo direct sum of neutrosophic linear subalgebras of V over the group G.

Now we proceed onto define yet another new algebraic structure in neutrosophic group linear algebras over a group G.

**DEFINITION 3.2.16:** *Let V be a neutrosophic group linear algebra over the group G. Let W $\subseteq$ V be a proper subgroup of V. Suppose H $\subseteq$ G be a proper subsemigroup of G.*

*If W is a neutrosophic semigroup linear algebra over the semigroup H then we call W to be a pseudo neutrosophic semigroup linear subalgebra of the neutrosophic group linear algebra V.*

We will illustrate this situation by some examples.

*Example 3.2.38:* Let

$$V = \left\{ \begin{pmatrix} a & b & c \\ d & e & f \\ g & h & i \\ k & l & m \end{pmatrix} \middle| a,b,...,m \in QI \right\}$$

be a neutrosophic group linear algebra over the group G = ZI.

$$W = \left\{ \begin{pmatrix} a & a & a \\ a & a & a \\ a & a & a \\ a & a & a \end{pmatrix} \middle| a \in QI \right\}$$

is a pseudo neutrosophic semigroup linear subalgebra of the neutrosophic group linear algebra over the semigroup $Z^+I \cup \{0\}$.

*Example 3.2.39:* Let V = {QI × QI × QI × QI × QI × QI} be a neutrosophic group linear algebra over the group G = ZI. Take W = {QI × {0} × QI × {0} × QI × {0}} $\subseteq$ V; W is a pseudo neutrosophic semigroup linear subalgebra of V over the semigroup S = $3Z^+I \cup \{0\} \subseteq$ G.

It may so happen that at times we may have neutrosophic group linear algebra V over a group G but may not have pseudo neutrosophic semigroup linear subalgebras of V. This is given by these classes of neutrosophic group linear algebras.

*Example 3.2.40:* Let $V = Z_pI \times Z_pI \times \ldots \times Z_pI$ be a neutrosophic group linear algebra over the group $G = Z_pI$ (p a prime); V has no pseudo neutrosophic semigroup linear subalgebras.

*Example 3.2.41:* Let $V = \{(a_{ij})_{m \times n} \mid a_{ij} \in Z_pI;$ p a prime$\}$ be a neutrosophic group linear algebra over the group $G = Z_pI$. V has no pseudo neutrosophic semigroup linear subalgebra as $G = Z_pI$ has no proper subset which is a semigroup under addition.

*Example 3.2.42:* Let $V = \{Z_pI[x] \mid$ p is a prime and $Z_pI[x]$ is a collection of polynomials in the variable x with coefficient from $Z_pI\}$ be a neutrosophic group linear algebra over the group $G = Z_pI$. V has no pseudo neutrosophic semigroup linear subalgebra as $G = Z_pI$ has no proper subset which is a semigorup.

Now we proceed onto define yet another new algebraic structure of the neutrosophic group linear algebra.

**DEFINITION 3.2.17:** *Let V be a neutrosophic group linear algebra over the group G. Let P be a proper neutrosophic subset of V. P is just a set and it is not a closed structure with respect to addition. If P is a neutrosophic group vector space over G then we call P to be a neutrosophic pseudo group vector subspace of V over G.*

We will illustrate this by some simple examples.

*Example 3.2.43:* Let
$$V = \left\{ \begin{pmatrix} a_1 & a_2 \\ 0 & 0 \end{pmatrix}, \begin{pmatrix} a_1 & a_2 \\ b_1 & b_2 \end{pmatrix}, \begin{pmatrix} 0 & 0 \\ b_1 & b_2 \end{pmatrix} \middle| a_1, a_2, b_1, b_2 \in QI \right\}$$

be a neutrosophic group linear algebra over the group $G = ZI$. Let

$$W = \left\{ \begin{pmatrix} a_1 & a_2 \\ 0 & 0 \end{pmatrix}, \begin{pmatrix} 0 & 0 \\ b_1 & b_2 \end{pmatrix} \middle| a_1, a_2, b_1, b_2 \in QI \right\} \subseteq V;$$

W is only a neutrosophic group vector space over the group G. Thus W is a pseudo neutrosophic group vector subspace of V over the group G.

*Example 3.2.44:* Let

$$V = \left\{ \begin{pmatrix} a_1 & a_2 \\ a_3 & 0 \\ 0 & 0 \end{pmatrix}, \begin{pmatrix} 0 & 0 \\ 0 & a_4 \\ 0 & a_5 \end{pmatrix}, \begin{pmatrix} 0 & 0 \\ 0 & 0 \\ a_6 & 0 \end{pmatrix}, \begin{pmatrix} a_1 & a_2 \\ a_3 & a_4 \\ a_6 & a_5 \end{pmatrix} \middle| \begin{array}{l} a_i \in QI \\ 1 \le i \le 6 \end{array} \right\}$$

be a neutrosophic group linear algebra over the group $G = QI$. Take

$$W = \left\{ \begin{pmatrix} 0 & 0 \\ 0 & 0 \\ a_6 & 0 \end{pmatrix}, \begin{pmatrix} 0 & 0 \\ 0 & a_4 \\ 0 & a_5 \end{pmatrix} \middle| a_6, a_4, a_5 \in QI \right\} \subseteq V,$$

W is a pseudo neutrosophic group vector subspace of V over the group G.

*Example 3.2.45:* Let

$$V = \left\{ \begin{pmatrix} a & b \\ c & d \end{pmatrix}, \begin{pmatrix} a & 0 \\ 0 & b \end{pmatrix}, \begin{pmatrix} 0 & c \\ d & 0 \end{pmatrix} \middle| a, b, c, d \in Z_7 I \right\}$$

be a neutrosophic group linear algebra over $Z_7 I = G$.

$$W = \left\{ \begin{pmatrix} a & 0 \\ 0 & b \end{pmatrix}, \begin{pmatrix} 0 & c \\ d & 0 \end{pmatrix} \middle| a, b, c, d \in Z_7 I \right\} \subseteq V$$

is a pseudo neutrosophic group vector space over the group G.

**Chapter Four**

# NEUTROSOPHIC FUZZY SET LINEAR ALGEBRA

In this chapter we introduce the new notion of neutrosophic set fuzzy linear algebra, neutrosophic semigroup fuzzy linear algebra and neutrosophic group fuzzy linear algebra.

Recall, as fuzzy vector space $(V, \eta)$ or $\eta V$ is an ordinary vector space $V$ over a field $F$ with a map $\eta: V \to [0,1]$ satisfying of following conditions:

(1) $\eta(a + b) \geq \min\{\eta(a), \eta(b)\}$
(2) $\eta(-a) = \eta(a)$
(3) $\eta(0) = 1$
(4) $\eta(ra) \geq \eta(a)$

for all $a, b \in V$ and $r \in F$ where $F$ is a field.

We now define the notion of neutrosophic set fuzzy linear algebra.

**DEFINITION 4.1:** *Let V be a neutrosophic set linear algebra over the set S. We say V with a map $\eta$ is a neutrosophic fuzzy set linear algebra if*
$$\eta : V \to \langle [0, 1] \cup [0, I] \rangle = N([0, 1])$$

*(Here $N([0, 1]) = \{a + bI \mid a, b \in [0, I]\}$). ($N([0, 1])$ will be known as the fuzzy neutrosophic set or neutrosophic fuzzy set), such that $\eta (a + b) \geq \min(\eta (a), \eta(b))$ for all $a, b \in V$ and $\eta (I) = I$ and is denoted by $V\eta$ or $\eta V$ or $V(\eta)$.*

*Since we known in the neutrosophic set vector space V merely we take V to be a set but in case of neutrosophic set linear algebra we assume V is closed with respect to some operation usually denoted as '+' so the additional condition $\eta(a + b) \geq \min(\eta(a), \eta(b))$ is essential for every $a, b \in V$.*

We will illustrate this situation by some examples.

***Example 4.1:*** Let $V = Q^+I$ be a neutrosophic set linear algebra over the set $S = Z^+I$.
Define $\eta: V \to N([0, 1])$

$$\eta(x) = \begin{cases} I & \text{if } x = aI \\ I + \dfrac{1}{a+b} & \text{if } x = a + bI \text{ and } a + b > 1 \\ I + 1 & \text{if } x = a + bI \text{ and } a + b \leq 1 \end{cases}$$

Clearly $V\eta$ is a neutrosophic set fuzzy linear algebra.

***Example 4.2:*** Let

$$V = \left\{ \begin{pmatrix} aI & bI \\ cI & dI \end{pmatrix} \middle| a, b, c, d \in Z^+I \right\}$$

be a neutrosophic set linear algebra over the set $S = 10Z^+I$.
Define $\eta: V \to N([0, 1])$ as follows:

$$\eta \begin{pmatrix} aI & bI \\ cI & dI \end{pmatrix} = \begin{cases} I + \dfrac{1}{a} & \text{if } a \geq b \\ I + \dfrac{1}{b} & \text{if } b \geq d \\ I + \dfrac{1}{c} & \text{if } c \geq a \\ I + \dfrac{1}{d} & \text{if } a \geq d \end{cases}$$

We see $V_\eta$ is a neutrosophic set fuzzy linear algebra.

**DEFINITION 4.2:** *Let V be a neutrosophic set vector space over the set S. Let $W \subseteq V$ be a neutrosophic set vector subspace of V defined over the set S. If $\eta : W \to N([0, 1])$ then $W_\eta$ is called the neutrosophic fuzzy set vector subspace of V.*

We now proceed onto define the notion of neutrosophic fuzzy set linear subalgebra.

**DEFINITION 4.3:** *Let V be a neutrosophic set linear algebra over the set S. Suppose W is a neutrosophic set linear subalgebra of V over S. Let $\eta: W \to N([0, 1])$. $\eta W$ is a neutrosophic set fuzzy linear subalgebra if $\eta(a+b) \geq \min \{\eta(a), \eta(b)\}$ for $a, b \in W$.*

Now we proceed onto define the notion of neutrosophic fuzzy semigroup vector spaces.

**DEFINITION 4.4:** *A neutrosophic fuzzy semigroup vector space or a fuzzy neutrosophic semigroup vector space $(V, \eta)$ or $V\eta$ where V is an ordinary neutrosophic semigroup vector space over the semigroup S; with a map $\eta : V \to N([0, 1])$ satisfying the following condition; $\eta(ra) \geq \eta(a)$ for all $a \in V$ and $r \in S$.*

We will illustrate this by some examples.

***Example 4.3:*** Let $V = \{(a_1I, a_2I, a_3I, a_4I, a_5I, a_6I, a_7I) \mid a_i \in Z^+I, 1 \leq i \leq 7\}$ be a neutrosophic set vector space over the set $S = 5Z^+I$.

Define $\eta : V \to N([0, 1])$ by

$$\eta(a_1I, a_2I, a_3I, a_4I, a_5I, a_6I, a_7I) = \frac{1}{a_1 + \ldots + a_7} + I$$

for every $(a_1I, a_2I, a_3I, a_4I, a_5I, a_6I, a_7I)$ in V. $V_\eta$ is a neutrosophic set fuzzy vector space.

***Example 4.4:*** Let

$$V = \left\{ \begin{bmatrix} a_1I & a_2I \\ a_3I & a_4I \\ a_5I & a_6I \\ a_7I & a_8I \\ a_9I & a_{10}I \end{bmatrix} \mid a_i \in QI, 1 \leq i \leq 10 \right\}$$

be a neutrosophic set vector space over a set $S = Z^+I$.
Define $\eta : V \to N([0, 1])$ by

$$\eta \begin{bmatrix} a_1I & a_2I \\ a_3I & a_4I \\ a_5I & a_6I \\ a_7I & a_8I \\ a_9I & a_{10}I \end{bmatrix} = \begin{cases} \frac{1}{|a_i|} + I & \text{if } |a_i| \in Z \\ I+1 & \text{otherwise} \\ 0 & \text{if } a_i = 0; i = 1, 2, \ldots, 10 \end{cases}$$

$V_\eta$ is a neutrosophic set fuzzy vector space.

***Example 4.5:*** Let

$$V = \left\{ \begin{pmatrix} aI & bI \\ cI & dI \end{pmatrix} \mid aI, bI, cI, dI \in Z^+I \cup \{0\} \right\}$$

be a neutrosophic set linear algebra over the set $S = 3Z^+I \cup \{0\}$. Let $\eta : V \to N([0, 1])$. Define

$$\eta \begin{pmatrix} aI & bI \\ cI & dI \end{pmatrix} = \begin{cases} I + \dfrac{1}{a+d} & \text{if } a+d \neq 0 \\ I + \dfrac{1}{b+c} & \text{if } b+c \neq 0 \\ 1 & \text{if } a+d = 0 \text{ and if } b+c = 0 \end{cases}$$

*Example 4.6:* Let $V = \{(a_1I, a_2I, a_3I, a_4I, a_5I, a_6I) \mid a_i I \in QI; 1 \leq i \leq 6\}$ be a neutrosophic set linear algebra over the set $S = Z^+I \cup \{0\}$. Let $W = \{(a_1I, a_2I, a_3I, a_4I, a_5I, a_6I) \mid a_i I \in Z^+I \cup \{0\}\}$ be a neutrosophic set linear subalgebra of V over S.

Define $\eta: W \to N([0, 1])$ by

$$\eta (a_1I, a_2I, \ldots, a_6I) = \begin{cases} \dfrac{I+1}{5} & \text{if atleast one } a_i \neq 0 \text{ or } 1 \\ I & \text{if all } a_i \text{'s are } 1 \\ 1 & \text{if all } a_i \text{'s are } 0; 1 \leq i \leq 6 \end{cases}$$

It is easily verified that $W_\eta$ or $\eta (W)$ is a neutrosophic fuzzy set linear subalgebra.

*Example 4.7:* Let

$$V = \left\{ \begin{pmatrix} a_1I & a_2I & a_3I \\ a_4I & a_5I & a_6I \end{pmatrix} \middle| a_iI \in RI; 1 \leq i \leq 5 \right\}$$

be a neutrosophic set linear algebra over the set $S = Z^+I \cup \{0\}$. Let

$$W = \left\{ \begin{pmatrix} a_1I & a_2I & a_3I \\ a_4I & a_5I & a_6I \end{pmatrix} \middle| a_iI \in Z^+I \cup \{0\}; 1 \leq i \leq 6 \right\} \subseteq V$$

be a neutrosophic set linear subalgebra over S.

Define $\eta : W \to N([0, 1])$ by

$$\eta\begin{pmatrix} a_1I & a_2I & a_3I \\ a_4I & a_5I & a_6I \end{pmatrix} = \begin{cases} I + \dfrac{1}{a_1 + a_2 + a_3} & \text{if } a_1 + a_2 + a_3 \neq 0 \\ I + \dfrac{1}{a_4 + a_5 + a_6} & \text{if } a_4 + a_5 + a_6 \neq 0 \\ 1 & \text{if } a_i = 0; 1 \leq i \leq 6 \end{cases}$$

$W_\eta$ is a neutrosophic set fuzzy linear subalgebra.

*Example 4.8:* Let

$$V = \left\{(a_1I, a_2I, a_3I), \begin{bmatrix} b_1I & b_2I \\ b_3I & b_4I \\ b_5I & b_6I \end{bmatrix} \mid a_i, b_j \in Z^+I \cup \{0\};\right.$$

$1 \leq i \leq 3; 1 \leq j \leq 6\}$ be a neutrosophic set vector space over the set $S = 3Z^+I \cup \{0\}$. Let $W = \{(a_1I, a_2I, a_3I) \mid a_iI \in Z^+I \cup \{0\}; 1 \leq i \leq 3\}$. Define $\eta: W \to N([0,1])$

$$\eta(a_1I, a_2I, a_3I) = I + \begin{cases} I + \dfrac{1}{a_i} & \text{if } a_i \neq 0; 1 \leq i \leq 3 \\ 1 & \text{if } a_i = 0; 1 \leq i \leq 3 \end{cases}$$

$W_\eta$ is a neutrosophic set fuzzy vector subspace.

*Example 4.9:* Let

$$V = \left\{\begin{pmatrix} a_1I & a_2I & a_3I \\ 0 & 0 & 0 \end{pmatrix}, \begin{pmatrix} 0 & 0 & 0 \\ b_1I & b_2I & b_3I \end{pmatrix},\right.$$

$(a_1I, a_2I, a_3I, a_4I) \mid b_jI, a_iI \in Z^+I \cup \{0\}\}$

be a neutrosophic set vector space over the set $S = 5Z^+I \cup \{0\}$. Let

$$W = \left\{ \begin{pmatrix} 0 & 0 & 0 \\ b_1I & b_2I & b_3I \end{pmatrix}, (a_1I, a_2I, a_3I, a_4I) \mid \right.$$

$$\left. b_jI, a_iI \in Z^+I \cup \{0\}, 1 \le i \le 4 \text{ and } 1 \le j \le 3 \right\} \subseteq V$$

be a neutrosophic set vector subspace of V over S. Define

$$\eta \begin{pmatrix} 0 & 0 & 0 \\ b_1I & b_2I & b_3I \end{pmatrix} = \begin{cases} I + \dfrac{1}{b_i} & \text{if } b_i \ne 0; 1 \le i \le 3 \\ 1 & \text{if } b_i = 0; 1 \le i \le 3 \end{cases}$$

and

$$\eta(a_1I, a_2I, a_3I, a_4I) = \begin{cases} I + \dfrac{1}{\sum\limits_i a_i} & \text{if } \sum\limits_i a_i \ne 0 \\ 1 & \text{if } \sum\limits_i a_i = 0 \end{cases}$$

$W_\eta$ is a neutrosophic set fuzzy vector subspace.

*Example 4.10:* Let

$$V = \left\{ \begin{pmatrix} a_1I & a_2I \\ 0 & a_3I \end{pmatrix}, \begin{pmatrix} 0 & b_1I \\ b_2I & 0 \end{pmatrix} \mid a_i, b_j \in Z^+I \cup \{0\}; \right.$$

$\left. 1 \le i \le 3 \text{ and } 1 \le j \le 2 \right\}$ be a neutrosophic semigroup vector space over the semigroup $S = 3Z^+I \cup \{0\}$.

Let $\eta : V \to N([0,1])$ be defined such that

$$\eta \begin{pmatrix} a_1I & a_2I \\ 0 & a_3I \end{pmatrix} = \begin{cases} I + \dfrac{1}{a_1 + a_2 + a_3} & \text{if } a_1 + a_2 + a_3 \ne 0 \\ 1 & \text{if } a_1 + a_2 + a_3 = 0 \end{cases}$$

and

$$\eta \begin{pmatrix} 0 & b_1 I \\ b_2 I & 0 \end{pmatrix} = \begin{cases} I + \dfrac{1}{b_1 + b_2} & \text{if } b_1 + b_2 \neq 0 \\ 1 & \text{if } b_1 + b_2 = 0 \end{cases}$$

$V_\eta$ is a fuzzy neutrosophic semigroup vector space or neutrosophic fuzzy semigroup vector space or neutrosophic semigroup fuzzy vector space.

*Example 4.11:* Let

$$V = \left\{ \begin{pmatrix} a_1 I \\ a_2 I \\ a_3 I \\ a_4 I \end{pmatrix}, (b_1 I, b_2 I, b_3 I), \begin{pmatrix} c_1 I & c_2 I & c_3 I \\ 0 & c_4 I & 0 \end{pmatrix} \right|$$

$a_i I, b_j I, c_k I \in Z^+ I \cup \{0\}$; $1 \leq i \leq 4$, $1 \leq j \leq 3$ and $1 \leq k \leq 4\}$

be a neutrosophic semigroup vector space over the semigroup $S = Z^+ I \cup \{0\}$. Define $\eta : V \to N([0, 1])$ as

$$\eta \begin{pmatrix} a_1 I \\ a_2 I \\ a_3 I \\ a_4 I \end{pmatrix} = \begin{cases} I + \dfrac{1}{a_1 + a_2} & \text{if } a_1 + a_2 \neq 0 \\ I + \dfrac{1}{a_3} & \text{if } a_3 \neq 0 \\ I + \dfrac{1}{a_4} & \text{if } a_4 \neq 0 \\ 1 & \text{if } a_i = 0; 1 \leq i \leq 4 \end{cases}$$

$$\eta(b_1 I, b_2 I, b_3 I) = \begin{cases} I + \dfrac{1}{b_i} & \text{if } b_i \neq 0; 1 \leq i \leq 3 \\ 1 & \text{if } b_i = 0; 1 \leq i \leq 3 \end{cases}$$

$$\eta \begin{bmatrix} c_1I & c_2I & c_3I \\ 0 & c_4I & 0 \end{bmatrix} = \begin{cases} \dfrac{1}{\sum_i c_i} + I & \text{if } \sum_i c_i \neq 0 \\ 1 & \text{if } \sum_i c_i = 0 \end{cases}$$

$V_\eta$ is a neutrosophic semigroup fuzzy vector space.

Now we proceed onto define the notion of neutrosophic semigroup fuzzy vector subspace.

**DEFINITION 4.5:** *Let V be a neutrosophic semigroup vector space over the semigroup S. Let W $\subseteq$ V be a neutrosophic semigroup vector subspace of V over S. We say $W_\eta$ is a neutrosophic semigroup fuzzy vector subspace if*
$$\eta : W \to N([0,1])$$
*such that*
$$\eta(I) = I$$
$$\eta(rx) \geq \eta(x)$$
*for all x, y $\in$ W and r $\in$ S.*

We will illustrate this situation by some simple examples.

*Example 4.12:* Let

$$V = \left\{ \begin{pmatrix} a_1 & a_2 & a_3 & a_4 \\ 0 & 0 & 0 & 0 \end{pmatrix}, \begin{pmatrix} 0 & 0 & 0 & 0 \\ b_1 & b_2 & b_3 & b_4 \end{pmatrix} \mid \right.$$
$$a_i, b_j \in QI; \ 1 \leq i \leq 4, \ 1 \leq j \leq 4 \}$$

be a neutrosophic semigroup vector space over the semigroup S $= Z^+I \cup \{0\}$. Let

$$W = \left\{ \begin{pmatrix} 0 & 0 & 0 & 0 \\ b_1 & b_2 & b_3 & b_4 \end{pmatrix} \mid b_i \in QI, \ 1 \leq i \leq 4 \right\} \subseteq V$$

be a neutrosophic semigroup vector subspace of V over the semigroup S.

Define $\eta: W \to N([0,1])$ by

$$\eta\begin{pmatrix} 0 & 0 & 0 & 0 \\ b_1 & b_2 & b_3 & b_4 \end{pmatrix} = \begin{cases} I + \dfrac{1}{b_i} & \text{if } b_i \neq 0;\ 1 \leq i \leq 4 \\ 1 & \text{if } b_i = 0;\ 1 \leq i \leq 4 \end{cases}$$

$W_\eta$ is a neutrosophic semigroup fuzzy vector subspace.

*Example 4.13:* Let

$$V = \left\{ \sum_{i=1}^{n} a_i I x_i \,\middle|\, a_i I \in Z^+ I \cup \{0\};\ 1 \leq i \leq n \right\}$$

be a neutrosophic semigroup vector space over the semigroup $S = Z^+ I \cup \{0\}$.
Let

$$W = \left\{ \sum_{i=1}^{n} a_i I x^{2i} \,\middle|\, a_i I \in 2Z^+ I \cup \{0\} \right\} \subseteq V.$$

W is a neutrosophic semigroup vector subspace of V over the semigroup S.
Define $\eta : W \to N([0,1])$ by

$$\eta\left( \sum_{i=1}^{n} a_i I x^{2i} \right) = \begin{cases} \dfrac{1}{\sum_{i=1}^{n} a_i} + I & \text{if } \sum_i a_i \neq 0 \\ 1 & \text{if } \sum_i a_i = 0 \end{cases}$$

$W_\eta$ is a neutrosophic semigroup fuzzy vector subspace.

Now we define neutrosophic semigroup fuzzy linear algebra.

**DEFINITION 4.6:** *Let V be a neutrosophic semigroup linear algebra over the semigroup S. We say $V_\eta$ or $\eta V$ or $V(\eta)$ is a*

*neutrosophic semigroup fuzzy linear algebra if $\eta : V \to N([0,1])$; such that*

$$\eta(x+y) \geq min(\eta(x), \eta(y));$$
$$\eta(rx) \geq \eta(x)$$

*for $r \in S$ and $y, x \in V$.*

**Example 4.14:** Let $V = \{(Z^+I \cup \{0\}) \times (Z^+I \cup \{0\}) \times (Z^+I \cup \{0\}) \times (Z^+I \cup \{0\})\}$ be a neutrosophic semigroup linear algebra over the semigroup $S = (2Z^+I \cup \{0\})$. Define $\eta : V \to N([0,1])$ where

$$\eta(a_1I, a_2I, a_3I, a_4I) = \begin{cases} I + \dfrac{1}{a_i} & \text{if } a_i \neq 0; 1 \leq i \leq 4 \\ 1 & \text{if } a_i = 0; 1 \leq i \leq 4 \end{cases}$$

$V_\eta$ is a neutrosophic semigroup fuzzy linear algebra.

**Example 4.15:** Let $V = \{(a_1I, \ldots, a_{10}I) | a_iI \in Z_{11}I; 1 \leq i \leq 10\}$ be a neutrosophic semigroup linear algebra over the semigroup $S = Z_{11}I$. Define $\eta : V \to N([0,1])$

$$\eta(a_1I, \ldots, a_{10}I) = \begin{cases} \dfrac{1}{a_i} + I & \text{if } a_i \neq 0; 1 \leq i \leq 10 \\ 1 & \text{if } a_i = 0; 1 \leq i \leq 10 \end{cases}$$

$V_\eta$ is a neutrosophic semigroup fuzzy linear algebra.

Now we can define neutrosophic semigroup fuzzy linear subalgebra as in case of neutrosophic semigroup fuzzy vector subspaces. We leave this task to the interested reader. However we will illustrate that by some examples.

**Example 4.16:** Let

$$V = \left\{ \begin{pmatrix} a_1I & a_2I \\ 0 & a_3I \end{pmatrix}, \begin{pmatrix} 0 & a_4I \\ a_5I & 0 \end{pmatrix}, \begin{pmatrix} b_1I & b_2I \\ b_3I & b_4I \end{pmatrix} \right| $$
$a_iI, b_jI \in QI$ for $1 \leq i \leq 5; 1 \leq j \leq 4\}$

be a neutrosophic semigroup linear algebra over the semigroup $S = Z^+I \cup \{0\}$.

Take

$$W = \left\{ \begin{pmatrix} a_1I & a_2I \\ 0 & a_3I \end{pmatrix} \middle| a_iI \in QI; 1 \le i \le 3 \right\} \subseteq V.$$

W is a neutrosophic semigroup linear subalgebra of V over the semigroup S.

Define $\eta : W \to N([0,1])$ by

$$\eta \begin{pmatrix} a_1I & a_2I \\ 0 & a_3I \end{pmatrix} = \begin{cases} \dfrac{1}{\sum a_i} + I & \text{if } \sum a_i \ge 0 \\ I & \text{if } \sum a_i < 0 \\ 1 & \text{if } \sum a_i = 0 \end{cases}$$

W$\eta$ is a neutrosophic semigroup fuzzy linear subalgebra.

*Example 4.17:* Let $V = \{QI [x] \mid$ all polynomials in the variable x with coefficients from QI$\}$ be a neutrosophic semigroup linear algebra over the semigroup $S = Z^+I \cup \{0\}$.

Let $W = \{ZI [x] \mid$ all polynomials in the variable x with coefficients from ZI$\} \subseteq V$.

Define $\eta : V \to N([0,1])$ by

$$\eta(p(x)) = \begin{cases} I + \dfrac{1}{\deg(p(x))} & \text{if } \deg(p(x)) \ge 1 \\ I & \text{if } \deg(p(x)) = \text{constant} \\ 1 & \text{if } p(x) = 0 \end{cases}$$

where constant is a neutrosophic number.

$W_\eta$ is a neutrosophic semigroup fuzzy linear subalgebra.

Now we proceed onto define the notion of neutrosophic group fuzzy linear algebra. Just we recall the definition of group fuzzy linear algebra.

**DEFINITION 4.7:** *Let V be a group linear algebra over the group G. Let $\eta : V \to [0,1]$ be such that*
$$\eta (a + b) \geq min( \eta(a), \eta (b))$$
$$\eta(-a) = \eta (a)$$
$$\eta(0) = 1$$
$$\eta(ra) = \eta (a)$$
*for all $a, b \in V$ and $r \in G$; we call $V_\eta$ the group fuzzy linear algebra.*

**DEFINITION 4.8:** *Let V be a neutrosophic group vector space over the group G. Let $\eta : V \to N([0,1])$ be such that $\eta (ra) \geq \eta (a)$ for all $a \in V$ and $\eta (I) = I$ for all $a \in V$. We call $V_\eta$ or $\eta V$ or $V(\eta)$ to be the neutrosophic group fuzzy vector space.*

We will illustrate this situation by some examples.

*Example 4.18:* Let

$$V = \left\{ \begin{pmatrix} a_1I \\ b_1I \\ c_1I \end{pmatrix}, \begin{pmatrix} a_1I & b_1I \\ 0 & c_1I \end{pmatrix}, (a_1I, b_1I, c_1I) \mid a_1I, b_1I, c_1I \in QI \right\}$$

be a neutrosophic group vector space over the group $G = ZI$.

Define for $x \in V$

$$\eta (x) = \begin{cases} I + \dfrac{1}{a+b+c} & \text{if } a+b+c \geq 1 \\ I & \text{if } a+b+c < 1 \\ 1 & \text{if } x \text{ has no neutrosophic component} \end{cases}$$

$V_\eta$ is a neutrosophic group fuzzy vector space.

*Example 4.19:* Let
$$V = \left\{\sum_{i=1}^{n} a_i Ix^{2i}, \sum_{i=1}^{m} a_i Ix^{5i} \mid a_i \in ZI\right\}$$

be a neutrosophic group vector space over the group $G = ZI$.

Define
$$\eta(p(x)) = \begin{cases} I + \dfrac{1}{\deg(p(x))} & \text{if } \deg(p(x)) \geq 1 \\ I & \text{if } \deg(p(x)) \text{ is neutrosophic constant} \\ 1 & \text{if } \deg(p(x)) \text{ is an interger} \end{cases}$$

$V_\eta$ is a neutrosophic group fuzzy vector space.

*Example 4.20:* Let
$$V = \left\{ \begin{bmatrix} a_1 I & a_2 I \\ a_3 I & a_4 I \\ a_5 I & a_6 I \\ a_7 I & a_8 I \end{bmatrix} \middle| a_i I \in QI; 1 \leq i \leq 8 \right\}$$

be a neutrosophic group linear algebra over the group $G = ZI$.
Define $\eta : V \to N([0,1])$ by

$$\eta \begin{bmatrix} a_1 I & a_2 I \\ a_3 I & a_4 I \\ a_5 I & a_6 I \\ a_7 I & a_8 I \end{bmatrix} =$$

$$\begin{cases} I + \dfrac{1}{2} & \text{if at least one of the } a_i \text{ is non zero } 1 \leq i \leq 8 \\ 1 & \text{of all the entrires in the matrix is zero} \end{cases}$$

$V_\eta$ is a neutrosophic group fuzzy linear algebra.

*Example 4.21:* Let

$$V = \left\{\sum_{i=1}^{n} a_i I x^i \mid a_i I \in QI;\ 1 \leq i \leq n\right\}$$

be a neutrosophic group linear algebra over the group $G = ZI$. Define $\eta : V \to N([0,1])$

$$\eta\left(\sum_{i=1}^{n} a_i I x^i\right) = \begin{cases} I + \dfrac{1}{\deg(p(x))} & \text{if } \deg(p(x)) \geq 1 \\ I & \text{if } \deg(p(x)) \text{ is a neutrosophic constant} \\ 1 & \text{if } \deg(p(x)) \text{ is a constant, i.e., } p(x) = 0 \end{cases}$$

$V_\eta$ or $\eta V$ is a neutrosophic group fuzzy linear algebra.

Next we proceed onto define fuzzy substructures.

**DEFINITION 4.9:** *Let V be a neutrosophic group vector space over the group G. Let $W \subseteq V$ be a neutrosophic group vector subspace of V over G. Define $\eta : W \to N([0,1])$ as $\eta(ra) \geq \eta(a)$ for all $r \in G$ and $a \in W$. We call $W_\eta$ or $\eta W$ to be a neutrosophic group fuzzy vector subspace.*

We will illustrate this situation by some simple examples.

*Example 4.22:* Let

$$V = \left\{ \begin{pmatrix} a_1 I & a_2 I & a_3 I \\ 0 & a_4 I & 0 \end{pmatrix}, \begin{pmatrix} 0 & a_5 I & 0 \\ a_6 I & 0 & a_7 I \end{pmatrix} \mid a_i I \in QI;\ 1 \leq i \leq 7 \right\}$$

be a neutrosophic group vector space over the group $G = ZI$.

Let

$$W = \left\{ \begin{pmatrix} 0 & a_5I & 0 \\ a_6I & 0 & a_7I \end{pmatrix} \mid a_iI \in QI, i = 5, 6, 7 \right\} \subseteq V$$

be a neutrosophic group vector subspace of V over G.
Define $\eta : W \to N([0,1])$ by

$$\eta \begin{pmatrix} 0 & a_5I & 0 \\ a_6I & 0 & a_7I \end{pmatrix} = \begin{cases} I + \dfrac{1}{5} & \text{if at least one } a_i \text{ is non zero; } i = 5, 6, 7 \\ 1 & \text{if } a_i\text{'s are zero; } i = 5, 6, 7 \end{cases}$$

$W_\eta$ is a neutrosophic group fuzzy vector subspace of V.

*Example 4.23:* Let

$$V = \left\{ (a_1I, a_2I, a_3I), \begin{bmatrix} a_1I & a_2I \\ a_3I & a_4I \\ a_5I & a_6I \\ a_7I & a_8I \end{bmatrix}, \begin{bmatrix} a_1I & a_2I & a_3I \\ a_4I & a_5I & a_6I \end{bmatrix} \mid \right.$$
$$a_iI \in QI; 1 \le i \le 8 \}$$

be a neutrosophic group vector space over the group G = ZI. Let

$$W = \left\{ \begin{bmatrix} a_1I & a_2I & a_3I \\ a_4I & a_5I & a_6I \end{bmatrix} \right\} \subseteq V$$

be a subspace of V.
Define $\eta: W \to N([0, 1])$ by

$$\begin{bmatrix} a_1I & a_2I & a_3I \\ a_4I & a_5I & a_6I \end{bmatrix} = \begin{cases} I + \dfrac{1}{2} & \text{if } a_i \ne 0 \text{ for } i = 1, 2, \ldots, 6 \\ 1 & \text{if } a_i = 0; 1 \le i \le 6 \end{cases}$$

$W\eta$ is a neutrosophic group fuzzy vector subspace. Now as in case of neutrosophic group fuzzy vector subspace we can define neutrosophic group fuzzy linear subalgebra.

We will however give some examples of neutrosophic group fuzzy linear subalgebras.

*Example 4.24:* Let

$$V = \left\{ \begin{pmatrix} a_1I & a_2I & a_3I \\ a_4I & a_5I & a_6I \\ a_7I & a_8I & a_9I \end{pmatrix} \middle| a_iI \in QI; 1 \leq i \leq 9 \right\}$$

be a neutrosophic group linear algebra over the group $G = ZI$. Let

$$W = \left\{ \begin{pmatrix} a_1I & a_2I & a_3I \\ a_4I & a_5I & a_6I \\ a_7I & a_8I & a_9I \end{pmatrix} \middle| a_iI \in ZI; 1 \leq i \leq 9 \right\} \subseteq V$$

be a neutrosophic group linear subalgebra over the group $G = ZI$.
Define $\eta : W \to N([0,1])$ by

$$\eta \begin{pmatrix} a_1I & a_2I & a_3I \\ a_4I & a_5I & a_6I \\ a_7I & a_8I & a_9I \end{pmatrix} = \begin{cases} I + \dfrac{1}{|a_i|} & \text{if } a_i \neq 0; 1 \leq i \leq 9 \\ 1 & \text{if } a_i = 0; 1 \leq i \leq 9 \end{cases}$$

$W_\eta$ is a neutrosophic group fuzzy linear subalgebra.

*Example 4.25:* Let

$$V = \left\{ \sum_{i=0}^{n} a_i I x^i \middle| a_i I \in QI \right\}$$

be a neutrosophic group linear over the group G = ZI. Let

$$W = \left\{ \sum_{i=0}^{n} a_i I x^i \mid a_i I \in ZI \right\} \subseteq V$$

be a neutrosophic group linear subalgebra of V.
Define $\eta : W \to N([0,1])$ by

$$\eta \left( \sum_{i=1}^{n} a_i I x^i \right) = \eta(p(x)) =$$

$$\begin{cases} I + \dfrac{1}{\deg(p(x))} & \text{if } p(x) \text{ is not a constant} \\ 1 & \text{if } p(x) \text{ is zero} \\ I & \text{if } p(x) \text{ is a neutrosophic integer} \end{cases}$$

$W_\eta$ is a neutrosophic group fuzzy linear subalgebra.

**DEFINITION 4.10:** *Let V be a neutrosophic group linear algebra over the group G. Let $W \subseteq V$, where W is a subgroup of V and $H \subseteq G$ be a proper subgroup of G; so that W is a neutrosophic subgroup linear subalgebra of V over the subgroup H of G. Let $\eta : W \to N([0,1])$ if $W_\eta$ is a neutrosophic group fuzzy linear algebra then we call $W\eta$ to be a neutrosophic subgroup fuzzy linear subalgebra.*

We will illustrate this by some simple example.

*Example 4.26:* Let

$$V = \left\{ \begin{pmatrix} a_1 I & a_2 I \\ a_3 I & 0 \end{pmatrix} \mid a_i \in QI, 1 \leq i \leq 3 \right\}$$

be a neutrosophic group linear algebra over the group G = ZI.
Let

$$W = \left\{ \begin{pmatrix} a_1I & a_2I \\ a_3I & 0 \end{pmatrix} \,\middle|\, a_i \in ZI,\, 1 \leq i \leq 3 \right\} \subseteq V$$

be a neutrosophic subgroup linear subalgebra of V over the subgroup $H = 3ZI \subseteq G$.
Define $\eta : W \to N([0,1])$

$$\eta \begin{pmatrix} a_1I & a_2I \\ a_3I & 0 \end{pmatrix} = \begin{cases} I + \dfrac{1}{2} & \text{if } a_i \neq 0;\ 1 \leq i \leq 3 \\ 1 & \text{if } a_i = 0;\ 1 \leq i \leq 3 \end{cases}$$

W$\eta$ is a neutrosophic subgroup fuzzy linear subalgebra.

*Example 4.27:* Let

$$V = \left\{ \begin{bmatrix} a_1I & a_2I & 0 \\ a_3I & 0 & a_4I \end{bmatrix} \,\middle|\, a_i I \in ZI,\, 1 \leq i \leq 4 \right\}$$

be a neutrosophic group linear algebra over the group $G = ZI$.
Let

$$W = \left\{ \begin{bmatrix} a_1I & a_2I & 0 \\ a_3I & 0 & a_4I \end{bmatrix} \,\middle|\, a_iI \in 5ZI,\, 1 \leq i \leq 4 \right\} \subseteq V$$

be a neutrosophic subgroup linear subalgebra over the subgroup $H = 10ZI \subseteq ZI$.

Define $\eta : W \to N([0,1])$

$$\eta \begin{bmatrix} a_1I & a_2I & 0 \\ a_3I & 0 & a_4I \end{bmatrix} = \begin{cases} I + \dfrac{1}{5} & \text{if } a_i \neq 0;\ 1 \leq i \leq 4 \\ 1 & \text{if } a_i = 0;\ 1 \leq i \leq 4 \end{cases}$$

W$\eta$ is a neutrosophic subgroup fuzzy sublinear algebra.

The importance of this structure is we do not demand neutrosophic field over which these structures are defined. Even a neutrosophic set is sufficient we know when we define fuzzy vector space or fuzzy linear algebra the field over which they are defined do not play any prominent role.

Another advantage of working with these fuzzy neutrosophic vector spaces is in most of the cases they become fuzzy equivalent.

Further our transformation to fuzzy set up demands only values from the set $N[(0,1)] = \{a + bI \mid a, b \in [0,1]\}$.

Further if we go for neutrosophic Markov process or Markov chains the probability matrix is a square matrix with positive entries from $N([0,1])$.

**Chapter Five**

# NEUTROSOPHIC SET BIVECTOR SPACES

In this chapter we introduce the notion of neutrosophic set bivector spaces and neutrosophic group bivector spaces. We enumerate some of their properties. These are useful on the study of mathematical models.

**DEFINITION 5.1:** *Let $V = V_1 \cup V_2$ where $V_1$ and $V_2$ are two distinct neutrosophic set vector spaces defined over the same set S. That is $V_1 \not\subset V_2$ and $V_2 \not\subset V_1$, we may have $V_1 \cap V_2 = \phi$ or non empty. We call V to be the neutrosophic set vector bispace or neutrosophic set bivector space over the set S.*

We will illustrate this by some examples.

*Example 5.1:* Let $V = V_1 \cup V_2$ where

$$V_1 = \left\{ \begin{pmatrix} a_1I & a_2I \\ 0 & a_3I \end{pmatrix}, \begin{pmatrix} a_4I & 0 \\ a_5I & 0 \end{pmatrix} \middle| a_iI \in Q^+I \cup \{0\}; 1 \leq i \leq 5 \right\}$$

and

$$V_2 = \left\{ (a_1I \quad a_2I \quad a_3I), \begin{pmatrix} a_1I \\ a_2I \\ a_3I \\ a_4I \end{pmatrix} \middle| a_iI \in Z^+I \right\}$$

be neutrosophic set vector spaces over the set $S = 5Z^+I$. V is a neutrosophic set bivector space over the set S.

*Example 5.2:* Let

$$V = \left\{ \begin{pmatrix} a_1I & a_2I \\ a_3I & a_4I \\ a_5I & a_6I \end{pmatrix}, \begin{pmatrix} b_1I & b_2I & b_3I \\ b_4I & b_5I & b_6I \\ b_7I & b_8I & b_9I \end{pmatrix} \middle| a_iI, b_jI \in Z_{12}I \right\} \cup$$

$\{(a_1I, a_2I, a_3I, a_4I), 0, Z^+I \mid a_iI \in Z^+I \cup \{0\}, 1 \leq i \leq 4\} = V_1 \cup V_2$ be such that $V_1$ is a neutrosophic set vector space over the set $S = \{0, 1\}$ and $V_2$ is a neutrosophic set vector space over the set $S = \{0, 1\}$. Thus V is a neutrosophic set bivector space over the set $S = \{0, 1\}$.

Now we proceed onto define substructure in neutrosophic set bivector spaces.

**DEFINITION 5.2:** *Let $V = V_1 \cup V_2$ be a neutrosophic set bivector space defined over the set S. A proper biset $W = W_1 \cup W_2 \subseteq V_1 \cup V_2 = V$, ($W_1 \subseteq V_1$ and $W_2 \subseteq V_2$) such that $W_1$ and $W_2$ are distinct and each $W_i$ is a neutrosophic set vector subspace of $V_i$ over the set S; $1 \leq i \leq 2$ is called the neutrosophic set bivector subspace of V over the set S.*

We will illustrate this definition by some examples.

*Example 5.3:* Let
$$V = V_1 \cup V_2 =$$

$$\left\{ \begin{pmatrix} a_1I \\ a_2I \\ a_3I \end{pmatrix}, (b_1I, b_2I, b_3I, b_4I, b_5I) \;\middle|\; a_iI, b_jI \in ZI; 1 \le i \le 3; 1 \le j \le 5 \right\}$$

$$\cup \left\{ \begin{pmatrix} a_1I & a_2I \\ a_3I & 0 \end{pmatrix}, ZI \;\middle|\; a_iI \in QI; 1 \le i \le 3 \right\}$$

be a neutrosophic set bivector space over the set $S = \{0, 1\}$. Take $W = W_1 \cup W_2 = \{[b_1I, b_2I, b_3I, b_4I, b_5I| b_iI \in ZI; 1 \le i \le 5\} \cup \{ZI\} \subseteq V_1 \cup V_2$; W is a neutrosophic set bivector subspace of V over the set $S = \{0, 1\}$.

*Example 5.4:* Let
$$V = V_1 \cup V_2 = \{3Z^+I, 5ZI\} \cup$$

$$\left\{ \begin{pmatrix} a_1I & a_2I & a_3I \\ a_4I & a_5I & a_6I \end{pmatrix}, (b_1I, b_2I) \;\middle|\; a_iI, b_jI \in 5Z^+I \cup \{0\}; 1 \le i \le 6; 1 \le j \le 2 \right\}$$

be a neutrosophic set bivector space over the set $S = 10Z^+I$. Take $W = W_1 \cup W_2 = \{3Z^+I\} \cup \{(b_1I, b_2I) \mid b_1I, b_2I \in Z^+I\} \subseteq V_1 \cup V_2$, W is a neutrosophic set bivector subspace of V over the set $S = 10Z^+I$.

**DEFINITION 5.3:** *Let $V = V_1 \cup V_2$ be a neutrosophic set bivector space over the set S. Let $X = X_1 \cup X_2 \subseteq V_1 \cup V_2 = V$, we say X is a bigenerating biset of V if $X_1$ is the generating neutrosophic set vector space of $V_1$ over S and $X_2$ is the generating neutrosophic set vector space of $V_2$ over S. The number of elements in $X = X_1 \cup X_2$ is the bidimension of V and is denoted by $|X| = (|X_1|, |X_2|)$ or $|X_1| \cup |X_2|$.*

We will illustrate this by some examples.

*Example 5.5:* Let
$$V = V_1 \cup V_2 =$$

$$\{(a, b, c), \begin{bmatrix} x_1 & x_2 & x_3 \\ x_4 & x_5 & x_6 \end{bmatrix} \mid a, b, c, x_i \in N(Q); 1 \leq i \leq 6\}$$

$$\cup \{(a, a, a, a, a), \begin{bmatrix} b \\ b \\ b \\ b \\ b \end{bmatrix} \mid a, b \in N(Z)\}$$

be a neutrosophic set bivector space over the set $S = N(Z)$.

$$X = \{(a, b, c), \begin{bmatrix} x_1 & x_2 & x_3 \\ x_4 & x_5 & x_6 \end{bmatrix} \mid a, b, c, x_i \in N(Q); 1 \leq i \leq 6\}$$

$$\cup \{(1, 1, 1, 1, 1), \begin{bmatrix} 1 \\ 1 \\ 1 \\ 1 \\ 1 \end{bmatrix}\}$$

$$= X_1 \cup X_2$$

is the bigenerating biset of V over S. Clearly $|X| = (\infty, 2)$.

*Example 5.6:* Let

$$V = \left\{ \begin{pmatrix} 1 & 1 \\ 0 & 0 \end{pmatrix}, \begin{pmatrix} 0 & 0 \\ 1 & 1 \end{pmatrix}, \begin{pmatrix} 1 & 0 \\ 0 & 0 \end{pmatrix}, \begin{pmatrix} 0 & 1 \\ 0 & 0 \end{pmatrix}, \begin{pmatrix} 1 & 0 \\ 1 & 0 \end{pmatrix}, \begin{pmatrix} 0 & 1 \\ 0 & 1 \end{pmatrix}, \right.$$

$$\left\{ \begin{pmatrix} 1 & 0 \\ 0 & 1 \end{pmatrix}, \begin{pmatrix} 0 & 1 \\ 1 & 0 \end{pmatrix}, \begin{pmatrix} 0 & 0 \\ 0 & 0 \end{pmatrix} \right\} \cup$$

$$\{(0\ 0\ 1), (0\ 0\ 0), (1\ 1\ 0\ 0), (0\ 0\ 0\ 0),$$
$$(1\ 1\ 1), (1\ 0\ 1), (0\ 0\ 1\ 1), (1\ 0\ 1\ 0)\}$$
$$= V_1 \cup V_2$$

be a neutrosophic set bivector space over the set $S = \{0, 1\}$.

$$X = X_1 \cup X_2 =$$

$$\left\{ \begin{pmatrix} 1 & 1 \\ 0 & 0 \end{pmatrix}, \begin{pmatrix} 0 & 0 \\ 1 & 1 \end{pmatrix}, \begin{pmatrix} 1 & 0 \\ 0 & 0 \end{pmatrix}, \begin{pmatrix} 0 & 1 \\ 0 & 0 \end{pmatrix}, \begin{pmatrix} 1 & 0 \\ 1 & 0 \end{pmatrix}, \begin{pmatrix} 0 & 1 \\ 0 & 1 \end{pmatrix}, \right.$$

$$\left. \begin{pmatrix} 1 & 0 \\ 0 & 1 \end{pmatrix}, \begin{pmatrix} 0 & 1 \\ 1 & 0 \end{pmatrix} \right\} \cup$$

$$\{(0\ 0\ 1), (1\ 1\ 0\ 0), (1\ 1\ 1), (0\ 0\ 1\ 1), (1\ 0\ 1)\ (1\ 0\ 1\ 0)\}$$

is the bigenerating biset of V over S. $|X| = (|X_1|, |X_2|) = (8, 6)$; thus bidimension of V is finite.

Now we have special substructures which we define in the following:

**DEFINITION 5.4:** *Let $V = V_1 \cup V_2$ be a neutrosophic set bivector space over the set S. Suppose $W = W_1 \cup W_2 \subseteq V_1 \cup V_2 = V$ is such that W is only a set bivector space over the set S, then we call W to be a pseudo neutrosophic set bivector subspace of V over S.*

**DEFINITION 5.5:** *Let $V = V_1 \cup V_2$ be a neutrosophic set bivector space over the set S. Suppose $W = W_1 \cup W_2 \subseteq V = V_1 \cup V_2$ is such that $W_1$ is a neutrosophic set vector subspace of $V_1$ and $W_2$ is just a set vector subspace of $V_2$ then we call W to a quasi neutrosophic set bivector subspace of V over S.*

We will illustrate the two definitions by some examples.

*Example 5.7:* Let

$$V = V_1 \cup V_2 = \left\{ \begin{pmatrix} a & b \\ 0 & 0 \end{pmatrix}, \begin{pmatrix} 0 & 0 \\ c & d \end{pmatrix} \middle| a,b,c,d \in N(Q) \right\}$$

$$\cup \left\{ (x_1, x_2, x_3, x_4, x_5), \begin{pmatrix} y_1 \\ y_2 \\ y_3 \end{pmatrix} \middle| \begin{array}{l} x_i, y_j \in N(Z); \\ 1 \leq i \leq 5; \\ 1 \leq j \leq 3 \end{array} \right\}$$

be a neutrosophic set bivector space over the set $S = \{0, 1\}$. Take

$$W = W_1 \cup W_2 = \left\{ \begin{pmatrix} a & b \\ 0 & 0 \end{pmatrix}, \begin{pmatrix} 0 & 0 \\ c & d \end{pmatrix} \middle| a,b,c,d \in Q \right\} \cup$$

$$\left\{ (x_1, x_2, x_3, x_4, x_5), \begin{pmatrix} y_1 \\ y_2 \\ y_3 \end{pmatrix} \middle| \begin{array}{l} x_i, y_j \in N(Z); \\ 1 \leq i \leq 5; \\ 1 \leq j \leq 3 \end{array} \right\}$$

$$\subseteq V_1 \cup V_2 = V;$$

W is a pseudo neutrosophic set bivector subspace of V over S.

*Example 5.8:* Let $V = V_1 \cup V_2 = \{N(Q)[x];$ that is $N(Q)[x]$ is the set of all polynomials in the variable x with coefficients from $N(Q)\} \cup$

$$\left\{ (x, y, z), \begin{bmatrix} x \\ y \\ z \\ w \end{bmatrix} \middle| x, y, z, w \in N(Z_7) \right\}$$

be a neutrosophic set vector bispace over the set S = {0, 1}. Take
$$W = W_1 \cup W_2 = \{Q[x]\} \cup \{(x, y, z) \mid x, y, z \in Z_7\}$$
$$\subseteq V_1 \cup V_2;$$
W is a pseudo neutrosophic set bivector subspace over the set S = {0, 1}.

At this juncture it is essential to make the following observations.

**THEOREM 5.1:** *Let $V = V_1 \cup V_2$ be a neutrosophic set bivector space over the set S (where $S = ZI$ or $mZ^+I$, $m \in N$ or $QI$ or $Q^+I$ or $RI$, $R^+I$, or $CI$ or $Z_mI$; $m \leq \infty$ and $m \in \{0, 1, ..., n \mid n \leq \infty\}$). Then V has no pseudo neutrosophic set bivector subspace.*

*Proof:* Given the set over which the neutrosophic set bivector space $V = V_1 \cup V_2$ is defined is a pure neutrosophic set then by the very definition of neutrosophic set bivector space both $V_1$ and $V_2$ cannot have a set vector subspace hence the claim.

We will illustrate this by some examples.

*Example 5.9:* Let $V = V_1 \cup V_2 = \{N(Q)[x];$ all polynomials in the variable x with coefficients from the neutrosophic set N(Q)} $\cup$

$$\left\{ \begin{pmatrix} a_1I & a_2I & a_5I \\ a_3I & a_4I & a_6I \end{pmatrix}, (aI, aI, aI, aI, aI) \mid a_i, a \in QI; 1 \leq i \leq 6 \right\}$$

be a neutrosophic set bivector space over the set $S = Z^+I$. Clearly we have in $V_1$ a subset $W_1 = \{Q[x]\}$ which is a proper subset of $V_1$ but however $Q[x]$ is not a set vector space over the set $S = Z^+I$. Thus we see V does not have a pseudo neutrosophic set bivector subspace.

*Example 5.10:* Let $V = V_1 \cup V_2 = \{M_{5 \times 6} = (m_{ij})$ is the collection of all $5 \times 6$ matrices with entries from QI} $\cup \{(a_1I, a_2I, ..., a_{11}I)$ $| a_i \in ZI; 1 \leq i \leq 11\}$ be a neutrosophic set bivector space over the set $S = Z^+I$. Clearly V has no pseudo neutrosophic set bivector subspace.

**DEFINITION 5.6:** *Let $V = V_1 \cup V_2$ where $V_1$ is a neutrosophic set vector space over the set S and $V_2$ is just a set vector space over the same set S. We call V to be a quasi neutrosophic set bivector space over S.*

*Note:* It is important to note the set S in the definition can only be a subset of reals or complex or rationals or integers and is never a neutrosophic subset. Thus for quasi neutrosophic set bivector space to be defined the set over which it is defined is only an ordinary set.

We will illustrate this situation by some examples.

*Example 5.11:* Let

$$V = V_1 \cup V_2 = \left\{ \begin{pmatrix} a & b \\ c & d \end{pmatrix} \middle| a,b,c,d \in Q \right\} \cup$$

$$\left\{ \begin{pmatrix} aI & bI & cI \\ dI & 0 & fI \end{pmatrix} \middle| aI, bI, cI, dI, fI \in ZI \right\}$$

be a quasi neutrosophic set bivector space over the set $S = Z$.

*Example 5.12:* Let

$$V = V_1 \cup V_2$$

$$= \left\{ \begin{pmatrix} a \\ b \\ c \end{pmatrix}, (x_1 \quad x_2 \quad x_3 \quad x_4 \quad x_5) \middle| \begin{array}{l} a,b,c,x_i \in Q; \\ 1 \leq i \leq 5 \end{array} \right\} \cup$$

$$\left\{ \begin{pmatrix} a_1I & a_2I & a_3I \\ 0 & 0 & 0 \end{pmatrix}, \begin{pmatrix} 0 & 0 & 0 \\ b_1I & 0 & b_2I \end{pmatrix} \middle| \begin{array}{l} a_iI, b_jI \in ZI; \\ 1 \leq i \leq 4, 1 \leq j \leq 2 \end{array} \right\}$$

be a quasi neutrosophic set bivector space over the set $S = 3Z^+ \cup \{0\}$.

It is important to note that quasi neutrosophic set bivector spaces can either have quasi neutrosophic set bivector subspaces or pseudo neutrosophic set bivector subspaces. Clearly quasi neutrosophic set bivector spaces do not contain neutrosophic set bivector subspaces.

**DEFINITION 5.7:** *Let $V = V_1 \cup V_2$ be a quasi neutrosophic bivector space over the set S. Suppose $W = W_1 \cup W_2 \subseteq V_1 \cup V_2$ is a proper subset of V and W is also a quasi neutrosophic bivector space over S then we define W to be a quasi neutrosophic bivector subspace of V over the set S.*

We will illustrate this by some examples.

*Example 5.13:* Let
$$V = V_1 \cup V_2$$

$$= \left\{ \begin{pmatrix} a_1 & a_2 & a_3 \\ a_4 & a_5 & a_6 \end{pmatrix}, \begin{pmatrix} b_1 & b_2 \\ b_3 & b_4 \\ b_5 & b_6 \\ b_7 & b_8 \end{pmatrix} \middle| \begin{array}{l} a_i, b_j \in Q; \\ 1 \leq i \leq 6 \\ 1 \leq j \leq 8 \end{array} \right\}$$

$$\cup \{(5ZI \times 5ZI \times 3ZI), (7ZI \times 13ZI \times 11ZI \times 17ZI)\}$$

be a quasi neutrosophic vector space over the set $S = Z$. Let

$$W = W_1 \cup W_2$$
$$= \left\{ \begin{pmatrix} a_1 & a_2 & a_3 \\ a_4 & a_5 & a_6 \end{pmatrix} \right\} \cup \{(5ZI \times 5ZI \times 3ZI)\}$$
$$\subseteq V_1 \cup V_2,$$

W is a quasi neutrosophic bivector subspace of V over the set S.

*Example 5.14:* Let

$$V = V_1 \cup V_2 = \left\{ (Q \times Q), \begin{pmatrix} a_1 \\ a_2 \\ a_3 \\ a_4 \\ a_5 \end{pmatrix} \middle| a_i \in Z; 1 \le i \le 5 \right\} \cup$$

$$\left\{ \begin{pmatrix} a_1 I & a_2 I \\ a_3 I & a_4 I \end{pmatrix}, \begin{pmatrix} a_1 I & a_2 I & a_3 I \\ 0 & a_4 I & 0 \end{pmatrix} \middle| \begin{array}{l} a_i I \in QI; \\ 1 \le i \le 4 \end{array} \right\}$$

be a quasi neutrosophic bivector space over the set $S = 5Z$. Take

$$W = W_1 \cup W_2$$

$$= \left\{ \begin{pmatrix} a_1 \\ a_2 \\ a_3 \\ a_4 \\ a_5 \end{pmatrix} \middle| a_i \in Z; 1 \le i \le 5 \right\} \cup \left\{ \begin{pmatrix} a_1 I & a_2 I \\ a_3 I & a_4 I \end{pmatrix} \middle| a_i I \in QI; 1 \le i \le 4 \right\}$$

$$\subseteq V_1 \cup V_2,$$

W is a quasi neutrosophic bivector subspace of V over the set S.

**DEFINITION 5.8:** *Let $V = V_1 \cup V_2$ be a quasi neutrosophic set bivector space over the set S. Choose $W = W_1 \cup W_2 \subseteq V_1 \cup V_2$ such that W is a set bivector space over the set S then we call W to be a pseudo quasi neutrosophic set bivector subspace of V over the set S.*

We will illustrate this situation by some examples.

*Example 5.15:* Let

$$V = V_1 \cup V_2$$

$$= \left\{ \begin{pmatrix} a & b \\ c & 0 \end{pmatrix}, \begin{pmatrix} 0 & a \\ b & c \end{pmatrix}, \begin{pmatrix} a & 0 \\ b & c \end{pmatrix}, \begin{pmatrix} a & b \\ 0 & c \end{pmatrix} \middle| a,b,c \in Q \right\} \cup$$

$$\left\{ \begin{pmatrix} aI & bI \\ 0 & 0 \end{pmatrix}, \begin{pmatrix} aI & 0 \\ bI & c \end{pmatrix}, \begin{pmatrix} aI & 0 \\ bI & c \end{pmatrix}, \begin{pmatrix} aI & 0 \\ 0 & bI \end{pmatrix}, \begin{bmatrix} x \\ y \end{bmatrix}, \begin{pmatrix} 0 & aI \\ bI & 0 \end{pmatrix} \middle| \begin{array}{l} aI, bI \in ZI \\ x, y \in 5Z \end{array} \right\}$$

be a quasi neutrosophic set bivector space over the set $S = 5Z^+ \cup \{0\}$. Take

$$W = W_1 \cup W_2$$

$$= \left\{ \begin{pmatrix} a & b \\ c & 0 \end{pmatrix}, \begin{pmatrix} 0 & a \\ b & c \end{pmatrix} \middle| a,b,c \in Z \right\} \cup \left\{ \begin{bmatrix} x \\ y \end{bmatrix} \middle| x,y \in 5Z \right\}$$

$$\subseteq V_1 \cup V_2;$$

W is pseudo quasi neutrosophic set bivector space over the set S.

*Example 5.16:* Let

$$V = V_1 \cup V_2$$

$$= \left\{ \begin{pmatrix} a & b & c \\ d & 0 & 0 \end{pmatrix}, \begin{bmatrix} a & b \\ c & d \\ 0 & 0 \end{bmatrix}, (a,b,c,d), \begin{bmatrix} a \\ b \\ c \\ d \end{bmatrix} \middle| a,b,c,d \in Z \right\}$$

$$\cup \left\{ \begin{pmatrix} a & b \\ c & d \end{pmatrix} \middle| a,b,c,d \in N(Q) \right\}$$

be a quasi neutrosophic set bivector space over the set $S = Z$. Take

$$W = W_1 \cup W_2 =$$

$$\left\{ \begin{pmatrix} a & b & c \\ d & 0 & 0 \end{pmatrix}, (a,b,c,d) \,\middle|\, a,b,c,d \in 5N \right\} \cup$$

$$\left\{ \begin{pmatrix} a & b \\ c & d \end{pmatrix} \,\middle|\, a,b,c,d \in Q \right\}$$

$$\subseteq V_1 \cup V_2$$

is a pseudo quasi neutrosophic set bivector subspace of V over the set S.

**DEFINITION 5.9:** *Let $V = V_1 \cup V_2$ be a quasi neutrosophic set bivector space over the set S. Suppose $W = W_1 \cup W_2 \subseteq V_1 \cup V_2$ and W is a quasi neutrosophic set bivector space over the set T, T a proper subset of S, then we call W to be a quasi neutrosophic subset bivector subspace of V over the subset T of S.*

We will illustrate this situation by some examples.

*Example 5.17:* Let

$$V = V_1 \cup V_2 = \left\{ \begin{pmatrix} a & b \\ c & d \end{pmatrix} \,\middle|\, a,b,c,d \in N(Z_{20}) \right\} \cup$$

$$\left\{ Z_{20} \times Z_{20} \times Z_{20}, \begin{bmatrix} a \\ b \\ c \end{bmatrix} \,\middle|\, a,b,c \in Z_{20} \right\}$$

be a quasi neutrosophic set bivector space over the set $S = Z_{20}$. Take

$$W = W_1 \cup W_2$$

$$= \left\{ \begin{pmatrix} aI & bI \\ cI & dI \end{pmatrix} \,\middle|\, aI, bI, cI, dI \in Z_{20}I \right\} \cup \left\{ \begin{pmatrix} a \\ b \\ c \end{pmatrix} \,\middle|\, a,b,c \in Z_{20} \right\}$$

$$\subseteq V_1 \cup V_2;$$

W is a quasi neutrosophic subset bivector subspace of V over the set T = {0, 5, 10, 15} ⊆ S.

*Example 5.18:* Let
$$V = V_1 \cup V_2$$
$$= \left\{ \begin{pmatrix} a & b & c \\ d & e & f \end{pmatrix}, (a,b,c,d,e,f) \,\middle|\, a,b,c,d,e,f \in N(Q) \right\} \cup$$

$$\left\{ (Z \times Z \times Z \times Z), \begin{bmatrix} a \\ b \\ c \\ d \\ e \end{bmatrix} \,\middle|\, a, b, c, d, e \in Z \right\}$$

be a quasi neutrosophic set bivector space over the set $2Z = S$. Take

$$W = W_1 \cup W_2 =$$
$$\left\{ \begin{pmatrix} a & b & c \\ d & e & f \end{pmatrix}, (a,b,c,d,e,f) \,\middle|\, a,b,c,d,e,f \in N(Q) \right\}$$
$$\cup \{(Z \times Z \times Z \times Z)\} \subseteq V_1 \cup V_2,$$

W is a quasi neutrosophic subset bivector subspace of V over the set $T = 10Z \subseteq S$.

Now we proceed onto define the notion of neutrosophic set bilinear algebras.

**DEFINITION 5.10:** *Let $V = V_1 \cup V_2$ be such that $V_1$ is a neutrosophic set linear algebra over the set S and $V_2$ is also a neutrosophic set linear algebra over the same set S, then we call V to be a neutrosophic set bilinear algebra over the set S only if $V_1 \not\subseteq V_2$ and $V_2 \not\subseteq V_1$ or $V_1 \neq V_2$.*

We will illustrate this situation by some examples.

*Example 5.19:* Let

$$V = V_1 \cup V_2 = \left\{ \begin{pmatrix} a & b \\ c & d \end{pmatrix} \middle| a,b,c,d \in N(Z_{15}) \right\} \cup \{Z_{15}I \times Z_{15}I \times Z_{15}I \times Z_{15}I\}$$

be a neutrosophic set bilinear algebra over the set $S = Z_{15}$.

*Example 5.20:* Let

$$V = V_1 \cup V_2$$

$$= \{(a_1I, a_2I, a_3I, a_4I, a_5I, a_6I) \mid a_iI \in QI; 1 \le i \le 6)\}$$
$$\cup \left\{ \begin{pmatrix} a & b \\ 0 & c \end{pmatrix} \middle| a,b,c \in QI \right\}$$

be a neutrosophic set bilinear algebra over the set $S = QI$.

Now we proceed onto define substructure of neutrosophic set bilinear algebras.

**DEFINITION 5.11:** *Let $V = V_1 \cup V_2$ be a neutrosophic set bilinear algebra over the set S. Suppose $W = W_1 \cup W_2 \subseteq V_1 \cup V_2$ is a proper bisubset of V and W itself is a neutrosophic set bilinear algebra over the set S then we call W to be a neutrosophic set bilinear subalgebra of V over the set S.*

We will illustrate this by some examples.

*Example 5.21:* Let

$$V = V_1 \cup V_2$$
$$= \left\{ \begin{pmatrix} a & b \\ c & d \end{pmatrix} \middle| a,b,c,d \in N(Z) \right\} \cup \{ZI \times ZI \times ZI \times ZI\}$$
$$(= \{(aI, bI, cI, dI) \mid aI, bI, cI, dI \in ZI\})$$

be a neutrosophic set bilinear algebra over the set $S = 5Z$. Take

$$W = W_1 \cup W_2$$
$$= \left\{ \begin{pmatrix} aI & bI \\ cI & dI \end{pmatrix} \middle| aI, bI, cI, dI \in ZI \right\} \cup \{ZI \times ZI \times \{0\} \times \{0\}\}$$
$$\subseteq V_1 \cup V_2;$$

W is a neutrosophic set bilinear subalgebra of V over the set S = 5Z.

*Example 5.22:* Let
$$V = V_1 \cup V_2$$
$$= \left\{ \begin{pmatrix} a & b & c \\ d & e & f \end{pmatrix} \middle| a,b,c,d,e,f \in N(Z_{22}) \right\}$$
$$\cup \left\{ \begin{pmatrix} a & b \\ c & d \end{pmatrix} \middle| a,b,c,d \in N(Z_{22}) \right\}$$

be a neutrosophic set bilinear algebra over the set $S = Z_{22}$.
Take
$$W = \left\{ \begin{pmatrix} a & b & c \\ 0 & d & 0 \end{pmatrix} \middle| a,b,c,d \in N(Z_{22}) \right\} \cup$$
$$\left\{ \begin{pmatrix} a & b \\ 0 & 0 \end{pmatrix} \middle| a,b \in N(Z_{22}) \right\}$$
$$= W_1 \cup W_2 \subseteq V_1 \cup V_2;$$

W is a neutrosophic set bilinear subalgebra of V over the set S.

We have the following result the proof of which is left as an exercise to the reader.

**THEOREM 5.2:** *Every neutrosophic set bilinear algebra V over a set S is always a neutrosophic set bivector space over the set S but however in general every neutrosophic set bivector space over a set S is not a neutrosophic set bilinear algebra over S.*

**DEFINITION 5.12:** *Let $V = V_1 \cup V_2$ be a neutrosophic set bilinear algebra over the set S. Suppose $X = X_1 \cup X_2 \subseteq V_1 \cup V_2 = V$ is such that $X_1$ is a generating set of $V_1$ over S and $X_2$ is a generating set of $V_2$ over S then we call $X = X_1 \cup X_2$ to be the generating bisubset of V over S. The bidimension of V is the cardinality of $X = X_1 \cup X_2$ denoted by $|X| = (|X_1|, |X_2|)$.*

We will illustrate this situation by some examples.

*Example 5.23:* Let

$$V = V_1 \cup V_2 = \left\{ \begin{pmatrix} a & a \\ a & a \end{pmatrix} \middle| a \in NI \right\} \cup \{(a, a, a, a, a, a) \mid a \in ZI\}$$

be a neutrosophic set bilinear algebra over the set $S = ZI$. Clearly

$$X = \left\{ \begin{pmatrix} I & I \\ I & I \end{pmatrix} \right\} \cup \{(I, I, I, I, I, I)\} = X_1 \cup X_2$$

bigenerates V over ZI. The bidimension of V is (1, 1).

*Example 5.24:* Let

$$V = V_1 \cup V_2 = ZI[x] \cup \left\{ \begin{pmatrix} x_1 & x_2 & x_3 \\ x_4 & x_5 & x_6 \end{pmatrix} \middle| x_i \in ZI \right\}$$

be a neutrosophic set bilinear algebra over the set $S = ZI$.

$$X = \{I, Ix, Ix^2, \ldots\} \cup \begin{pmatrix} I & 0 & 0 \\ 0 & 0 & 0 \end{pmatrix}, \begin{pmatrix} 0 & I & 0 \\ 0 & 0 & 0 \end{pmatrix}, \begin{pmatrix} 0 & 0 & I \\ 0 & 0 & 0 \end{pmatrix},$$

$$\begin{pmatrix} 0 & 0 & 0 \\ I & 0 & 0 \end{pmatrix}, \begin{pmatrix} 0 & 0 & 0 \\ 0 & I & 0 \end{pmatrix}, \begin{pmatrix} 0 & 0 & 0 \\ 0 & 0 & I \end{pmatrix} \right\} \subseteq V_1 \cup V_2$$

is a bigenerator of V over the set S. Clearly bidimension of V is infinite as $|X| = (|X_1|, |X_2|) = (\infty, 6)$.

**DEFINITION 5.13:** *Let $V = V_1 \cup V_2$ be a neutrosophic set bilinear algebra over the set S. Suppose $W = W_1 \cup W_2 \subseteq V_1 \cup V_2$ is such that $W = W_1 \cup W_2$ is only a neutrosophic set bivector space over the set S then we call W to be a pseudo neutrosophic set bivector subspace of V over the set S.*

We will illustrate this by some examples.

*Example 5.25:* Let

$$V = \left\{ \begin{pmatrix} a & b \\ c & d \end{pmatrix} \middle| a,b,c,d \in ZI \right\} \cup \{ZI \times ZI \times ZI\} = V_1 \cup V_2$$

be a neutrosophic set bilinear algebra over the set $S = ZI$.
Let

$$W = \left\{ \begin{pmatrix} a & 0 \\ 0 & 0 \end{pmatrix}, \begin{pmatrix} 0 & b \\ 0 & 0 \end{pmatrix} \middle| a,b \in ZI \right\} \cup$$

$$\{5ZI \times 3ZI \times \{0\}, \{0\} \times \{0\} \times ZI\} \subseteq V_1 \cup V_2$$

Then, W is only a neutrosophic set bivector space over the set S = ZI. Thus W is a pseudo neutrosophic set bilinear subalgebra over S.

*Example 5.26:* Let

$$V = \left\{ \begin{pmatrix} a_1 & a_2 & a_3 \\ a_4 & a_5 & a_6 \end{pmatrix} \middle| a_i \in QI; 1 \leq i \leq 6 \right\} \cup$$

{QI[x]; all polynomials in the variable x with coefficients from QI} = $V_1 \cup V_2$ be a neutrosophic set bilinear algebra over the set S = ZI. Choose

$$W = W_1 \cup W_2 =$$

$$\left\{ \begin{pmatrix} a_1 & 0 & 0 \\ 0 & 0 & a_6 \end{pmatrix} \begin{pmatrix} 0 & a_2 & 0 \\ a_4 & 0 & a_6 \end{pmatrix} \middle| a_1, a_6, a_2, a_4 \in QI \right\} \cup$$

{5ZI[x], 3ZI[x] where 5ZI [x] denotes the set of all polynomials in the variable x with coefficients from 5ZI and 3ZI [x] denotes the set of all polynomials in the variable x with coefficients for 3ZI} $\subseteq V_1 \cup V_2$. W is a pseudo neutrosophic set bivector space over the set S = ZI.

**DEFINITION 5.14:** *Let $V = V_1 \cup V_2$ be a neutrosophic set bilinear algebra over the set S. Let $W = W_1 \cup W_2$ be such that $W_1$ is a neutrosophic set linear subalgebra of $V_1$ over S and $W_2$ is only a neutrosophic set vector subspace of $V_2$ then we call $W = W_1 \cup W_2$ to be a quasi neutrosophic set bilinear subalgebra of V over the set S.*

We will illustrate this by some simple examples.

*Example 5.27:* Let

$$V = V_1 \cup V_2 = \left\{ \begin{pmatrix} a & b \\ c & d \end{pmatrix} \middle| a,b,c,d \in QI \right\} \cup \left\{ \begin{pmatrix} a_1 & a_2 & a_3 \\ a_4 & a_5 & a_6 \end{pmatrix} \right\}$$

such that $a_i \in QI$; $1 \leq i \leq 6$} be neutrosophic set bilinear algebra over the set S = ZI.
Consider

$$W = W_1 \cup W_2$$

$$= \left\{ \begin{pmatrix} a & 0 \\ 0 & 0 \end{pmatrix}, \begin{pmatrix} 0 & b \\ 0 & 0 \end{pmatrix} \middle| a,b \in QI \right\} \cup$$

$$\left\{ \begin{pmatrix} a_1 & a_2 & a_3 \\ a_4 & a_5 & a_6 \end{pmatrix} \middle| a_i \in ZI; 1 \leq i \leq 6 \right\}$$

$$\subseteq V_1 \cup V_2;$$

W is only a quasi neutrosophic set bilinear subalgebra of V over the set S = ZI.

*Example 5.28:* Let

$$V = \left\{ \begin{pmatrix} a_1 & a_2 \\ a_3 & a_4 \\ a_5 & a_6 \\ a_7 & a_8 \end{pmatrix} \middle| a_i \in QI; 1 \le i \le 8 \right\} \cup$$

$$\left\{ \begin{pmatrix} a_1 & a_2 & a_3 & a_4 & a_5 \\ a_6 & a_7 & a_8 & a_9 & a_{10} \end{pmatrix} \middle| a_i \in ZI; 1 \le i \le 10 \right\}$$

$$= V_1 \cup V_2$$

be a neutrosophic set bilinear algebra over the set S = ZI. Take

$$W = W_1 \cup W_2$$

$$= \left\{ \begin{pmatrix} a_1 & a_2 \\ 0 & 0 \\ 0 & 0 \\ a_7 & a_8 \end{pmatrix}, \begin{pmatrix} 0 & 0 \\ a_3 & a_4 \\ a_5 & a_6 \\ 0 & 0 \end{pmatrix} \middle| a_i \in QI; 1 \le i \le 8 \right\} \cup$$

$$\left\{ \begin{pmatrix} a_1 & a_2 & a_3 & a_4 & a_5 \\ a_6 & a_7 & a_8 & a_9 & a_{10} \end{pmatrix} \middle| a_i \in 5ZI; 1 \le i \le 10 \right\}$$

$$\subseteq V_1 \cup V_2;$$

W is a quasi neutrosophic set bilinear subalgebra of V over the set S.

**DEFINITION 5.15:** *Let $V = V_1 \cup V_2$ be such that $V_1$ is a neutrosophic set linear algebra over the set S and $V_2$ is a only a set linear algebra over a set S then we call V to be a quasi neutrosophic set bilinear algebra over the set S.*

We shall illustrate this situation by some examples.

*Example 5.29:* Let

$$V = V_1 \cup V_2 = \left\{ \begin{pmatrix} a & b \\ c & d \end{pmatrix} \middle| a, b, c, d \in QI \right\} \cup$$

$$\left\{ \begin{pmatrix} a_1 & a_2 & a_3 \\ a_4 & a_5 & a_6 \end{pmatrix} \middle| a_i \in ZI; 1 \le i \le 6 \right\}$$

where $V_1$ is a neutrosophic set linear algebra over the set $S = Z^+ \cup \{0\}$ and $V_2$ is a set linear algebra over the same set $S = Z^+ \cup \{0\}$. We see $V = V_1 \cup V_2$ is a quasi neutrosophic set bilinear algebra over the set S.

*Example 5.30:* Let $V = V_1 \cup V_2$ where $V_1 = \{Z[x];$ all polynomials in the variable x with coefficients from the ring $Z\}$ is a set linear algebra over the set $S = Z$ and $V_2 = \{QI \times QI \times QI \times QI \times Q \times Q \times QI\}$ be a neutrosophic set linear algebra over the set $S = Z$. $V = V_1 \cup V_2$ is a quasi neutrosophic set bilinear algebra over the set S.

**DEFINITION 5.16:** *Let $V = V_1 \cup V_2$ where $V_1$ is a neutrosophic set linear algebra over the set S and $V_2$ is a neutrosophic set vector space over the same set S then we call V to be a pseudo neutrosophic set bilinear algebra over the set S.*

We will illustrate this by some examples.

*Example 5.31:* Let
$$V = V_1 \cup V_2$$

$$= \left\{ \begin{pmatrix} a_1 & a_2 & a_3 \\ 0 & a_4 & 0 \end{pmatrix}, \begin{pmatrix} a_1 & 0 & 0 \\ a_5 & 0 & a_6 \end{pmatrix} \middle| a_i \in ZI; 1 \leq i \leq 6 \right\} \cup$$

$$\left\{ \begin{pmatrix} a & b \\ c & d \end{pmatrix} \middle| a,b,c,d \in QI \right\}$$

be a pseudo neutrosophic set bilinear algebra over the set $S = ZI$.

*Example 5.32:* Let

$$V = V_1 \cup V_2$$

$$= \left\{ \begin{pmatrix} a & a & a & a & a \\ a & a & a & a & a \end{pmatrix}, \begin{pmatrix} a & a & a \\ a & a & a \end{pmatrix} \middle| a \in QI \right\}$$

$$\cup \left\{ \begin{pmatrix} a & a \\ a & a \end{pmatrix} \middle| a \in ZI \right\}$$

be a pseudo neutrosophic set bilinear algebra over the set $S = Z^+I$; here $V_1$ is only a neutrosophic set vector space over $S = Z^+I$ and $V_2$ is a neutrosophic set linear algebra over the set $S = Z^+I$.

We define neutrosophic set bilinear transformation of these in the following.

**DEFINITION 5.17:** *Let $V = V_1 \cup V_2$ and $W = W_1 \cup W_2$ be any two neutrosophic set bivector spaces defined over the same set S. A map $T = T_1 \cup T_2$ from $V = V_1 \cup V_2$ to $W_1 \cup W_2$ is a neutrosophic set bilinear transformation of the neutrosophic set bivector spaces if $T_1 : V_1 \to W_1$ and $T_2 : V_2 \to W_2$ are neutrosophic set linear transformations of the neutrosophic set vector spaces.*

Here $T = T_1 \cup T_2 : V = V_1 \cup V_2 \to W = W_1 \cup W_2$ where $T_1 : V_1 \to W_1 \cup T_2 : V_2 \to W_2 = T:V \to W = T_1 \cup T_2 : V_1 \cup V_2 \to$

$W_1 \cup W_2$. On similar lines one can define neutrosophic set bilinear transformations of neutrosophic set bilinear algebras.

If $W = W_1 \cup W_2$ is replaced by $V = V_1 \cup V_2$ we call such neutrosophic set bilinear transformations as neutrosophic set bilinear operators on V.

We denote the collection of all neutrosophic set bilinear transformations of V to W by $NHom_S$ ($V = V_1 \cup V_2$, $W = W_1 \cup W_2$) and that of the neutrosophic set bilinear operators of V to V by $NHom_S$ ($V = V_1 \cup V_2$), $V = V_1 \cup V_2$).

We will illustrate this situation by some examples.

*Example 5.33:* Let
$$V = V_1 \cup V_2$$

$$= \left\{ \begin{pmatrix} a & 0 \\ b & 0 \end{pmatrix}, \begin{pmatrix} 0 & b \\ 0 & d \end{pmatrix} \middle| a,b,c,d \in QI \right\} \cup$$

$$\left\{ \begin{pmatrix} a_1 & a_2 & a_3 \\ 0 & 0 & a_4 \end{pmatrix}, \begin{pmatrix} 0 & 0 & 0 & 0 \\ a_5 & a_6 & 0 & 0 \end{pmatrix} \middle| a_i \in QI; 1 \le i \le 6 \right\}$$

and

$$W = W_1 \cup W_2$$

$$= \{QI \times QI \times QI \times QI\}$$

$$\cup \left\{ \begin{pmatrix} a_1 & a_2 \\ a_3 & a_4 \\ 0 & 0 \end{pmatrix}, \begin{pmatrix} 0 & 0 \\ a_5 & a_6 \end{pmatrix} \middle| a_i \in QI; 1 \le i \le 6 \right\}$$

be neutrosophic set bivector spaces defined over the set $S = ZI$.
Define $\eta = \eta_1 \cup \eta_2 = V = V_1 \cup V_2 \to W = W_1 \cup W_2$ where $\eta_1 : V_1 \to W_1$ and $\eta_2 : V_2 \to W_2$ as

$$\eta_1 \begin{pmatrix} a & 0 \\ b & 0 \end{pmatrix} = \{(a, b, 0, 0)\},$$

$$\eta_1 \begin{pmatrix} 0 & b \\ 0 & d \end{pmatrix} = \{(0, 0, b, d)\},$$

$$\eta_2 \begin{pmatrix} a_1 & a_2 & a_3 \\ 0 & 0 & a_4 \end{pmatrix} = \begin{pmatrix} a_1 & a_2 \\ a_3 & a_4 \\ 0 & 0 \end{pmatrix},$$

and

$$\eta_2 \begin{pmatrix} 0 & 0 & 0 & 0 \\ a_5 & a_6 & 0 & 0 \end{pmatrix} = \begin{pmatrix} 0 & 0 \\ a_5 & a_6 \end{pmatrix}.$$

It is easily verified that $\eta : V \to W$ is a neutrosophic set bilinear transformation of V to W.

*Example 5.34:* Let

$$V = \left\{ \begin{pmatrix} a & b \\ 0 & c \end{pmatrix}, \begin{pmatrix} 0 & a \\ b & c \end{pmatrix}, \begin{pmatrix} a & 0 \\ b & 0 \end{pmatrix}, \begin{pmatrix} 0 & a \\ 0 & b \end{pmatrix} \middle| a, b, c \in QI \right\}$$

$\cup$ {(a b c 0), (0 0 a), (0 0 a b 0 c) | a, b, c $\in$ ZI} be a neutrosophic set bivector space over the set S = ZI.
Define $\eta : V \to V$ as
$$\eta = \eta_1 \cup \eta_2 : V_1 \cup V_2 = V \to V_1 \cup V_2 = V$$
$\eta_1 : V_1 \to V_1$ is such that

$$\eta_1 \begin{pmatrix} a & b \\ 0 & c \end{pmatrix} = \begin{pmatrix} 0 & a \\ b & c \end{pmatrix},$$

$$\eta_1 \begin{pmatrix} 0 & a \\ b & c \end{pmatrix} = \begin{pmatrix} a & b \\ 0 & c \end{pmatrix},$$

$$\eta_1 \begin{pmatrix} a & 0 \\ b & 0 \end{pmatrix} = \begin{pmatrix} 0 & a \\ 0 & b \end{pmatrix}$$

and

$$\eta_1 \begin{pmatrix} 0 & a \\ 0 & b \end{pmatrix} = \begin{pmatrix} a & 0 \\ b & 0 \end{pmatrix}.$$

$\eta_2 : V_2 \to V_2$ is such that
$$\eta_2 (a, b, c, 0) = (0, 0, a, b, 0, c)$$
$$\eta_2 (0\ 0\ a) = (0\ 0\ a)$$
and
$$\eta_2 (0, 0, a, b, 0, c) = (a, b, c, 0).$$

It is easily verified that $\eta : V \to V$ $(\eta_1 : V_1 \to V_1) \cup (\eta_2 : V_2 \to V_2)$ is a neutrosophic set bilinear operator on V.

Interested reader can study the structure of $NHom_S (V, W)$ and $NHom_S (V, V)$ where V and W are neutrosophic set bivector space over the set S.

Next we proceed onto give examples of neutrosophic set bilinear algebra transformations.

***Example 5.35:*** Let

$$V = V_1 \cup V_2 = \left\{ \begin{pmatrix} a & b \\ c & d \end{pmatrix} \middle| a, b, c, d \in QI \right\} \cup \left\{ \sum_{i=1}^{n} a_i x^i = p(x) \right\}$$

be the collection of all polynomials in the variable x with coefficients from QI} and

$$W = W_1 \cup W_2 = \{(a, b, c, d) \mid a, b, c, d \in QI\}$$
$$\cup \left\{ \sum_{i=1}^{n} a_i x^{2i} = q(x) \right\}$$

be the collection of polynomials in the variable x of even degree with coefficients from QI} be two neutrosophic set bilinear algebra over the set S = ZI.

Define $\eta = \eta_1 \cup \eta_2 : V_1 \cup V_2 = V \to W_1 \cup W_2 = W$ as $\eta_1 : V_1 \to W_1$ and $\eta_2 : V_2 \to W_2$ such that

$$\eta_1 \begin{pmatrix} a & b \\ c & d \end{pmatrix} = \{a, b, c, d\}$$

and

$$\eta_2 \left( \sum_{i=1}^{n} a_i x^i \right) = \sum_{i=1}^{n} a_i x^{2i}.$$

It is easily verified that η is a neutrosophic set bilinear transformation from V to W.

*Example 5.36:* Let

$$V = V_1 \cup V_2 =$$

$$\left\{ \begin{pmatrix} a_1 & a_2 & a_3 \\ a_4 & a_5 & a_6 \end{pmatrix} \middle| a_i \in QI; 1 \le i \le 6 \right\} \cup$$

$$\left\{ \begin{pmatrix} a_1 & a_2 & a_3 \\ a_4 & a_5 & a_6 \\ a_7 & a_8 & a_9 \end{pmatrix} \middle| a_i \in QI; 1 \le i \le 9 \right\}$$

be a neutrosophic set bilinear algebra over the set $S = ZI$. Define $\eta = \eta_1 \cup \eta_2 : V = V_1 \cup V_2 \to V = V_1 \cup V_2$ by

$$\eta_1 : V_1 \to V_1 \text{ and } \eta_2 : V_2 \to V_2$$

$$\eta_1 \begin{pmatrix} a_1 & a_2 & a_3 \\ a_4 & a_5 & a_6 \end{pmatrix} = \begin{pmatrix} a_4 & a_5 & a_6 \\ a_1 & a_2 & a_3 \end{pmatrix}$$

and

$$\eta_2 \begin{pmatrix} a_1 & a_2 & a_3 \\ a_4 & a_5 & a_6 \\ a_7 & a_8 & a_9 \end{pmatrix} = \begin{pmatrix} a_7 & a_8 & a_9 \\ a_4 & a_5 & a_6 \\ a_1 & a_2 & a_3 \end{pmatrix}.$$

It is easily verified that η is a neutrosophic set bilinear operator on V.

We proceed onto define the notion of neutrosophic biset bivector spaces.

**DEFINITION 5.18:** *Let $V = V_1 \cup V_2$ be such that $V_1$ is a neutrosophic set vector space over the set $S_1$ and $V_2$ is a neutrosophic set vector space over the set $S_2$; $S_1 \neq S_2$; $S_1 \not\subset S_2$. We define $V = V_1 \cup V_2$ to be the neutrosophic biset bivector space over the biset $S = S_1 \cup S_2$.*

We will illustrate this by some examples.

***Example 5.37:*** Let $V = V_1 \cup V_2$ where .

$$V_1 = \left\{ (a, b, c), \begin{bmatrix} a \\ b \\ c \\ d \end{bmatrix} \,\middle|\, a, b, c, d \in QI \right\}$$

be a neutrosophic set vector space over the set $S_1 = 5ZI$ and

$$V_2 = \left\{ \begin{pmatrix} a & 0 \\ b & 0 \end{pmatrix}, \begin{pmatrix} a_1 & a_2 & a_3 \\ a_4 & a_5 & a_6 \end{pmatrix} \,\middle|\, a, b, a_i \in ZI; 1 \leq i \leq 6 \right\}$$

be a neutrosophic set vector space over the set $S_2 = 7Z^+I \cup \{0\}$. $V = V_1 \cup V_2$ is a neutrosophic biset bivector space over the biset $S = S_1 \cup S_2$.

***Example 5.38:*** Let $V = V_1 \cup V_2$ where $V_1 = \{ZI[x]$; i.e., all polynomials in the variable x with coefficients from $ZI\}$ is a neutrosophic set vector space over the set $S_1 = 3Z^+I \cup \{0\}$ and

$$V_2 = \left\{ (a, b, c), \begin{pmatrix} a & b \\ c & d \end{pmatrix} \,\middle|\, a, b, c, d \in ZI \right\}$$

be a neutrosophic set vector space over the set $S_2 = 5Z^+ I \cup \{0\}$. $V = V_1 \cup V_2$ is a neutrosophic biset bivector space over the set $S = S_1 \cup S_2$.

These types of algebraic structures will find their applications in mathematical models.

*Example 5.39:* Let $V = V_1 \cup V_2$ where $V_1 = \{(0\ 0\ 0), (I\ 0\ 0), (0\ 1\ 1), (0\ 0\ 0\ 0\ 0), (0\ 0\ 0\ I), (0\ 0\ 0\ 0), (I\ 1\ I\ 1), (1\ 0\ I\ 0)\ (1, I, 0, 0, I), (0\ 0\ I\ 1\ 1)\}$ is a neutrosophic set vector space over the set $S_1 = Z_2 = \{0, 1\}$ and

$$V_2 = \left\{ \begin{pmatrix} a & b & c \\ d & 0 & e \\ 0 & 0 & f \end{pmatrix}, \begin{pmatrix} a & 0 & 0 \\ b & 0 & 0 \\ e & 0 & g \end{pmatrix} \middle| a,b,c,d,e,f,g \in Z_5 I \right\}$$

is a neutrosophic set vector space over the set $S_2 = Z_5\ I$. Thus $V = V_1 \cup V_2$ is a neutrosophic biset bivector space over the biset $S = S_1 \cup S_2$.

**DEFINITION 5.19:** *Let $V = V_1 \cup V_2$ be a neutrosophic biset bivector space over the biset $S = S_1 \cup S_2$. Let $W = W_1 \cup W_2 \subseteq V_1 \cup V_2$; if W is a neutrosophic biset bivector space over the set $S = S_1 \cup S_2$ then we call W to be a neutrosophic biset bivector subspace of V over the biset $S = S_1 \cup S_2$.*

*Example 5.40:* Let $V = V_1 \cup V_2$ where

$$V_1 = \left\{ \begin{pmatrix} a & a \\ a & a \end{pmatrix} \middle| a \in Z_{12}I \right\}$$

is a neutrosophic set vector space over the set $S_1 = Z_{12}I$ and

$$V_2 = \left\{ \begin{pmatrix} a_1 & a_2 & a_3 \\ 0 & 0 & 0 \end{pmatrix}, \begin{pmatrix} 0 & 0 & 0 \\ a_4 & a_5 & a_6 \end{pmatrix} \middle| a_i \in Z_{15}I \right\}$$

is a neutrosophic set vector space over the set $S_2 = Z_{15}I$. $V = V_1 \cup V_2$ is a neutrosophic biset bivector space over $S = S_1 \cup S_2$. Let $W = W_1 \cup W_2$ where

$$W_1 = \left\{ \begin{pmatrix} a & a \\ a & a \end{pmatrix} \middle| a \in \{0, 3I, 6I, 9I\} \right\} \subseteq V_1$$

be a neutrosophic set vector subspace of $V_1$ over $S_1$ and

$$W_2 = \left\{ \begin{pmatrix} a_1 & a_2 & a_3 \\ 0 & 0 & 0 \end{pmatrix} \middle| a_i \in Z_{15}I \right\} \subseteq V_2$$

be the neutrosophic set vector subspace of $V_2$ over $S_2$ $W = W_1 \cup W_2 \subseteq V_1 \cup V_2$ is a neutrosophic biset bivector subspace of V over the biset $S = S_1 \cup S_2$.

*Example 5.41:* Let $V = V_1 \cup V_2$ where

$$V_1 = \left\{ \begin{pmatrix} a_1 \\ a_2 \\ a_3 \end{pmatrix}, (a_1, a_2, a_3) \middle| a_i \in ZI; 1 \le i \le 3 \right\}$$

is a neutrosophic set vector space over the set $S_1 = ZI$ and

$$V_2 = \left\{ \begin{pmatrix} a_1 & 0 & a_2 & 0 \\ 0 & a_3 & 0 & a_4 \end{pmatrix}, \begin{pmatrix} a_1 & 0 \\ 0 & a_3 \\ a_2 & 0 \\ 0 & a_4 \end{pmatrix} \middle| a_i \in Z_{25}I; 1 \le i \le 4 \right\}$$

be a neutrosophic set vector space over the set $S_2 = Z_{25}I$. $V = V_1 \cup V_2$ is a neutrosophic biset bivector space over the biset $S = S_1 \cup S_2$.

Take $W = W_1 \cup W_2 \subseteq V_1 \cup V_2$ where

$$W_1 = \{(a_1, a_2, a_3) \mid a_i \in ZI; 1 \le i \le 3\} \subseteq V_1$$

be the neutrosophic set vector subspace of $V_1$ over $S_1$ and

$$W_2 = \left\{ \begin{pmatrix} a_1 & 0 & a_2 & 0 \\ 0 & a_3 & 0 & a_4 \end{pmatrix} \middle| a_i \in Z_{25}I; 1 \leq i \leq 4 \right\} \subseteq V_2$$

be the neutrosophic set vector subspace of $V_2$ over $S_2$. Then W = $W_1 \cup W_2$ is a neutrosophic biset bivector subspace of V over the biset $S = S_1 \cup S_2$.

**DEFINITION 5.20:** *Let $V = V_1 \cup V_2$ be a neutrosophic biset bivector space over the biset $S = S_1 \cup S_2$.*

*If $X = X_1 \cup X_2 \subseteq V_1 \cup V_2$ is such that $X_1$ generates $V_1$ over $S_1$ and $X_2$ generates $V_2$ over $S_2$ then we say $X = X_1 \cup X_2$ is the bigenerator of the neutrosophic biset bivector space $V = V_1 \cup V_2$ over the biset $S = S_1 \cup S_2$.*

*The bicardinality of $X = X_1 \cup X_2$ denoted by $|X| = |X_1 \cup X_2| = |X_1| \cup |X_2|$ or $(|X_1|, |X_2|)$, gives the bidimension of V over S.*

We will illustrate this situation by some examples.

*Example 5.42:* Let

$$V = \left\{ (a\ a\ a\ a), \begin{bmatrix} a \\ a \\ a \\ a \\ a \\ a \end{bmatrix} \middle| a \in ZI \right\} \cup$$

$$\left\{ (a, a, a, a, a, a, a), \begin{pmatrix} a & a \\ a & a \end{pmatrix} \middle| a \in Z_{12}I \right\}.$$

$V_1 \cup V_2$ be a neutrosophic biset bivector space over the biset $S = ZI \cup Z_{15}I = S_1 \cup Z_{12}$.

Take

$$X = \{(1\ 1\ 1\ 1), \begin{bmatrix} 1 \\ 1 \\ 1 \\ 1 \\ 1 \end{bmatrix}\} \cup \{(1, 1, 1, 1, 1, 1, 1), \begin{pmatrix} 1 & 1 \\ 1 & 1 \end{pmatrix}\}$$

$$= X_1 \cup X_2 \subseteq V_1 \cup V_2;$$

clearly X bigenerates V and bidimension of V is (2, 2).

***Example 5.43:*** Let $V = V_1 \cup V_2$, where $V_1 = \{ZI[x]$; all polynomials in the variable x with coefficients from ZI$\}$ is a neutrosophic set vector space over the set $S_1 = ZI$ and

$$V_2 = \left\{ \begin{pmatrix} a & b \\ c & d \end{pmatrix} \middle| a, b, c, d \in Z_3 I \right\}$$

be a neutrosophic set vector space over the set $S_2 = Z_3 I$. $V = V_1 \cup V_2$ is a neutrosophic biset bivector space over the biset $S = S_1 \cup S_2$. V is infinitely generated by any $X = X_1 \cup X_2$ over $S = S_1 \cup S_2$.

Now we proceed onto define neutrosophic biset bilinear algebra and enumerate few of their properties.

**DEFINITION 5.21:** *Let $V = V_1 \cup V_2$ where $V_1$ is a set linear algebra over the set $S_1$ and $V_2$ another neutrosophic set linear algebra over the set $S_2$ which is distinct and different from $V_1$; further $S_1 \neq S_2$ or $S_1 \not\subset S_2$ or $S_2 \not\subset S_1$; then we call V to be the biset bilinear algebra over the biset $S = S_1 \cup S_2$.*

We will illustrate this by some examples.

***Example 5.44:*** Let $V = V_1 \cup V_2$ where

$$V_1 = \left\{ \begin{pmatrix} a & b \\ c & d \end{pmatrix} \middle| a,b,c,d \in Z_5 I \right\}$$

is a neutrosophic set linear algebra over the set $S_1 = Z_5 I$ and $V_2 = \{QI \times QI \times QI \times QI\}$ is a neutrosophic set linear algebra over the set $S_2 = QI$. Now $V = V_1 \cup V_2$ is a neutrosophic biset bilinear algebra over the biset $S = S_1 \cup S_2$.
If we take

$$X = \left\{ \begin{pmatrix} I & 0 \\ 0 & 0 \end{pmatrix}, \begin{pmatrix} 0 & I \\ 0 & 0 \end{pmatrix}, \begin{pmatrix} 0 & 0 \\ I & 0 \end{pmatrix}, \begin{pmatrix} 0 & 0 \\ 0 & I \end{pmatrix} \right\}$$

$$\cup \{(I, 0, 0, 0), (0, I, 0, 0), (0, 0, I, 0), (0, 0, 0, I)\}$$
$$= X_1 \cup X_2 \subseteq V_1 \cup V_2;$$

X bigenerates V as a neutrosophic biset bilinear algebra over the biset $S = S_1 \cup S_2$. The bidimensin of V is (4, 4).

***Example 5.45:*** Let $V = V_1 \cup V_2$ where

$$V_1 = \left\{ \begin{pmatrix} a_1 & a_2 & a_3 \\ a_4 & a_5 & a_6 \end{pmatrix} \middle| a_i \in QI; 1 \le i \le 6 \right\}$$

is a neutrosophic set linear algebra over the set $S_1 = ZI$ and

$$V_2 = \left\{ \begin{pmatrix} a & b \\ c & d \end{pmatrix} \middle| a,b,c,d \in RI \right\}$$

be a neutrosophic set linear algebra over the set $S_2 = Q^+ I$. $V = V_1 \cup V_2$ is a neutrosophic biset bilinear algebra over the biset $S = S_1 \cup S_2$. Clearly both $V_1$ and $V_2$ are infinitely generated as neutrosophic set linear algebras over the sets $S_1$ and $S_2$ respectively.

We now proceed onto define the substructures in neutrosophic biset bilinear algebras.

**DEFINITION 5.22:** *Let $V = V_1 \cup V_2$ be a neutrosophic biset bilinear algebra over the biset $S = S_1 \cup S_2$. Let $W = W_1 \cup W_2 \subseteq V_1 \cup V_2$; if W is a neutrosophic biset bilinear algebra over the biset $S = S_1 \cup S_2$ then we call W to be a neutrosophic biset bilinear subalgebra of V over the biset $S = S_1 \cup S_2$.*

*Example 5.46:* Let $V = V_1 \cup V_2$ where $V_1 = \{ZI \times ZI \times ZI \times ZI\}$ is a neutrosophic set linear over the set $S_1 = 5Z^+I$ and

$$V_2 = \left\{ \begin{pmatrix} a_1 & a_2 \\ a_4 & a_3 \end{pmatrix} \middle| a_i \in QI; 1 \leq i \leq 4 \right\}$$

be a neutrosophic set linear algebra over the set $S_2 = 7Z^+I$. $V = V_1 \cup V_2$ is a neutrosophic biset bilinear algebra over the biset $S = S_1 \cup S_2$. Take $W = W_1 \cup W_2 \subseteq V_1 \cup V_2$ where

$$W_1 = \{ZI \times \{0\} \times \{0\} \times ZI\} \subseteq V_1$$

and

$$W_2 = \left\{ \begin{pmatrix} a_1 & a_2 \\ a_4 & a_3 \end{pmatrix} \middle| a_i \in ZI; 1 \leq i \leq 4 \right\} \subseteq V_2;$$

$W = W_1 \cup W_2$ is a neutrosophic biset bilinear subalgebra of V over the biset $S = S_1 \cup S_2$.

*Example 5.47:* Let $V = V_1 \cup V_2$ where

$$V_1 = \left\{ \begin{pmatrix} a_1 & a_2 & a_3 \\ 0 & a_4 & a_5 \\ 0 & 0 & a_6 \end{pmatrix} \middle| a_i \in ZI; 1 \leq i \leq 6 \right\}$$

is neutrosophic set linear algebra over the set $S_1 = 11Z^+I \cup \{0\}$ and

$$V_2 = \{(a_1 \quad a_2 \quad a_3 \quad a_4 \quad a_5) \mid a_i \in QI; 1 \leq i \leq 5\}$$

be a neutrosophic set linear algebra over the set $S_2 = 17Z^+I \cup \{0\}$. $V = V_1 \cup V_2$ is a neutrosophic biset bilinear algebra over the biset $S = S_1 \cup S_2$.

Take $W = W_1 \cup W_2 \subseteq V_1 \cup V_2$ where

$$W_1 = \left\{ \begin{pmatrix} a_1 & a_2 & a_3 \\ 0 & a_4 & a_5 \\ 0 & 0 & a_6 \end{pmatrix} \middle| a_i \in Z^+I; 1 \leq i \leq 6 \right\} \subseteq V_1$$

and

$$W_2 = \{(a_1 \quad a_2 \quad a_3 \quad a_4 \quad a_5) \mid a_i \in ZI; 1 \leq i \leq 5\} \subseteq V_2;$$

$W$ is a neutrosophic biset bilinear subalgebra of $V$ over the biset $S = S_1 \cup S_2$.

Now we proceed onto define the notion of quasi neutrosophic biset bilinear algebra.

**DEFINITION 5.23:** *Let $V = V_1 \cup V_2$ where $V_1$ is a neutrosophic set vector space over the set $S_1$ and $V_2$ is a neutrosophic set linear algebra over the set $S_2$ where $V_1 \neq V_2$; $V_1 \not\subset V_2$, $V_2 \not\subset V_1$ and $S_1 \neq S_2$; $S_1 \not\subset S_2$ and $S_2 \not\subset S_1$. We define $V = V_1 \cup V_2$ to be a quasi neutrosophic biset bilinear algebra over the biset $S = S_1 \cup S_2$.*

We will illustrate this by some examples.

*Example 5.48:* Let $V = V_1 \cup V_2$ where

$$V_1 = \left\{ (a_1, a_2, a_3, a_4), \begin{bmatrix} a_1 \\ a_2 \\ a_3 \\ a_4 \\ a_5 \end{bmatrix} \middle| a_i \in ZI; 1 \leq i \leq 5 \right\}$$

is a neutrosophic set vector space over the set $S_1 = 3Z^+ I \cup \{0\}$ and

$$V_2 = \left\{ \begin{pmatrix} a_1 & a_2 & a_3 \\ a_5 & a_6 & a_7 \\ a_8 & a_4 & a_9 \end{pmatrix} \middle| a_i \in QI; 1 \le i \le 9 \right\}$$

is a neutrosophic set linear algebra over the set $S_2 = 5Z^+ \cup \{0\}$. $V = V_1 \cup V_2$ is a quasi neutrosophic biset bilinear algebra over the biset $S = S_1 \cup S_2$.

*Example 5.49:* Let $V = V_1 \cup V_2$ where

$$V_1 = \left\{ \begin{pmatrix} a_1 & a_2 \\ a_3 & a_4 \\ a_5 & a_6 \end{pmatrix}, \begin{pmatrix} a_1 & a_2 & a_3 & a_4 \\ a_5 & a_6 & a_7 & a_8 \\ a_9 & a_{10} & a_{11} & a_{12} \end{pmatrix} \right.$$

where $a_i \in QI$; $1 \le i \le 12\}$ is a neutrosophic set vector space over the set $S_1 = ZI$ and

$$V_2 = \left\{ \begin{pmatrix} a_1 & a_2 & a_3 \\ a_4 & a_5 & a_6 \\ a_7 & a_8 & a_9 \end{pmatrix} \middle| a_i \in Q^+I \cup \{0\}; 1 \le i \le 9 \right\}$$

is a neutrosophic set linear algebra over the set $S_2 = Q^+I \cup \{0\}$. Thus $V = V_1 \cup V_2$ is a quasi neutrosophic biset bilinear algebra over the biset $S = S_1 \cup S_2$.

It is important to mention here that we have substructure defined on them.

**DEFINITION 5.24:** *Let $V = V_1 \cup V_2$ where V is a quasi neutrosophic biset bilinear algebra over the biset $S = S_1 \cup S_2$. Take $W = W_1 \cup W_2 \subseteq V_1 \cup V_2$ if W is a quasi neutrosophic biset bilinear algebra over the biset $S = S_1 \cup S_2$, then we define W to* *be a quasi neutrosophic biset bilinear subalgebra of V over the biset $S = S_1 \cup S_2$.*

We will illustrate this situation by some examples.

***Example 5.50:*** Let $V = V_1 \cup V_2$ where

$$V_1 = \left\{ \begin{pmatrix} a_1 & a_2 & a_3 \\ 0 & a_4 & a_5 \end{pmatrix}, \begin{pmatrix} a_1 & a_2 \\ a_3 & a_4 \\ a_5 & 0 \end{pmatrix} \middle| a_i \in ZI; 1 \le i \le 5 \right\}$$

is a neutrosophic set vector space over the set $S_1 = 3Z^+I \cup \{0\}$ and $V_2 = \{(a_1, a_2, a_3, a_4, a_5, a_6, a_7) \mid a_i \in QI, 1 \le i \le 7\}$ is a neutrosophic set linear algebra over the set $S_2 = 13Z^+I \cup \{0\}$. $V = V_1 \cup V_2$ is a quasi neutrosophic biset bilinear algebra over the biset $S = S_1 \cup S_2$.
Take $W = W_1 \cup W_2 \subseteq V_1 \cup V_2$ where

$$W_1 = \left\{ \begin{pmatrix} a_1 & a_2 & a_3 \\ 0 & a_4 & a_5 \end{pmatrix}, \begin{pmatrix} a_1 & a_2 \\ a_3 & a_4 \\ a_5 & 0 \end{pmatrix} \middle| a_i \in 3ZI; 1 \le i \le 5 \right\}$$

and $W_2 = \{(a_1, 0, a_3, 0, a_5, 0, a_7) \mid a_1, a_3, a_5, a_7 \in QI\}$; W is a quasi neutrosophic biset bilinear subalgebra of V over the biset $S = S_1 \cup S_2$.

***Example 5.51:*** Let $V = V_1 \cup V_2$ where $V_1 = \{Z^+I[x];$ all polynomials in the variable x with coefficients from $Z^+I \cup \{0\}\}$ is a neutrosophic set linear algebra over the set $S_1 = 3Z^+I \cup \{0\}$ and $V_2 = \{3ZI[x]$ and $7ZI[x]$; that all polynomials in the variable x with coefficients from the 3ZI and 7ZI respectively$\}$ is a neutrosophic set vector space over the set $S_2 = 8Z^+I \cup \{0\}$. Thus $V = V_1 \cup V_2$ is a quasi neutrosophic biset bilinear algebra over the biset $S = S_1 \cup S_2$.
Take $W = W_1 \cup W_2$ where $W_1 = \{$The set of all even degree polynomials in the variable x with coefficients from $Z^+I \cup \{0\}\}$

$\subseteq V_1$ and $W_2 = \{3Z^+I\ [x]$ and $7Z^+I\ [x]$; collection of all polynomials in the variable x with coefficients from $3Z^+I \cup \{0\}$ and $7Z^+ \cup \{0\}$ respectively$\} \subseteq V_2$; $W = W_1 \cup W_2$ is a quasi neutrosophic biset bilinear subalgebra of V over the biset $S = S_1 \cup S_2$.

It may so happen that a neutrosophic biset bilinear algebra contain quasi neutrosophic bilinear subalgebra. We now define this concept.

**DEFINITION 5.25:** *Let $V = V_1 \cup V_2$ be a neutrosophic biset bilinear algebra over the biset $S = S_1 \cup S_2$. Suppose $W = W_1 \cup W_2 \subseteq V_1 \cup V_2$ where $W_1$ is just a proper subset of $V_1$ but $W_1$ is only a neutrosophic set vector space over the set $S_1$ and $W_2$ is proper subset of $V_2$ and $W_2$ is a neutrosophic set linear algebra over the set $S_2$. We call $W = W_1 \cup W_2 \subseteq V_1 \cup V_2$ to be a quasi neutrosophic biset bilinear subalgebra of V over the biset $S = S_1 \cup S_2$.*

We will illustrate this situation by some simple examples.

*Example 5.52:* Let
$$V = V_1 \cup V_2$$
$$= \{(a_1\ a_2\ a_3\ a_4\ a_5) | a_i \in QI; 1 \le i \le 5\}$$
$$\cup \left\{ \begin{pmatrix} a_1 & a_2 \\ a_3 & a_4 \end{pmatrix} \middle| a_i \in ZI; 1 \le i \le 4 \right\}$$

be a neutrosophic biset bilinear algebra over the biset $S = S_1 \cup S_2$ where $S_1 = Q^+I \cup \{0\}$ and $S_2 = 3ZI$. Take $W = W_1 \cup W_2 \subseteq V_1 \cup V_2$ where

$$W_1 = \{(a_1\ a_2\ 0\ a_3\ a_4) | a_i \in QI; 1 \le i \le 4\} \subseteq V_1$$

is only a neutrosophic set vector space over the set $S_1 = Q^+I \cup \{0\}$.

$$W_2 = \left\{ \begin{pmatrix} 0 & a_2 \\ a_3 & 0 \end{pmatrix} \middle| a_i \in 3Z^+I; i = 2,3 \right\} \subseteq V_2$$

be a neutrosophic set vector space over the set $S_2 = 3ZI$. $W = W_1 \cup W_2 \subseteq V_1 \cup V_2$ is a quasi neutrosophic biset bilinear subalgebra of V over the biset $S = S_1 \cup S_2$.

*Example 5.53:* Let

$$V = V_1 \cup V_2 = \left\{ \begin{pmatrix} a_1 & a_2 \\ a_4 & a_3 \end{pmatrix} \middle| a_i \in ZI; 1 \leq i \leq 4 \right\} \cup$$

{ZI[x]; all polynomials in the variable x with coefficients from ZI} be a neutrosophic biset bilinear algebra over the biset $S = S_1 \cup S_2 = 3ZI \cup 7ZI$.
Take

$$W = W_1 \cup W_2 = \left\{ \begin{pmatrix} a & b \\ 0 & c \end{pmatrix} \middle| a, b, c \in ZI \right\} \cup$$

{3ZI[x] and 7ZI[x]; all polynomials in the variable x with coefficients from 3ZI and 7ZI respectively} $\subseteq V_1 \cup V_2$; W is only a quasi neutrosophic biset bilinear subalgebra over the biset S.

Now we proceed onto define yet another new substructure.

**DEFINITION 5.26:** *Let $V = V_1 \cup V_2$ be a neutrosophic biset bivector space (bilinear algebra) over the biset $S = S_1 \cup S_2$. Let $W = W_1 \cup W_2 \subseteq V_1 \cup V_2$ be a proper subbiset of $V_1 \cup V_2$. Let $T = T_1 \cup T_2 \subseteq S_1 \cup S_2 = S$, where T is also a proper subbiset of $S = S_1 \cup S_2$. If W is a neutrosophic biset bivector space (bilinear algebra) over the biset T then we call W to be a neutrosophic subbiset bivector subspace (bilinear subalgebra) of V over the subbiset T of S.*

We will illustrate this situation by some examples.

**Example 5.54:** Let $V = V_1 \cup V_2$ where

$$V_1 = \left\{ \begin{pmatrix} a & b & c \\ d & e & f \end{pmatrix}, \begin{pmatrix} a & b \\ c & d \end{pmatrix} \middle| a,b,c,d,e,f \in QI \right\}$$

be a neutrosophic set vector space over the set $S_1 = ZI$ and

$$V_2 = \left\{ \begin{pmatrix} a & b & c & d \\ e & f & g & h \end{pmatrix}, \begin{pmatrix} a \\ b \\ c \\ d \\ e \end{pmatrix} \middle| a,b,c,d,e,f,g,h \in QI \right\}$$

be a neutrosophic set vector space over the set $S_2 = Q^+I$. $V = V_1 \cup V_2$ is a neutrosophic biset bivector space over the biset $S = S_1 \cup S_2$.
Take

$$W = W_1 \cup W_2$$

$$= \left\{ \begin{pmatrix} a & b & c \\ d & e & f \end{pmatrix} \middle| a,b,c,d,e,f \in ZI \right\} \cup \left\{ \begin{pmatrix} a \\ b \\ c \\ d \\ e \end{pmatrix} \middle| a,b,c,d,e \in QI \right\}$$

$$\subseteq V_1 \cup V_2$$

and $T = 3ZI \cup 7ZI \subseteq S_1 \cup S_2$; W is a neutrosophic subbiset bivector subspace of V over the subbiset $T = T_1 \cup T_2 = 3ZI \cup 7ZI \subseteq S_1 \cup S_2 = S$.

**Example 5.55:** Let $V = V_1 \cup V_2 = \{QI[x]$, the set of all polynomials in the variable x with coefficients from $QI\} \cup$

$$\left\{ \begin{pmatrix} a & b & c \\ d & e & f \\ g & h & i \end{pmatrix} \middle| a,b,c,d,e,f,g,h,i \in RI \right\}$$

be the neutrosophic biset bilinear algebra over the biset $S = S_1 \cup S_2 = ZI \cup Q^+I$.

Take $W = W_1 \cup W_2 = \{ZI[x]$; that is the set of all polynomials in the variable x with coefficients from $ZI\} \cup$

$$\left\{ \begin{pmatrix} a & b & c \\ d & e & f \\ g & h & i \end{pmatrix} \middle| a,b,c,d,e,f,g,h,i \in ZI \right\}$$

W is a neutrosophic subbiset bilinear subalgebra of V over the subbiset $T = T_1 \cup T_2 = 3ZI \cup Z^+I \subseteq S_1 \cup S_2$.

Now we proceed onto define the new notion of neutrosophic semigroup bivector space over the semigroup S.

**DEFINITION 5.27:** *Let $V = V_1 \cup V_2$ where ($V_1 \neq V_2$, $V_1 \not\subset V_2$ and $V_2 \not\subset V_1$) $V_1$ is a neutrosophic semigroup vector space over the semigroup S and $V_2$ is a neutrosophic semigroup vector space over the semigroup S, then we say $V = V_1 \cup V_2$ to be the neutrosophic semigroup bivector space over the semigroup S.*

We will illustrate this situation by some examples.

***Example 5.56:*** Let $V = V_1 \cup V_2 = \{ZI \times ZI \times Z^+I \times 3Z^+I\} \cup$

$$\left\{ \begin{pmatrix} a & b & c \\ d & e & f \end{pmatrix} \middle| a,b,c,d,e,f \in QI \right\}$$

be a neutrosophic semigroup bivector space over the semigroup $S = Z^+I$.

*Example 5.57:* Let $V = V_1 \cup V_2 = \{(a, b, c), (a, b, c, d, e) \mid a, b, c, d, e \in Z_7I\} \cup$

$$\left\{ \begin{pmatrix} a & b \\ c & d \\ e & f \end{pmatrix}, \begin{pmatrix} a \\ b \\ c \\ d \end{pmatrix} \middle| a, b, c, d \in Z_7I \right\}$$

be a neutrosophic semigroup bivector space over the semigroup $= Z_7I$.

Now we proceed onto define the notion of neutrosophic semigroup bivector subspace.

**DEFINITION 5.28:** *Let $V = V_1 \cup V_2$ be a neutrosophic semigroup bivector space over the semigroup S. Let $W = W_1 \cup W_2 \subseteq V_1 \cup V_2 = V$ be a proper biset of V; if W is a neutrosophic semigroup bivector space over the semigroup S then we call W to be a neutrosophic semigroup bivector subspace of V over the semigroup S.*

We will illustrate this definition by some examples.

*Example 5.58:* Let

$$V = V_1 \cup V_2 = \left\{ \begin{pmatrix} a & a & a \\ a & a & a \end{pmatrix} \middle| a \in Z^+I \cup \{0\} \right\} \cup$$

$$\left\{ \begin{pmatrix} a & a \\ a & a \\ b & b \\ b & b \end{pmatrix}, \begin{bmatrix} a \\ b \\ c \\ d \end{bmatrix} \middle| a, b, c, d \in Z^+I \cup \{0\} \right\}$$

be a neutrosophic semigroup bivector space over the semigroup $S = 3Z^+I \cup \{0\}$.

Take

$$W = W_1 \cup W_2 = \left\{ \begin{pmatrix} a & a & a \\ a & a & a \end{pmatrix} \middle| a \in 5Z^+I \cup \{0\} \right\} \cup$$

$$\left\{ \begin{bmatrix} a \\ b \\ c \\ d \end{bmatrix} \middle| a,b,c,d \in Z^+I \cup \{0\} \right\} \subseteq V_1 \cup V_2,$$

W is a neutrosophic semigroup bivector subspace of V over the semigroup S.

***Example 5.59:*** Let $V = V_1 \cup V_2 = \{(0\ 0\ 0), (1\ 1\ I), (0\ 0\ 0\ 0), (I\ 1\ 1\ I), (1\ 1\ 0\ I), (0\ 0\ I\ 0), (0\ 0\ 0\ 0), (I\ 1\ 1\ 1\ I)\} \cup \{Z_2I \times Z_2I \times Z_2I \times Z_2I \times Z_2I\}$ be a neutrosophic semigroup bivector space over the semigroup $S = \{0, 1\} = Z_2 =$ (addition modulo 2). Consider $W = W_1 \cup W_2 = \{(0\ 0\ 0\ 0), (I\ 1\ 1\ I), (1\ 1\ 0\ I), (0\ 0\ I\ 0)\} \cup \{Z_2I \times Z_2I \times \{0\} \times \{0\} \times \{0\}\} \subseteq V_1 \cup V_2$; W is a neutrosophic semigroup bivector subspace of V over the semigroup $S = Z_2$.

**DEFINITION 5.29:** *Let $V = V_1 \cup V_2$ be such that $V_1$ is a neutrosophic semigroup linear algebra over the semigroup S and $V_2$ is a neutrosophic semigroup linear algebra over the semigroup S; ($V_1 \neq V_2$; $V_1 \not\subseteq V_2$ and $V_2 \not\subseteq V_1$). $V = V_1 \cup V_2$ is defined as the neutrosophic semigroup bilinear algebra over the semigroup S.*

We will illustrate this definition by some examples.

***Example 5.60:*** Let
$$V = V_1 \cup V_2$$
$$= \left\{ \begin{bmatrix} a & a & a \\ a & a & a \\ a & a & a \end{bmatrix} \middle| a \in Z_{12}I \right\} \cup \{(a\ a\ a\ a\ a\ a\ a) \mid a \in Z_{12}\ I\}$$

be a neutrosophic semigroup bilinear algebra over the semigroup $S = Z_{12}I$.

*Example 5.61:* Let $V = V_1 \cup V_2 = \{Z_5I [x]$; all polynomials in the variable x with coefficients from $Z_5I\} \cup \{Z_5I \times Z_5I \times Z_5I\}$ be a neutrosophic semigroup bilinear algebra over the semigroup $S = Z_5I$.

It is important at this juncture to mention that every neutrosophic semigroup bilinear algebra is a neutrosophic semigroup bivector space but a neutrosophic semigroup bivector space in general is not a neutrosophic semigroup bilinear algebra.
   The interested reader can prove the above statement.

Now we proceed onto define the notion of bigenerator and bidimension of this algebraic structure.

**DEFINITION 5.30:** *Let $V = V_1 \cup V_2$ be a neutrosophic semigroup bivector space (bilinear algebra) over the semigroup S. Let $X = X_1 \cup X_2 \subseteq V_1 \cup V_2$, if $X_1$ generates $V_1$ as a neutrosophic semigroup vector space (linear algebra) over the semigroup S and if $X_2$ generates $V_2$ as a neutrosophic semigroup vector space (linear algebra) over the semigroup S then, $X = X_1 \cup X_2$ is called the bigenerator of the neutrosophic semigroup bivector space (bilinear algebra) over the semigroup S.*

   *The bidimension of V is $|X| = (|X_1| \cup |X_2|)$ or $(|X_1|, |X_2|)$ over the semigroup S. Even if one of $|X_1|$ or $|X_2|$ are infinite we say the bidimension of V is infinite only when both $|X_1|$ and $|X_2|$ is finite we say bidimension of V is finite.*

We will illustrate this by some simple examples.

*Example 5.62:* Let $V = V_1 \cup V_2 = \{(I, 0, 0), (1, 0, 0), (0\ 0\ 0), (0, 0), (I, I), (1, 1), (1, I), (I, 1)\} \cup \{(a\ a\ a) \mid a \in Z_2I\}$ be a neutrosophic semigroup bivector space over the semigroup $S = Z_2$. $X = X_1 \cup X_2 = \{(I, 0, 0), (1, 0, 0)\ (I, I)\ (1, 1), (1, I), (I, 1)\} \cup$

{(1 1 1), (I I I)} $\subseteq V_1 \cup V_2$ is the bigenerator of V and the bidimension of V is (6, 2).

*Example 5.63:* Let
$$V = V_1 \cup V_2$$
$$= \left\{ \begin{bmatrix} a & a \\ a & a \end{bmatrix} \middle| a \in ZI \right\} \cup \{(a\ a\ a\ a\ a\ a) \mid a \in ZI\}$$

be a neutrosophic semigroup linear algebra over the semigroup $S = ZI$.
Take
$$X = \left\{ \begin{bmatrix} 1 & 1 \\ 1 & 1 \end{bmatrix} \right\} \cup \{(1\ 1\ 1\ 1\ 1\ 1)\} = X_1 \cup X_2 \subseteq V_1 \cup V_2;$$

X bigenerates V and bidimension of V is (1, 1).

**DEFINITION 5.31:** *Let $V = V_1 \cup V_2$ where $V_1$ is a neutrosophic semigroup vector space over the semigroup S and $V_2$ is a neutrosophic semigroup linear algebra over the semigroup S, then we call V to be a quasi neutrosophic semigroup bilinear algebra over the semigroup S.*

We will illustrate this by some examples.

*Example 5.64:* Let
$$V = V_1 \cup V_2$$
$$= \{(a\ a\ a),\ (a\ a\ a\ a),\ (a\ a),\ (a\ a\ a\ a\ a\ a) \mid a \in Z_3I\} \cup$$
$$\left\{ \begin{pmatrix} a & a & a & a \\ a & a & a & a \\ a & a & a & a \end{pmatrix} \middle| a \in Z_3I \right\}$$

be a neutrosophic semigroup quasi bilinear algebra over the semigroup $S = Z_3I$.

*Example 5.65:* Let

$$V = V_1 \cup V_2$$

$$= \left\{ \begin{pmatrix} a & a \\ a & a \\ b & b \\ b & b \end{pmatrix} \middle| a, b \in ZI \right\} \cup \{(a\ a\ b\ b) \mid a, b \in ZI\}$$

be a quasi neutrosophic semigroup bilinear algebra over the semigroup $S = ZI$.

Now we proceed onto define another substructure.

**DEFINITION 5.32:** *Let $V = V_1 \cup V_2$ be a quasi neutrosophic semigroup bilinear algebra over the semigroup S. Let $W = W_1 \cup W_2 \subseteq V_1 \cup V_2$ be a proper subset of V such that W is a quasi neutrosophic semigroup bilinear algebra over the semigroup S, then we call W to be a quasi neutrosophic semigroup bilinear subalgebra of V over the semigroup S.*

*Example 5.66:* Let

$$V = V_1 \cup V_2 = \left\{ \begin{pmatrix} a & a & a & a & a \\ a & a & a & a & a \end{pmatrix} \middle| a \in ZI \right\}$$

$$\cup \{(a\ a\ a\ a\ a\ a) \mid a \in ZI\}$$

be a quasi neutrosophic semigroup bilinear algebra over the semigroup $S = ZI$.
Take

$$W = W_1 \cup W_2$$

$$= \left\{ \begin{pmatrix} a & a & a & a & a \\ a & a & a & a & a \end{pmatrix} \middle| a \in 5ZI \right\} \cup \{(a\ a\ a\ a\ a) \mid a \in 5ZI\}$$

$$\subseteq V_1 \cup V_2,$$

W is a quasi neutrosophic semigroup bilinear subalgebra of V over the semigroup S = ZI.

*Example 5.67:* Let

$$V = V_1 \cup V_2 = \left\{ \begin{bmatrix} a & b \\ c & d \end{bmatrix} \middle| a,b,c,d \in ZI \right\} \cup$$

$$\left\{ (a, a), \begin{bmatrix} a \\ a \\ a \\ a \\ a \end{bmatrix}, \begin{bmatrix} a & a \\ a & a \end{bmatrix} \middle| a \in ZI \right\}$$

be a quasi neutrosophic semigroup bilinear algebra over the semigroup S = ZI. Take

$$W = W_1 \cup W_2 = \left\{ \begin{pmatrix} a & a \\ a & a \end{pmatrix} \middle| a \in ZI \right\} \cup \left\{ (a, a), \begin{bmatrix} a & a \\ a & a \end{bmatrix} \right\}$$

$\subseteq V_1 \cup V_2$; W is a quasi neutrosophic semigroup bilinear subalgebra of V over the semigroup S = ZI.

Now we proceed onto define the new notion of neutrosophic group bivector space and neutrosophic group bilinear algebras.

**DEFINITION 5.33:** *Let $V = V_1 \cup V_2$ be such that $V_1 \neq V_2$, $V_1 \not\subseteq V_2$ and $V_2 \not\subseteq V_1$, $V_1$ and $V_2$ are neutrosophic group vector spaces over the same group G then we call V to be a neutrosophic group bivector space defined over the group G.*

*Note:* The group G can be an ordinary group or a neutrosophic group.

We will illustrate this by some examples.

*Example 5.68:* Let

$$V = V_1 \cup V_2 = \left\{ \begin{pmatrix} a & a \\ a & a \end{pmatrix} \middle| a \in Z_3I \right\} \cup \{(a\ a\ a\ a) \mid a \in Z_3I\}.$$

$V = V_1 \cup V_2$ is a neutrosophic group bivector space over the group $G = Z_3I$. In fact $V = V_1 \cup V_2$ is also a neutrosophic group bivector space over the group $G = Z_3$.

*Example 5.69:* Let

$V = V_1 \cup V_2 = \{(0, 0, 0)\ (I, I, I),\ (0, I, I)\ (I, 0, I)\ (0, 0, 0, 0, 0)\ (I, I, 0, I, I)\ (I, 0, I, 0, I)\ (0, 0, I, I, 0)\ (0, 0, 0, 0)\ (I, I, I, I)\ (I, I, 0, 0)\ (0, 0, I, I)\} \cup$

$$\left\{ \begin{pmatrix} a & a \\ a & a \end{pmatrix} \begin{pmatrix} a & a & a & a \\ a & a & a & a \end{pmatrix} \middle| a \in Z_2I \right\}$$

be a neutrosophic group bivector space over the group $G = Z_2I$.

Now we proceed onto define the neutrosophic group bivector subspace.

**DEFINITION 5.34:** *Let $V = V_1 \cup V_2$ be a neutrosophic group bivector space over the group G. $W = W_1 \cup W_2 \subseteq V_1 \cup V_2$ is said to be a neutrosophic group bivector subspace of V over G if W itself is a neutrosophic group bivector space over the group G.*

We will illustrate this by some examples.

*Example 5.70:* Let

$$V = V_1 \cup V_2 = \left\{ \begin{pmatrix} a & a & a \\ a & a & a \end{pmatrix}, \begin{pmatrix} a & a \\ a & a \\ a & a \\ a & a \\ a & a \end{pmatrix} \middle| a \in Z_5I \right\} \cup$$

$$\{(a, a, a)\ (a, a)\ (a, a, a, a, a)\ |\ a \in Z_5I\}$$

be a neutrosophic group bivector space over the group $G = Z_5I$. Let

$$W = W_1 \cup W_2 \subseteq \left\{ \begin{pmatrix} a & a & a \\ a & a & a \end{pmatrix} \middle|\ a \in Z_5I \right\} \cup$$

$$\{(a, a)\ (a, a, a)\ |\ a \in Z_5I\} \subseteq V_1 \cup V_2,$$

W is a neutrosophic group bivector subspace of V over the group $G = Z_5I$.

*Example 5.71:* Let

$$V = V_1 \cup V_2 = \left\{ \begin{pmatrix} a & b \\ c & d \end{pmatrix}, \begin{pmatrix} a & a & a & a \\ a & a & a & a \end{pmatrix} \middle|\ a, b, c, d \in QI \right\} \cup$$

$$\left\{ (a\ a\ a\ a\ a\ a\ a),\ (a\ a\ a),\ (a, a)\ \begin{bmatrix} a \\ a \\ a \\ a \\ a \end{bmatrix}, \begin{bmatrix} a \\ a \\ a \\ a \end{bmatrix}, \begin{bmatrix} a \\ a \\ a \\ a \\ a \\ a \end{bmatrix}\ \middle|\ a \in QI \right\}$$

be a neutrosophic group bivector space over the group $G = ZI$.

$$W = W_1 \cup W_2 = \left\{ \begin{pmatrix} a & b \\ c & d \end{pmatrix} \middle|\ a, b, c, d \in QI \right\} \cup$$

$$\left\{ (a\ a\ a)\ \begin{bmatrix} a \\ a \\ a \end{bmatrix}, \begin{bmatrix} a \\ a \\ a \\ a \\ a \end{bmatrix},\ (a, a)\ \middle|\ a \in QI \right\}$$

⊆ $V_1 \cup V_2$ is a neutrosophic group bivector subspace of V over the group G = ZI.

Now we proceed onto define the notion of pseudo neutrosophic semigroup bivector subspace of V.

**DEFINITION 5.35:** *Let $V = V_1 \cup V_2$ be a neutrosophic group bivector space over the group G. Let $W = W_1 \cup W_2 \subseteq V_1 \cup V_2$ and $H \subseteq G$ be a semigroup of the group G. If W is a neutrosophic semigroup bivector space over the semigroup H then we call W to be a neutrosophic pseudo semigroup bivector subspace of V.*

We will illustrate this by some examples.

*Example 5.72:* Let $V = V_1 \cup V_2 = \{ZI[x]\} \cup \{(a\ b\ c) \mid a, b, c \in ZI\}$ be a neutrosophic group bivector space over the group G = ZI. Take $W = W_1 \cup W_2 = \{Z^+I\,[x]\} \cup \{(a\ b\ c) \mid a, b, c \in Z^+I\} \subseteq V_1 \cup V_2$. W is a neutrosophic semigroup bivector space over the semigroup $Z^+I \subseteq ZI$. W is a pseudo neutrosophic semigroup bivector subspace of V over the semigroup $Z^+I \subseteq Z$.

*Example 5.73:* Let

$$V = V_1 \cup V_2 = \left\{ \begin{pmatrix} a & b & c \\ d & e & f \end{pmatrix} \middle| a,b,c,d,e,f \in ZI \right\}$$

$\cup \{(a, b, c, d, e) \mid a, b, c, d, e \in ZI\}$ be a neutrosophic group bivector space over the group ZI.
Take

$$W = W_1 \cup W_2 = \left\{ \begin{pmatrix} a & b & 0 \\ c & d & 0 \end{pmatrix} \middle| a,b,c,d \in Z^+I \cup \{0\} \right\}$$

$\cup \{(a, a, a, a, a) \mid a \in ZI\} \subseteq V_1 \cup V_2$.

W is a neutrosophic semigroup bivector subspace of V over the semigroup $2Z^+I \cup \{0\}$.

**DEFINITION 5.36:** *Let $V = V_1 \cup V_2$ be a neutrosophic group bivector space over the group G. Let $W = W_1 \cup W_2 \subseteq V_1 \cup V_2 = V$; take S a proper subset of G. If W is a neutrosophic set bivector space over the set S then we call W to be the pseudo neutrosophic set bivector subspace of V over the set S.*

We will illustrate this situation by some simple examples.

*Example 5.74:* Let
$$V = V_1 \cup V_2$$

$$= \{(a\ b, c, d) \mid a, b, c, d \in 2ZI\} \cup \left\{ \begin{pmatrix} a & b \\ c & d \end{pmatrix} \middle| a, b, c, d \in ZI \right\}$$

be neutrosophic group bivector space over the group $G = ZI$. Take

$$W = W_1 \cup W_2$$
$$= \{(a\ b, c, d) \mid a, b, c, d \in 2Z^+I \cup \{0\}\} \cup$$

$$\left\{ \begin{pmatrix} a & b \\ c & d \end{pmatrix} \middle| a, b, c, d \in 4Z^+I \cup \{0\} \right\}.$$

W is a neutrosophic pseudo set bivector subspace of V over the set $S = \{0, 2I, 2^2I, 2^3I, \ldots, 2^nI \mid n \in N\}$.

*Example 5.75:* Let

$$V = V_1 \cup V_2 = \{QI \times ZI \times QI \times QI\} \cup$$
$$\left\{ \begin{pmatrix} a & b \\ c & d \end{pmatrix} \middle| a, b, c, d \in 3ZI \right\}$$

be a neutrosophic group bivector space over the group $G = ZI$. Take

$$W = W_1 \cup W_2 = \{Q^+I \times \{0\} \times \{0\} \times Q^+I\} \cup$$
$$\left\{ \begin{pmatrix} a & a \\ a & a \end{pmatrix} \middle| a,b,c,d \in 9Z^+I \right\}$$

$\subseteq V_1 \cup V_2$; W is a neutrosophic pseudo set bivector subspace of V over the set $S = \{0, 3I, 3^2I, \ldots, 3^nI \mid n \in N\}$.

Now we proceed onto define yet another new notion, viz. neutrosophic bisemigroup bivector space.

**DEFINITION 5.37:** *Let $V = V_1 \cup V_2$ where $V_1$ is a neutrosophic semigroup vector space over the semigroup $S_1$ and $V_2$ is also neutrosophic semigroup vector space over the semigroup $S_2$. ($S_1 \neq S_2$, $S_1 \not\subseteq S_2$ and $S_2 \not\subseteq S_1$ and $V_1 \neq V_2$, $V_1 \not\subseteq V_2$ and $V_2 \not\subseteq V_1$) we call V to be the neutrosophic bisemigroup bivector space over the bisemigroup $S = S_1 \cup S_2$.*

We illustrate this by some examples.

***Example 5.76:*** Let $V = V_1 \cup V_2$ where

$$V_1 = \left\{ \begin{pmatrix} a \\ a \end{pmatrix}, (a, a, a, a, a) \mid a \in Z_7I \right\}$$

and

$$V_2 = \left\{ \begin{pmatrix} a & a \\ a & a \\ a & a \\ a & a \end{pmatrix}, \begin{pmatrix} a \\ a \\ a \end{pmatrix} \middle| a \in Z^+I \right\}.$$

V is a neutrosophic bisemigroup bivector space over the bisemigroup $S = S_1 \cup S_2 = Z_7I \cup Z^+I$.

**DEFINITION 5.38:** *Let $V = V_1 \cup V_2$ where $V_1$ is a neutrosophic semigroup linear algebra over the semigroup $S_1$ and $V_2$ is a neutrosophic semigroup linear algebra over the semigroup $S_2$* ($V_1 \neq V_2$, $V_1 \not\subseteq V_2$ and $V_2 \not\subseteq V_1$ and $S_1 \neq S_2$, $S_2 \not\subseteq S_1$, $S_1 \not\subseteq S_2$). V is a neutrosophic bisemigroup bilinear algebra over the bisemigroup $S = S_1 \cup S_2$.

We will illustrate this by some examples.

*Example 5.77:* Let

$$V = V_1 \cup V_2 = \left\{ \begin{pmatrix} a & b & c \\ d & e & f \end{pmatrix} \middle| a,b,c,d,e,f \in QI \right\} \cup$$

{ZI[x]; all polynomials in the variable x with coefficients from ZI} be a neutrosophic bisemigroup bilinear algebra over the bisemigroup $S = S_1 \cup S_2 = Q^+I \cup ZI$.

*Example 5.78:* Let

$$V = V_1 \cup V_2 = \{(x, y) \mid x, y \in Z_{17}I\} \cup$$

$$\left\{ \begin{pmatrix} a & a & a \\ a & 0 & a \end{pmatrix} \middle| a \in Z_{12}I \right\}$$

be a neutrosophic bisemigroup bilinear algebra over the bisemigroup $S = S_1 \cup S_2 = Z_{17}I \cup Z_{12}I$.

**DEFINITION 5.39:** *Let $V = V_1 \cup V_2$ be a neutrosophic bisemigroup bivector space over the bisemigroup $S = S_1 \cup S_2$. $W = W_1 \cup W_2 \subseteq V_1 \cup V_2$ is called a neutrosophic bisemigroup bivector subspace of V over the bisemigroup S if W is a neutrosophic bisemigroup bivector space over the bisemigroup S.*

We give examples of the definition.

*Example 5.79:* Let
$V = V_1 \cup V_2 = \{(x, y, z), (a, a, a, a) \mid x, y, z, a \in Z_{12}I\} \cup$

$$\left\{ \begin{pmatrix} a \\ a \\ a \\ a \\ a \end{pmatrix}, (a,b) \,\middle|\, a,b \in Z_{21}I \right\}$$

be a neutrosophic bisemigroup bivector space over the bisemigroup $S = S_1 \cup S_2 = Z_{12}I \cup Z_{21}I$.
Take
$$W = W_1 \cup W_2$$

$$= \{(x, y, z) \mid x, y, z \in Z_{12}I\} \cup \left\{ \begin{pmatrix} a \\ a \\ a \\ a \\ a \end{pmatrix} \,\middle|\, a \in Z_{21}I \right\} \subseteq V_1 \cup V_2;$$

W is a neutrosophic bisemigroup bivector subspace of V over the bisemigroup S.

*Example 5.80:* Let
$$V = V_1 \cup V_2$$
$$= \{QI \times QI, ZI \times ZI \times ZI\} \cup$$

$$\left\{ \begin{pmatrix} a & a \\ a & a \end{pmatrix}, \begin{pmatrix} a & b \\ c & d \\ e & f \\ g & h \end{pmatrix} \,\middle|\, a,b,c,d,e,f,g,h \in Z_{20}I \right\}$$

be a neutrosophic bisemigroup bivector space over the bisemigroup $S = S_1 \cup S_2 = ZI \cup Z_{20}I$.
Take
$$W = W_1 \cup W_2$$

$$= \{QI \times QI\} \cup \left\{ \begin{pmatrix} a & a \\ a & a \end{pmatrix} \right\} \subseteq V_1 \cup V_2,$$

W is a neutrosophic bisemigroup bivector subspace of V over the bisemigroup S.

**DEFINITION 5.40:** *Let $V = V_1 \cup V_2$ be a neutrosophic bisemigroup bilinear algebra over the bisemigroup $S = S_1 \cup S_2$. $W = W_1 \cup W_2 \subseteq V_1 \cup V_2$ is defined as the neutrosophic bisemigroup bilinear subalgebra of V over $S_1 \cup S_2 = S$ if W is a neutrosophic bisemigroup bilinear algebra over the bisemigroup $S = S_1 \cup S_2$.*

We will give an example.

*Example 5.81:* Let
$$V = V_1 \cup V_2$$
$$= \left\{ \begin{pmatrix} a & b & c \\ d & e & f \end{pmatrix} \middle| a,b,c,d,e,f \in QI \right\} \cup \left\{ \begin{pmatrix} a \\ a \\ a \end{pmatrix} \middle| a \in Z_{12}I \right\}$$

be a neutrosophic bisemigroup bilinear algebra over the bisemigroup $S = S_1 \cup S_2 = ZI \cup Z_{12}I$.
Take
$$W = W_1 \cup W_2$$
$$= \left\{ \begin{pmatrix} a & a & a \\ a & a & a \end{pmatrix} \middle| a \in QI \right\} \cup \left\{ \begin{pmatrix} a \\ a \\ a \end{pmatrix} \middle| a \in \{0, 2I, 4I, 6I, 8I, 10I\} \right\}$$

$$\subseteq V_1 \cup V_2;$$

W is a neutrosophic bisemigroup bilinear subalgebra of V over the bisemigroup S.

**DEFINITION 5.41:** *Let $V = V_1 \cup V_2$ be a neutrosophic bisemigroup bivector space (bilinear algebra) over the bisemigroup $S = S_1 \cup S_2$. Take $W = W_1 \cup W_2 \subseteq (V_1 \cup V_2)$ ($W_1 \neq W_2$, $W_1 \not\subseteq W_2$, $W_2 \not\subseteq W_1$) and $T = T_1 \cup T_2 \subseteq S_1 \cup S_2$ ($T_1 \neq T_2$, $T_1 \not\subseteq T_2$ and $T_2 \not\subseteq T_1$) be a subbisemigroup of $S = S_1 \cup S_2$. If W is a neutrosophic bisemigroup bivector space (bilinear algebra) over the bisemigroup $T = T_1 \cup T_2$ then we call W to be a neutrosophic subbisemigroup subbivector space (subbilinar algebra) of V over the subbisemigroup T of $S = S_1 \cup S_2$.*

We will illustrate this by some simple examples.

*Example 5.82:* Let
$$V = V_1 \cup V_2$$

$$= \left\{ \begin{pmatrix} a & b \\ c & d \end{pmatrix}, (a,a,a,a,a) \,\middle|\, a,b,c,d \in Z_{12}I \right\} \cup$$

$$\left\{ \begin{pmatrix} a & a & a & a \\ a & a & a & a \end{pmatrix}, (x,y) \,\middle|\, a, x, y \in ZI \right\}$$

be a neutrosophic bisemigroup bivector space over the bisemigroup $S = S_1 \cup S_2 = Z_{12}I \cup ZI$. Take $W = W_1 \cup W_2$

$$= \left\{ \begin{pmatrix} a & b \\ c & d \end{pmatrix} \,\middle|\, a,b,c,d \in Z_{12}I \right\} \cup \{(x, y) \mid x, y \in ZI\}$$

$$\subseteq V_1 \cup V_2;$$

W is a neutrosophic subbisemigroup bivector subspace of V over the subbisemigroup $T = T_1 \cup T_2 = \{0, 2I, 4I, 6I, 8I, 10I\} \cup \{Z^+I \cup \{0\}\} \subseteq S_1 \cup S_2$ of the bisemigroup S.

*Example 5.83:* Let
$$V = V_1 \cup V_2$$

$$= \left\{ \begin{pmatrix} a & a & a & a & a \\ a & a & a & a & a \end{pmatrix} \middle| a \in QI \right\} \cup$$

$$\left\{ \begin{pmatrix} a & b \\ c & d \\ e & f \\ g & h \end{pmatrix} \middle| a, b, c, d, e, f, g, h \in Z_{24}I \right\}$$

be a neutrosophic bisemigroup bilinear algebra over the bisemigroup $S = S_1 \cup S_2 = QI \cup Z_{24}I$.
Take

$$W = W_1 \cup W_2$$

$$= \left\{ \begin{pmatrix} a & a & a & a & a \\ a & a & a & a & a \end{pmatrix} \middle| a \in ZI \right\} \cup \left\{ \begin{pmatrix} a & a \\ a & a \\ a & a \\ a & a \end{pmatrix} \middle| a \in Z_{24}I \right\}$$

$$\subseteq V_1 \cup V_2,$$

W is a neutrosophic subbisemigroup bilinear subalgebra of V over the subbisemigroup $T = T_1 \cup T_2 = ZI \cup 2Z_{24}I$ of the bisemigroup $S = S_1 \cup S_2 = QI \cup Z_{24}I$.

Now we proceed onto define the bidimension.

**DEFINITION 5.42:** *Let $V = V_1 \cup V_2$ be a neutrosophic bisemigroup bivector space (bilinear algebra) over the bisemigroup $S = S_1 \cup S_2$. Take $X = X_1 \cup X_2 \subseteq V_1 \cup V_2$; if $X_1$ generates $V_1$ as a neutrosophic semigroup vector space (linear algebra) over the semigroup $S_1$ and $X_2$ generates $V_2$ as a neutrosophic semigroup vector space (linear algebra) over $S_2$ then we say $X = X_1 \cup X_2$ bigenerates V over the bisemigroup $S = S_1 \cup S_2$.*

*The cardinality of $X = X_1 \cup X_2$ denoted by $|X_1| \cup |X_2|$ or $(|X_1|, |X_2|)$ is called the bidimension of the neutrosophic*

*bisemigroup bivector space (bilinear algebra) $V = V_1 \cup V_2$. Even if one of $X_1$ or $X_2$ is of infinite dimension then we say the bidimension of V is infinite, only when both $X_1$ and $X_2$ are of finite cardinality we say V is of finite bidimension over the bisemigroup $S = S_1 \cup S_2$.*

We will illustrate this by some simple examples.

***Example 5.84:*** Let $V = V_1 \cup V_2 = \{(0, 0), (I, I), (0, 0, 0), (I, I, I) (0\ 0\ 0\ 0\ 0), (I, 0, I, 0, I)\} \cup \{(5I, 5I), (0, 0), (10I, 10I), (15I, 15I), (20I, 20I), (5I, 0), (10I, 0), (20I, 0), (15I, 0)\}$ be a neutrosophic bisemigroup bivector space over the bisemigroup $S = S_1 \cup S_2 = Z_2I \cup Z_{25}I$.

Take $X = X_1 \cup X_2$ $\{(I, I), (I, I, I), (I, 0, I, 0, I)\} \cup (5I, 5I), (5I, 0)\} \subseteq V_1 \cup V_2$. X is a bigenerator of V and the bidimension of V is (3, 2)

***Example 5.85:*** Let

$$V = V_1 \cup V_2$$

$$= \left\{ \begin{pmatrix} a & b & c \\ a & b & c \\ a & b & c \end{pmatrix} \middle| a, b, c \in ZI \right\} \cup$$

$$\left\{ \begin{pmatrix} a & a & a & a & a \\ b & b & b & b & b \\ c & c & c & c & c \end{pmatrix} \middle| a, b, c \in Q^+I \right\}$$

be a neutrosophic bisemigroup bilinear algebra over the bisemigroup $S = S_1 \cup S_2 = ZI \cup Q^+I$.

Take $X = X_1 \cup X_2$

$$= \left\{ \begin{pmatrix} 1 & 0 & 0 \\ 1 & 0 & 0 \\ 1 & 0 & 0 \end{pmatrix}, \begin{pmatrix} 0 & 1 & 0 \\ 0 & 1 & 0 \\ 0 & 1 & 0 \end{pmatrix}, \begin{pmatrix} 0 & 0 & 1 \\ 0 & 0 & 1 \\ 0 & 0 & 1 \end{pmatrix} \right\}$$

$$\cup \left\{ \begin{pmatrix} 1 & 1 & 1 & 1 & 1 \\ 0 & 0 & 0 & 0 & 0 \\ 0 & 0 & 0 & 0 & 0 \end{pmatrix}, \begin{pmatrix} 0 & 0 & 0 & 0 & 0 \\ 1 & 1 & 1 & 1 & 1 \\ 0 & 0 & 0 & 0 & 0 \end{pmatrix}, \begin{pmatrix} 0 & 0 & 0 & 0 & 0 \\ 0 & 0 & 0 & 0 & 0 \\ 1 & 1 & 1 & 1 & 1 \end{pmatrix} \right\}$$

$\subseteq V_1 \cup V_2$, X bigenerates V over the bisemigroup $S = S_1 \cup S_2$. The bidimension of the neutrosophic bisemigroup bilinear algebra $V = V_1 \cup V_2$ over the bisemigroup $S = S_1 \cup S_2$ is (3, 3)

Now we proceed onto define the notion of neutrosophic bigroup bivector space over the bigroup.

**DEFINITION 5.43:** *Let $V = V_1 \cup V_2$ be such that $V_1 \neq V_2$, $V_1 \not\subseteq V_2$ and $V_2 \not\subseteq V_1$. If $V_1$ is a neutrosophic group vector space over the group $G_1$ and $V_2$ is a neutrosophic group vector space over the group $G_2$ and $G_1 \neq G_2$, $G_1 \not\subseteq G_2$ and $G_2 \not\subseteq G_1$ then we call $V = V_1 \cup V_2$ to be a neutrosophic bigroup bivector space over the bigroup $G = G_1 \cup G_2$.*

We will illustrate this by some examples.

*Example 5.86:* Let $V = V_1 \cup V_2 = \{(0, 0), (I, I), (1, 1), (0, 0, 0), (I, I, I), (I, 0, I, 0, I, 0 I) (0, 0, 0, 0, 0, 0, 0), (1, 1, 1, 1, 1), (I, I, I, I, I), (0, 0, 0, 0, 0)\} \cup$

$$\left\{ \begin{pmatrix} a & a \\ a & a \end{pmatrix}, \begin{pmatrix} a \\ a \end{pmatrix}, (a, a), (a, a, a) \right\}$$

such that $a \in Z_{12}I$ be a neutrosophic bigroup bivector space over the bigroup $G = G_1 \cup G_2 = \{N(Z_2) \cup Z_{12}I\}$.

*Example 5.87:* Let

$$V = V_1 \cup V_2 = \left\{ Z_2I \times Z_2I \times Z_2I, \begin{bmatrix} x \\ y \\ z \end{bmatrix} \mid x, y, z \in Z_2I \right\} \cup$$

$$\left\{ \begin{pmatrix} a & b \\ c & d \end{pmatrix}, \begin{pmatrix} a \\ b \end{pmatrix}, (c, d) \mid a, b, c, d \in ZI \right\}$$

be a neutrosophic bigroup bivector space over the bigroup $G = G_1 \cup G_2 = Z_2I \cup ZI$.

We now proceed onto define substructures in the neutrosophic bigroup bivector spaces.

**DEFINITION 5.44:** *Let $V = V_1 \cup V_2$ be a neutrosophic bigroup bivector space over the bigroup $G = G_1 \cup G_2$. $W = W_1 \cup W_2 \subseteq V_1 \cup V_2$ is defined to be a neutrosophic bigroup bivector subspace of $V = V_1 \cup V_2$ over the bigroup $G = G_1 \cup G_2$ if $W = W_1 \cup W_2$ is itself a neutrosophic bigroup bivector space over the bigroup $G = G_1 \cup G_2$.*

**DEFINITION 5.45:** *Let $V = V_1 \cup V_2$ be a neutrosophic bigroup bivector space over the bigroup $G = G_1 \cup G_2$. Let $W = W_1 \cup W_2 \subseteq V_1 \cup V_2$ and $H = H_1 \cup H_2 \subseteq G_1 \cup G_2$ be such that $H_1$ is a proper subset of $G_1$ and is a semigroup under the operations of $G_1$ and $H_2$ is also a proper subset of $G_2$ and is a semigroup of $G_2$. If $W = W_1 \cup W_2$ is a neutrosophic bisemigroup bivector space over the bisemigroup $H = H_1 \cup H_2$ then we define $W = W_1 \cup W_2 \subseteq V_1 \cup V_2$ ($W_1 \neq W_2$, $W_1 \not\subseteq W_2$, and $W_2 \not\subseteq W_1$ with $H_1 \neq H_2$, $H_1 \not\subseteq H_2$ and $H_2 \not\subseteq H_1$) to be a pseudo neutrosophic bisemigroup bivector subspace of $V$ over the bisemigroup $H = H_1 \cup H_2$ contained in $G = G_1 \cup G_2$.*

**DEFINITION 5.46:** *Let $V = V_1 \cup V_2$ be a neutrosophic bigroup bivector space over the bigroup $G = G_1 \cup G_2$. Let $W = W_1 \cup W_2 \subseteq V_1 \cup V_2$ and $P = P_1 \cup P_2 \subseteq G_1 \cup G_2$ is such that $P_i$ is a proper subgroup of $G_i$, $i=1, 2$ if $W$ is a neutrosophic bigroup bivector space over the bigroup $P = P_1 \cup P_2$ ($P_1 \neq P_2$, $P_1 \not\subseteq P_2$ and $P_2 \not\subseteq P_1$) then we call $W$ to be a neutrosophic subbigroup bivector subspace of $V$ over the subgroup $P$ of the bigroup $G$.*

We will illustrate these definitions by some examples.

*Example 5.88:* Let

$$V = V_1 \cup V_2$$

$$= \left\{ \begin{pmatrix} a & b \\ c & d \end{pmatrix}, \begin{bmatrix} a \\ a \\ a \end{bmatrix} \middle| a,b,c,d \in QI \right\} \cup$$

$$\left\{ (a,a,a,a,a), \begin{pmatrix} a \\ a \\ a \\ a \\ a \\ a \end{pmatrix} \middle| a \in Z_{12}I \right\}$$

be a neutrosophic bigroup bivector space over the bigroup $G = QI \cup Z_{12}I$. Take $P = P_1 \cup P_2 = ZI \cup \{0, 2I, 4I, 6I, 8I, 10I\}$ be a subbigroup of G.

$$W = W_1 \cup W_2$$

$$= \left\{ \begin{pmatrix} a & b \\ c & d \end{pmatrix} \middle| a,b,c,d \in QI \right\} \cup \left\{ (a,a,a,a,a) \middle| a \in Z_{12}I \right\}$$

$$\subseteq V_1 \cup V_2,$$

W is a neutrosophic subbigroup bivector subspace of V over the subbigroup $P = P_1 \cup P_2 \subseteq G_1 \cup G_2 = QI \cup Z_{12}I$.

*Example 5.89:* Let

$$V = V_1 \cup V_2$$

$$= \left\{ \begin{pmatrix} a & b & c \\ d & e & f \end{pmatrix}, \begin{bmatrix} a \\ a \\ a \\ a \\ a \end{bmatrix} \middle| a,b,c,d,e,f \in QI \right\} \cup$$

$$\left\{ \begin{pmatrix} a & b \\ c & d \\ e & f \end{pmatrix}, (a,a,a,a,a,a) \middle| a,b,c,d,e,f \in Z_{20}I \right\}$$

be a neutrosophic bigroup bivector space over the bigroup $G = QI \cup Z_{20}I = G_1 \cup G_2$.
Take

$$W = W_1 \cup W_2 = \left\{ \begin{bmatrix} a \\ a \\ a \\ a \\ a \end{bmatrix} \middle| a \in QI \right\}$$

$\cup$ {(a, a, a, a, a, a)| a $\in$ $Z_{20}I$} $\subseteq$ $V_1 \cup V_2$ and $P = P_1 \cup P_2 = ZI \cup \{0, 10I\} \subseteq G_1 \cup G_2 = QI \cup Z_{20}I$. W is a neutrosophic subbigroup subbivector space over the subbigroup $P = P_1 \cup P_2 \subseteq G_1 \cup G_2$.

*Example 5.90:* Let $V = V_1 \cup V_2$

$$= \left\{ \begin{pmatrix} a & b \\ c & d \end{pmatrix}, \begin{pmatrix} a \\ a \\ a \end{pmatrix}, [a,a,a,a,a] \middle| a,b,c,d \in QI \right\} \cup$$

$$\left\{ \begin{pmatrix} a_1 & a_2 & a_3 & a_4 \\ a_5 & 0 & 0 & a_6 \end{pmatrix}, \begin{pmatrix} 0 & 0 & 0 & a_4 \\ 0 & a_5 & a_6 & 0 \end{pmatrix} \middle| a_i \in Z_{16}I; 1 \le i \le 6 \right\}$$

be a neutrosophic bigroup bivector space over the bigroup $G = G_1 \cup G_2 = \{QI\} \cup Z_{16}I$. Take $W = W_1 \cup W_2$

$$= \left\{ \begin{pmatrix} a & b \\ c & d \end{pmatrix}, \begin{bmatrix} a \\ a \\ a \end{bmatrix} \middle| a,b,c,d \in QI \right\} \cup$$

$$\left\{ \begin{pmatrix} a_1 & a_2 & a_3 & a_4 \\ a_5 & 0 & 0 & a_6 \end{pmatrix} \middle| a_i \in Z_{16}I; 1 \leq i \leq 6 \right\}$$

$$\subseteq V_1 \cup V_2 = V,$$

W is a neutrosophic bigroup subbivector space over the bigroup $G = G_1 \cup G_2$.

*Example 5.91:* Let
$$V = V_1 \cup V_2 =$$

$$\left\{ (a, a), \begin{bmatrix} a \\ a \\ a \\ a \\ a \\ a \end{bmatrix} \middle| a \in Z_{25}I \right\} \cup \left\{ \begin{pmatrix} a & a & a \\ a & a & a \end{pmatrix}, \begin{pmatrix} a \\ a \\ a \end{pmatrix} \middle| a \in QI \right\}$$

be a neutrosophic bigroup bivector space over the bigroup $G = G_1 \cup G_2 = Z_{25}I \cup QI$.
Take

$$W = W_1 \cup W_2$$

$$= \left\{ \begin{bmatrix} a \\ a \\ a \\ a \\ a \\ a \end{bmatrix} \middle| a \in Z_{25}I \right\} \cup \left\{ \begin{pmatrix} a & a & a \\ a & a & a \end{pmatrix} \middle| a \in QI \right\} \subseteq V_1 \cup V_2,$$

W is a neutrosophic bigroup subbivector space over the bigroup $G = G_1 \cup G_2$.

*Example 5.92:* Let
$$V = V_1 \cup V_2 = \left\{ \begin{pmatrix} a & b & c \\ 0 & d & e \\ 0 & 0 & f \end{pmatrix}, \begin{bmatrix} a & b \\ c & d \end{bmatrix} \middle| a, b, c, d, e, f \in QI \right\}$$

$$\cup \left\{ \begin{bmatrix} a \\ a \\ a \end{bmatrix}, (x, y, z) \,\middle|\, a, x, y, z \in ZI \right\}$$

be a neutrosophic bigroup bivector space over the bigroup $G = G_1 \cup G_2 = 3ZI \cup 8ZI$.

Choose
$$W = W_1 \cup W_2$$

$$= \left\{ \begin{pmatrix} a & b & c \\ 0 & d & e \\ 0 & 0 & f \end{pmatrix} \,\middle|\, a, b, c, d, e, f \in QI \right\}$$

$\cup \{(x, y, z) \mid (x, y, z) \in Z_{36}I\} \subseteq V_1 \cup V_2$;

W is a pseudo neutrosophic subbisemigroup bivector subspace over the subbisemigroup $H = H_1 \cup H_2 = \{3Z^+I \cup 48Z^+I\} \subseteq G_1 \cup G_2$.

*Example 5.93:* Let
$$V = V_1 \cup V_2$$

$$= \left\{ \begin{pmatrix} a & b & c \\ d & e & f \end{pmatrix}, \begin{pmatrix} a \\ b \\ c \\ d \end{pmatrix} \,\middle|\, a, b, c, d, e, f \in ZI \right\} \cup$$

$$\left\{ \begin{pmatrix} a & b \\ c & d \\ e & f \end{pmatrix}, (a, b, c, d, e) \,\middle|\, a, b, c, d, e, f \in QI \right\}$$

be a neutrosophic bigroup bivector space over the bigroup $G = G_1 \cup G_2 = 5ZI \cup 7ZI$.

Consider
$$W = W_1 \cup W_2$$

$$= \left\{ \begin{pmatrix} a \\ b \\ c \\ d \end{pmatrix} \middle| a,b,c,d \in ZI \right\} \cup \{[a, b, c, d, e] \mid a, b, c, d, e \in QI\}$$

$$\subseteq V_1 \cup V_2;$$

W is a pseudo neutrosophic subbisemigroup bivector subspace of V over the subsemigroup

$$5Z^+I \cup 7Z^+I \subseteq 5ZI \cup 7ZI.$$

Now we proceed onto define neutrosophic pseudo bigroup bivector space.

**DEFINITION 5.47:** *Let $V = V_1 \cup V_2$ where $V_1$ is a neutrosophic semigroup vector space over the semigroup $S_1$ and $V_2$ is a neutrosophic group vector space over the group $G_1$ ($S_1 \not\subseteq G_1$, $G_1 \not\subseteq S_1$, $S_1 \neq G_1$, $V_1 \neq V_2$, $V_1 \not\subseteq V_2$ and $V_2 \not\subseteq V_1$). We call $V = V_1 \cup V_2$ to be a neutrosophic pseudo bigroup bivector space over the pseudo bigroup $G = S_1 \cup G_1$.*

We will illustrate this situation by some examples.

*Example 5.94:* Let
$$V = V_1 \cup V_2$$

$$= \left\{ \begin{bmatrix} a & b \\ c & d \end{bmatrix} \middle| a,b,c,d \in Z^+I \cup \{0\} \right\} \cup \left\{ \begin{bmatrix} a & a & a \\ a & a & a \end{bmatrix} \middle| a \in QI \right\}$$

is a pseudo neutrosophic bigroup bivector space over the pseudo bigroup
$$G = (Z^+I \cup \{0\}) \cup QI = G_1 \cup G_2.$$

*Example 5.95:* Let
$$V = V_1 \cup V_2 =$$

$$\left\{ \begin{pmatrix} a \\ a \\ a \end{pmatrix}, (a, a, a, a, a) \;\middle|\; a \in 3Z^+I \cup \{0\} \right\} \cup$$

$$\left\{ \begin{pmatrix} a & b \\ c & d \end{pmatrix}, \begin{pmatrix} a & b & c \\ d & e & f \end{pmatrix} \;\middle|\; a,b,c,d,e,f \in QI \right\}$$

be a pseudo neutrosophic bigroup bivector space over the pseudo bigroup $G = S_1 \cup G_1 = \{3Z^+I \cup \{0\}\} \cup QI$.

**Chapter Six**

# NEUTROSOPHIC FUZZY GROUP BILINEAR ALGEBRA

In this chapter we proceed onto define the neutrosophic fuzzy analogue of the definition given in chapter four of this book. Through out this book $N([0, 1]) = \{a + b_i \mid a, b \in [0, I]\}$ is the fuzzy neutrosophic set or neutrosophic fuzzy set.

**DEFINITION 6.1:** *Let $V = V_1 \cup V_2$ be a neutrosophic set bivector space over the set S. We say the neutrosophic set bivector space $V = V_1 \cup V_2$ with the bimap $\eta = \eta_1 \cup \eta_2$ is a neutrosophic fuzzy set bivector space or neutrosophic set fuzzy bivector space if $\eta : \eta_1 \cup \eta_2 : V = V_1 \cup V_2 \to [0, 1]$ and*
$$\eta(r_1 a_1 \cup r_2 a_2) = (\eta_1 \cup \eta_2)(r_1 a_1 \cup r_2 a_2)$$
$$= \eta_1(r_1 a_1) \cup \eta_2(r_2 a_2) \geq \eta_1(a_1) \cup \eta_2(a_2)$$
*for all $a_i \in V_i$; $i = 1, 2$ and $r_1, r_2 \in S$.*

*We call* $V_\eta = V_\eta = (V_1 \cup V_2)_{\eta_1 \cup \eta_2}$ *or* $\eta V = (\eta_1 \cup \eta_2)(V_1 \cup V_2)$ *to be the neutrosophic fuzzy set bivector space over the set S.*

We will illustrate this by some examples.

*Example 6.1:* Let
$$V = V_1 \cup V_2$$

$$= \left\{ \begin{pmatrix} a & a \\ a & 0 \end{pmatrix}, \begin{pmatrix} a & a & a \\ 0 & a & 0 \end{pmatrix} \middle| a \in QI \right\} \cup \left\{ (a, a, a, a), \begin{bmatrix} a & a \\ a & a \\ 0 & a \\ a & 0 \end{bmatrix} \middle| a \in ZI \right\}$$

be a neutrosophic set bivector space over the set $S = \{0, 2, 2^2, \ldots, 2^n \mid n \in N\}$.
Define $\eta : V \to N([0, 1])$
i.e., $\eta = \eta_1 \cup \eta_2 : V_1 \cup V_2 \to N([0, 1])$ as

$$\eta_1 \begin{pmatrix} a & a \\ a & 0 \end{pmatrix} = \begin{cases} I + \dfrac{1}{2} & \text{if } a \neq 0 \\ 1 & \text{if } a = 0 \end{cases}$$

$$\eta_1 \begin{pmatrix} a & a & a \\ 0 & a & 0 \end{pmatrix} = \begin{cases} I + \dfrac{1}{4} & \text{if } a \neq 0 \\ 1 & \text{if } a = 0 \end{cases}$$

and

$$\eta_2 (a, a, a, a) = \begin{cases} I + \dfrac{1}{5} & \text{if } a \neq 0 \\ 1 & \text{if } a = 0 \end{cases}$$

$$\eta_2 \begin{bmatrix} a & a \\ a & a \\ 0 & a \\ a & 0 \end{bmatrix} = \begin{cases} I + \dfrac{1}{8} & \text{if } a \neq 0 \\ 1 & \text{if } a = 0 \end{cases}.$$

$V_\eta = (V_1 \cup V_2)_{\eta_1 \cup \eta_2}$ is a neutrosophic fuzzy set bivector space.

*Example 6.2:* Let
$$V = V_1 \cup V_2$$

$$= \left\{ \begin{pmatrix} a & a \\ b & b \end{pmatrix}, \begin{pmatrix} a & b \\ a & b \\ a & b \end{pmatrix} \middle| a, b \in Z_{12}I \right\} \cup$$

$$\left\{ \begin{pmatrix} a & a \\ b & b \\ c & c \\ d & d \end{pmatrix}, \begin{pmatrix} a \\ b \\ c \\ d \end{pmatrix}, (a,b,c,d) \middle| a,b,c,d \in Z_{17}I \right\}$$

be a neutrosophic set bivector space over the set $S = \{0, 1\}$.
Define $\eta = \eta_1 \cup \eta_2 : V = V_1 \cup V_2 \to N([0, 1])$ where $\eta_1 : V_1 \to N([0, 1])$ and $\eta_2 : V_2 \to N([0, 1])$.
Defined by

$$\eta_1 \begin{pmatrix} a & a \\ b & b \end{pmatrix} = \begin{cases} I + \dfrac{1}{4} & \text{if } a \neq 0 \text{ or } b \neq 0 \\ 1 & \text{if } a = b = 0 \end{cases},$$

$$\eta_1 \begin{pmatrix} a & b \\ a & b \\ a & b \end{pmatrix} = \begin{cases} I + \dfrac{1}{6} & \text{if } a \neq 0 \text{ or } b \neq 0 \\ 1 & \text{if } a = b = 0 \end{cases},$$

$$\eta_2 \begin{pmatrix} a & a \\ b & b \\ c & c \\ d & d \end{pmatrix} = \begin{cases} I + \dfrac{1}{8} & \text{if } a \neq 0 \text{ or } b \neq 0 \text{ or } c \neq 0 \text{ or } d \neq 0 \\ 1 & \text{if } a = b = c = d = 0 \end{cases},$$

$$\eta_2 \begin{pmatrix} a \\ b \\ c \\ d \end{pmatrix} = \begin{cases} I + \dfrac{1}{4} & \text{if } \begin{bmatrix} a \\ b \\ c \\ d \end{bmatrix} \neq \begin{bmatrix} 0 \\ 0 \\ 0 \\ 0 \end{bmatrix} \\ 1 & \text{if } a = b = c = d = 0 \end{cases}$$

and

$$\eta_2 (a, b, c, d) = \begin{cases} I + \dfrac{1}{18} & \text{if } a \neq 0 \text{ or } b \neq 0 \text{ or } c \neq 0 \text{ or } d \neq 0 \\ 1 & \text{if } a = b = c = d = 0 \end{cases}.$$

$V_\eta = (V_1 \cup V_2)_{\eta_1 \cup \eta_2}$ is a neutrosophic set fuzzy bivector space.

Now we proceed onto define the notion of neutrosophic set fuzzy bilinear algebra.

**DEFINITION 6.2:** *Let $V = V_1 \cup V_2$ be a neutrosophic set bilinear algebra over the set S. A neutrosophic set fuzzy bilinear algebra or neutrosophic fuzzy set bilinear algebra $(V, \eta) = (V_1 \cup V_2, \eta_1 \cup \eta_2)$ or $(\eta_1 \cup \eta_2) (V_1 \cup V_2)$ is a bimap $\eta = \eta_1 \cup \eta_2 : V_1 \cup V_2 \to N[(0, 1)]$ such that*
$$\eta_i (a_i + b_i) \geq \min( \eta_i (a_i), \eta_i (b_i))$$
$$\eta_i (r_i a_i) \geq \eta_i (a_i)$$
*for $a_i, b_i \in V_i$, $r_i \in S$. $i = 1, 2$.*

We illustrate this by some examples.

*Example 6.3:* Let
$$V = V_1 \cup V_2$$
$$= \{ZI \times ZI \times ZI\} \cup \left\{ \begin{pmatrix} a & a \\ a & a \end{pmatrix} \middle| a \in QI \right\}$$

be a neutrosophic set bilinear algebra over the set $S = Z^+I$. Define $\eta = \eta_1 \cup \eta_2 : V = V_1 \cup V_2 \to N([0, 1])$ where $\eta_1 : V_1 \to N([0, 1])$ and $\eta_2 : V_2 \to N([0, 1])$ such that

$$\eta_1(a, b, c) = \begin{cases} I + \dfrac{1}{3} & \text{if } (a,b,c) \neq (0,0,0) \\ \dfrac{1}{3} & \text{if } (a,b,c) = (0,0,0) \end{cases}$$

and

$$\eta_2 \begin{pmatrix} a & a \\ a & a \end{pmatrix} = \begin{cases} I + \dfrac{1}{4} & \text{if } \begin{pmatrix} a & a \\ a & a \end{pmatrix} \neq \begin{pmatrix} 0 & 0 \\ 0 & 0 \end{pmatrix} \\ \dfrac{1}{4} & \text{if } \begin{pmatrix} a & a \\ a & a \end{pmatrix} = \begin{pmatrix} 0 & 0 \\ 0 & 0 \end{pmatrix} \end{cases}$$

$V_\eta = (V_1 \cup V_2)_{\eta_1 \cup \eta_2}$ is a neutrosophic set fuzzy bilinear algebra.

*Example 6.4:* Let

$$V = V_1 \cup V_2$$

$$= \left\{ \begin{pmatrix} a & b & c & d \\ e & f & g & h \end{pmatrix} \middle| a,b,c,d,e,f \in Q^+I \cup \{0\} \right\}$$

$$\cup \left\{ \begin{pmatrix} a \\ b \\ c \\ d \\ e \\ f \end{pmatrix} \middle| a,b,c,d,e,f \in Q^+I \cup \{0\} \right\}$$

be a neutrosophic set bilinear algebra over the set $S = Z^+$.
Define $\eta: V \to N([0, 1])$
$$\eta = \eta_1 \cup \eta_2 : V = V_1 \cup V_2 \to N([0, 1])$$
by
$$\eta_1 : V_1 \to N([0, 1])$$
and
$$\eta_2 : V_2 \to N([0, 1])$$
as

$$\eta_1 \begin{pmatrix} a & b & c & d \\ e & f & g & h \end{pmatrix} = \begin{cases} I + \dfrac{1}{8} & \text{if } \begin{pmatrix} a & b & c & d \\ e & f & g & h \end{pmatrix} \neq \begin{pmatrix} 0 & 0 & 0 & 0 \\ 0 & 0 & 0 & 0 \end{pmatrix} \\ \dfrac{1}{8} & \text{if } \begin{pmatrix} a & b & c & d \\ e & f & g & h \end{pmatrix} = \begin{pmatrix} 0 & 0 & 0 & 0 \\ 0 & 0 & 0 & 0 \end{pmatrix} \end{cases}$$

$$\eta_2 \begin{pmatrix} a \\ b \\ c \\ d \\ e \\ f \end{pmatrix} = \begin{cases} I + \dfrac{1}{9} & \text{if } \begin{pmatrix} a \\ b \\ c \\ d \\ e \\ f \end{pmatrix} \neq \begin{pmatrix} 0 \\ 0 \\ 0 \\ 0 \\ 0 \\ 0 \end{pmatrix} \\ \dfrac{1}{9} & \text{if } \begin{pmatrix} a \\ b \\ c \\ d \\ e \\ f \end{pmatrix} = \begin{pmatrix} 0 \\ 0 \\ 0 \\ 0 \\ 0 \\ 0 \end{pmatrix} \end{cases}$$

$V_\eta = (V_1 \cup V_2)_{\eta_1 \cup \eta_2}$ is a neutrosophic set fuzzy bilinear algebra.

Now we proceed onto define the notion of neutrosophic set fuzzy bivector subspace.

**DEFINITION 6.3:** *Let $V = V_1 \cup V_2$ be a neutrosophic setbivector space over the set S and $W = W_1 \cup W_2 \subseteq V_1 \cup V_2$ be the neutrosophic set bivector subspace of V over the set S. Define $\eta = \eta_1 \cup \eta_2: W = W_1 \cup W_2 \to N[(0, 1)]$ then $W_\eta = (W_1 \cup W_2)_{\eta_1 \cup \eta_2}$ is called the neutrosophic set fuzzy bivector subspace of V.*

We will illustrate this situation by some examples.

*Example 6.5:* Let $V = V_1 \cup V_2$

$$= \left\{ (Z^+I \times Z^+I \times Z^+I \times Z^+I), \begin{pmatrix} a \\ b \\ c \\ d \end{pmatrix} \middle| a, b, c, d \in Z^+I \right\}$$

$$\cup \left\{ \begin{pmatrix} a & b \\ 0 & c \end{pmatrix}, \begin{pmatrix} a & a & a & a \\ a & a & a & a \end{pmatrix} \middle| a, b, c \in Z^+I \right\}$$

be a neutrosophic set bivector space over the set $S = Z^+I$.
Let
$$W = W_1 \cup W_2$$

$$= \left\{ \begin{pmatrix} a \\ b \\ c \\ d \end{pmatrix} \middle| a, b, c, d \in Z^+I \right\} \cup \left\{ \begin{pmatrix} a & b \\ 0 & c \end{pmatrix} \middle| a, b, c \in Z^+I \right\}$$

$$\subseteq V_1 \cup V_2$$

be a neutrosophic set bivector subspace of V over the set S.
Define
$$\eta = \eta_1 \cup \eta_2 : W_1 \cup W_2 \to N([0, 1])$$
by
$$\eta_1 : W_1 \to N([0, 1])$$
such that
$$\eta_1 \begin{pmatrix} a \\ b \\ c \\ d \end{pmatrix} = I + \frac{1}{3}$$

$$\eta_2 \begin{pmatrix} a & b \\ 0 & c \end{pmatrix} = I + \frac{1}{4},$$

then $W_\eta = (W_1 \cup W_2)_{\eta_1 \cup \eta_2}$ is a neutrosophic set fuzzy bivector subspace of V.

*Example 6.6:* Let

$$V = V_1 \cup V_2 = \{Z^+I [x], (a, a, a, a) \mid a \in Z^+I\} \cup$$

$$\left\{ \begin{bmatrix} a & b \\ c & d \\ e & f \\ g & h \end{bmatrix}, (a, a, a, a, a), \right.$$

$$\begin{bmatrix} a & a & a & a & a \\ b & b & b & b & b \end{bmatrix} \mid a, b, c, d, e, f, g, h \in Q^+I\}$$

be a neutrosophic set bivector space over the set $S = 3Z^+I$. Take

$$W = W_1 \cup W_2$$
$$= \{(a, a, a, a) \mid a \in Z^+I\} \cup \left\{ \begin{bmatrix} a & a & a & a & a \\ b & b & b & b & b \end{bmatrix} \mid a, b \in Z^+I \right\}$$
$$\subseteq V_1 \cup V_2$$

to be a neutrosophic set bivector subspace of V over the set S. Define $\eta = \eta_1 \cup \eta_2: W_1 \cup W_2 \to N([0, 1])$ where $\eta_1 : W_1 \to N([0, 1])$ and $\eta_2 : W_2 \to N([0, 1])$;

$$\eta_1 (a, a, a, a) = I + \frac{1}{20}$$

and

$$\eta_2 \begin{bmatrix} a & a & a & a & a \\ b & b & b & b & b \end{bmatrix} = I + \frac{1}{8}$$

$W_\eta = (W_1 \cup W_2)_{\eta_1 \cup \eta_2}$ is a neutrosophic set fuzzy bivector subspace of V.

Now we proceed onto define the notion of neutrosophic fuzzy semigroup bivector space.

**DEFINITION 6.4:** *A neutrosophic semigroup fuzzy bivector space or a neutrosophic fuzzy semigroup bivector space $V\eta$ (or $\eta V$ or $(V_1 \cup V_2)_{\eta_1 \cup \eta_2}$) is a neutrosophic semigroup bivector space $V = V_1 \cup V_2$ over the semigroup S, with a bimap $\eta = \eta_1 \cup \eta_2 : V = V_1 \cup V_2 \to N([0, 1])$ satisfying the following conditions;*

$$\eta(ra) = \eta(r_1 a_1 \cup r_2 a_2)$$
$$= (\eta_1 \cup \eta_2)(r_1 a_1 \cup r_2 a_2)$$
$$= \eta_1(r_1 a_1) \cup \eta_2(r_2 a_2) \geq \eta_1(a_1) \cup \eta_2(a_2);$$
$$\text{i.e., } \eta_1(r_1, a_1) \geq \eta_1(a_1) \text{ and } \eta_2(r_2 a_2) \geq \eta_2(a_2)$$

*for all $a_1 \in V_1$, $a_2 \in V_2$ and $r_1, r_2 \in S$.*

We will illustrate this situation by some examples.

*Example 6.7:* Let

$$V = V_1 \cup V_2$$
$$= \{ZI \times ZI \times ZI \times ZI, Z^+I \times Z^+I \times Z^+I\} \cup$$

$$\left\{ \begin{pmatrix} a & d \\ b & e \\ c & f \end{pmatrix}, \begin{bmatrix} a & a & a & a & a \\ b & b & b & b & b \end{bmatrix} \middle| a,b,c,d,e,f \in Z^+I \right\}$$

be a neutrosophic semigroup bivector space over the semigroup $S = Z^+I$. Define $\eta = \eta_1 \cup \eta_2 : V = V_1 \cup V_2 \to N[(0, 1)]$ where $\eta_1 : V_1 \to N([0, 1])$ and $\eta_2 : V_2 \to N([0, 1])$ with

$$\eta_1(a, b, c, d) = \begin{cases} I + \dfrac{1}{4} & \text{if } (a\ b\ c\ d) \\ 1 & \text{if } a = b = c = d = 0 \end{cases}$$

$$\eta_1(a, b, c) = I + \frac{1}{3}$$

and

$$\eta_2 \begin{pmatrix} a & d \\ b & e \\ c & f \end{pmatrix} = 0.5I + \frac{1}{8}$$

$$\eta_2 \begin{bmatrix} a & a & a & a & a \\ b & b & b & b & b \end{bmatrix} = 0.8I + \frac{1}{6}.$$

$V_\eta = (V_1 \cup V_2)_{\eta_1 \cup \eta_2}$ is a neutrosophic semigroup fuzzy bivector space.

*Example 6.8:* Let

$$V = V_1 \cup V_2$$

$$= \left\{ \begin{pmatrix} a & b & c \\ 0 & d & e \\ 0 & 0 & f \end{pmatrix} \middle| a,b,c,d,e,f \in Z_5I \right\} \cup \{Z_5I \times Z_5I \times Z_5I \times Z_5I\}$$

be a neutrosophic semigroup bivector space over the semigroup $S = Z_6$. Define $\eta : V \to N([0, 1])$ that is $\eta = \eta_1 \cup \eta_2 : V_1 \cup V_2 \to N[(0, 1)]$ where

$$\eta_1 : V_1 \to N([0, 1])$$

and

$$\eta_2 : V_2 \to N([0, 1])$$

by

$$\eta_1 \begin{pmatrix} a & b & c \\ 0 & d & e \\ 0 & 0 & f \end{pmatrix} = \begin{cases} 0.9I + \dfrac{1}{5} & \text{if } \begin{pmatrix} a & b & c \\ 0 & d & e \\ 0 & 0 & f \end{pmatrix} \neq \begin{pmatrix} 0 & 0 & 0 \\ 0 & 0 & 0 \\ 0 & 0 & 0 \end{pmatrix} \\ 1 & \text{if } \begin{pmatrix} a & b & c \\ 0 & d & e \\ 0 & 0 & f \end{pmatrix} = \begin{pmatrix} 0 & 0 & 0 \\ 0 & 0 & 0 \\ 0 & 0 & 0 \end{pmatrix} \end{cases}$$

and

$$\eta_2(a, b, c, d) = \begin{cases} I + \dfrac{1}{5} & \text{if } (a,b,c,d) \neq (0,0,0,0) \\ 1 & \text{if } (a,b,c,d) = (0,0,0,0) \end{cases}.$$

$V_\eta = (V_1 \cup V_2)_{\eta_1 \cup \eta_2}$ is a neutrosophic semigroup fuzzy bivector space.

Now we proceed onto define the notion of neutrosophic semigroup fuzzy bivector subspace. This is defined in two ways.

**DEFINITION 6.5:** *Let $V = V_1 \cup V_2$ be a neutrosophic semigroup bivector space over the semigroup S. $V_\eta = (V_1 \cup V_2)_{\eta_1 \cup \eta_2}$ be a neutrosophic semigroup fuzzy bivector space. Suppose $W = W_1 \cup W_2 \subseteq V_1 \cup V_2$ be a neutrosophic semigroup bivector subspace of V, then we define $W_{\bar{\eta}} = (W_1 \cup W_2)_{\bar{\eta}_1 \cup \bar{\eta}_2}$ to be the neutrosophic semigroup fuzzy bivector subspace of $V_\eta$, where $\bar{\eta}: W \to N([0, 1])$ is such that $\bar{\eta} = \bar{\eta}_1 \cup \bar{\eta}_2$ is the restriction of $\eta$ to W i.e $\bar{\eta}_i : W_i \to N([0, 1])$ where $\eta_i : V_i \to N([0, 1])$. $i = 1, 2$.*

**DEFINITION 6.6:** *Let $V = V_1 \cup V_2$ be a neutrosophic semigroup bivector space defined over the semigroup S. Let $W = W_1 \cup W_2 \subseteq V_1 \cup V_2$ be a neutrosophic semigroup bivector subspace of V over the semigroup S. Let $= \eta_1 \cup \eta_2 : W = W_1 \cup W_2 \to N([0, 1])$ be a bimap such that $W\eta$ is a neutrosophic semigroup fuzzy bivector space; then we call $W_\eta = (W_1 \cup W_2)_{\eta_1 \cup \eta_2}$ to be a neutrosophic semigroup fuzzy bivector subspace of V.*

It is important to note the following. In general the two definitions are not equivalent for $\eta : W = W_1 \cup W_2 \to N([0, 1])$, $\eta$ may not be defined on $V \setminus W = V_1 \setminus W_1 \cup V_2 \setminus W_2$, where as $\bar{\eta} : W \to N([0, 1])$ is only a restriction bimap.

Here we give examples of both the definitions.

***Example 6.9:*** Let

$$V = V_1 \cup V_2$$

$$= \left\{ (a, b, c, d, e), \begin{bmatrix} a \\ a \\ a \end{bmatrix} \mid a, b, c, d, e \in Z^+I \right\} \cup$$

$$\left\{ Z^+I \times Z^+I \times Z^+I, \left\{ \begin{pmatrix} a & b \\ c & d \end{pmatrix} \mid a, b, c, d \in Z^+I \right\} \right\}$$

be a neutrosophic semigroup bivector space over the semigroup $S = Z^+I$. Let $\eta = \eta_1 \cup \eta_2 : V = V_1 \cup V_2 \to N([0, 1])$ where $\eta_1 : V_1 \to N([0, 1])$ and $\eta_2 : V_2 \to N([0, 1])$ are defined by

$$\eta_1(a, b, c, d, e) = I + \frac{1}{7}$$

$$\eta_1 \left( \begin{bmatrix} a \\ a \\ a \end{bmatrix} \right) = I + \frac{1}{3}, \quad \eta_2(a, a, a) = I + \frac{1}{26}$$

$$\eta_2 \begin{pmatrix} a & b \\ c & d \end{pmatrix} = I + \frac{1}{4}$$

$V_\eta = (V_1 \cup V_2)_{\eta_1 \cup \eta_2}$ is a neutrosophic semigroup fuzzy bivector space. Let

$$W = W_1 \cup W_2$$

$$= \left\{ \begin{pmatrix} a \\ a \\ a \end{pmatrix} \mid a \in Z^+I \right\} \cup \left\{ \begin{pmatrix} a & b \\ c & d \end{pmatrix} \mid a, b, c, d \in Z^+I \right\}$$

$\subseteq V_1 \cup V_2$.

Define $\overline{\eta}: W = W_1 \cup W_2 \to N([0, 1])$ as follows, $\overline{\eta} = \overline{\eta}_1 \cup \overline{\eta}_2 :$
$W_1 \cup W_2 \to N([0, 1])$

$$\overline{\eta}_1 : W_1 \to N([0, 1])$$

and

$$\overline{\eta}_2 : W_2 \to N([0, 1])$$

is such that

$$\overline{\eta}_1 \begin{pmatrix} a \\ a \\ a \end{pmatrix} = I + \frac{1}{3}$$

and

$$\overline{\eta}_2 \begin{pmatrix} a & b \\ c & d \end{pmatrix} = I + \frac{1}{4}$$

$W_\eta = (W_1 \cup W_2)_{\eta_1 \cup \eta_2}$ is the neutrosophic semigroup fuzzy bivector space of $V_\eta$. $\overline{\eta}$ is clearly the restriction bimap of $\eta$ on W.

*Example 6.10:* Let
$$V = V_1 \cup V_2$$

$$= \left\{ \begin{pmatrix} a & a & a & a \\ b & b & b & b \end{pmatrix}, \begin{pmatrix} a & b \\ a & b \\ a & b \\ a & b \\ a & b \\ a & b \\ a & b \end{pmatrix} \middle| a, b \in Z^+ I \right\} \cup$$

$$\left\{ (a, a, a, a, a), \begin{bmatrix} a & b & c \\ d & e & f \end{bmatrix} \middle| a, b, c, d, e, f \in Q^+ I \right\}$$

be a neutrosophic semigroup bivector space over the semigroup $S = Z^+ I$.

Define $\eta = \eta_1 \cup \eta_2 = V = V_1 \cup V_2 \to N([0, 1])$ where $\eta_1 : V_1 \to N([0, 1])$ and $\eta_2 : V_2 \to N([0, 1])$ are such that

$$\eta_1 \begin{pmatrix} a & a & a & a \\ b & b & b & b \end{pmatrix} = I + \frac{1}{8}$$

$$\eta_1 \begin{pmatrix} a & b \\ a & b \\ a & b \\ a & b \\ a & b \\ a & b \end{pmatrix} = I + \frac{1}{12}$$

and

$$\eta_2 (a, a, a, a, a) = I + \frac{1}{5}$$

$$\eta_2 \begin{bmatrix} a & b & c \\ d & e & f \end{bmatrix} = I + \frac{1}{6}$$

$V_\eta = (V_1 \cup V_2)_{\eta_1 \cup \eta_2}$ is a neutrosophic semigroup fuzzy bivector space. Let

$$W = W_1 \cup W_2$$

$$= \left\{ \begin{pmatrix} a & a & a & a \\ a & a & a & a \end{pmatrix} \middle| a \in Z^+I \right\} \cup \{(a, a, a, a, a) \mid a \in Q^+I\}$$

$$\subseteq V_1 \cup V_2$$

be a neutrosophic semigroup bivector subspace of V over the semigroup S.
Define
$$\overline{\eta} = \overline{\eta}_1 \cup \overline{\eta}_2 : W_1 \cup W_2 \to N([0, 1])$$
where
$$\overline{\eta}_1 : W_1 \to N([0, 1])$$

and
$$\overline{\eta}_2 : W_2 \to N([0, 1])$$
are such that
$$\overline{\eta}_1 \begin{pmatrix} a & a & a & a \\ a & a & a & a \end{pmatrix} = I + \frac{1}{8}$$
and
$$\overline{\eta}_2 \ (a, a, a, a, a) = I + \frac{1}{5}.$$

$W_\eta = W_{\overline{\eta}} = (W_1 \cup W_2)_{\overline{\eta}_1 \cup \overline{\eta}_2}$ is a neutrosophic semigroup fuzzy bivector space of $V_\eta$.

*Example 6.11:* Let
$$V = V_1 \cup V_2$$
$$= \{(a, b, c, d, e)\,, \begin{pmatrix} a \\ a \\ a \end{pmatrix} \mid a, b, c, d, e, \in Z^+I\}$$

$$\cup \ \{Z^+I \times Z^+I \times Z^+I, \begin{pmatrix} a & b \\ c & d \end{pmatrix} \mid a, b, c, d, \in Z^+I\}$$

be a neutrosophic semigroup bivector space over the semigroup $Z^+I$. Let
$$W = W_1 \cup W_2$$
$$= \left\{ \begin{pmatrix} a \\ a \\ a \end{pmatrix} \middle| a \in Z^+I \right\} \cup \left\{ \begin{pmatrix} a & b \\ c & d \end{pmatrix} \middle| a, b, c, d \in Z^+I \right\} \subseteq V_1 \cup V_2$$

be a neutrosophic semigroup bivector subspace of V over the semigroup S. Define $\eta : W \to N([0, 1])$ i.e.,
$$\eta = \eta_1 \cup \eta_2 = W_1 \cup W_2 \to N([0, 1])$$
where

$$\eta_1 : W_1 \to N([0, 1])$$
and
$$\eta_2 : W_2 \to N([0, 1])$$
are such that
$$\eta_1 \begin{pmatrix} a \\ a \\ a \end{pmatrix} = I + \frac{1}{100}$$
and
$$\eta_2 \begin{pmatrix} a & b \\ c & d \end{pmatrix} = I + \frac{1}{200}$$

$W_\eta = (W_1 \cup W_2)_{\eta_1 \cup \eta_2}$ is a neutrosophic semigroup fuzzy bivector space of V. Clearly $W_{\overline{\eta}}$ and $W_\eta$ are distinct given in examples 6.9 and 6.11 respectively.

*Example 6.12:* Let
$$V = V_1 \cup V_2$$

$$= \left\{ \begin{pmatrix} a & a & a & a \\ b & b & b & b \end{pmatrix}, \begin{pmatrix} a & b \\ a & b \\ a & b \\ a & b \\ a & b \\ a & b \end{pmatrix} \middle| a, b \in Z^+I \right\}$$

$$\cup \left\{ (a,a,a,a,a), \begin{bmatrix} a & b & c \\ d & e & f \end{bmatrix} \middle| a,b,c,d,e,f \in Q^+I \right\}$$

be a neutrosophic semigroup bivector space over the semigroup $S = Z^+I$.
Let
$$W = W_1 \cup W_2$$

$$= \left\{ \begin{pmatrix} a & a & a & a \\ a & a & a & a \end{pmatrix} \middle| a \in Z^+I \right\} \cup \{(a, a, a, a, a) \mid a \in Q^+I\}$$

$$\subseteq V_1 \cup V_2$$

be a neutrosophic semigroup bivector subspace of V over the semigroup $S = Z^+I$.

Define $\eta = \eta_1 \cup \eta_2 = W_1 \cup W_2 = W \to N([0, 1])$ where

$$\eta_1 : W_1 \to N([0, 1])$$

and

$$\eta_2 : W_2 \to N([0, 1])$$

with

$$\eta_1 \begin{pmatrix} a & a & a & a \\ a & a & a & a \end{pmatrix} = I + 1$$

and

$$\eta_2 (a, a, a, a, a) = I + 0.5.$$

$W_\eta = (W_1 \cup W_2)_{\eta_1 \cup \eta_2}$ is a neutrosophic semigroup fuzzy bivector space of V. Clearly $W_{\bar\eta}$ and $W_\eta$ are different in examples 6.10 and 6.12 respectively.

Now we proceed onto define the notion of neutrosophic semigroup fuzzy bilinear algebra.

**DEFINITION 6.7:** *Let $V = V_1 \cup V_2$ be a neutrosophic semigroup bilinear algebra over the semigroup S. We say $V_\eta = (V_1 \cup V_2)_{\eta_1 \cup \eta_2}$ or $(\eta_1 \cup \eta_2)(V_1 \cup V_2)$ is a neutrosophic semigroup fuzzy bilinear algebra if $\eta = \eta_1 \cup \eta_2 : V = V_1 \cup V_2 \to N([0, 1])$ is such that $\eta_1 : V_1 \to N([0, 1])$ and $\eta_2 : V_2 \to N([0, 1])$ satisfy the conditions $\eta_i (x_i + y_i) \geq \min (\eta_i (x_i), \eta_i (y_i)); \eta_i (rx_i) = r\eta_i (x_i); i=1, 2$ for every $r \in S$ and $x_i, y_i \in V_i ; i=1, 2.$*

We will illustrate this situation by some examples.

***Example 6.13:*** Let $V = V_1 \cup V_2$
$$= \{(a, a, a, a, a, a, a) \mid a \in Q^+I\} \cup$$

$$\left\{ \begin{pmatrix} a & b \\ c & d \end{pmatrix} \middle| a, b, c, d \in Z^+I \right\}$$

be a neutrosophic semigroup bilinear algebra over the semigroup $S = Z^+I$.

Define $\eta = \eta_1 \cup \eta_2 = V = V_1 \cup V_2 \to N([0, 1])$ where $\eta_1 : V_1 \to N([0, 1])$ and $\eta_2 : V_2 \to N([0, 1])$ are such that

$$\eta_1 (a, a, a, a, a, a, a) = I + \frac{1}{8}$$

$$\eta_2 \begin{pmatrix} a & b \\ c & d \end{pmatrix} = I + \frac{1}{9}$$

$V_\eta = (V_1 \cup V_2)_{\eta_1 \cup \eta_2}$ is a neutrosophic semigroup fuzzy bilinear algebra.

*Example 6.14:* Let

$$V = V_1 \cup V_2$$

$$= \left\{ \begin{pmatrix} a_1I & a_2I & a_3I \\ a_4I & a_5I & a_6I \end{pmatrix} \middle| a_iI \in Z_2I; 1 \leq i \leq 6 \right\}$$

and

$$V_2 = \{(a_1I, a_2I, a_3I, a_4I) \mid a_iI \in Z_2I; 1 \leq i \leq 4\}$$

be a neutrosophic semigroup bilinear algebra over the semigroup $S = Z_2$. Define $\eta = \eta_1 \cup \eta_2 : V_1 \cup V_2 \to N([0, 1])$ where $\eta_1 : V_1 \to N([0, 1])$ and $\eta_2 : V_2 \to N([0, 1])$; such that

$$\eta_1 \begin{pmatrix} a_1I & a_2I & a_3I \\ a_4I & a_5I & a_6I \end{pmatrix} = \begin{cases} I + \frac{1}{6} & \text{if } a_i \neq 0 \text{ for some } i, 1 \leq i \leq 6 \\ 1 & \text{if } a_i = 0 \text{ for } i = 1, 2, 3, 4, 5, 6 \end{cases}$$

and

$$\eta_2(a_1I, a_2I, a_3I, a_4I) = \begin{cases} I + \dfrac{1}{4} & \text{if } a_i \neq 0 \text{ for some } i = 1, 2, 3, 4 \\ 1 & \text{if } a_i = 0; 1 \leq i \leq 4 \end{cases}$$

$V_\eta = (V_1 \cup V_2)_{\eta_1 \cup \eta_2}$ is a neutrosophic semigroup fuzzy bilinear algebra.

The reader is expected to define the two types of neutrosophic semigroup fuzzy bilinear subalgebra as in case of neutrosophic semigroup bivector spaces.

**DEFINITION 6.8:** *Let $V = V_1 \cup V_2$ be a neutrosophic group bivector space over the group G. Define $\eta = \eta_1 \cup \eta_2 : V = V_1 \cup V_2 \to N([0, 1])$ a bimap such that*
$$\eta_1 : V_1 \to N([0, 1])$$
*and*
$$\eta_2 : V_2 \to N([0, 1])$$
*where*
$$\eta_i (a_i + b_i) \geq \min (\eta_i (a_i), \eta_i (b_i));$$
$$\eta_i (a_i) = \eta_i (- a_i)$$
$$\eta_i (0) = 1$$
$$\eta_i (ra_i) \geq \eta (a_i)$$

*for all $a_i, b_i \in V_i$, $r \in G$ for $i=1, 2$. We call $V_\eta = (V_1 \cup V_2)_{\eta_1 \cup \eta_2}$ to be the neutrosophic group fuzzy bivector space.*

It is pertinent to mention here that the concept of neutrosophic group fuzzy bilinear algebra and neutrosophic group fuzzy bivector spaces and fuzzy equivalent.
We will illustrate this by some examples.

*Example 6.15:* Let
$$V = V_1 \cup V_2$$
$$= \left\{ \begin{pmatrix} aI & bI & cI \\ dI & eI & fI \end{pmatrix} \middle| aI, bI, cI, dI, eI, fI \in ZI \right\} \cup \{ZI \times ZI\}$$

be a neutrosophic group bivector space over the group Z.

Let
$$\eta = \eta_1 \cup \eta_2 : V = V_1 \cup V_2 \to N([0, 1])$$
be defined as follows:

$$\eta_1 : V_1 \to N([0, 1])$$

and

$$\eta_2 : V_2 \to N([0, 1])$$

$$\eta_1 \begin{pmatrix} aI & bI & cI \\ dI & eI & fI \end{pmatrix} = \begin{cases} I + \dfrac{1}{|a|} & \text{if } a \neq 0 \\ I + \dfrac{1}{|b|} & \text{if } b \neq 0 \\ I + \dfrac{1}{|c|} & \text{if } a = 0 = b \text{ and } c \neq 0 \\ I + \dfrac{1}{|d|} & \text{if } a = b = c = 0 \text{ and } d \neq 0 \\ I + \dfrac{1}{|e|} & \text{if } a = b = c = d = 0;\ e \neq 0 \\ I + \dfrac{1}{|f|} & \text{if } a = b = c = d = e = 0;\ f \neq 0 \\ 1 & \text{if } a = b = c = d = e = f = 0 \end{cases}$$

$$\eta_2 (x, y) = \begin{cases} I + \dfrac{1}{2} & \text{if } (x, y) \neq (0,0) \\ 0 & \text{if } (x, y) = (0,0) \end{cases}$$

$V\eta = (V_1 \cup V_2)_{\eta_1 \cup \eta_2}$ is a neutrosophic group fuzzy bivector space.

*Example 6.16:* Let
$$V = V_1 \cup V_2$$

$$= \left\{ \begin{pmatrix} a & b \\ c & d \end{pmatrix} \middle| a,b,c,d \in QI \right\} \cup \{QI \times QI \times QI\}$$

be a neutrosophic group bilinear algebra over the group $G = QI$. Define $\eta = \eta_1 \cup \eta_2: V = V_1 \cup V_2 \to N([0, 1])$ where

$$\eta_1 : V_1 \to N([0, 1])$$

and

$$\eta_2 : V_2 \to N([0, 1])$$

with

$$\eta_1 \begin{pmatrix} a & b \\ c & d \end{pmatrix} = \begin{cases} I + \dfrac{1}{5} & \text{if } |ad = bc| \neq 0 \\ 0 & \text{if } |ad = bc| = 0 \end{cases}$$

and

$$\eta_2(a, b, c) = \begin{cases} I + \dfrac{1}{5} & \text{if atleast one of } a \text{ or } b \text{ or } c \text{ is non zero} \\ 0 & \text{if } a = b = c = 0 \end{cases}$$

$V_\eta = (V_1 \cup V_2)_{\eta_1 \cup \eta_2}$ is a neutrosophic group fuzzy bilinear algebra. Now we proceed onto define neutrosophic group fuzzy bilinear subalgebra or neutrosophic group fuzzy bivector subspace as both are fuzzy equivalent.

**DEFINITION 6.9:** *Let $V = V_1 \cup V_2$ a neutrosophic group bivector space over the group G. Let $V_\eta = (V_1 \cup V_2)_{\eta_1 \cup \eta_2}$ be the neutrosophic group fuzzy bivector space of V. Let $W = W_1 \cup W_2 \subseteq V_1 \cup V_2$ be a proper neutrosophic group bivector subspace of V $W_{\overline{\eta}} = (W_1 \cup W_2)_{\overline{\eta}_1 \cup \overline{\eta}_2}$ is a neutrosophic group fuzzy bivector subspace of $V_\eta$ if $\overline{\eta} = \overline{\eta}_1 \cup \overline{\eta}_2 : W_1 \cup W_2 \to N([0, 1])$ is the restriction bimap of $\eta : \eta_1 \cup \eta_2 : V_1 \cup V_2 \to N([0, 1])$ i.e., $\overline{\eta}_i : W_i \to N([0, 1])$ is the restriction map, $\eta_i : V_i \to N([0, 1])$ for $i = 1, 2$.*

We will illustrate this by an example.

*Example 6.17:* Let
$$V = V_1 \cup V_2$$
$$= \left\{ \begin{pmatrix} a & b \\ c & d \end{pmatrix} \middle| a, b, c, d \in QI \right\} \cup \{QI \times QI \times QI\}$$

be a neutrosophic group bivector space over the group $G = QI$.
Define $\eta = \eta_1 \cup \eta_2 : V_1 \cup V_2 \to N([0, 1])$ where
$$\eta_1 : V_1 \to N([0, 1])$$
and
$$\eta_2 : V_2 \to N([0, 1])$$
Define

$$\eta_1 \begin{pmatrix} a & b \\ c & d \end{pmatrix} = \begin{cases} I + \dfrac{1}{4} & \text{if } a \neq 0, b \neq 0, c \neq 0, d \neq 0 \\ I + \dfrac{1}{3} & \text{if one of } a, b, \text{or } c \text{ or } d \text{ is zero} \\ I + \dfrac{1}{2} & \text{if two of } a, b, c, \text{or } d \text{ is zero} \\ I + 1 & \text{if only one of } a, b, c \text{ or } d \text{ is non zero} \\ 1 & \text{if } a = b = c = d = 0 \end{cases}$$

$$\eta_2(a, b, c) = \begin{cases} I + \dfrac{1}{3} & \text{if } a \neq 0, b \neq 0, c \neq 0 \\ I + \dfrac{1}{2} & \text{if one of } a, b, \text{or } c \text{ is zero} \\ I + 1 & \text{if one of } a \text{ or } b \text{ or } c \text{ is non zero} \\ 1 & \text{if } a = b = c = 0 \end{cases}$$

$V\eta = (V_1 \cup V_2)_{\eta_1 \cup \eta_2}$ is a neutrosophic group fuzzy bivector space.

Let
$$W = W_1 \cup W_2$$

$$= \left\{ \begin{pmatrix} a & b \\ 0 & 0 \end{pmatrix} \middle| a, b \in QI \right\} \cup (QI \times QI \times \{0\}) \}$$

$$\subseteq V_1 \cup V_2$$

be a neutrosophic group bivector space over the group $G = QI$. Define $\bar{\eta} = \bar{\eta}_1 \cup \bar{\eta}_2 : W_1 \cup W_2 \to N([0, 1])$ by; where

$$\bar{\eta}_1 : W_1 \to N([0, 1])$$

and

$$\bar{\eta}_2 : W_2 \to N([0, 1])$$

is such that

$$\bar{\eta}_1 \begin{pmatrix} a & b \\ 0 & 0 \end{pmatrix} = \begin{cases} I+1 & \text{if one of a or b is non zero} \\ I+\dfrac{1}{2} & \text{if } a \neq 0, b \neq 0 \\ 1 & \text{if } a = 0 = b \end{cases}$$

$$\bar{\eta}_2 (a, b, 0) = \begin{cases} I+\dfrac{1}{2} & \text{if } a \neq 0, b \neq 0 \\ I+1 & \text{if } a = 0 \text{ or } b = 0 \\ 1 & \text{if } a = b = 0 \end{cases}$$

$W_{\bar{\eta}} = (W_1 \cup W_2)_{\bar{\eta}_1 \cup \bar{\eta}_2}$ is a neutrosophic group fuzzy bivector subspace of $V\eta$.

**DEFINITION 6.10:** *Let $V = V_1 \cup V_2$ be a neutrosophic group bivector space over the group G. Let $W = W_1 \cup W_2 \subseteq V_1 \cup V_2$ be a neutrosophic group bivector subspace of V over G. Define a bimap $\eta = \eta_1 \cup \eta_2$; $W_1 \cup W_2 \to N([0, 1])$ so that $W_\eta = (W_1 \cup W_2)_{\eta_1 \cup \eta_2}$ is a neutrosophic group fuzzy bivector space then we call $W\eta$ to be the neutrosophic group fuzzy bivector subspace of V.*

We will illustrate this definition by an example.

*Example 6.18:* Let
$$V = V_1 \cup V_2$$
$$= \{ZI \times ZI \times ZI \times ZI\} \cup \left\{ \begin{pmatrix} a & b \\ c & d \end{pmatrix} \middle| a,b,c,d \in QI \right\}$$

be a neutrosophic group bivector space over the group $G = ZI$. Let
$$W = W_1 \cup W_2$$
$$= \{ZI \times ZI \times ZI \times \{0\}\} \cup \left\{ \begin{pmatrix} a & b \\ 0 & d \end{pmatrix} \middle| a,b,d \in QI \right\}$$
$$\subseteq V_1 \cup V_2$$

be a neutrosophic group bivector space over the group G.

Define $\eta : W \to N([0, 1])$; that is $\eta = \eta_1 \cup \eta_2 : W = W_1 \cup W_2 \to N([0, 1])$ where $\eta_1 : W_1 \to N([0, 1])$ and $\eta_2 : W_2 \to N([0, 1])$.

$$\eta_1(x, y, z, 0) = \begin{cases} I + \dfrac{1}{3} & \text{if atleast one of } x, y, z \text{ is non zero} \\ 1 & \text{if } x = y = z = 0 \end{cases}$$

$$\eta_2 \begin{pmatrix} a & b \\ 0 & d \end{pmatrix} = \begin{cases} I + \dfrac{1}{4} & \text{if atleast one of } a, b, d \text{ is non zero} \\ 1 & \text{if } a = b = d = 0 \end{cases}$$

$W\eta = (W_1 \cup W_2)_{\eta_1 \cup \eta_2}$ is a neutrosophic group fuzzy bivector subspace of V.

**DEFINITION 6.11:** *Let $V = V_1 \cup V_2$ be a neutrosophic biset bivector space over the biset $S = S_1 \cup S_2$. Let $\eta = \eta_1 \cup \eta_2 : V_1 \cup V_2 \to N([0, 1])$ where $\eta_1 : V_1 \to N([0, 1])$ and $\eta_2 : V_2 \to N([0, 1])$ such that $\eta_1(r_1 a_1) \geq \eta_1(a_1)$ for all $r_1 \in S_1$ and $a_1 \in V_1$ and*

$\eta_2 (r_2 \, a_2) \geq \eta_2(a_2)$ for all $r_2 \in S_2$ and $a_2 \in V_2$. $V_\eta = (V_1 \cup V_2)_{\eta_1 \cup \eta_2}$ is a neutrosophic biset fuzzy bivector space.

We will give an example of this definition.

*Example 6.19:* Let
$$V = V_1 \cup V_2$$

$$= \{(Z^+I \times Z^+ \times Z^+I)\} \cup \left\{ \begin{pmatrix} a & a & a & a \\ a & a & a & a \end{pmatrix} \middle| a \in Z^+I \right\}$$

be a neutrosophic biset bivector space over the biset $3Z^+ \cup Z^+I = S = S_1 \cup S_2$.
Define $\eta = \eta_1 \cup \eta_2 : V_1 \cup V_2 \to N([0, 1])$ where $\eta_1 : V_1 \to N([0, 1])$ and $\eta_2 : V_2 \to N([0, 1])$ given by

$$\eta_1 (a, b, c) = I + \frac{1}{b}$$

$$\eta_2 \begin{pmatrix} a & a & a & a \\ a & a & a & a \end{pmatrix} = I + \frac{1}{8}$$

$V_\eta = (V_1 \cup V_2)_{\eta_1 \cup \eta_2}$ is a neutrosophic biset fuzzy bivector space.

Now the notion of neutrosophic biset fuzzy bilinear algebra is left as an exercise for the reader to define and give examples of it. The two types of neutrosophic biset fuzzy bivector subspaces can also defined analogous to earlier definitions.

**DEFINITION 6.12:** *Let $V = V_1 \cup V_2$ be a neutrosophic bisemigroup bivector space over the bisemigroup $S = S_1 \cup S_2$. Let $\eta = \eta_1 \cup \eta_2 : V = V_1 \cup V_2 \to N([0, 1])$ be a bimap; if $\eta_i (a_i + b_i) \geq min (\eta_i(a_i), \eta_i(b_i)) : \eta_i(r_i a_i) \geq r_i \eta(a_i)$ for all $r_i \in S_i$ and $a_i, b_i \in V_i$; $i = 1, 2$, then we call $V_\eta = (V_1 \cup V_2)_{\eta_1 \cup \eta_2}$ to be a neutrosophic bisemigroup fuzzy bivector space.*

We will illustrate this by an example.

*Example 6.20:* Let
$$V = V_1 \cup V_2$$
$$= \{(a, a, a) \mid a \in Z_6I\} \cup \left\{ \begin{pmatrix} a \\ a \\ a \\ a \\ a \end{pmatrix} \middle| a \in Z_7I \right\}$$

be a neutrosophic bisemigroup bivector space over the bisemigroup $S = S_1 \cup S_2 = Z_6 \cup Z_7$.

Define
$$\eta = \eta_1 \cup \eta_2 : V_1 \cup V_2 \rightarrow N([0, 1])$$
where
$$\eta_1: V_1 \rightarrow N([0, 1])$$
and
$$\eta_2: V_2 \rightarrow N([0, 1])$$
by

$$\eta_1(a, a, a) = \begin{cases} I + \dfrac{1}{3} & \text{if } a \neq 0 \\ 1 & \text{if } a = 0 \end{cases}$$

and

$$\eta_2 \begin{pmatrix} a \\ a \\ a \\ a \\ a \end{pmatrix} = \begin{cases} I + \dfrac{1}{4} & \text{if } a \neq 0 \\ 1 & \text{if } a = 0 \end{cases}.$$

$V_\eta = (V_1 \cup V_2)_{\eta_1 \cup \eta_2}$ is a neutrosophic bisemigroup fuzzy bivector space.

Next we proceed onto define the notion of neutrosophic bigroup fuzzy bivector spaces.

**DEFINITION 6.13:** *Let $V = V_1 \cup V_2$ be a neutrosophic bigroup bivector space over the bigroup $G = G_1 \cup G_2$. Let $\eta = \eta_1 \cup \eta_2 : V_1 \cup V_2 \to N([0, 1])$ be a bimap where $\eta_1 : V_1 \to N([0, 1])$ and $\eta_2 : V_2 \to N([0, 1])$ are such that $V_1 \eta_1$ and $V_2 \eta_2$ are neutrosophic group fuzzy vector spaces then $V_\eta = (V_1 \cup V_2)_{\eta_1 \cup \eta_2} = V_1\eta_1 \cup V_2\eta_2$ is a neutrosophic bigroup fuzzy bivector space.*

We will illustrate this situation by some example.

*Example 6.21:* Let
$$V = V_1 \cup V_2$$

$$= \left\{ \begin{pmatrix} a & b \\ c & d \end{pmatrix} \middle| a,b,c,d \in Z_{20}I \right\}$$

$$\cup \{Z^+I \times Z^+I \times Z^+I \times Z^+I \times Z^+I\}$$

be a neutrosophic bigroup fuzzy bivector space over the bigroup $G = Z_{20} \cup Z^+I$.

Define
$$\eta : V \to N([0, 1])$$
where
$$\eta = \eta_1 \cup \eta_2 : V_1 \cup V_2 \to N([0, 1]);$$

$$\eta_1 : V_1 \to N([0, 1])$$
and
$$\eta_2 : V_2 \to N([0, 1]).$$

$$\eta_1 \begin{pmatrix} a & b \\ c & d \end{pmatrix} = \begin{cases} I + \dfrac{1}{8} & \text{if } |ad - bc| \neq 0 \\ 1 & \text{if } ad - bc = 0 \end{cases}$$

and

$$\eta_2\,(a, b, c, d, e) = \begin{cases} I + \dfrac{1}{5} & \text{if } (a,b,c,d,e) \neq (0,0,0,0,0) \\ 1 & \text{if } (a,b,c,d,e) = (0,0,0,0,0) \end{cases}$$

$V_\eta = (V_1 \cup V_2)_{\eta_1 \cup \eta_2}$ is a neutrosophic bigroup fuzzy bivector space.

As in case of neutrosophic biset vector space we can define two types of neutrosophic bigroup fuzzy bivector subspaces. This task is left as an exercise for the reader.

Chapter Seven

# SUGGESTED PROBLEMS

(1) Find all mixed neutrosophic integer set subspace of the mixed neutrosophic integer set vector space $V = \{20I, 0, 40 - I, 25I - 5, 45I - 3, 28I - 4, 12 + 4I, 3 + 4I, 40 - 4I, 10I, 201, 28I, 54\} \subseteq N(Z)$ over the set $S = \{0, 1\} \subseteq Z$.

   a. Does V have pseudo pure neutrosophic integer set vector subspace?

   b. Does V contain any pseudo set integer set vector subspace?

(2) Find all pure neutrosophic integer set vector subspaces of the pure neutrosophic integer set vector space $V = \{0, I, 28I, 42 + I, -79 + 3I, 442, 89I, 200I + 4002, 42I + 381, 451I\}$ over the set $S = \{0, 1\} \subseteq Z$.

(3) Let $V = \{3nI, 5mI + 2n, 0, 2m \mid m, n \in Z^+\}$ be the mixed neutrosophic integer set vector space over $S = \{0, 1\} \subseteq Z$.

a. Find 3 mixed neutrosophic integer set vector subspace of V.

b. Find 3 pseudo pure neutrosophic integer set vector subspaces of V.

c. Find 3 pseudo set integer set vector subspace of V.

d. Can V have mixed neutrosophic integer subset vector subspace? Justify your claim.

(4) Let $V = \left\{ \begin{pmatrix} aI & b \\ c & dI \end{pmatrix}, \begin{pmatrix} 0 & 0 \\ 3nI & 0 \end{pmatrix}, \begin{pmatrix} 0 & nZ \\ 0 & pzI \end{pmatrix} \middle| n, p, a, b, c, d \in Z^+ \right\}$

be a pure neutrosophic integer set vector space over $Z^+ \subseteq Z$.

a. Find atleast 5 pure neutrosophic integer set vector subspace of V over $Z^+$.

b. Find atleast 5 pure neutrosophic integer subset vector subspace of V.

(5) Let V = {2I, 0, 3I + 2, 27, 38 – 3I, 54I – 47, 280I, 249} and W = {41I, 156I, 31 – I, 48I, 56 + 47I, 56, 0, 27 + 48I, 56I + 21} be two mixed neutrosophic integer set vector spaces over the set S = {0, 1} $\subseteq$ Z.

a. Find a neutrosophic integer set linear transformation of V to W.

b. Find atleast one neutrosophic integer set subspace preserving linear transformation.

c. Find one neutrosophic integer set pseudo set subspace preserving linear transformation.

d. Find one neutrosophic integer pseudo pure subspace preserving linear transformation.

(6) Let V = {0, 2I, 25I, 3I – 2, 48 – 5I, 28I + 4} $\subseteq$ PN(Z) be a pure neutrosophic integer set vector space over the set S = {0, 1} $\subseteq$ Z.

Find all neutrosophic set linear operators on V. How many of them preserve the pure neutrosophic integer set subspaces of V.

(7) Let V = {0, 2I + 8, 27 – 3I, 4I, –51, 48I – 3, 27 – 8I, 9I + 8, – 5I, 49 – IM 28M – 26, 42, – 21} ⊆ N(Z) be the mixed neutrosophic integer set vector space over S = {0, 1}.

   a. How many neutrosophic integer set linear operators can be defined on V?

   b. Find those neutrosophic integer set linear operators on V which preserves subspaces of V.

(8) Let V = {I, 2I, 0, 19I – 3, 27 – I, I + 4, 28, 48 – 31I, 151I} be a neutrosophic integer set vector space over the set S = {0, 1}.

   a. Find all neutrosophic integer set vector subspace of V over S.

   b. Find atleast 3 neutrosophic integer linear operators on V which preserves all types of subspaces.

(9) Let V = {2ZI, 5Z, m + nI | m, n ∈ Z} be the neutrosophic integer set linear algebra over the set Z.

   a. Find at least 5 neutrosophic integer subset linear subalgebras of V over S ⊆ Z.

   b. Find 5 pseudo neutrosophic integer set vector subspaces of V.

   c. Find 5 neutrosophic integer set linear subalgebras of V.

   d. Find at least 5 neutrosophic integer set linear operators which preserves atleast one of the three substructures.

   e. Does there exist any neutrosophic integer set linear operator on V which preserves simultaneously all the three substructures mentioned in the problem.

(10) Let V = {3ZI, 25ZI, 41Z, 9Z + 32ZI} be a neutrosophic integer set vector space over the set Z.

   a. Find neutrosophic integer set subvector spaces of V over Z.

   b. Find neutrosophic integer set linear operators on V.

  c. Find neutrosophic integer subset subvector spaces of V.

(11) Let $V = \{2I, 0, 41I, 21 - I, 4, 43I, 25 + 3I, 27 - I, I - 42, 3, 48I + 90\} \subseteq N(Z)$ be a neutrosophic integer set vector space over the set $\{0, 1\}$. Find a neutrosophic integer set generator of V.

(12) Let $V = \{3ZI, 41Z, 28ZI + 31Z\} \subseteq N(Z)$ be a neutrosophic integer set vector space of V over $Z^+$.

  a. Find a neutrosophic integer set generator of V over $Z^+$.

  b. If V is defined over Z what will be the neutrosophic integer set generator of V over Z.

  c. If V is defined over $3Z^+$ what will be the neutrosophic integer set generator of V over $3Z^+$.

  d. If V is defined over $S = \{0, 1\}$ what will be the neutrosophic integer set generator of V over S.

  e. If V is defined over the set $S = \{-1, 2, 0, 1\}$, what will be the neutrosophic integer set generator of V over S.

(13) Let $V = \{3ZI, 2Z, m + nI \mid m, n \in Z\}$ be a neutrosophic integer set linear algebra over Z.

  a. Find the neutrosophic integer set generator of V over Z.

  b. Find the neutrosophic integer set generator of V when the same V is defined over the set $S = \{0, 1\}$.

  c. Suppose the neutrosophic integer set linear algebra V is defined over $3Z^+$ what will be the neutrosophic integer set generator of V over the set $3Z^+$.

  d. What will be the neutrosophic integer set generator of V if V is defined over the set $S = \{0, 1, 2, -1, -2\}$?

(14) Give some interesting properties about neutrosophic integer set linear algebras.

(15) Give an example of a neutrosophic integer set linear algebra which has no neutrosophic integer set linear subalgebra.

(16) Give an example of a neutrosophic integer set linear algebra which has no neutrosophic integer subset linear subalgebra.

(17) Give some interesting results on substructures of neutrosophic integer set vector spaces.

(18) Obtain some interesting properties about n-n set vector spaces.

(19) Find all n-n set vector subspaces of the n-n set vector space $V = \{7 \pm 7I\, 27I, 48 \pm 48I, 56, 42I, 7 \pm 70I, -28, 56 - 56I, -28 \pm 28I, 0\}$ defined over the set $S = \{0, 1, 1 - I\}$.

(20) Let V be a neutrosophic integer set vector space over a set $S \subseteq Z$. Let $\text{NHom}_S(V, V)$ denote the collection of all linear operators of V, what is the algebraic structure enjoyed by $\text{NHom}_S(V, V)$?

(21) Let W and V be neutrosophic integer set vector space over a set $S \subseteq Z$. Let $\text{NHom}_S(V, W)$ denote the set of all linear transformation of V into W. Find the algebraic structure enjoyed by $\text{NHom}_S(V, W)$.

(22) Give some interesting properties about n-n set linear algebra.

(23) Let $V = \{24I, 22 - I, 90 + 4I, 22 + 2I, 0, 21 + 9I, 21I, 30I, 94I\} \subseteq N(Z)$ be a n-n set vector space over the set $S = \{0, 1, I\}$. Find all n-n set vector subspaces of V over S.

(24) Let $V = \{m - mI \mid m \in Z\}$ be a n-n set linear algebra over $S = \{0, 1, 1 - I\} \subseteq N(Z)$. Find the set of all n-n set linear operators on V.

(25) Let $V = \{2I, 3I - 2, 0, 44I - 20, 27I + 9, 22 + I, I, 24I, 36I, 23I\} \subseteq N(Z)$ and $W = \{5 - 5I, 20 - 20I, 30 - 30I, I, 44 - 44I, 26I, 72I\} \subseteq N(U)$ be n-n set vector space over the set $S = \{1, 0, I\}$. Find atleast 5 distinct n-n set linear transformation from V to W which preserves n-n set vector subspaces of V and one n-n set linear transformation from V to W which does not preserve subspaces!

(26) Find the algebraic structure enjoyed by $NHom_S(V, W)$ where V and W are n-n set vector spaces defined over S; where $NHom_S(V, W)$ is the collection of all n-n set linear transformations of V to W.

(27) What is the algebraic structure of $NHom_S(V, V)$; where $NHom_S(V, V)$ is the collection of all n-n set linear operators of a n-n set vector space over the set $S \subseteq N(Z)$?

(28) When V and W in problem 26 is replaced by n-n set linear algebra, what is the structure of $NHom_S(V, W)$?

(29) When V in the problem 27 is replaced by n-n set linear algebra, what is the structure of $NHom_S(V, V)$?

(30) Let $V = \left\{ \begin{pmatrix} m-mI & 0 \\ 0 & m-mI \\ m-mI & m-mI \end{pmatrix} \middle| m \in Z^+ \right\} \subseteq N(Z)$ be the n-n set linear algebra over the set $S = \{0, 1, 1-I\} \subseteq N(Z)$.

    a. Find $NHom_S(V, V)$.

    b. Find n-n subset linear subalgebras of V over F.

    c. Does V have pseudo n-n set vector subspaces?

    d. Does there exist a n-n set linear operator on V which does not preserve any n-n set linear subalgebras?

(31) Obtain some interesting properties about set neutrosophic integer vector spaces.

(32) Let $V = \{ZI\} \subseteq N(Z)$ be a set neutrosophic integer linear algebra over the set $S = Z^+$.

    a. Find at least 3 set neutrosophic sublinear algebras of V.

    b. Does V have pseudo set neutrosophic integer subvector spaces?

    c. Find some subset neutrosophic integer sublinear algebras of V.

(33) Let V = {3 ± 3I, 5 ± 5I, 0, 3I, 8 ± 8I, 7I, 42 ± 42I} ⊆ N(Z (I)); V is a set neutrosophic integer vector space over the set S = {0, 1, I, 1 – I} ⊆ N(Z) or equivalently N(Z(I)).

    a. Find atleast 3 set neutrosophic integer vector subspaces over the set S.

    b. Does there exist subset neutrosophic integer vector subspace in V?

(34) Obtain some interesting properties about set neutrosophic integer linear algebras.

(35) Let V ⊆ N($Z_n$) be a set neutrosophic modulo integer vector space over the set S = {0, 1}; find some interesting properties about V.

(36) Let V = {I, 0, 3I, 14I, 5I, 25I, 10, 10I, 3, 14} ⊆ N($Z_{26}$) be a set neutrosophic modulo integer vector space over the set S = {0, 1, I, 1 + 25I, 25 + I}.

    a. Find a set neutrosophic modulo integer vector linear operator T on V which preserves all substructures.

    b. Find NS ($Hom_S$ (V, V)). What is the | $NSHom_S$ (V, V)| ?

(37) Prove V = {0, I, 2I, …., 22I} ⊆ N($Z_{23}$) is doubly simple neutrosophic modulo integer linear algebra over the set S = {0, I, 1 + 23I}.

(38) Obtain some interesting properties about V = {0, I + (p – 1), (p – 1) I + 1, 2I + (p – 2), (p – 2)I + 2, (p – 3) + 3I, (p – 3) I + 3, … , $\frac{p-1}{2}$ I + $\frac{p+1}{2}$, $\frac{p+1}{2}$ I + $\frac{p-1}{2}$} ⊆ N($Z_p$); p a prime, a set neutrosophic modulo integer algebra over the set S = V ∪ {0, 1, I}.

(39) Obtain some interesting properties about set neutrosophic real matrix vector spaces.

(40) Can there exist a set neutrosophic real matrix linear algebra which has finite cardinality?

(41) Give an example of a neutrosophic group linear algebra of infinite cardinality.

(42) Does there exist a neutrosophic set linear algebra which is simple?

(43) Give an example of a neutrosophic group linear algebra of finite dimension.

(44) Obtain some interesting results on neutrosophic group vector space.

(45) Define a linear transformation from the neutrosophic semigroup vector space V to W where

$$V = \left\{ \begin{pmatrix} a_1I & a_2I & a_3I \\ a_4I & a_5I & a_6I \end{pmatrix} \middle| a_iI \in QI; 1 \leq i \leq 6 \right\}$$

and

$$W = \left\{ \begin{pmatrix} b_1I & b_2I & b_7I \\ b_3I & b_4I & b_8I \\ b_5I & b_6I & b_9I \end{pmatrix} \middle| b_iI \in QI; 1 \leq i \leq 9 \right\}.$$

Does there exists a linear transformation T from V to W such that $T^{-1}$ exists?

(46) Does there exists neutrosophic set vector spaces V and W defined over the set S such that there does not exist any T : V → W such that $T^{-1}$ exist? Justify your claim.

(47) If G is a simple group, can we say if V is a neutrosophic group vector space defined over G also is simple?

(48) Give an example of a simple neutrosophic group linear algebra.

(49) Let $V = \left\{ \begin{pmatrix} aI & bI \\ cI & dI \end{pmatrix} \middle| aI, bI, cI, dI \in QI \right\}$ be a neutrosophic group linear algebra defined over the group G = pZI (p a prime). Find proper neutrosophic group linear subalgebras of V.

Can V contain proper neutrosophic subgroup linear subalgebras? Justify your answer!

(50) Give some interesting properties about neutrosophic semigroup linear algebras.

(51) Give an example of a simple neutrosophic semigroup linear algebra.

(52) Does the neutrosophic semigroup linear algebra $V = \left\{ \begin{pmatrix} aI & bI \\ cI & dI \\ eI & fI \end{pmatrix} \middle| aI, bI, cI, dI, eI, fI \in Z_2I \right\}$ defined over $S = Z_2I$ have proper neutrosophic subsemigroup linear subalgebra? Justify your claim. Can V have proper neutrosophic semigroup linear subalgebras?

(53) Can a neutrosophic semigroup linear algebra have pseudo neutrosophic semigroup linear subalgebra? If so give examples of them?

(54) Let $V = \{Z_7I [x]$ where $Z_7I [x]$ consists of all polynomials with coefficients from $Z_7I\}$ be a neutrosophic semigroup linear algebra over the semigroup $S = Z_7I$. Can V have pseudo neutrosophic semigroup linear subalgebras? Justify your claim.

(55) Prove if V is any neutrosophic semigroup linear algebra over a semigroup $S = ZI$ (or $Z^+I$ or $Q^+I$ or $R^+I$ or QI or RI) then V cannot contain pseudo neutrosophic semigroup linear subalgebras.

(56) Give some interesting properties about the substructures of a neutrosophic group vector space.

(57) Give an example of a neutrosophic group vector space which is simple.

(58) Give an example of a neutrosophic group linear algebra which has no proper pseudo neutrosophic group linear subalgebras.

(59) Give an example of a simple neutrosophic group linear algebra defined over a group G where the group G is not simple.

(60) Let $V = \left\{ \begin{pmatrix} aI & aI \\ aI & aI \end{pmatrix} \middle| aI \in Z_7I \right\}$ be a neutrosophic group linear algebra over the group $G = N(Z_7)$. Is V simple? Justify your claim. Can this problem be generalized for any prime p?

(61) Find all neutrosophic set vector subspaces of $V = \left\{ \begin{pmatrix} aI & 0 \\ bI & 0 \end{pmatrix}, \begin{pmatrix} 0 & aI \\ bI & cI \end{pmatrix}, (aI, bI, cI), \begin{pmatrix} aI \\ bI \end{pmatrix}, Z_{11}I \times Z_{11}I \middle| aI, bI, cI \in Z_{11}I \right\}$ over the set $S = Z_{11}I$.

(62) Find all neutrosophic subset vector subspaces of $V = \left\{ \begin{pmatrix} aI & aI \\ 0 & aI \end{pmatrix}, \begin{pmatrix} aI & 0 \\ aI & aI \end{pmatrix}, (aI, aI, aI) \middle| aI \in Z_{20}I \right\}$ defined over the set $S = \{0, 5, 10\}$.

(63) Find all neutrosophic subset vector subspaces of $V = \{(aI, bI, cI) \mid aI, bI, cI \in Z_{20}I\}$ defined over the set $N(Z_{20})$.

(64) Can we say if V is a neutrosophic set vector space over the set S and if V has a proper neutrosophic set vector subspaces over the set S then V has a proper neutrosophic subset vector subspace?

(65) Is the claim if V a neutrosophic set vector space over the set S has proper neutrosophic subset vector subspace then V has proper neutrosophic set vector subspace? Justify your answer.

(66) Given $V = \left\{ \begin{pmatrix} aI & bI & cI \\ 0 & dI & 0 \end{pmatrix}, \begin{pmatrix} a_1I & a_2I & a_3I \\ a_4I & a_5I & a_6I \\ a_7I & a_8I & a_9I \end{pmatrix}, (aI, aI, aI, aI) \middle| \begin{array}{l} a,b,c,c,d,a_i \in Z_{25}I \\ 1 \le i \le 9 \end{array} \right\}$

is a neutrosophic set vector space over the set $S = N(Z_{25})$.

a. Find all proper neutrosophic subset vector subspaces of V.

b. Find all proper neutrosophic set vector subspaces of V.

c. Does V have pseudo neutrosophic set linear subalgebra?

d. Can V have pseudo neutrosophic set vector subspaces?

(67) Is $Z_{13}I$ a simple neutrosophic semigroup?

a. Construct a neutrosophic set vector space V over $Z_{13}I$ which is a simple neutrosophic set vector space (over $Z_{13}I$).

b. Construct a neutrosophic set vector space V over $N(Z_{13})$ which is a simple neutrosophic set vector space (over $N(Z_{13})$).

(68) Let $V = \left\{ \begin{pmatrix} aI & aI & aI & aI \\ aI & aI & aI & aI \end{pmatrix} \middle| a \in Z_{19}I \right\}$ be a neutrosophic semigroup vector space over the semigroup $S = Z_{19}I$. Is V simple? Justify your claim.

(69) Let V in problem (68) be taken as a neutrosophic group vector space over $S = N(Z_{19})$. Is V simple?

(70) Is V defined in problems (68) and (69) neutrosophic semigroup linear algebra over $S = Z_{19}I$ and neutrosophic group linear algebra over the group $S = N(Z_{19})$ respectively?

(71) Give some interesting properties about the collection of linear transformation operators of a neutrosophic group linear algebras V over G. Does the collection of such linear operators of V form a neutrosophic group linear algebra over G?

(72) Given $V = \left\{ \begin{pmatrix} aI & aI \\ aI & aI \end{pmatrix} \middle| aI \in Z_{23}I \right\}$ be a neutrosophic group linear algebra over the group $G = N(Z_{23})$. Find the set of all

linear operators of V to V. Does the collection of linear operators of V to V form a neutrosophic group linear algebra over the group $G = N(Z_{23})$?

(73) Let V be a neutrosophic set vector space over the set S. Let $N(Hom_S (V, V))$ denote the set of all linear operators of V to V. What is the algebraic structure of $NHom_S (V, V)$?

(74) Suppose V is a neutrosophic semigroup linear algebra over the semigroup S. Let $N(Hom_S (V, V))$ denote the collection of all linear operators of V to V. Is $N(Hom_S (V, V))$ a neutrosophic semigroup linear algebra over the semigroup S?

(75) Let V and W be neutrosophic semigroup vector spaces over the semigroup S. Suppose $N(Hom_S (V, W))$ denote the collection of all linear transformations of V to W, what is the algebraic structure enjoyed by the collection $N(Hom_S (V, W))$? Suppose $IN(Hom_S(V, W))$ denote the collection of all invertible linear transformations of V to W what is the algebraic structure of $IN(Hom_S (V, W))$?

If $N(Hom_S (V, W)$ denote the collection of all linear transformations of W to V can we find any relation between the collections $N(Hom_S (V, W))$ and $N(Hom_S(W, V))$?

Does these exist any relation between $I(NHom_S(V, W))$ and $INHom_S (W, V)$?

(76) Give some important properties about $I(NHom_G(V, W)$ (ii) $INHom_G(W, V)$ where V and W are neutrosophic group linear algebras defined over the group G.

(77) What is the difference between the algebraic structures of $N(Hom_S(V, W))$ and $(NHom_G(V, W)$? (Here V and W are neutrosophic semigroup linear algebras defined over the semigroup S and V and W are neutrosophic group linear algebras defined over the group G respectively).

(78) If V and W are finite dimensional neutrosophic set vector spaces over the set S, what can be said about the dimension of $NHom_S(V, W)$? Is $N(Hom_S(V, W)$ finite dimensional?

(79) Does there exist any relation between the dimensions of V and W and that of the dimension of N(Hom$_S$ (V, W))?

(80) Does the structures V, W and N(Hom$_S$(V, W)) enjoy any common algebraic properties? (Here V and W are neutrosophic set vector spaces defined over the set S).

What happens if V and W are neutrosophic set linear algebras defined over S?

(81) Let V be a neutrosophic group linear algebra over the group G of dimension say N(n < ∞). Suppose N(H$_G$(V, V)) denote the set of all linear operators on V. CaN(NHom$_G$(V, V)) have any form of dimension associated with it?

(82) Let $V = \left\{ \begin{pmatrix} aI & aI \\ aI & aI \end{pmatrix} \middle| aI \in Z_{11}I \right\}$ be a neutrosophic group vector space over the group G = $Z_{11}I$. Find N(Hom$_G$(V, V)). Suppose V is realized as a neutrosophic semigroup linear algebra what can be said about N(Hom$_G$(V, V))? Suppose V is realized only as a neutrosophic set linear algebra over G, what can we say about the algebraic structure of N(Hom$_G$(V, V))?

(83) Obtain some interesting properties about V = $V_1 \cup V_2$ where V is a neutrosophic set bivector space.

(84) Give an example of a neutrosophic set bivector space.

(85) Let $V = V_1 \cup V_2 = \left\{ \begin{pmatrix} aI & bI \\ cI & 0 \end{pmatrix} \middle| aI, bI, cI \in ZI \right\} \cup \{(aI, aI, aI,$ aI, aI) | aI ∈ ZI} be a neutrosophic set bivector space over the set S = $3Z^+I$

   a. Find neutrosophic set bivector subspaces of V.

   b. Find pseudo neutrosophic set bilinear sub algebras of V.

   c. Does V have pseudo neutrosophic set bilinear subalgebras?

(86) Obtain some interesting properties about neutrosophic biset bivector spaces.

(87) Prove a neutrosophic biset bivector space in general is not a neutrosophic biset bilinear algebra.

(88) Let $V = V_1 \cup V_2$ be a neutrosophic biset bivector space over the biset $S = Z_{25}I \cup ZI$ where $V_1 = \left\{ \begin{bmatrix} a \\ a \\ a \\ a \\ a \end{bmatrix}, Z_{25}I \times Z_{25}I \times Z_{25}I \right\}$ and $V_2 = \{(a, b, c, d, e), \begin{pmatrix} a & b & c \\ d & e & f \end{pmatrix} \mid a, b, c, d, e, f \in ZI\}$.

   a. Find atleast 5 neutrosophic biset bivector subspaces of V.
   b. Find three distinct neutrosophic biset bivector operators on V.
   c. Find three distinct neutrosophic subset bilinear bivector subspaces of V.
   d. What is the algebraic structure enjoyed by $N(\text{Hom}_S(V, V)) = \{$set of all bilinear operators on V$\}$?
   e. Can V have pseudo neutrosophic biset bilinear subalgebras?

(89) Let $V = V_1 \cup V_2$

$= \left\{ \begin{pmatrix} a & a & a \\ a & a & a \\ a & a & a \end{pmatrix} \middle| a \in QI \right\} \cup \left\{ \begin{pmatrix} a & b \\ c & d \end{pmatrix} \middle| a, b, c, d \in QI \right\}$

be a neutrosophic group bivector space over the group $Z = G$. Find the bidimension of V over Z. What is the dimension of $S(\text{Hom}_G(V, V))$ a neutrosophic group vector space over G or a neutrosophic group linear algebra over G or a neutrosophic group bivector space over G? Justify your claim!

(90) Let $V = V_1 \cup V_2$ be a neutrosophic semigroup bivector space over a semigroup S. Prove in general V is not a neutrosophic group bivector space even if S is a group!

(91) Let $V = V_1 \cup V_2 = \left\{ \begin{pmatrix} a & b & c \\ d & e & f \end{pmatrix} \middle| a,b,c,d,e,f \in Z_{12}I \right\} \cup \left\{ \begin{pmatrix} a & a \\ b & b \\ c & c \\ d & d \end{pmatrix} \middle| a,b,c,d \in Z_{27}I \right\}$ be a neutrosophic bigroup bivector over the bigroup $G = Z_{12} \cup Z_{27}$. Find a neutrosophic bigroup bivector subspace W of V. Define a neutrosophic linear operator T which preserves this subspace; that is $T(W) \subseteq W$. What is the bidimension of $S(Hom_G(V, V))$? Find the algebraic structure enjoyed by $S(Hom_G(V, V))$.

(92) Give some interesting properties of neutrosophic semigroup bivector spaces.

(93) What is the difference between neutrosophic bisemigroup bivector space and neutrosophic semigroup bivector space?

(94) Obtain some interesting features especially enjoyed by neutrosophic group bilinear algebra.

(95) Give an example of a neutrosophic bigroup bivector space of bidimension (7, 5) over a bigroup $G = G_1 \cup G_2$.

(96) Can one claim the neutrosophic bidimension in general is reduced if we consider neutrosophic bigroup bilinear algebra instead of neutrosophic bigroup bivector spaces?

(97) What is the advantage of neutrosophic set linear algebras over a neutrosophic linear algebras?

(98) What is the benefit of using neutrosophic biset bivector spaces instead of a neutrosophic set bivector space?

(99) Give some nice applications of neutrosophic bisemigroup bilinear algebras.

(100) Which bispace is useful in neutrosophic modeling; neutrosophic bigroup bivector space or neutrosophic group bivector space?

(101) What is the advantage of using neutrosophic bigroup bivector space instead of using a neutrosophic group bivector space?

(102) Is their any benefit in using a neutrosophic set bivector space instead of a neutrosophic set vector space?

(103) What is the advantage of studying neutrosophic bigroup bivector space in the place of neutrosophic bigroup bilinear algebra?

(104) Does there exist any generalized neutrosophic bistructures other than the neutrosophic set bivector spaces?

(105) Give an example of a (11, 19) bidimensional neutrosophic bigroup bivector space.

(106) Can one say a neutrosophic bigroup bivector space of bidimension $(p_1, p_2)$ ($p_1$ and $p_2$, two distinct primes) is always bisimple?

(107) Can a neutrosophic group bivector space of bidimension say (12, 15) be simple? Justify your claim!

(108) Give an example of a neutrosophic group bilinear algebra of finite bidimension, which is simple.

(109) Give an example of a neutrosophic group bilinear algebra of infinite bidimension, which is simple.

(110) Can one say a neutrosophic group bilinear algebra of bidimension (1, 1) is always simple? Justify your claim.

(111) Give an example of a neutrosophic group bilinear algebra of bidimension (1, 1).

(112) What is the difference between a simple neutrosophic group bivector space and a simple neutrosophic group bilinear algebra?

(113) Suppose $V = V_1 \cup V_2$ is a neutrosophic group bivector space defined over the group $G = \{(0, 1, 2, \ldots, p - 1)$ such that p is a prime$\}$. Is V a simple neutrosophic bivector space over G?

(114) Find a bigenerating subset of the neutrosophic group bivector space $V = \left\{ \begin{pmatrix} a & b & c \\ d & e & f \end{pmatrix} \middle| a,b,c,d,e,f \in Z_{12}I \right\} \cup \left\{ \begin{pmatrix} a & a & a & a \\ a & a & a & a \end{pmatrix} \middle| a \in Z_{12}I \right\}$ over the group $G = Z_{12}I$.

(115) Let $V = V_1 \cup V_2 = \left\{ \begin{pmatrix} a & a & a \\ a & a & a \end{pmatrix} \middle| a \in Z_{11}I \right\} \cup \left\{ \begin{pmatrix} a & a \\ a & a \\ a & a \\ a & a \end{pmatrix} \middle| a \in Z_{13}I \right\}$ be a neutrosophic bigroup bilinear algebra over the bigroup $G = G_1 \cup G_2 = Z_{11} \cup Z_{13}$. Find a bigenerating subset of V. If $V = V_1 \cup V_2$ is a neutrosophic bigroup bilinear algebra over the bigroup $G = G_1 \cup G_2 = Z_{11}I \cup Z_{13}I$, what is the bidimension of the bigenerating subset of V? If $V = V_1 \cup V_2$ is a neutrosophic bigroup bilinear algebra over the bigroup $G = G_1 \cup G_2 = N(Z_{11}) \cup N(Z_{13})$, what is the bigenerating subset of V?

When is V simple?

(116) Give an example of a (7, 14) bidimension neutrosophic semigroup bivector space.

(117) Give an example of a neutrosophic semigroup bilinear algebra of bidimension (7, 14).

(118) Compare the algebraic structure in problems (116) and (117).

(119) Give by an real model that neutrosophic set vector space is useful than a neutrosophic vector space.

(120) Obtain some interesting properties about neutrosophic set fuzzy bilinear algebras.

(121) Let $V = V_1 \cup V_2$

$$= \left\{ \begin{pmatrix} a & b & c \\ d & e & f \end{pmatrix}, \begin{pmatrix} a & b \\ c & d \\ e & f \end{pmatrix} \middle| a,b,c,d,e,f \in QI \right\} \cup$$

$$\left\{ (a,a,a,a,a), \begin{pmatrix} a & a \\ b & b \\ c & c \\ d & d \end{pmatrix} \middle| a,b \in Z^+I \right\}$$

be a neutrosophic set bivector space over the set $S = Z^+I$.

Let $W = W_1 \cup W_2 = \left\{ \begin{pmatrix} a & b & c \\ d & e & f \end{pmatrix} \middle| a,b,c,d,e,f \in Z^+I \right\} \cup$

$\{(a, a, a, a, a) \mid a \in Z^+I\} \subseteq V_1 \cup V_2$ be a neutrosophic set bivector subspace of V over the set S. Find neutrosophic linear bioperator on V which preserves W. Find one neutrosophic linear bioperator of V which does not preserve W. Find $V_\eta$ and $W_{\bar\eta}$.

(122) Give some interesting properties about neutrosophic semigroup fuzzy bivector spaces.

(123) Prove the notion of neutrosophic semigroup fuzzy bilinear algebra and neutrosophic semigroup fuzzy bivector space are fuzzy equivalent.

(124) Let $V = V_1 \cup V_2$

$$= \left\{ \begin{pmatrix} a & a & a \\ a & a & a \\ b & b & b \end{pmatrix}, \begin{pmatrix} a & b \\ a & b \\ a & b \\ a & b \end{pmatrix}, (a,b,a,b,a,b) \middle| a,b \in Z^+I \right\} \cup$$

$$\left\{ \begin{pmatrix} a & a & a \\ 0 & b & b \\ 0 & 0 & a \end{pmatrix}, \begin{pmatrix} a & a & a & a \\ 0 & b & b & b \\ 0 & 0 & a & a \\ 0 & 0 & 0 & b \end{pmatrix} \middle| a, b \in Z^+ I \right\}$$

be a neutrosophic semigroup bivector space over the semigroup $S = Z^+ I \cup \{0\}$. Define $\eta : V \to N([0, 1])$ so that $V_\eta = (V_1 \cup V_2)_{\eta_1 \cup \eta_2}$ is a neutrosophic semigroup fuzzy bivector space. Define a $\eta$ so that $W_{\overline{\eta}}$ does not exist for any proper neutrosophic semigroup bivector subspace W of V where $\overline{\eta}$ is the restriction of $\eta$ to W.

(125) Let $V = V_1 \cup V_2$

$$= \left\{ \begin{pmatrix} a & a \\ a & a \\ a & a \end{pmatrix}, (a\ a\ a\ a\ a) \middle| a \in Z_{17} I \right\}$$

$$\cup \left\{ \begin{pmatrix} a & a & a \\ a & a & a \\ a & a & a \end{pmatrix}, \begin{pmatrix} a \\ a \\ a \\ a \\ a \end{pmatrix} \middle| a \in Z_{11} I \right\}$$

be a neutrosophic bigroup bivector space over the bigroup $G = Z_{17} I \cup Z_{11} I$. Find a $\eta : V \to N[Z_0, 1]$ so that $\eta$ is a substructure preserving neutrosophic bigroup fuzzy bivector space.

(126) Obtain some interesting properties about neutrosophic bigroup fuzzy bivector spaces over a bigroup $G = G_1 \cup G_2$.

(127) Let $V = V_1 \cup V_2 =$

$$\left\{ \begin{pmatrix} a & a & a \\ a & a & a \\ a & a & a \end{pmatrix} \middle| a \in Z_{23} I \right\} \cup$$

$$\left\{ \begin{pmatrix} a & a & a & a & a & a \\ a & a & a & a & a & a \end{pmatrix} \middle| a \in Z_{19}I \right\}$$

be a neutrosophic bigroup bilinear algebra over the bigroup $G = G_1 \cup G_2 = Z_{23}I \cup Z_{19}I$.

   a. Find a bigenerator of V.

   b. If G is replaced by $Z_{23} \cup Z_{19}$, what is the bigenerator of V?

   c. Define $\eta:V \to N([0, 1])$ so that $V_\eta = (V_1 \cup V_2)_{\eta_1 \cup \eta_2}$ is a neutrosophic bigroup fuzzy bilinear algebra.

   d. Can $V = V_1 \cup V_2$ have any proper neutrosophic bigroup bilinear subalgebras?

   e. Is V simple?

   f. Define two neutrosophic bigroup fuzzy bilinear algebras say $\eta$ and $\theta$ so that $\eta$ and $\theta$ do not agree on any element on V.

(128) Let V and W be any two distinct neutrosophic bisemigroup bivector spaces over the same bisemigroup $S = S_1 \cup S_2$. Find $N(\text{Hom}_S(V, W))$.

(129) Obtain some interesting properties about the neutrosophic bigroup bilinear algebras.

(130) Let

$$V = \left\{ \begin{pmatrix} a & a & a & a \\ b & b & b & b \\ c & c & c & c \end{pmatrix}, \begin{bmatrix} a & b \\ a & b \\ a & b \\ a & b \\ a & b \end{bmatrix} \middle| a,b,c \in ZI \right\} \cup$$

$$\left\{ (x,y), \begin{pmatrix} a & a & a & a & a \\ b & b & b & b & b \\ c & c & c & c & c \end{pmatrix}, [a,a,a,a] \middle| x,y,a,b,c \in ZI \right\}$$

and

$$W = \left\{ \begin{pmatrix} a & b & c \\ 0 & a & b \\ 0 & 0 & c \end{pmatrix}, \begin{bmatrix} a & 0 & 0 \\ a & b & 0 \\ a & b & c \end{bmatrix} \middle| a,b,c \in Z^+I \right\} \cup$$

$$\left\{ \begin{pmatrix} a & b & c & d \\ e & f & g & h \\ i & j & k & l \end{pmatrix}, (a,b,c,d,e,f) \middle| a,b,c,d,e,f,g,h,i,j,k,l \in ZI \right\}$$

be any two neutrosophic set bivector spaces over the set $S = Z^+I$. Find $N(\text{Hom}_S(V, W))$. What will happen if $S$ is replaced by $S_1 = \{0, 1, 2, 3, 4\}$.

  a. Will this affect the bigenerating subset of V and W?

  b. What is bidimension of V and W are neutrosophic set bivector spaces over the set S?

  c. What is the bidimension of V and W as neutrosophic set bivector spaces over the set $S_1$?

(131) Can one prove bidimension depends on the set over which the neutrosophic set bivector space is defined?

(132) Prove or disprove the bigenerator of a neutrosophic semigroup bilinear algebra V is dependent on the semigroup over which V is defined.

(133) Obtain some interesting properties about bidimension and bigenerators of the neutrosophic set bivector spaces.

(134) Suppose $V = V_1 \cup V_2$ is such that V can be treated as a neutrosophic set bivector space as well as neutrosophic set bilinear algebra over the same set S. Will they have different sets of bigenerators when $V = V_1 \cup V_2$ is just a neutrosophic set bivector space and another set of bigenerator when the same V is a neutrosophic set bilinear algebra over the same S.? Justify your claim.

(135) Let $V = V_1 \cup V_2$ be any neutrosophic bigroup bivector space over the bigroup G.

a. Is it possible for V to have more than one bigenerating bisubset of V?

b. Will the bidimension of V be the same if V has more than one bigenerating bisubset?

(136) Let $V = V_1 \cup V_2 =$

$$\left\{ \begin{pmatrix} a & b & c \\ d & e & f \end{pmatrix} \middle| a,b,c,d,e,f \in Z_2I \right\} \cup$$

$$\left\{ \begin{pmatrix} a & b \\ c & d \end{pmatrix} \middle| a,b,c,d \in Z_3I \right\}$$

be a neutrosophic bigroup bivector space over the bigroup $G = Z_2I \cup Z_3I$.

a. What is the bidimension of V?

b. How many bigenerating bisubset V has?

(137) Characterize those neutrosophic bigroup bivector spaces which have a unique generating subbiset!

(138) Give an example of a neutrosophic bigroup bivector space which has more than one bigenerating bisubset.

(139) Let $V = V_1 \cup V_2 =$

$$\left\{ \begin{pmatrix} a & a & a \\ b & b & b \\ c & c & c \end{pmatrix} \middle| a,b,c \in Z_{12}I \right\} \cup \left\{ \begin{pmatrix} a & a & a \\ a & a & a \\ a & a & a \end{pmatrix} \middle| a \in Z_{25}I \right\}$$

be a neutrosophic bigroup bilinear algebra over the bigroup $G = Z_{12}I \cup Z_{25}I$.

a. Find a bigenerating bisubset of V.

b. How many sets of bigenerating bisubset of V exist?

c. What is the bidimension of V?

(140) Let $V = V_1 \cup V_2 =$

$$\left\{ (x,y,z), \begin{bmatrix} x \\ y \\ z \\ w \end{bmatrix}, \begin{bmatrix} x & y \\ z & w \end{bmatrix}, \begin{bmatrix} x & y & z \\ 0 & y & 0 \end{bmatrix} \middle| x,y,z,w \in ZI \right\} \cup$$

$$\left\{ \begin{pmatrix} a & b & c & g \\ d & e & f & h \end{pmatrix}, (a,b,c,d), \begin{bmatrix} a \\ b \\ c \\ d \\ e \end{bmatrix}, \begin{bmatrix} a & b & c & d \\ 0 & e & f & g \\ 0 & 0 & h & i \\ 0 & 0 & 0 & f \end{bmatrix} \middle| a,b,c,d,e,f,g,h,i \in ZI \right\}$$

be a neutrosophic group bivector space over the group G = ZI.

a. What is the bidimension of V?

b. Find a bigenerator of V?

c. Find all proper neutrosophic group bivector subspaces of V. Is that collection finite or infinite?

d. Define $\eta : V \to N([0, 1])$ which can preserve all proper neutrosophic group bivector subspaces; $V_\eta$ is the neutrosophic group fuzzy bivector space.

e. Find a $\eta : V \to N([0, 1])$ which does not yield even a single $W_{\bar{\eta}}$ where W is a proper neutrosophic group bivector subspace of V; $\bar{\eta}$ the extension of $\eta$ on V.

f. Can V have pseudo neutrosophic semigroup bivector subspaces?

g. If group ZI is replaced by pZI, p a prime will $V = V_1 \cup V_2$ have different bidimension and bigenerator? Justify your answer.

(141) Give an example of neutrosophic bigroup bilinear algebra which is simple but of infinite bidimension.

(142) Give an example of a neutrosophic group bivector space which is simple and of finite bidimension (n, m), n and m are not prime.

(143) Suppose $V = V_1 \cup V_2$ is a neutrosophic bigroup bilinear algebra over a bigroup $G = G_1 \cup G_2$, where V is of finite bidimension, can one conclude V has only a finite number of neutrosophic bigroup bilinear subalgebras? Justify your claim.

(144) Let $V = V_1 \cup V_2 =$

$$\left\{ \begin{bmatrix} a & a & a \\ 0 & a & a \\ 0 & 0 & a \end{bmatrix} \middle| a \in Z_5 I \right\} \cup \left\{ \begin{bmatrix} a & 0 & 0 & 0 \\ a & a & 0 & 0 \\ a & a & a & 0 \\ a & a & a & a \end{bmatrix} \middle| a \in Z_{12} I \right\}$$

be a neutrosophic bigroup bilinear algebra over the bigroup $G = G_1 \cup G_2 = Z_5 I \cup Z_{12} I$.

    a. Define $\eta : V \to N[(0, 1)]$ so that $V_\eta$ is a neutrosophic bigroup fuzzy bilinear algebra.

    b. Can V have neutrosophic bigroup bilinear subalgebras?

    c. Can V have pseudo neutrosophic subbisemigroup bilinear subalgebras?

Hence or otherwise if $V = V_1 \cup V_2$ is a neutrosophic bigroup bilinear algebra over a bigroup $G = Z_p I \cup Z_n I$ (p a prime and n not a prime) of the form;

$$V = \left\{ \begin{bmatrix} a & a & a \\ 0 & a & a \\ 0 & 0 & a \end{bmatrix} \middle| a \in Z_p I \right\} \cup \left\{ \begin{bmatrix} a & 0 & 0 & 0 \\ a & a & 0 & 0 \\ a & a & a & 0 \\ a & a & a & a \end{bmatrix} \middle| a \in Z_n I \right\}$$

answer the above three questions.

What will happen if

$$V = V_1 \cup V_2$$

$$= \left\{ \begin{bmatrix} a & b & c \\ 0 & d & e \\ 0 & 0 & f \end{bmatrix} \middle| a,b,c,d,e,f \in Z_pI \right\} \cup$$

$$\left\{ \begin{bmatrix} a & 0 & 0 & 0 \\ b & c & 0 & 0 \\ d & e & f & 0 \\ g & h & i & j \end{bmatrix} \middle| a,b,c,d,e,f,g,h,i \in Z_nI \right\}$$

is a neutrosophic bigroup bilinear algebra over the bigroup;

a. if $G = Z_pI \cup Z_nI$?

b. If over $G = Z_p \cup Z_n$?

(145) Give an example of a neutrosophic bigroup bilinear algebra over a bigroup which has entries from $N(\mathbb{C})$.

(146) Let $V = V_1 \cup V_2$ be a neutrosophic biset bilinear algebra over the biset $S = S_1 \cup S_2$. Let $N(M_G(V, V))$ denote the set of all neutrosophic biset linear operators from V to V.

What is the algebraic structure of $N(M_G(V, V))$?

(147) Does the bidimension of a neutrosophic bisemigroup bilinear algebra depend on the bisemigroup over which it is defined?

(148) Let $V = V_1 \cup V_2 = \{Z_{12}I \times Z_{12}I \times Z_{12}I \times Z_{12}I\} \cup \{ Z_{12}I \times Z_{12}I \times Z_{12}I \times Z_{12}I \times Z_{12}I \times Z_{12}I, \begin{bmatrix} a & a \\ a & a \end{bmatrix} \middle| a \in Z_{12}I \}$ be a neutrosophic set bivector space over $S = Z_{12}$.

What is the bidimension of V? Find a bigenerating bisubset of V. What is the bidimension of V if $S = Z_{12}I$?

a. Will the bigenerating bisubset be different?

b. What is the bidimension of V if S if replaced by $S_1 = \{0, 6I\}$?

(149) Let $V = V_1 \cup V_2 = Z_5 I [x] \cup Z_2 I [x]$ be a neutrosophic bigroup bilinear algebra over the bigroup $G = Z_5 \cup Z_2$.

  a. What is the bidimension of V over G?

  b. If G is replaced by $G_1 = Z_5I \cup Z_2I$; what is the bidimension of V?

  c. If G is replaced by $G_2 = Z_5 \cup Z_2I$; what is the bidimension?

  d. If G is replaced by $H = Z_5I \cup Z_2$; what is the bidimension of V?

(150) Let $V = V_1 \cup V_2 =$

$$\left\{ \begin{bmatrix} a & a & a \\ a & a & a \end{bmatrix} \middle| a \in ZI \right\} \cup \{ZI \times ZI \times ZI\}$$

and

$$W = \left\{ \begin{bmatrix} a & b \\ c & d \\ e & f \end{bmatrix} \middle| a,b,c,d,e,f \in ZI \right\} \cup$$

$$\left\{ \begin{bmatrix} a & a & a & a \\ b & b & b & b \\ c & c & c & c \end{bmatrix} \middle| a,b,c \in ZI \right\}$$

be neutrosophic group linear algebras over the group $G = ZI$.

  a. Find bidimensions of V and W.

  b. What is the bidimension of $N(M_{ZI}(V, W))$?

  c. If G is replaced by pZI, p a prime find the bidimensions of V and W.

  d. What is the bidimension of $N(M_{pZI}(V, W))$?

e. Find Vη a neutrosophic group fuzzy bilinear algebra.

(151) Let $V = V_1 \cup V_2 =$

$$\left\{ \begin{bmatrix} a & b & c \\ d & e & f \end{bmatrix} \bigg| a, b, c, d, e, f \in Z_{30}I \right\} \cup$$

$$\left\{ \begin{bmatrix} a & a & a \\ a & a & a \\ a & a & a \end{bmatrix} \bigg| a \in Z_{30}I \right\}$$

be a neutrosophic group bilinear algebra over the group $G = Z_{30}I$.

a. Can V have pseudo neutrosophic semigroup bilinear algebras over $G = Z_{30}I$?

b. What is the bidimension of V over G?

c. What is the bidimension of $N(M_G(V, V))$?

(152) Let $V = V_1 \cup V_2 = \{Z_8 \times Z_8 \times Z_8 \times Z_8I\} \cup \{Z_8I \times Z_8 I \times Z_8 I \times Z_8I \times Z_8 I \times Z_8 I\}$ be a neutrosophic group bilinear algebra over the group $G = \{0, 4\}$ addition under modulo 8.

a. Find the bidimension of V over $G = \{0, 4\}$.

b. What is the bigenerating set of V over $G = \{0, 4\}$?

c. If G is replaced by $H = \{0, 2, 4, 6\}$, what is the bidimension of V and the bigenerating set of V over H?

d. What is the bidimension of $N(M_H(V, V))$?

(153) Let $V = V_1 \cup V_2$ be a neutrosophic group bilinear algebra over the group G. Suppose $H_1, H_2, \ldots, H_n$ be n distinct subgroups of G. Suppose $V = V_1 \cup V_2$ is also a neutrosophic group bilinear algebra over each of the subgroups $H_1, H_2, \ldots, H_n$; compare them.

(154) Let $V = \left\{ \begin{bmatrix} a & b \\ c & d \end{bmatrix} \middle| a,b,c,d \in Z_{18}I \right\} \cup \{(a\ a\ a\ a) \mid a \in Z_{18}\ I\}$

be a neutrosophic semigroup bilinear algebra over the semigroup $S = Z_{18}$.

  a. Define $\eta = \eta_1 \cup \eta_2 : V \to N([0, 1])$ so that $V_\eta$ is a neutrosophic semigroup fuzzy bilinear algebra.

  b. Is every bimap $\eta : V \to N([0, 1])$ is such that, $V\eta$ is a neutrosophic semigroup fuzzy bilinear algebra?

(155) Let $V = V_1 \cup V_2 = \{ZI \times Z \times ZI \times ZI \times Z\} \cup$

$\left\{ \begin{bmatrix} a & b & c \\ d & e & f \end{bmatrix} \middle| a,b,c \in ZI, d,e,f, \in Z \right\}$

be a neutrosophic group bilinear algebra over the group $G = Z$.

Find $\eta : V \to N([0, 1])$ so that $V\eta$ is a neutrosophic group fuzzy bilinear algebra.

(156) Give some interesting properties about neutrosophic group fuzzy bivector spaces.

(157) Let $V = V_1 \cup V_2$ be any neutrosophic quasi semigroup bilinear algebra over the semigroup S.

Find some interesting properties about this algebraic structure.

(158) Let $V = V_1 \cup V_2 = Z_{12}I[x] \cup Z^+I[x]$ be a neutrosophic bisemigroup bivector space defined over the bisemigroup $S = Z_{12}I \cup Z^+I$.

  a. Find a neutrosophic bisemigroup bivector subspace of V.

  b. Find a bigenerating bisubset of V.

  c. What is the bidimension of V?

(159) Let $V = V_1 \cup V_2 =$

$$\left\{ \begin{bmatrix} a & b \\ c & d \end{bmatrix} \middle| a,b,c,d \in N(Z_{12}I) \right\} \cup \{N(Z_{12})I\,[x]\}$$

be a neutrosophic group bivector space over the group $G = Z_{12}$.

  a. Find pseudo group bivector subspace of V.
  b. Find the bidimension of V.
  c. If G is replaced by $G_1 = Z_{12}I$ can V have pseudo group bivector subspaces?
  d. Define $\eta : V \to N([0, 1])$ such that $V\eta$ is a neutrosophic group fuzzy bivector space.
  e. Can V have pseudo neutrosophic group bilinear subalgebras?

(160) Let $V = V_1 \cup V_2 = \{(Z_7I)^{12}\} \cup \{(Z_6I)^8\}$ be a neutrosophic bigroup bivector space over the bigroup $G = Z_7 \cup Z_6$.

  a. What is bidimension of V over G?
  b. If G is replaced by $G_1 = Z_7I \cup Z_6I$ what is the bidimension of V over G?
  c. Does V have a pseudo neutrosophic bisemigroup bivector subspace W of V over G?

(161) Let $V = V_1 \cup V_2 =$

$$\left\{ \begin{bmatrix} a & b \\ c & d \end{bmatrix} \middle| a,b,c,d \in 3ZI \right\} \cup$$

$$\left\{ \begin{bmatrix} a & b \\ c & d \end{bmatrix} \middle| a,b,c,d \in 5ZI \right\}$$

be a neutrosophic bigroup bivector space over the bigroup $G = 3ZI \cup 5ZI$.

If
$$W = W_1 \cup W_2 =$$

$$\left\{ \begin{bmatrix} a & a \\ 0 & 0 \end{bmatrix} \middle| a \in 3ZI \right\} \cup \left\{ \begin{bmatrix} a & a \\ 0 & 0 \end{bmatrix} \middle| a \in 5ZI \right\}$$

$\subseteq V_1 \cup V_2$ is a neutrosophic bigroup bivector subspace of V, find the bidimension of W over G.

# FURTHER READING

# INDEX

# ABOUT THE AUTHORS

**Dr.W.B.Vasantha Kandasamy** is an Associate Professor in the Department of Mathematics, Indian Institute of Technology Madras, Chennai. In the past decade she has guided 12 Ph.D. scholars in the different fields of non-associative algebras, algebraic coding theory, transportation theory, fuzzy groups, and applications of fuzzy theory of the problems faced in chemical industries and cement industries.

She has to her credit 646 research papers. She has guided over 68 M.Sc. and M.Tech. projects. She has worked in collaboration projects with the Indian Space Research Organization and with the Tamil Nadu State AIDS Control Society. This is her 45$^{th}$ book.

On India's 60th Independence Day, Dr.Vasantha was conferred the Kalpana Chawla Award for Courage and Daring Enterprise by the State Government of Tamil Nadu in recognition of her sustained fight for social justice in the Indian Institute of Technology (IIT) Madras and for her contribution to mathematics. (The award, instituted in the memory of Indian-American astronaut Kalpana Chawla who died aboard Space Shuttle Columbia). The award carried a cash prize of five lakh rupees (the highest prize-money for any Indian award) and a gold medal.
She can be contacted at vasanthakandasamy@gmail.com
Web Site: http://mat.iitm.ac.in/home/wbv/public_html

**Dr. Florentin Smarandache** is a Professor of Mathematics at the University of New Mexico in USA. He published over 75 books and 150 articles and notes in mathematics, physics, philosophy, psychology, rebus, literature.

In mathematics his research is in number theory, non-Euclidean geometry, synthetic geometry, algebraic structures, statistics, neutrosophic logic and set (generalizations of fuzzy logic and set respectively), neutrosophic probability (generalization of classical and imprecise probability). Also, small contributions to nuclear and particle physics, information fusion, neutrosophy (a generalization of dialectics), law of sensations and stimuli, etc. He can be contacted at smarand@unm.edu

**K. Ilanthenral** is the editor of The Maths Tiger, Quarterly Journal of Maths. She can be contacted at ilanthenral@gmail.com